\documentclass[12pt,a4paper]{article} 

\usepackage{amsmath,amssymb}
\usepackage{amssymb}
 \usepackage{amsfonts} \usepackage{amsmath}

 \usepackage{hyperref}
\usepackage{tikz,fullpage}
\usetikzlibrary{backgrounds}
\usetikzlibrary{calc}

\newtheorem{lemma}{Lemma}[section]
\newtheorem{proposition}{Proposition}[section]

\newtheorem{conjecture}{Conjecture}[section]

\newtheorem{example}{Example}[section]
\newcommand{\sy}{  \mathbf{y}_1  }
\newcommand{\py}{  \mathbf{y}_2  }

\newcommand{\bee}[1]{  {b}_{#1}  }
\newcommand{\BEE}[1]{  \mathbf{b}_{#1}  }

\newcommand{\hesi}{ \|\hspace{-0.27 ex}\|\hspace{-0.27 ex}\| }

\newcommand{\kos}[2]{[\frac{#1}{#2}] }

\newcommand{\Mir}{\mathrm{Mir}}
\newcommand{\mrMir}{\mathrm{^{\star}Mir}}
\newcommand{\mr}[1]{\mathrm{#1}}
\newcommand{\mi}[1]{\mathit{#1}}

\newcommand{\imia}[1]{\stackrel{\Diamond}{#1}}
\newcommand{\imi}[1]{\stackrel{\wedge}{#1}}
\newcommand{\ima}[1]{\stackrel{\vee}{#1}}

\newcommand{\smia}[1]{\stackrel{\mr{0}}{#1}}
\newcommand{\smi}[1]{\stackrel{\frown}{#1}}
\newcommand{\sma}[1]{\stackrel{\smile}{#1}}

\newcommand{\gbar}{{-\hspace{-1.5 ex}g}}
\newcommand{\fbar}{{-\hspace{-1.5 ex}f}}
\newcommand{\abar}{{-\hspace{-1.5 ex}a}}
\newcommand{\bbar}{{-\hspace{-1.5 ex}b}}

\newcommand{\elliast}{{{\ell i}^\ast\!}}
\newcommand{\mg}[1]{{\text{\boldmath$#1$}}}

\newcommand{\mga}{{\mathbf{a}}}
\newcommand{\mgb}{{\mathbf{b}}}
\newcommand{\mgc}{{\mathbf{c}}}
\newcommand{\mgd}{{\mathbf{d}}}  \newcommand{\mgD}{{\mathbf{D}}}

\newcommand{\mgf}{{\mathbf{f}}}

\newcommand{\mgp}{{\mathbf{p}}}
\newcommand{\mgq}{{\mathbf{q}}}

\newcommand{\bu}{\bullet}

\newcommand{\caB}{\mathcal{B}}

\newcommand{\caI}{\mathcal{I}}
\newcommand{\caJ}{\mathcal{J}}
\newcommand{\caK}{\mathcal{K}}
\newcommand{\caL}{\mathcal{L}}
\newcommand{\caM}{\mathcal{M}}

\newcommand{\caP}{\mathcal{P}}
\newcommand{\caQ}{\mathcal{Q}}
\newcommand{\caR}{\mathcal{R}}

\newcommand{\caW}{\mathcal{W}}

\newcommand{\doC}{\mathbb{C}}
\newcommand{\doD}{\mathbb{D}}

\newcommand{\doH}{\mathbb{H}}
\newcommand{\doI}{\mathbb{I}}
\newcommand{\doJ}{\mathbb{J}}

\newcommand{\doN}{\mathbb{N}}

\newcommand{\doQ}{\mathbb{Q}}
\newcommand{\doR}{\mathbb{R}}
\newcommand{\doS}{\mathbb{S}}

\newcommand{\doZ}{\mathbb{Z}}

\begin{document}

%%%%%%%%%%%%%%%%%%%%%%%%%%%%%%%%%%%%%%%%%%%%%%%%%%%%%%%%%%%%%%%%%%%%%%%%%%%%%%%%%%%%%%%%%%%%
%%%%%%%%%%%%%%%%%%%%%%%%%%%%%%%%%%%%%%%%%%%%%%%%%%%%%%%%%%%%%%%%%%%%%%%%%%%%%%%%%%%%%%%%%%%%
%%%%%%%%%%%%%%%%%%%%%%%%%%%%%%%%%%%%%%%%%%%%%%%%%%%%%%%%%%%%%%%%%%%%%%%%%%%%%%%%%%%%%%%%%%%%
\title{Power series with sum-product Taylor coefficients and their resurgence algebra. }

\maketitle

\author{Jean Ecalle (Orsay, CNRS) , Shweta Sharma (Orsay)}
 \\
 \\
 \\
 \noindent    
{\bf Abstract.}{\it  The present paper is devoted to power series of SP type, i.e. with coefficients
that are syntactically sum-product combinations. Apart from their applications to analytic knot
theory and the so-called ``Volume Conjecture'', SP-series are interesting in their own right, on
at least four counts\,: (i) they generate quite distinctive resurgence algebras (ii) they are one
of those relatively rare instances when the resurgence properties have to be derived directly
from the Taylor coefficients (iii) some of them produce singularities that unexpectedly verify  finite-order
differential equations  (iv) all of them are best handled with the help of two remarkable,
infinite-order integral-differential transforms, \textup{mir} and \textup{nir}. }
\tableofcontents  
   
%

%%%%%%%%%%%%%%%%%%%%%%%%%%%%%%%%%%%%%%%%%%%%%%%%%%%%%%%%%%
%% Some resurgence properties of knot-related functions.
%%%%%%%%%%%%%%%%%%%%%%%%%%%%%%%%%%%%%%%%%%%%%%%%%%%%%%%%%%

%\documentclass[12pt,a4paper]{article}\input{SP_commands}\begin{document}

%%%%%%%%%%%%%%%%%%%%%%%%%%%%%%%%%%%%%%%%%%%%%%%%%%%%%%%%%%%%%%%%%%%%%%%%%%%%%%%%%%%%%%%%%%%%
%%%%%%%%%%%%%%%%%%%%%%%%%%%%%%%%%%%%%%%%%%%%%%%%%%%%%%%%%%%%%%%%%%%%%%%%%%%%%%%%%%%%%%%%%%%%
%%%%%%%%%%%%%%%%%%%%%%%%%%%%%%%%%%%%%%%%%%%%%%%%%%%%%%%%%%%%%%%%%%%%%%%%%%%%%%%%%%%%%%%%%%%%

\section{Introduction.}

%%%%%%%%%%%%%%%%%%%%%%%%%%%%%%%%%%%%%%%%%%%%%%%%%%%%%%%%%%%%%%%%%%%%%%%%%%%%%%%%%%%%%%%%%%%%
%%%%%%%%%%%%%%%%%%%%%%%%%%%%%%%%%%%%%%%%%%%%%%%%%%%%%%%%%%%%%%%%%%%%%%%%%%%%%%%%%%%%%%%%%%%%
\subsection{Power series with coefficients of sum-product type.}

{\bf The notion of SP series.}\\
  {\it Sum-product} series (or {\it SP}-series for short) are Taylor series\,:
\begin{equation}\label{a1}
j(\zeta):=\sum_{n\geq 0} \;J(n)\; \zeta^n
\end{equation} 
whose coefficients are syntactically of sum-product (SP) type\,: 
\begin{equation}\label{a2}
J(n):=\!\!\sum_{\epsilon\leq m < n}
\prod_{\epsilon \leq k \leq m}F(\frac{k}{n})=
\!\!\!\sum_{\epsilon \leq m < n }\!\!\exp\Big(\!\!-\!\!\sum_{\epsilon \leq k
\leq m} f(\frac{k}{n})\Big)
\quad (\epsilon \in\{0,1\})
\end{equation}
Summation starts at $\epsilon=0$ unless $F(0)\in\{0,\infty\}$,
in which case it starts at $\epsilon=1$.
It always ends at $n-1$, not $n$.\footnote{
This choice is to ensure near-invariance
under the change $F(x)\mapsto 1/F(1-x)$. See \S3.5}
The two driving functions $F$ and $f$ are connected under $F\equiv \exp(-f)$.
Unless stated otherwise, $F$ will be assumed to be meromorphic,
and special attention shall be paid to the case when $F$ has neither
zeros nor poles, i.e. when $f$ is holomorphic.

The importance of SP-series comes from their analytic properties (isolated
singularities of a quite distinctive type) and their frequent occurence
in various fields of mathematics (ODEs, knot theory etc). 

As for the above definition, it is less arbitrary than may seem at first sight.
Indeed, none of the following changes\,:

(i) changing the grid $\{k/n\}$ to $\{\mi{Const}\,k/n\}$

(ii) changing the lower summation bounds from 0 or 1 to 2,3 \dots

(iii) changing the upper summation bound from $n-1$  to $n$ or
$n\!-\!2, n\!-\!3 $ etc or a multiple thereof 

(iv) replacing the 0-accumulating products $\prod F(\frac{k}{n})$
by 1-accumulating products $\prod F(\frac{n-k}{n})$
\\
- none of these changes, we claim, would make much difference or even (allowing for minor adjustments)
take us beyond the class of SP-series.
\\

\noindent
{\bf Special cases of SP series.} 
\\
For $F$ a polynomial or rational function (resp. a trigonometric polynomial)
and for Taylor coefficients $J(n)$ defined by pure products $\prod$ (rather than
sum-products $\sum\prod$)
the series $j(\zeta)$ would be
of hypergeometric (resp. $q$-hypergeometric) type. Thus the theory
of SP-series extends \---
and bridges \--- two important fields. But it covers wider ground. In fact, the
main impulse for developping it came from knot theory, 
and we didn't get involved in the subject until Stavros Garoufalidis\footnote{
an expert in knot theory who visited Orsay in the fall of 2006.
}
and Ovidiu Costin\footnote{
an analyst who together with S. Garoufalidis has been pursuing an approach to the subject 
parallel to ours, but distinct \,: for a comparison, see \S12.1.
} 
drew our attention to its potential.
\\

\noindent
{\bf Overview.}
\\ In this first paper, halfway between survey and full treatment\footnote{Two follow-up investigations [SS1],[SS2] are being planned. },
we shall attempt five things\,:
\\

(i) bring out the {\it main analytic features} of SP-series, such as the dichotomy
between their two types of singularities ({\it outer/inner}), and produce
{\it complete systems of resurgence equations}, which encode in compact form
the whole Riemann surface structure.
\\

(ii)  {\it localise} and {\it formalise} the problem, i.e. break it down into
the separate study of a number of {\it local} singularities, each of which is
produced by a specific {\it non-linear functional transform} capable of
a full analytical description, which reduces everything to {\it formal} 
manipulations on power series.
\\

(iii) sketch the {\it general picture} for arbitrary driving functions $F$ and
$f$\, \--- pending a future, detailed investigation.
\\

(iv) show that in many instances ($f$ polynomial, $F$ monomial or even just rational) our local singularities
satisfy ordinary differential equations, but of a very distinctive type, which accounts for the `rigidity' of their
resurgence equations, i.e. the occurrence in them of essentially {\it discrete} Stokes constants.\footnote{contrary
to the usual situation, where these Stokes or resurgence constants are free to vary continuously.}
\\

(v) sketch numerous examples and then give a careful treatment, theoretical {\it and} numerical,
of {\it one special case} chosen for its didactic value (it illustrates all
the main SP-phenomena) and its practical relevance to knot theory
(specifically, to the knot $4_1$).

%%%%%%%%%%%%%%%%%%%%%%%%%%%%%%%%%%%%%%%%%%%%%%%%%%%%%%%%%%%%%%%%%%%%%%%%%%%%%%%%%%%%%%%%%%%%
%%%%%%%%%%%%%%%%%%%%%%%%%%%%%%%%%%%%%%%%%%%%%%%%%%%%%%%%%%%%%%%%%%%%%%%%%%%%%%%%%%%%%%%%%%%%
\subsection{The outer/inner dichotomy and the ingress factor.}
{\bf The outer/inner dichotomy.}
\\ 
Under analytic continuation, SP-series give
rise to two distinct types of {\it singularities}, also referred to
as {\it `generators'}, since under alien derivation they generate the resurgence algebra
of our SP-series. On the one hand, we  have the
{\it outer generators}, so-called because they never recur under alien derivation (but
produce inner generators), and on the other hand we have the {\it inner generators},
so-called because they recur indefinitely under
alien derivation (but never re-produce the outer generators). These two are,
by any account, the main types of generators, but for completeness we add two further
classes\,: the {\it original generators} (i.e. the SP-series themselves) and
the {\it exceptional generators}, which don't occur naturally, but
can prove useful as auxiliary adjuncts.
 \\

\noindent
{\bf A gratifying surprise\,: the {\it mir}-transform.}
\\
We shall see that {\it outer}  generators can be viewed as infinite sums
of {\it inner generators}, and that the latter 
 can be constructed quite explicitely
by subjecting the driving function $F$
to a chain of nine local transforms, all of which are elementary,
save for one crucial step\,:
the {\it mir}-transform. Furthermore, this 
 {\it mir}-transform , though resulting from
an unpromising mix of complex operations\footnote{
two Laplace transforms, direct and inverse, with a few violently
non-linear operations thrown in.
},
will turn out to be an {\it integro-differential operator}, of infinite order
but with a transparent expression that sheds much light on its analytic properties.
We regard this fascinating {\it mir}-transform, popping out of nowhere yet highly
helpful, as the centre-piece of this investigation.
\\

\noindent
{\bf The ingress factor and the cleansing of SP-series.}
\\
Actually, rather than directly considering the SP-series $j(\zeta)$
with coefficients $J(n)$, it shall prove expedient to study the
slightly modified series  $j^\#(\zeta)$
with coefficients $J^\#(n)$ obtained after division
by a suitably defined {\it `ingress factor'} $I\!g_F(n)$ of strictly local character\,:
\begin{equation}\label{a3}
j^\#(\zeta)=\sum J^\#(n) \,\zeta^n\quad\quad \mi{with}\quad\quad J^\#(n):=J(n)/I\!g_F(n)
\end{equation}
This purely technical trick involves no loss of information\footnote{
since information about $j^\#$ immediately translates into information about $j$,
and vice versa.} 
and achieves two things\,:
\\

(i) the various {\it outer} and {\it inner} generators will now appear as
 purely {\it local} transforms of the driving function $F$, taken at suitable
base points.
\\

(ii) distinct series $j_{F_i}(\zeta)$ relative to distinct base points $x_i$
(or, put another way, to distinct translates $F_i(x):=F(x+x_i)$ of the same
driving function) will lead to exactly the same {\it inner generators}
and so to the same {\it inner algebra}\;\--- which wouldn't be the case
but for the pre-emptive removal of $I\!g_F$. 
\\

In any case, as we shall see, the ingress factor is a relatively innocuous
function and (even when it is divergent-resurgent, as may happen) the effect
not only of {\it removing it} but also, if we so wish, of {\it putting it back} can
be completely mastered.
%%%%%%%%%%%%%%%%%%%%%%%%%%%%%%%%%%%%%%%%%%%%%%%%%%%%%%%%%%%%%%%%%%%%%%%%%%%%%%%%%%%%%%%%%%%%
%%%%%%%%%%%%%%%%%%%%%%%%%%%%%%%%%%%%%%%%%%%%%%%%%%%%%%%%%%%%%%%%%%%%%%%%%%%%%%%%%%%%%%%%%%%%
\subsection{The four gates to the inner algebra.}
%%%%%%%%%%%%%%%%%%%%%%%%%%%%%%%%%%%%%%%%%%%%%%%%%%%%%%%%%%%%%%%%%%%%%%%%%%%%%%%%%%%%%%%%%%%%
%%%%%%%%%%%%%%%%%%%%%%%%%%%%%%%%%%%%%%%%%%%%%%%%%%%%%%%%%%%%%%%%%%%%%%%%%%%%%%%%%%%%%%%%%%%%
We have just described the various types of singularities or `generators'
we are liable to encounter when analytically continuing a SP-series. Amongst these, as we saw,
the {\it inner generators} stand out. They span the {\it inner algebra}, which
is the problem's hard, invariant core. Let us now review the situation
once again, but from another angle, by asking\,: how many {\it gates}
are there for entering the unique {\it inner algebra} ? There are,
in effect, four types\,:
\\

\noindent
{\bf Gates of type 1\,: original generators}. We may of course enter through an {\it
original generator}, i.e. through a  {\it SP series}, relative to any base point
$x_0$ of our choosing. Provided we remove the corresponding ingress factor, we shall always
arrive at the same inner algebra.
\\ 

\noindent
{\bf Gates of type 2\,: outer generators}. We may enter through an {\it
outer generator}, i.e. through the mechanism of the nine-link chain of
section 5, again relative to any base point. Still, when $F$ does have zeros $x_i$, 
these qualify as privileged base points, since in that case we can
make do with the simpler four-link chain of section 5.
\\ 
\\
{\bf Gates of type 3\,: inner generators.} We may enter through an {\it inner generator}, i.e. 
 via the mechanism of the nine-link chain of section 4, but only from a base point $x_i$
where $f$ (not $F$ !) vanishes\,.    
By so doing,
 we do not properly speaking {\it enter} the inner algebra, but rather start {\it right
there}.  Due to the ping-pong phenomenon, this inner generator then generates all the
other ones. The method, though, has the drawback of introducing a jarring
dissymmetry, by giving precedence to {\it one} inner generator over all others.
\\
\\ 
{\bf Gates of type 4\,: exceptional generators.} We may enter through
 a {\it exceptional} or {\it `mobile' generator}, i.e. once again
 via the mechanism of the nine-link chain of section 4, but relative to
any base point $x_0$
where $f$ doesn't vanish\footnote{
so that the so-called {\it tangency order} is  $\kappa=0$,
whereas for the inner generators it is $\geq 1$ and generically $=1$. See \S4.}. It turns
out that any such ``exceptional generator" generates all
the inner generators (\--- and what's more, {\it symmetrically} so \---), but isn't
generated by them. In other words, it gracefully {\it self-eliminates}, thereby atoning
for its parasitical character. Exceptional
generators, being `mobile', have the added advantage that their base point $x_0$ can be taken
arbitrarily close to the base point $x_i$ of any given inner generator, which fact proves
quite helpful, computationally and also theoretically.

%\end{document}
%%%%%%%%%%%%%%%%%%%%%%%%%%%%%%%%%%%%%%%%%%%%%%%%%%%%%%%%%%%%%%%%%%%%%%%%%%%%%%%%%%%%%%%%%%%%
%%%%%%%%%%%%%%%%%%%%%%%%%%%%%%%%%%%%%%%%%%%%%%%%%%%%%%%%%%%%%%%%%%%%%%%%%%%%%%%%%%%%%%%%%%%%
%%%%%%%%%%%%%%%%%%%%%%%%%%%%%%%%%%%%%%%%%%%%%%%%%%%%%%%%%%%%%%%%%%%%%%%%%%%%%%%%%%%%%%%%%%%%
%%%%%%%%%%%%%%%%%%%%%%%%%%%%%%%%%%%%%%%%%%%%%%%%%%%%%%%%%%%%%%%%%%%%%%%%%%%%%%%%%%%%%%%%%%%%
%%%%%%%%%%%%%%%%%%%%%%%%%%%%%%%%%%%%%%%%%%%%%%%%%%%%%%%%%%%%%%%%%%%%%%%%%%%%%%%%%%%%%%%%%%%%
%%%%%%%%%%%%%%%%%%%%%%%%%%%%%%%%%%%%%%%%%%%%%%%%%%%%%%%%%%%%%%%%%%%%%%%%%%%%%%%%%%%%%%%%%%%%
%%%%%%%%%%%%%%%%%%%%%%%%%%%%%%%%%%%%%%%%%%%%%%%%%%%%%%%%%%%%%%%%%%%%%%%%%%%%%%%%%%%%%%%%%%%%
%%%%%%%%%%%%%%%%%%%%%%%%%%%%%%%%%%%%%%%%%%%%%%%%%%%%%%%%%%%%%%%%%%%%%%%%%%%%%%%%%%%%%%%%%%%%
%%%%%%%%%%%%%%%%%%%%%%%%%%%%%%%%%%%%%%%%%%%%%%%%%%%%%%%%%%%%%%%%%%%%%%%%%%%%%%%%%%%%%%%%%%%%

%

%%%%%%%%%%%%%%%%%%%%%%%%%%%%%%%%%%%%%%%%%%%%%%%%%%%%%%%%%%
%% Some resurgence properties of knot-related functions.
%%%%%%%%%%%%%%%%%%%%%%%%%%%%%%%%%%%%%%%%%%%%%%%%%%%%%%%%%%

%\documentclass[12pt,a4paper]{article}\input{SP_commands}\begin{document}

%%%%%%%%%%%%%%%%%%%%%%%%%%%%%%%%%%%%%%%%%%%%%%%%%%%%%%%%%%%%%%%%%%%%%%%%%%%%%%%%%%%%%%%%%%%%
%%%%%%%%%%%%%%%%%%%%%%%%%%%%%%%%%%%%%%%%%%%%%%%%%%%%%%%%%%%%%%%%%%%%%%%%%%%%%%%%%%%%%%%%%%%%
%%%%%%%%%%%%%%%%%%%%%%%%%%%%%%%%%%%%%%%%%%%%%%%%%%%%%%%%%%%%%%%%%%%%%%%%%%%%%%%%%%%%%%%%%%%%
 
%
\setcounter{section}{1}
%\addtocounter{2}
\section{Some resurgence background.}
%%%%%%%%%%%%%%%%%%%%%%%%%%%%%%%%%%%%%%%%%%%%%%%%%%%%%%%%%%%%%%%%%%%%%%%%%%%%%%%%%%%%%%%%%%%%
%%%%%%%%%%%%%%%%%%%%%%%%%%%%%%%%%%%%%%%%%%%%%%%%%%%%%%%%%%%%%%%%%%%%%%%%%%%%%%%%%%%%%%%%%%%%
\subsection{Resurgent functions and their three models.}
%%%%%%%%%%%%%%%%%%%%%%%%%%%%%%%%%%%%%%%%%%%%%%%%%%%%%%%%%%%%%%%%%%%%%%%%%
{\bf The four models\,: formal, geometric, upper/lower convolutive.}
\\
Resurgent `functions' exist in three/four types of models\,:
\\

(i) The {\it formal model}, consisting of formal power series $\tilde{\varphi}(z)$
of a variable $z\sim\infty$. The tilda points to the quality of being `formal',
i.e. possibly divergent. 
\\

(ii) The {\it geometric models} of direction $\theta$. They consist of 
sectorial analytic germs $\varphi_\theta(z)$ of the same variable $z\sim\infty$,
defined on sectors of aperture $>\pi$ and bisected by the axis $\arg(z^{-1})=\theta$.
\\

(iii) The {\it convolutive model}, consisting of {\it `global microfunctions'}
of a variable $\zeta\sim 0$. Each microfunction possesses one {\it minor} (exactly
defined, but with some {\it information missing}) and many {\it majors} (defined up to
regular germs at the origin, i.e. with some {\it redundant information}).
However, under a frequently fulfilled integrability condition at $\zeta\sim 0$,
the {\it minor} contains all the information, i.e. fully determines the
microfunction, in which case all calculations reduce to manipulations on
the sole minors.
As for the {\it `globalness'} of our microfunctions, it means that their minors
possess the property of ``endless analytic continuation"\,: they can be 
continued analytically in the $\zeta$-plane along any given (self-avoiding
or self-intersecting, whole or punctured\footnote{
the broken line may be punctured at a finite number  of 
singularity-carrying points $\zeta_1,\dots,\zeta_n$, in which case we demand analytic
continuability along {\it all} the $2^n$ paths that follow the broken line but
circumvent each  $\zeta_i$ to the right or to the left.
}) broken line starting from 0 and ending anywhere we like.

Usually, one makes do with a single convolutive model, but here it will be
convenient to adduce two of them\,: the {\it `upper'} and {\it `lower'}
models. In both, the minor-major relation has the same form\,;
\[\begin{array}{ccccccccc}
&\mi{minor}\quad&& \mi{major}
\\[1.5 ex]
\mi{upper}\quad&\smi{\varphi}(\zeta) &\equiv&
\frac{1}{2\pi i}\big(\sma{\varphi}(\zeta\,e^{-\pi i})-\sma{\varphi}(\zeta\,e^{\pi i})\big)
\\[1.5 ex]
\mi{lower}\quad&\imi{\varphi}(\zeta) &\equiv&
\frac{1}{2\pi i}\big(\ima{\varphi}(\zeta\,e^{-\pi i})-\ima{\varphi}(\zeta\,e^{\pi i})\big)
\end{array}\]
while the upper-lower correspondence goes like this\,:
\begin{equation}\label{a4}
\imi{\varphi}(\zeta)\equiv\, \partial_\zeta\smi{\varphi}(\zeta)
\quad \quad ; \quad\quad
\ima{\varphi}(\zeta)\equiv\, -\partial_\zeta\sma{\varphi}(\zeta)
\end{equation}
One of the points of resurgent analysis is to resum divergent
series of `natural origin', i.e. to go from the formal model to
the geometric one via one of the two convolutive models. 
Concretely, we go from {\it formal} to {\it convolutive}
by means of a formal or term-wise Borel transform\footnote{
thus {\it upper} (resp. {\it lower}) Borel takes 
$\sum a_n\,z^{-n} $ to 
$\sum a_n\,\frac{\zeta^{n}}{n!} $ 
(resp. $\sum a_n\,\frac{\zeta^{n-1}}{(n-1)!} $ ).
}
and from {\it convolutive} to {\it geometric} by means of a
$\theta$-polarised Laplace transform,
i.e. with integration along the half-axis $\arg(\zeta)=\theta$.

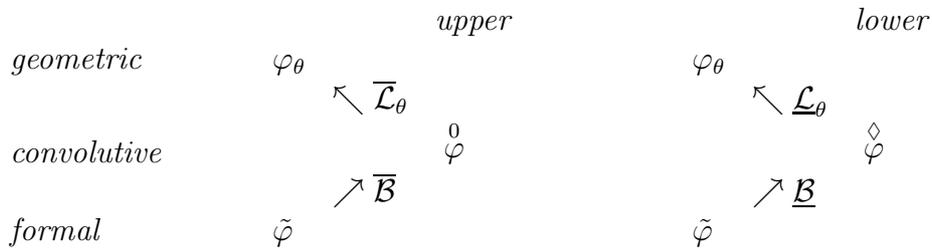
\[\begin{array}{llllllllllllllll}
 &&& \mi{upper}&&&& \mi{lower} 
\\
\mi{geometric}& \varphi_\theta &&&&\varphi_\theta
\\
&&\nwarrow \overline{\caL}_\theta &&&&\nwarrow \underline{\caL}_\theta 
\\
\mi{convolutive}\hspace{6. ex}&&&\;\smia{\varphi} &\quad\quad\quad\quad&&&\;\imia{\varphi}
\\
&&\nearrow \overline{\caB} &&&&\nearrow \underline{\caB}
\\
\mi{formal}&\tilde{\varphi} &&&&\tilde{\varphi}
\end{array}\]
%%%%%%%%%%%%%%%%%%%%%%%%%%%%%%%%%%%%%%%%%%%%%%%%%%%%%%%%%%%%%%%%%%%%%%%%%
{\bf Resurgent algebras\,: the multiplicative structure.}
\\
Resurgent functions are stable not just under {\it addition}
(which has the same form in all models) but also under a {\it product}
whose shape varies from model to model\,:
\\
(i) in the {\it formal} model, it is the ordinary
multiplication of power series.
\\
(ii) in the {\it geometric} model, it is the pointwise
multiplication of analytic germs.
\\
(iii) in the {\it convolutive} models, it is the {\it upper/lower}
convolution, with distinct expressions for minors\footnote{
minor convolution is possible only under a suitable
integrability condition at the origin. That condition
is automatically met when {\it one} (hence {\it all}) majors
verify $\sma{\varphi}\!\!\!(\zeta)\rightarrow 0$ or 
$\zeta\ima{\varphi}\!(\zeta)\rightarrow 0$
as $\zeta \rightarrow 0 $ uniformly on any sector of finite aperture.
}
 and majors\,: 
\[\begin{array}{ccccccccc}
\mi{minor\; convolution}&&&&
\\[1.5 ex]
\mi{upper\;}\;\overline{\ast}\quad&
(\smi{\varphi}_1 \overline{\ast} \smi{\varphi}_2)(\zeta)
&:=&\int_{0}^{\zeta}\;\smi{\varphi}_1\!(\zeta\!-\!\zeta_2)\; d\!\smi{\varphi}_2\!(\zeta_2)
\hspace{8.ex}
\\[1.5 ex]
\mi{lower\;}\;\underline{\ast}\quad&
(\imi{\varphi}_1 \underline{\ast} \imi{\varphi}_2)(\zeta)
&:=&\int_{0}^{\zeta}\;\imi{\varphi}_1\!(\zeta\!-\!\zeta_2)\imi{\varphi}_2\!(\zeta_2)\,d\zeta_2
\hspace{7.ex}
\end{array}\]
\[\begin{array}{ccccccccc}
\mi{major\; convolution}&&&&
\\[1.5 ex]
\mi{upper\;}\;\overline{\ast}\quad&
(\sma{\varphi}_1 \overline{\ast} \sma{\varphi}_2)(\zeta)
&:=&\frac{1}{2\pi i}\int_{c-i\,d}^{c+i\,d}\;\sma{\varphi}_1\!(\zeta\!-\!\zeta_2)\;
d\!\sma{\varphi}_2\!(\zeta_2) \,
\\[1.5 ex]
\mi{lower\; }\;\underline{\ast}\quad&
(\ima{\varphi}_1 \underline{\ast} \ima{\varphi}_2)(\zeta)
&:=&\frac{1}{2\pi i}\int_{c-i\,d}^{c+i\,d}\;
\ima{\varphi}_1\!(\zeta\!-\!\zeta_2)\ima{\varphi}_2\!(\zeta_2)\,d\zeta_2 \; 
\end{array}\]
%%%%%%%%%%%%%%%%%%%%%%%%%%%%%%%%%%%%%%%%%%%%%%%%%%%%%%%%%%%%%%%%%%%%%%%%%
{\bf The upper/lower Borel-Laplace transforms.}
\\
For simplicity, let us fix the polarisation $\theta=0$ and drop the index $\theta$
in the geometric model $\varphi_\theta(z)$. We get the familiar formulas,
reproduced here just for definiteness\,:
\[\begin{array}{ccccccccc}
\mi{multiplicative}&&&&\mi{convolutive} &&
\\[1.ex]
& \mi{upper}& \smia{\varphi}\!(\zeta)&\!=\!& \{\smi{\varphi}\!(\zeta),\sma{\varphi}\!(\zeta) \}
&\quad&
\\
&  \overline{\caB} \nearrow  \swarrow \overline{\caL}\quad & & &
\\
\varphi(z)& &&&\downarrow\partial_\zeta\;\;,\;\;\downarrow\!-\partial_\zeta &
\\
&  \underline{\caB} \searrow  \nwarrow \underline{\caL}\quad & & &
\\
& \mi{lower}& \imia{\varphi}\!(\zeta)&\!=\!&\{\imi{\varphi}\!(\zeta),\ima{\varphi}\!(\zeta) \}
&\quad&
\end{array}\]

\[\begin{array}{ccccccccc}
\mi{minors}&\mi{Borel}\quad&& \mi{Laplace}
\\[1.5 ex]
\mi{upper}\quad&\smi{\varphi}(\zeta)=\frac{1}{2\pi i}\int_{c-i\infty}^{c+i\infty}
e^{\zeta z}\varphi(z)\frac{dz}{z} 
&&
\varphi(z)=\int_{0}^{+\infty}
e^{-z\zeta}\,d\!\smi{\varphi}\!(\zeta)
\\[1.5 ex]
\mi{lower}\quad&\imi{\varphi}(\zeta)=\frac{1}{2\pi i}\int_{c-i\infty}^{c+i\infty}
e^{\zeta z}\varphi(z){dz} 
&&
\varphi(z)=\int_{0}^{+\infty}
e^{-z\zeta}\,\imi{\varphi}\!(\zeta)\,d\zeta
\end{array}\]
\[\begin{array}{ccccccccc}
\mi{majors}&\mi{Laplace}\quad&& \mi{Borel}
\\[1.5 ex]
\mi{upper}\quad&\sma{\varphi}(\zeta)=\int_{c}^{+\infty}
e^{-\zeta z}\varphi(z)\frac{dz}{z} 
&&
\varphi(z)=-\frac{1}{2\pi i}\int_{c-i\infty}^{c+i\infty}
e^{z\zeta}\,d\!\sma{\varphi}\!(\zeta)
\\[1.5 ex]
\mi{lower}\quad&\ima{\varphi}(\zeta)=\int_{c}^{+\infty}
e^{-\zeta z}\varphi(z){dz} 
&&
\varphi(z)=\frac{1}{2\pi i}\int_{c-i\infty}^{c+i\infty}
e^{z\zeta}\,\ima{\varphi}\!(\zeta)\,d\zeta
\end{array}\]
%%%%%%%%%%%%%%%%%%%%%%%%%%%%%%%%%%%%%%%%%%%%%%%%%%%%%%%%%%%%%%%%%%%%%%%%%
{\bf Monomials in all four models.}
\\
The following table covers not only the monomials $J_\sigma(z):=z^{-\sigma}$
but also\footnote{ after $\sigma$-differentiation} the whole range of
binomials $J_{\sigma,n}(z):=z^{-\sigma}\log^n(z)$
with $\sigma\in \doC, n \in \doN $.
\[\begin{array}{ccccccccc}
&&&\mi{minor} &&\mi{major} &&
\\[1.ex]
\sigma \not\in \doZ &\mi{upper}
 & &
\smi{J}_\sigma\!(\zeta)={\zeta^\sigma}/{\Gamma(1+\sigma)}
&,&
\sma{J}_\sigma\!(\zeta)=\zeta^\sigma\Gamma(-\sigma)
 &\quad&
\\
&  \nearrow  \swarrow\quad & & &
\\
J_\sigma(z)=z^{-\sigma}& &&\downarrow\partial_\zeta\;\;& &\;\;\downarrow\!-\partial_\zeta &
\\
&  \searrow  \nwarrow \quad & & &
\\
&\mi{lower}
 &
&\imi{J}_\sigma\!(\zeta)=\zeta^{\sigma-1}/\Gamma(\sigma)
&,&
\ima{J}_\sigma\!(\zeta)=\zeta^{\sigma-1}\Gamma(1-\sigma)
 &\quad&
\\[1. ex]
\dotfill & \dotfill & \dotfill & \dotfill & \dotfill & \dotfill & 
\\[1. ex]
s\in \doN^+ &\mi{upper}
 & &
\smi{J}_s\!(\zeta)=\frac{1}{s!}\;{\zeta^s}
&,&
\sma{J}_s\!(\zeta)=(-\zeta)^s\log(\frac{1}{\zeta})
 &\quad&
\\
&  \nearrow  \swarrow\quad & & &
\\
J_s(z)=z^{-s}& &&\downarrow\partial_\zeta\;\;& &\;\;\downarrow\!-\partial_\zeta &
\\
&  \searrow  \nwarrow \quad & & &
\\
&\mi{lower}
 &
&\imi{J}_s\!(\zeta)=\frac{1}{(s-1)!}\,\zeta^{s-1}
&,&
\ima{J}_s\!(\zeta)= (-\zeta)^{s-1}\log(\frac{1}{\zeta})
 &\quad&
\\[1. ex]
\dotfill & \dotfill & \dotfill & \dotfill & \dotfill & \dotfill & 
\\[1. ex]
&\mi{upper}
 & &
\smi{J}_0\!(\zeta)= 1
&,&
\sma{J}_s\!(\zeta)=\log(\frac{1}{\zeta})
 &\quad&
\\
&  \nearrow  \swarrow\quad & & &
\\
J_0(z)=1 & &&\downarrow\partial_\zeta\;\;& &\;\;\downarrow\!-\partial_\zeta &
\\
&  \searrow  \nwarrow \quad & & &
\\
&\mi{lower}
 &
&\imi{J}_0\!(\zeta)= 0
&,&
\ima{J}_s\!(\zeta)= \frac{1}{\zeta}
 &\quad&
\\[1. ex]
\dotfill & \dotfill & \dotfill & \dotfill & \dotfill & \dotfill & 
\\[1. ex]
s\in \doN^+ &\mi{upper}
 & &
\smi{J}_{\!-s}\!(\zeta)=0
&,&
\sma{J}_{\!-s}\!(\zeta)=(s-1)!\;\zeta^{-s}
 &\quad&
\\
&  \nearrow  \swarrow\quad & & &
\\
J_{\!-s}(z)=z^{s}& &&\downarrow\partial_\zeta\;\;& &\;\;\downarrow\!-\partial_\zeta &
\\
&  \searrow  \nwarrow \quad & & &
\\
&\mi{lower}
 &
&\imi{J}_{\!-s}\!(\zeta)=0
&,&
\ima{J}_{\!-s}\!(\zeta)= s!\;\zeta^{-s-1}
 &\quad&
\end{array}\]
%%%%%%%%%%%%%%%%%%%%%%%%%%%%%%%%%%%%%%%%%%%%%%%%%%%%%%%%%%%%%%%%%%%%%%%%%
{\bf The pros and cons of the upper/lower choices.}
\\

\noindent
{\it Advantages of the lower choice\,:}
\\

(i) $\{\underline{\caB},\underline{\caL}, \underline{\ast}\} $ 
are more usual/natural choices than
 $\{\overline{\caB},\overline{\caL}, \overline{\ast}\} $ 
\\

(ii) the operators $(\partial_z+\omega)^{-1}$ and 
$(e^{\omega\partial_z}-1)^{-1}$
which constantly occur  in the theory of {\it singular}
differential or difference equations and are ultimately responsible
for the frequent occurrence, in this theory, of both divergence and resurgence,
turn into {\it minor} multiplication by
$ (-\zeta+\omega)^{-1}$ 
or 
$(e^{-\omega\zeta}-1)^{-1}$
in the $\zeta$-plane\footnote{ or {\it major} multiplication by
$ (\zeta+\omega)^{-1}$ 
or 
$(e^{\omega\zeta}-1)^{-1}$
}, whereas with the upper choice we would be saddled with the more unwieldy operators
$ \partial_\zeta^{-1}(-\zeta+\omega)^{-1}\partial_\zeta$ 
and $\partial_\zeta^{-1}(e^{-\omega\zeta}-1)^{-1}\partial_\zeta$.
\\

\noindent
{\it Advantages of the upper choice\,:}
\\

(i) the `monomial' formulas for $J_\sigma$ ({\it supra}) assume a smoother shape, with
the simple sign change $-\sigma \mapsto \sigma$ instead of
$-\sigma \mapsto \sigma-1$.
\\

(ii) upper convolution $\overline{\ast}$ and pointwise multiplication 
(both in the $\zeta$-plane) have the same unit element, namely
$\smi{\varphi}_0(\zeta)\equiv 1 $, which is extremely
useful when studying  {\it dimorphy phenomena}\footnote{ 
i.e. the simultaneous stability of certain function rings
unter {\it two} unrelated products, like pointwise {\it multiplication}
and some form or other of {\it convolution}.
}, e.g. the dimorphy of poly- or hyperlogarithms.
\\
  
\noindent
{\it In the present investigation, we shall resort to both choices, because\,:}
\\

(i) the lower choice leads to simpler formulas when deriving a function's singularities
from its Taylor coefficient asymptotics (see \S2.3 {\it infra}).
\\

(ii) the upper choice is the one naturally favoured by the functional transforms
({\it nir/mir} and {\it nur/mur}) that lead to the {\it inner}
and {\it outer} generators of SP-series (see \S4.4-5 and \S5.3-5 {\it infra}).
%%%%%%%%%%%%%%%%%%%%%%%%%%%%%%%%%%%%%%%%%%%%%%%%%%%%%%%%%%%%%%%%%%%%%%%%%%%%%%%%%%%%%%%%%%%%
%%%%%%%%%%%%%%%%%%%%%%%%%%%%%%%%%%%%%%%%%%%%%%%%%%%%%%%%%%%%%%%%%%%%%%%%%%%%%%%%%%%%%%%%%%%%
\subsection{Alien derivations as a tool for Riemann surface description.}
Resurgent functions are acted upon by a huge range of exotic derivations,
the so-called {\it alien derivations} $\Delta_\omega$, with indices
$\omega$ ranging through the whole of $\doC_\bu:=\widetilde{\doC-\{0\}}$.
In other words, $\arg(\omega)$ is defined {\it exactly} rather than $\mi{mod}\; 2\pi $.
Together, these $\Delta_\omega$ generate a {\it free Lie algebra} on
$\doC_\bu$. Alien derivations, by pull-back, act on all three models.
There being no scope for confusion, the same symbols $\Delta_\omega$
can be used in each model. {\it Alien derivations}, however, are linear operators
which quantitatively measure the {\it singularities of minors} in the
$\zeta$-plane. To interpret or calculate {\it alien derivatives},
we must therefore go to (either of) the {\it convolutive models}, which in that sense enjoy
an undoubted primacy. However, for notational ease, it is often
convenient to write down {\it resurgence equations}\footnote{
i.e. any relation 
$ E(\varphi,\Delta_{\omega_1}\varphi,\dots,\Delta_{\omega_n}\varphi) =0 $,
linear or not, between a resurgent function $\varphi$ and one or several
of its alien derivatives.
} in the multplicative models (formal or geometric), the product there
being the more familiar multiplication.

For simplicity, in all the following definitions/identities the
indices $\omega$ are assumed to be on $\doR^+ \subset \doC_\bu$\,.
Adaptation to the general case is immediate. 
\\

\noindent
%%%%%%%%%%%%%%%%%%%%%%%%%%%%%%%%%%%%%%%%%%%%%%%%%%%%%%%%%%%%%%%%%%%%%%%%%
{\bf Definition of the operators $\Delta_\omega$ and $\Delta^\pm_\omega$.}
\[\begin{array}{cccccccc}
& \mi{multiplicative}\quad\quad & \mi{convolutive}&&\mi{convolutive}
\\[1.ex]
&&\smia{\varphi}:=\{\smi{\varphi},\sma{\varphi} \}
&\mapsto & 
\smia{\varphi}_\omega:=\{\smi{\varphi}_\omega,\sma{\varphi}_\omega \}
\\[-0.5 ex]
\Delta_\omega\,: & \varphi\;\mapsto\;  \varphi_\omega \quad\quad
\\[-0.5 ex]  
&&\imia{\varphi}:=\{\imi{\varphi},\ima{\varphi} \}
&\mapsto & 
\imia{\varphi}_\omega:=\{\imi{\varphi}_\omega,\ima{\varphi}_\omega \}
\end{array}\]
\[\begin{array}{cclrlccc}
\ima{\varphi}_{\omega}(\zeta)&:=&\sum_{\epsilon_1,\dots.\epsilon_{r-1}}&
\;\;\;\lambda_{\epsilon_1,\dots,\epsilon_{r-1}}&
\imi{\varphi}^{({\epsilon_1 \atop \omega_{1} }
{,\dots , \atop  ,\dots , }{\epsilon_{r-1} \atop \omega_{r-1} })}\!\!(\omega-\zeta)
\\[0.9 ex] 
\imi{\varphi}_{\omega}(\zeta)&:=&\sum_{\epsilon_1,\dots,\epsilon_{r}}&
\frac{\epsilon_r}{2\pi i}\; \lambda_{\epsilon_1,\dots,\epsilon_{r-1}}&
\imi{\varphi}^{({\epsilon_1 \atop \omega_{1} }
{,\dots , \atop ,\dots , }{\epsilon_r\;\; \atop \omega_r\;\; })}(\omega+\zeta)
\\[1.9 ex] 
\ima{\varphi}_{\omega}(\zeta)&:=&\sum_{\epsilon_0,\dots,\epsilon_{r-1}}&
\frac{\epsilon_0}{2\pi i}\; \lambda_{\epsilon_1,\dots,\epsilon_{r-1}}&
\ima{\varphi}^{({\epsilon_0 \atop \underline{\omega}_{0} }
{,\dots , \atop  ,\dots , }{\epsilon_{r-1} \atop \underline{\omega}_{r-1}
})}\!\!(\underline{\omega}+\zeta)
\\[0.9 ex] 
\imi{\varphi}_{\omega}(\zeta)&:=&\sum_{\epsilon_0,\dots,\epsilon_{r}}&
\frac{\epsilon_0\,\epsilon_r}{2\pi i} \;\lambda_{\epsilon_1,\dots,\epsilon_{r-1}}&
\ima{\varphi}^{({\epsilon_0 \atop \underline{\omega}_{0} }
{,\dots , \atop  ,\dots , }{\epsilon_r\;\; \atop \underline{\omega}_r\;\; })}
(\underline{\omega}-\zeta)
\end{array}\]
\begin{equation*}
0=:\omega_0 < \omega_1 <\omega_2 <\dots <\omega_{r-2}<\omega_{r-1}<\omega_{r}:=\omega
\quad\quad\quad  \Big(\;\underline{\omega}_i:=-\omega_i\; ,\; \forall i\;\Big)
\end{equation*}
The above relations should first be interpreted locally, i.e. for $\zeta/\omega <\!<1$,
and then extended globally by analytic continuation in $\zeta$. Here
$ \imi{\varphi}^{({\mg{\epsilon} \atop \mg{\omega}})}$
etc
denotes the branch corresponding to the left or right circumvention of
each intervening singularity $\omega_i$
if $\epsilon_i$ is $+$ or $-$, and to a branch weightage $\lambda$ that doesn't depend
on the increments $\omega_i$\,:
\begin{equation*}
\lambda_{\epsilon_1,\dots,\epsilon_{r-1}}:=\frac{p!\,q!}{r!}\quad \mi{with}
\quad p:=\sum_{\epsilon_i=+} 1
\;,\;
 q:=\sum_{\epsilon_i=-} 1
\quad\quad\quad 
\Big( \sum_{\epsilon_i} \lambda_{\epsilon_1,\dots,\epsilon_{r-1}}\equiv 1\Big)
\end{equation*}
The lateral operators $\Delta^\epsilon_\omega\,:\; \varphi\mapsto \varphi_{\omega_\epsilon}$ 
(with index
$\epsilon=\pm)$  are defined by the same formulas as above, but with
weights $\lambda_{\epsilon_1,...,\epsilon_{r-1}}$ replaced by the
much more elementary
$2\pi i\lambda^\epsilon_{\epsilon_1,...,\epsilon_{r-1}}$\,: 
\begin{eqnarray*}
\lambda^\epsilon_{\epsilon_1,...,\epsilon_{r-1}}&:=& 1\;\; \mi{if}\;\;
\epsilon_1=\epsilon_2=\dots=\epsilon_{r-1}=\epsilon\; \in {\pm 1}
\\ 
&:=& 0\;\;\; otherwise
\end{eqnarray*}
Thus, the minor-to-major and minor-to-minor formulas read\,:
\[\begin{array}{cclrlccc}
\ima{\varphi}_{\omega_\epsilon}(\zeta)&:=&&{2\pi i}
\;&
\imi{\varphi}^{({\epsilon \atop \omega_{1} }
{,\dots , \atop  ,\dots , }{\epsilon \atop \omega_{r-1} })}\!\!(\omega-\zeta)
\\[0.9 ex] 
\imi{\varphi}_{\omega_\epsilon}(\zeta)&:=&\sum_{\epsilon_{r}}&
{\epsilon_r}\; &
\imi{\varphi}^{({\epsilon \atop \omega_{1} }
{,\dots , \atop ,\dots , }
{\epsilon\;\; \atop \omega_{r-1}\;\; }
{, \atop  , }
{\epsilon_r\;\; \atop \omega_{r}\;\; }
)}
(\omega+\zeta)
\end{array}\]
This settles the action of alien operators in the {\it lower} convolutive model.
Their action in the {\it upper} model is exactly the same. Their action in the multiplicative models
is defined indirectly, by pull-back from the convolutive models (with the same notation  $\Delta_{\omega} $ holding for all models).
\\ 
\noindent
%%%%%%%%%%%%%%%%%%%%%%%%%%%%%%%%%%%%%%%%%%%%%%%%%%%%%%%%%%%%%%%%%%%%%%%%%
{\bf The operators $\Delta_\omega$ are derivations but the simpler  $\Delta^\pm_\omega$ are not.}
\\
Indeed, for any two test functions $\varphi_1, \varphi_2$ the identities hold\footnote{ 
For simplicity, we write the following identities in the {\it multiplicative} models.
When transposing them to the $\zeta$-plane, where they make more direct sense,
multiplication must of course be replaced by convolution.
}
\,:
\begin{eqnarray*}
\Delta_\omega ({\varphi_1}\,{\varphi_2})
&\equiv& 
(\Delta_\omega {\varphi_1})\;{\varphi_2}
+{\varphi_1}\;(\Delta_\omega {\varphi_2})
\\
\Delta^\pm_\omega ({\varphi_1}\,{\varphi_2})
&\equiv& 
(\Delta^\pm_\omega {\varphi_1})\;{\varphi_2}
+{\varphi_1}\;(\Delta^\pm_\omega {\varphi_2})
+\sum_{\omega_1+\omega_2=\omega  }
^{\frac{\omega_1}{\omega}>0,\frac{\omega_2}{\omega}>0 }
(\Delta^\pm_{\omega_1}\varphi_1)(\Delta^\pm_{\omega_2}\varphi_2)
\end{eqnarray*}
%%%%%%%%%%%%%%%%%%%%%%%%%%%%%%%%%%%%%%%%%%%%%%%%%%%%%%%%%%%%%%%%%%%%%%%%%
{\bf Lateral and median singularities.}
\\
The lateral and median operators are related by the following identities\,:
\begin{eqnarray}\label{a5}
1+\sum_{\omega>0} t^\omega \,\Delta_\omega^\pm
&=&
\exp\Big(\pm 2\pi i \sum_{\omega>0} t^\omega \,\Delta_\omega \Big)
\\ \label{a6}
2\pi i \sum_{\omega>0} t^\omega \,\Delta_\omega
&=&
\pm\log\Big(1+\sum_{\omega>0} t^\omega \,\Delta_\omega^\pm\Big)
\end{eqnarray}
Interpretation\,: we first expand {\it exp} and {\it log} the usual way, 
then equate the contributions of each power $t^\omega$ from the left- and
right-hand sides. Although the above formulas express each $\Delta_\omega^\pm$ as 
an infinite sum of (finite) $\Delta_\omega$-products, and {\it vice versa}, when applied
to any given test function $\varphi$ the infinite sums actually reduce
to a finite number of {\it non-vanishing} summands\footnote{due to the minors
having only isolated singularities.}.
\\

\noindent
%%%%%%%%%%%%%%%%%%%%%%%%%%%%%%%%%%%%%%%%%%%%%%%%%%%%%%%%%%%%%%%%%%%%%%%%%
{\bf Compact description of Riemann surfaces.}
\\
Knowing {\it all} the alien derivatives (of first and higher orders)
of a minor $\imi{\varphi}(\zeta)$ or  $\smi{\varphi}(\zeta)$ 
enables one to piece together that minor's behaviour on its {\it entire
Riemann surface} $\caR$ from the behaviour of its various alien derivatives
on their sole {\it holomorphy stars}, by means of the general formula\,:
\begin{equation}\label{a7}
\imi{\varphi}(\zeta_{_{\Gamma}})\equiv\, \imi{\varphi}(\zeta)
+\sum_{r\geq 1}\sum_{\omega_i\in\doC_\bu}
(2\pi i)^r
H^{\omega_1,\dots,\omega_r}_\Gamma
t^{\omega_1+...\omega_r}\Delta_{\omega_r}\dots\Delta_{\omega_1} \imi{\varphi}(\zeta)
\quad\quad
\end{equation}
Here, $\zeta_{_{\Gamma}}$
denotes any chosen point on $\caR$, reached from $0$ by following
a broken line $\Gamma$
in the $\zeta$-plane. Both sums $\sum_{r}$ and $\sum_{\omega_i}$
are finite\footnote{That is to say, when applied to any given resurgent function, they carry only finitely many non-vanishing terms.}. The coefficients $H^\bu_{\Gamma}$ are in $\doZ $.
Unlike in (\ref{a5})(\ref{a6}), $t^\omega$ in (\ref{a7}) should no longer
be viewed as the symbolic power of a free variable $t$, but as an shift
operator acting on functions of $\zeta$ and changing $\zeta$ into $\zeta+\omega$. \footnote{with $\zeta$ close to 0 and suitably positioned, to ensure that $\zeta+\omega$ be in the holomorphy star of the test function.}

To sum up\,:
\\
(i) alien derivations `uniformise' everything.
\\
(ii) a full knowledge of a minor's alien derivatives
(given for example by a {\it complete} system of resurgence equations)
implies a full knowlege of that minor's Riemann surface.
\\

\noindent
%%%%%%%%%%%%%%%%%%%%%%%%%%%%%%%%%%%%%%%%%%%%%%%%%%%%%%%%%%%%%%%%%%%%%%%%%
{\bf Strong versus weak resurgence.}
\\ 
``Proper" resurgence equations are relations of the form\,:
\begin{equation}\label{a8}
E(\varphi,\Delta_\omega\varphi)\equiv 0
\quad \quad \mi{or} \quad \quad
E(\varphi,\Delta_{\omega_1}\varphi,\dots,\Delta_{\omega_n}\varphi )\equiv 0
\end{equation}
with expressions $E$ that are typically non-linear (at least in $\varphi$)
and that may involve arbitrary scalar- or function-valued coefficients. Such equations
express unexpected {\it self-reproduction} properties \--- that is to say,
non-trivial relations between the minor (as a germ at $\zeta=0$) and 
its various singularities. Moreover, when the resurgent function $\varphi$,
in the multiplicative model, happens to be the formal solution
of some equation or system $S(\varphi)=0$ (think for example of a singular
differential, or difference, or functional, equation), the resurgence
of $\varphi$ as well as the exact shape of its resurgence equations (\ref{a8}),
can usually be derived almost without analysis, merely by letting each $\Delta_\omega$
act on $S(\varphi)$ in accordance with certain formal rules. Put another way\,:
we can deduce deep   {\it analytic} facts
from purely {\it formal-algebraic} manipulations. What we have here
 is {\it full-fledged
resurgence} \--- resurgence at its best and most useful.

But two types of situations may arise which lead to {\it watered-down}
forms of resurgence. 

One is the case when, due to severe constraints built into the
resurgence-generating problem, the coefficients inside $E$ 
are no longer free to vary continuously, but must assume
discrete, usually entire values\,: we then speak of {\it rigid resurgence}.

Another is the case when the expressions $E$ are linear or affine
functions of their arguments $\varphi$ and $\Delta_{\omega_i}\varphi$.
The self-reproduction aspect, to which resurgence owes its name, then
completely disappears, and makes way for a simple {\it exchange}
or {\it `ping-pong'} between singularities (in the linear case)
with possible `annihilations' (in the affine case).

Both restrictions entail a severe impoverishment of the resurgence
phenomenon. As it happens, and as we propose to show in this paper,
SP-series combine these two restrictions\,: they lead to fairly degenerate 
resurgence
patterns that are both {\it rigid} and {\it affine}. 
Furthermore, as a rule, SP-series verify no useful
equation or system $S(\varphi)=0$ that might give us a cue as to their
resurgence properties.
In cases
such as this, the resurgence apparatus (alien derivations etc)
ceases to be a vehicle for {\it proving things} and
retains only its (non-negligible\,!) {\it notational value}
(as a device for describing Riemann surfaces etc) while the onus
of proving the hard analytic facts falls on altogether different
tools, like {\it Taylor coefficient asymptotics}\footnote{see \S2.3}
and the {\it nir/mir}-transforms.\footnote{see \S4.4,\S4.5.}
\\

\noindent 
%%%%%%%%%%%%%%%%%%%%%%%%%%%%%%%%%%%%%%%%%%%%%%%%%%%%%%%%%%%%%%%%%%%%%%%%%
{\bf The pros and cons of the $2\pi i$ factor.}
\\
On balance, we gain more than we lose by inserting the $2\pi i $ 
factor into the definition (\ref{a7})
of the alien derivations. True, by removing it there
we would also eliminate it from
the identities relating minors to majors (see \S2.1), but the factor would sneak back into the  $J_\sigma $-identities {\it supra}, 
thus spoiling the whole set of `monomial' formulas. Worse still,
real-indexed derivations $\Delta_\omega$ acting on real-analytic
derivands $\varphi$ would no longer produce real-analytic derivatives
$\Delta_\omega \varphi$ \,\--- which would be particularly damaging
in ``all-real'' settings, e.g. when dealing with 
chirality 1 knots like $4_1$ (see \S9 {\it infra}).
%%%%%%%%%%%%%%%%%%%%%%%%%%%%%%%%%%%%%%%%%%%%%%%%%%%%%%%%%%%%%%%%%%%%%%%%%%%%%%%%%%%%%%%%%%%%
%%%%%%%%%%%%%%%%%%%%%%%%%%%%%%%%%%%%%%%%%%%%%%%%%%%%%%%%%%%%%%%%%%%%%%%%%%%%%%%%%%%%%%%%%%%%

\subsection{Retrieving the resurgence of a series from the resurgence of its Taylor
coefficients.}
SP-series are one of those rare instances where there is no shortcut for
calculating the singularities\,: we have no option but to deduce them from
a close examination of the asymptotics of the Taylor coefficients.\footnote{
The present section is based on a private communication (1992) by J.E. to Prof. G.K. Immink. An independent,
more detailed treatment was later given by O. Costin in [C2].
}
 
The better to respect the symmetry between our series $\varphi$ and its Taylor
coefficients $J$, we shall view them both as resurgent functions of
the variables $z$ resp. $n$ in the multiplicative models and
$\zeta$ resp. $\nu$ in the (lower) convolutive models. The aim then is to understand
the correspondence between the triplets\,:
\begin{equation}\label{a9}
\{ \tilde{\varphi}(z), \imia{\varphi}\!(\zeta),\varphi(z) \}
\longleftrightarrow 
\{ \tilde{J}(n), \imia{J}\!(\nu),J(n) \}
\end{equation}
and the alien derivatives attached to them.
\\
%%%%%%%%%%%%%%%%%%%%%%%%%%%%%%%%%%%%%%%%%%%%%%%%%%%%%%%%%%%%%%%%%%%%%%%%%%%%%%%%%%%%%%%%%%%%
%%%%%%%%%%%%%%%%%%%%%%%%%%%%%%%%%%%%%%%%%%%%%%%%%%%%%%%%%%%%%%%%%%%%%%%%%%%%%%%%%%%%%%%%%%%%

\noindent
{\bf Retrieving closest singularities.}\\
Let us start with the simplest case, when $\imi{\varphi}$
has a single singularity on the boundary of its disk of convergence, say at $\zeta_0$.
We can of course assume $\zeta_0$ to be real positive. 
\begin{equation*}
\tilde{\varphi}(z)=\sum_{0\leq n}(n\!+\!1)!\;J(n)\,z^{-n-1}
\,(\mi{dv^t})\;\stackrel{\underline{\caB}}{\mapsto}\;
\imi{\varphi}(\zeta)=\sum_{0\leq n}\,J(n)\,\zeta^{n}\;
({\it cv^t\,on}\,|\zeta|<\zeta_0)
\end{equation*} 
 In order to
deduce the closest singularity of  $\imi{\varphi}$ from the
closest singularity of $\imi{J}$,
 we first express $J(n)$ as a Cauchy integral on a circle $|\zeta|=|\zeta_0|-\epsilon$.
We then deform that circle to a contour $\Gamma$ which coincides with
the larger circle $|\zeta|=|\zeta_0|+\epsilon$ except for a slit along
the interval $[\zeta_0,\zeta_0+\epsilon]$ to avoid crossing the singularity at $\zeta_0$.
Lastly, we transform $\Gamma$ into $\Gamma_\ast$ 
(resp. $\underline{\Gamma_\ast}=-\Gamma_\ast $)
under the change 
$\zeta= \zeta_0 e^\nu$ .
\begin{eqnarray}  \label{a10}
J(n)&=&\frac{1}{2\pi i}\;\oint \imi{\varphi}(\zeta)\,\zeta^{-n-1}\,d\zeta
\\   \label{a11}
&=&\frac{1}{2\pi i}\;\int_\Gamma \imi{\varphi}(\zeta)\,\zeta^{-n-1}\,d\zeta
+o(\zeta_0^{-n})
\quad\quad\quad (\mi{contour\, deformation})
\\ \nonumber
&=&\frac{e^{-n\nu_{0}}}{2\pi i}\;\int_{\Gamma_\ast} 
\imi{\varphi}(\zeta_0 e^\nu)\, e^{-n\nu}\,d\nu
+o(e^{-n\nu_0})
\quad\; (\mi{setting\;}\zeta:=\zeta_0\,e^\nu=e^{\nu_0+\nu})
\\    \label{a12}
&=&\frac{e^{-n\nu_{0}}}{2\pi i}\;\int_{\Gamma_\ast} 
\ima{\varphi}_{\zeta_0}(\zeta_0 -\zeta_0 e^\nu)\, e^{-n\nu}\,d\nu
+o(e^{-n\nu_0})
\quad\quad\; (\mi{always})
\\    \label{a13}
&=&\frac{e^{-n\nu_{0}}}{2\pi i}\;\int_{\underline{\Gamma}_\ast} 
\ima{\varphi}_{\zeta_0}(\zeta_0 -\zeta_0 e^{-\nu} )\, e^{n\nu}\,d\nu
+o(e^{-n\nu_0})
\quad\quad\; (\mi{always})
\\    \label{a14}
&=&e^{-n\nu_0}\int_{0}^{c} 
\imi{\varphi}_{\zeta_0}(\zeta_0 e^\nu-\zeta_0)\, e^{-n\nu}\,d\nu
+o(e^{-n\nu_0})
\quad\; (\mi{if\,}\imia{\varphi}_{\zeta_0}\mi{\,integrable})
\end{eqnarray}
Therefore
\begin{equation}   \label{a15}
A(n)\equiv e^{-n\nu_0}\;A_{\nu_0}(n)+o(e^{-n \nu_0})
\quad\quad\quad (\mi{with\;}\nu_0:=\log(\zeta_0) )
\end{equation}
where $J_{\nu_0} $ denotes the (lower, and if need be, truncated) Borel transform of a
resurgent function
$\imia{J}_{\nu_0} $ linked to 
$\imia{\varphi}_{\zeta_0}:=\Delta_{\zeta_0}\imi{\varphi}$
by\,:
\[\begin{array}{llllllll}
\imia{J}_{\nu_0}\!(\nu)&\!\!=\!\!&\{\imi{J}_{\nu_0}\!(\nu),\ima{J}_{\nu_0}\!(\nu) \}
&\quad&
\imia{\varphi}_{\zeta_0}\!(\nu)&\!\!=\!\!&
\{\imi{\varphi}_{\zeta_0}\!(\nu),\ima{\varphi}_{\zeta_0}\!(\zeta)
\} &
\\[2. ex]
\imi{J}_{\nu_0}\!(\nu)&\!\!=\!\!&
\imi{\varphi}_{\zeta_0}\!(\zeta_0\, e^\nu-\zeta_0\;)
 &\quad&
\imi{\varphi}_{\zeta_0}\!(\zeta)&\!\!=\!\!&
\imi{J}_{\nu_0}\!(\;\log(1+\frac{\zeta}{\zeta_0}))
 &\quad(\mi{minors})
\\[2. ex]
\ima{J}_{\nu_0}\!(\nu)&\!\!=\!\!&
\ima{\varphi}_{\zeta_0}\!(\zeta_0 -\zeta_0\, e^{-\nu})
 &\quad&
\ima{\varphi}_{\zeta_0}\!(\zeta)&\!\!=\!\!&
\ima{J}_{\nu_0}\!(\!-\!\log(1\!-\!\frac{\zeta}{\zeta_0}))
 &\quad(\mi{majors})
\end{array}\]
%%%%%%%%%%%%%%%%%%%%%%%%%%%%%%%%%%%%%%%%%%%%%%%%%%%%%%%%%%%%%%%%%%%%%%%%%
{\bf Retrieving distant singularities.}\\
The procedure actually extends to farther-lying singularities. In fact,
if $\imi{J}$ is endlessly continuable, so is  $\imi{\varphi}$,
and the former's resurgence pattern neatly translates into the latter's, under a set of linear but non-trivial formulas
Here, however, we shall only require knowledge of those singularities
of $\imi{\varphi}$ which lie on its (0-centered, closed) star of holomorphy.
All the other singularities will follow under repeated alien differentiation.
%\end{document} 
%%%%%%%%%%%%%%%%%%%%%%%%%%%%%%%%%%%%%%%%%%%%%%%%%%%%%%%%%%%%%%%%%%%%%%%%%%%%%%%%%%%%%%%%%%%%
%%%%%%%%%%%%%%%%%%%%%%%%%%%%%%%%%%%%%%%%%%%%%%%%%%%%%%%%%%%%%%%%%%%%%%%%%%%%%%%%%%%%%%%%%%%%
%%%%%%%%%%%%%%%%%%%%%%%%%%%%%%%%%%%%%%%%%%%%%%%%%%%%%%%%%%%%%%%%%%%%%%%%%%%%%%%%%%%%%%%%%%%%
%%%%%%%%%%%%%%%%%%%%%%%%%%%%%%%%%%%%%%%%%%%%%%%%%%%%%%%%%%%%%%%%%%%%%%%%%%%%%%%%%%%%%%%%%%%%
%%%%%%%%%%%%%%%%%%%%%%%%%%%%%%%%%%%%%%%%%%%%%%%%%%%%%%%%%%%%%%%%%%%%%%%%%%%%%%%%%%%%%%%%%%%%
%%%%%%%%%%%%%%%%%%%%%%%%%%%%%%%%%%%%%%%%%%%%%%%%%%%%%%%%%%%%%%%%%%%%%%%%%%%%%%%%%%%%%%%%%%%%
%%%%%%%%%%%%%%%%%%%%%%%%%%%%%%%%%%%%%%%%%%%%%%%%%%%%%%%%%%%%%%%%%%%%%%%%%%%%%%%%%%%%%%%%%%%%
%%%%%%%%%%%%%%%%%%%%%%%%%%%%%%%%%%%%%%%%%%%%%%%%%%%%%%%%%%%%%%%%%%%%%%%%%%%%%%%%%%%%%%%%%%%%
%%%%%%%%%%%%%%%%%%%%%%%%%%%%%%%%%%%%%%%%%%%%%%%%%%%%%%%%%%%%%%%%%%%%%%%%%%%%%%%%%%%%%%%%%%%%

%

%%%%%%%%%%%%%%%%%%%%%%%%%%%%%%%%%%%%%%%%%%%%%%%%%%%%%%%%%%
%% Some resurgence properties of knot-related functions.
%%%%%%%%%%%%%%%%%%%%%%%%%%%%%%%%%%%%%%%%%%%%%%%%%%%%%%%%%%

%\documentclass[12pt,a4paper]{article}\input{SP_commands}\begin{document}
%%%%%%%%%%%%%%%%%%%%%%%%%%%%%%%%%%%%%%%%%%%%%%%%%%%%%%%%%%%%%%%%%%%%%%%%%%%%%%%%%%%%%%%%%%%%
%%%%%%%%%%%%%%%%%%%%%%%%%%%%%%%%%%%%%%%%%%%%%%%%%%%%%%%%%%%%%%%%%%%%%%%%%%%%%%%%%%%%%%%%%%%%

\section{The ingress factor.}
We must first describe the asymptotics of the ``product" part (for $m=n$)
of our ``sum-product" coefficients. This involves a trifactorisation\,:
\begin{equation}   \label{a16}
\prod_{0\leq k \leq n}F(\frac{k}{n})
=:P_F(n)
\sim \tilde{I\!g}_F(n)\,e^{-\nu_\ast\,n}\,\tilde{E}\!g_F(n)
\quad\mi{with}\quad\nu_\ast=\!\int_0^1\!\!f(x)dx\;
\end{equation}
with
\\
(i) an {\it ingress} factor ${I\!g}_{\!F}$, purely local at $x=0$.
\\
(ii) an {\it exponential} factor $e^{-\nu_0\,n} $, global on $[0,1]$
\\
(iii) an {\it egress} factor ${E\!g}_{\!F}$, purely local at $x=1$.

The non-trivial factors (ingress/egress) may be divergent-resurgent
(hence the tilda) but, at least for holomorphic data $F$, they 
always remain fairly elementary. They often vanish (when $F$ is {\it even}
at 0 or 1) and, even when divergent, they can always be resummed in a canonical way.
Lastly, as already hinted, it will prove technically convenient
to factor out the first of these (ingress), thereby replacing
the original SP-series $j(\zeta)$ by its `cleansed'
and more regular version $j^\#(\zeta)$.
%%%%%%%%%%%%%%%%%%%%%%%%%%%%%%%%%%%%%%%%%%%%%%%%%%%%%%%%%%%%%%%%%%%%%%%%%%%%%%%%%%%%%%%%%%%%
%%%%%%%%%%%%%%%%%%%%%%%%%%%%%%%%%%%%%%%%%%%%%%%%%%%%%%%%%%%%%%%%%%%%%%%%%%%%%%%%%%%%%%%%%%%%
\subsection{Bernoulli numbers and polynomials.}
 For future use, let us collect a few formulas about 
 two convenient variants of 
the classical Bernoulli numbers $B_k$ and Bernoulli polynomials $B_k(t)$.
\\

\noindent
{\bf The Bernoulli numbers and polynomials.}
\begin{eqnarray}     \label{a17}
\quad
\frak{b}_k&:=&
\frac{\frak{b}^\ast_k(0)}{k!}
=
\frac{B_{k+1}(1)}{(k\!+\!1)!}
\hspace{19.ex}
 (k\in -1-\doN)
\\[1.ex]     \label{a18}
\frak{b}(\tau) &:=&\! \frac{e^{\tau}}{e^{\tau}-1}\,=\,\sum_{k\geq
-1}\frak{b}_k\,\tau^{k}
=\tau^{-1}
+\frac{1}{2}
+\frac{1}{12}\tau-\frac{1}{720}\tau^3-\dots
\\[0.1 ex]    \label{a19}
\frak{b}_k^\ast(\tau)&:=&
\frak{b}(\partial_\tau)\,\tau^k
\hspace{26. ex}(k\in \doC\,,\,k\not=-1)
\\[1.ex]    \label{a20}
\frak{b}^{\ast\ast}(\tau,\zeta) \!\!&:=&\!\!
\sum_{k\geq 0}\frak{b}_k^\ast(\tau)\,\frac{\zeta^k}{k!}
\;=\;
\frak{b}(\partial_\tau)\,e^{\tau\zeta}
\;=\;
\frac{e^{\tau\zeta}e^{\zeta}}{e^{\zeta}- 1}-\frac{1}{\zeta}
\end{eqnarray}
For $k\in \doN$, we have $\frak{b}_k=\frac{B_{k+1}}{(1\!+\!k)} $
for the scalars, and 
the  series $\frak{b}^{\ast}_k(\tau)$ 
essentially coincide with the Bernoulli
polynomials. For all other values of $k$, the  scalars $\frak{b}_k$
are no longer defined and
the  $\frak{b}^{\ast}_k(\tau)$   
become divergent  series
in decreasing powers of $\tau$.
\begin{eqnarray}     \label{a21}
\frak{b}_k^\ast(\tau)&:=&
\sum_{s=-1}^{k}\,\frak{b}_s\,\tau^{k-s}\frac{k!}{(k\!-\!s)!}\;\;
=\;\;
\frac{B_{k+1}(\tau\!+\!1)}{k\!+\!1} 
\quad(\mi{if}\;k\in\doN)
\\[1. ex]    \label{a22}
\frak{b}_k^\ast(\tau)&:=&
\sum_{s=-1}^{+\infty}\,\frak{b}_s\,
\tau^{k-s}\frac{\Gamma(k\!+\!1)}{\Gamma(k\!+\!1\!-\!s)}
\hspace{14.ex}(\mi{if}\;k\in\doC-\doZ)
\\[1.ex]    \label{a23}
\frak{b}^{\ast}(\tau) \!\!&:=&\!\!
\frac{\tau^{k+1}}{k\!+\!1}+
\sum_{k\geq 0}\,(-1)^s\frak{b}_s\,\tau^{k-s}\frac{(s-k)!}{k!}
\hspace{4.ex}(\mi{if}\;k\in -2-\doN)
\end{eqnarray}

\noindent
{\bf The Euler-Bernoulli numbers and polynomials.}
\begin{eqnarray}    \label{a24}
\quad
\beta_k&:=&
\frac{\beta^\ast_k(0)}{k!}
=
\frac{B_{k+1}(\frac{1}{2})}{(k\!+\!1)!}
\hspace{19.ex}
 (k\in -1-\doN)
\\[1.ex]    \label{a25}
\beta(\tau) &:=&\! \frac{1}{e^{\tau/2}-e^{-\tau/2}}\,=\,\sum_{k\geq
-1}\beta_k\,\tau^{-k}
=\tau^{-1}-\frac{1}{24}\tau+\frac{7}{5760}\tau^3-\dots
\\[0.1 ex]    \label{a26}
\beta_k^\ast(\tau)&:=&
\beta(\partial_\tau)\,\tau^k
\hspace{26. ex}(k\in \doC\,,\,k\not=-1)
\\[1.ex]    \label{a27}
\beta^{\ast\ast}(\tau,\zeta) \!\!&:=&\!\!
\sum_{k\geq 0}\beta_k^\ast(\tau)\,\frac{\zeta^k}{k!}
\;=\;
\beta(\partial_\tau)\,e^{\tau\zeta}
\;=\;
\frac{e^{\tau\zeta}}{e^{\zeta/2}- e^{-\zeta/2}}-\frac{1}{\zeta}
\end{eqnarray}
For $k\in \doN$, the $\beta^\ast_k(\tau)$ essentially coincide with the 
Euler-Bernoulli
polynomials. For all other values of $k$, they are divergent series
in decreasing powers of $\tau$.
\begin{eqnarray}    \label{a28}
\beta_k^\ast(\tau)&:=&
\sum_{s=-1}^{k}\,\beta_s\,\tau^{k-s}\frac{k!}{(k\!-\!s)!}\;\;
=\;\;
\frac{B_{k+1}(\tau\!+\!\frac{1}{2})}{k\!+\!1} 
\quad(\mi{if}\;k\in\doN)
\\[1. ex]    \label{a29}
\beta_k^\ast(\tau)&:=&
\sum_{s=-1}^{+\infty}\,\beta_s\,
\tau^{k-s}\frac{\Gamma(k\!+\!1)}{\Gamma(k\!+\!1\!-\!s)}
\hspace{14.ex}(\mi{if}\;k\in\doC-\doZ)
\\[1.ex]    \label{a30}
\beta^{\ast}(\tau) \!\!&:=&\!\!
\frac{\tau^{k+1}}{k\!+\!1}+
\sum_{k\geq 0}\,(-1)^s\beta_s\,\tau^{k-s}\frac{(s-k)!}{k!}
\hspace{4.ex}(\mi{if}\;k\in -2-\doN)
\end{eqnarray}
For all $k\in \doN$ we have the parity relations
$\beta^\ast_{2k}= 0\,,\,\beta^\ast_{k}(-\tau)\equiv(-1)^{k+1}\beta_k(\tau) $
\\

\noindent
{\bf The Euler-MacLaurin formula.}
\\ 
We shall make constant use of the basic identities 
($\forall s \in \doN $)\,:
\begin{equation*}
\!\!\sum_{1\leq k \leq m} k^s 
\equiv \frak{b}^\ast_s(m)-\frak{b}^\ast_s(0)
\equiv {\beta}^\ast_s(m\!+\!\frac{1}{2})-{\beta}^\ast_s(\frac{1}{2})
\equiv \frac{B_{s+1}(m\!+\!1)-B_{s+1}(1)}{s+1}
\end{equation*}
and of these variants of the Euler-MacLaurin formula\,:
\begin{eqnarray}     \label{a31}
\!\!\!\!\!\sum_{k}^{0 \leq \frac{k}{n} \leq \bar{x}} f(\frac{k}{n})
\!&\!\!\sim\!\!&\!
n\!\int_{0}^{\bar{x}}\!\!\!f(x)dx\!+\!\frac{f(0)}{2}\!+\!\frac{f(\bar{x})}{2}
+\sum_{1\leq s\;\mi{odd}}\,\frac{\frak{b}_{s}}{n^s}
\big( {f^{(s)}(\bar{x})}-{f^{(s)}(0)}\big)\;\;\;\;\;\quad
\\   \label{a32}
\!\!\!\!\!\sum_{k}^{0 \leq \frac{k}{n} \leq \bar{x}} f(\frac{k}{n})
\!&\!\!\sim\!\!&\!
n\!\int_{0}^{\bar{x}}\!\!\!f(x)dx\!+\!\frac{f(0)}{2}\!+\!\frac{f(\bar{x})}{2}
+\sum_{1\leq s\; odd}\,\frac{\frak{b}^\ast_{s}(0)}{n^s}
\big( \bar{f}_{s}-f_{s} \big)
\end{eqnarray} 
where $f_{s}$ and $\bar{f}_{s}$ denote the Taylor coefficients of $f$ at $0$ and $\bar{x}$.
%%%%%%%%%%%%%%%%%%%%%%%%%%%%%%%%%%%%%%%%%%%%%%%%%%%%%%%%%%%%%%%%%%%%%%%%%%%%%%%%%%%%%%%%%%%%
%%%%%%%%%%%%%%%%%%%%%%%%%%%%%%%%%%%%%%%%%%%%%%%%%%%%%%%%%%%%%%%%%%%%%%%%%%%%%%%%%%%%%%%%%%%%
\subsection{Resurgence of the Gamma function.}

\begin{lemma}[Exact asymptotics of the Gamma function]{\!\!.}\\
The functions $ \Theta,\theta $ defined on $\{\Re{(n)}>0\}\subset \doC$ by
\begin{equation}   \label{a33}
\Theta(n)\equiv e^{\theta(n)}:=(2\pi)^{-\frac{1}{2}}\Gamma(n+1)\,n^{-n-\frac{1}{2}}\,e^n
\quad\quad \quad ( \theta(n) \;\mi{real\; if}\; n\; \mi{real})
\end{equation}
possess resurgent-resummable asymptotic expansions as $\Re{(n)}\rightarrow +\infty$\,:
\begin{equation}   \label{a34}
\Theta(n)=1+\sum_{1\leq s}\Theta_s\,n^{-s}
\quad; \quad
\theta(n)=\sum_{0\leq s}\theta_{1+2s}\,n^{-1-2s}
\quad \mi{(odd\;powers )}\quad
\end{equation}
with explicit lower/upper Borel transforms\,:
\begin{eqnarray}   \label{a35}
\imi{\theta}(\nu)&=&-\frac{1}{\nu^2}+\frac{1}{2\,\nu}\,\frac{1}{\tanh({\nu}/{2})}
\\   \label{a36}
\smi{\theta}(\nu)&=&+\frac{1}{\nu}+\frac{1}{2}\int_0^{\nu/2}\frac{1}{\tanh(t)}\,\frac{dt}{t}
\\ \nonumber
&=&
\frac{1}{12}\nu
-\frac{1}{2160}\nu^{3}
+\frac{1}{151200}\nu^{5}
-\frac{1}{8467200}\nu^{7}
+\frac{1}{431101440}\nu^{9}
\dots
\end{eqnarray}
\end{lemma}
This immediately follows from $\Gamma$'s functional equation. We get successively\,:
\begin{eqnarray*}
\frac{\Theta(n+\frac{1}{2})}{\Theta(n-\frac{1}{2})}
&\!\!=\!\!&
e\;\Big(\frac{n-\frac{1}{2}}{n+\frac{1}{2}}\Big)^n\;
\\[1. ex]
\theta(n+\frac{1}{2})-\theta(n-\frac{1}{2})
&\!\!=\!\!&
1+n\log(n-\frac{1}{2})-n\log(n+\frac{1}{2})
\\[1. ex]
\partial_n\,\frac{1}{n}\big(\theta(n+\frac{1}{2})-\theta(n+\frac{1}{2})\big)
&\!\!=\!\!&
-\frac{1}{n^2}
+\frac{1}{n-\frac{1}{2}}-\frac{1}{n+\frac{1}{2}}
\\[1.ex]
-\nu\,\partial_{\nu}^{-1}\Big((e^{-\nu/2}-e^{\nu/2})\imi{\theta}(\nu) \Big)
&\!\!=\!\!&
-\nu+e^{\nu/2}-e^{-\nu/2}
\\[1.ex]
\imi{\theta}\!(\nu) 
&\!\!=\!\!&
(e^{\nu/2}-e^{-\nu/2})^{-1}\;
\partial_\nu\;\,\nu^{-1}\;
(e^{\nu/2}-e^{-\nu/2})
\quad\quad
\\[1.ex]
\imi{\theta}\!(\nu) 
&\!\!=\!\!&
-\frac{1}{\nu^2}+\frac{1}{2\,\nu}\,\frac{1}{\tanh({\nu}/{2})}
\end{eqnarray*}
Laplace summation along $\arg(\nu)=0$ yields the exact values $\theta(n)$
and $\Theta(n)$. The only non-vanishing alien derivatives are
(in multiplicative notation)\,:
\begin{eqnarray}   \label{a37}
\Delta_\omega\;\tilde{\theta}&=&\frac{1}{\omega}
\hspace{12.ex}\forall\omega \in 2\pi i\doZ^\ast
\\   \label{a38}
\Delta_\omega\;\tilde{\Theta}&=&\frac{1}{\omega}\;\tilde{\Theta}
\hspace{10.ex}\forall\omega \in 2\pi i\doZ^\ast
\end{eqnarray}
Using formula (\ref{a38}) and its iterates  for crossing the vertical axis in the $\nu$-plane,
we can evaluate the quotient of the regular resummations of $\hat{\Theta}(\nu)$ along $\arg(\nu)=0$
and $\arg(\nu)=\pm \pi $, and the result of course agrees with the complement
formula\footnote{For details, see [E3] pp 243-244.}\,: 
\begin{equation}   \label{a39}
\frac{1}{\Gamma(n)\,\Gamma(1\!-\!n)}=\frac{\sin{\pi\,n}}{\pi}
\hspace{10.ex}\forall n \in \doC
\end{equation}

%%%%%%%%%%%%%%%%%%%%%%%%%%%%%%%%%%%%%%%%%%%%%%%%%%%%%%%%%%%%%%%%%%%%%%%%%%%%%%%%%%%%%%%%%%%%
%%%%%%%%%%%%%%%%%%%%%%%%%%%%%%%%%%%%%%%%%%%%%%%%%%%%%%%%%%%%%%%%%%%%%%%%%%%%%%%%%%%%%%%%%%%%
\subsection{Monomial/binomial/exponential factors.}
In view of definition  (\ref{a16}) and formula (\ref{a31}), for a generic input $F:=e^{-f}$ with $F(0),F(1)\not= 0,\infty$ we
get the asymptotic expansions\,:
\begin{eqnarray}   \label{a40}
\tilde{I\!g}_{_{F}}(n)&=&
\exp\Big(-\frac{1}{2}f(0)+\sum_{1\leq s \mi{odd}} \frac{\frak{b}_s}{n^s}f^{(s)}(0)\Big)
\hspace{22.ex}
\\[1.ex]   \label{a41}
\tilde{E\!g}_{_{F}}(n)&=&
\exp\Big(-\frac{1}{2}f(1)-\sum_{1\leq s \mi{odd}} \frac{\frak{b}_s}{n^s}f^{(s)}(1)\Big)
\hspace{22.ex}
\end{eqnarray} 
and the important parity relation\,\footnote{In \S5 it will account for the
relation between the two outer generators which are always present in the generic case
(i.e. when $F(0)$ and $F(1)\not=0,\infty$).}:
\begin{equation}   \label{a42}
\{F^\perp(x)=1/F(1\!-\!x)\} \Longrightarrow 
\{
1=\tilde{I\!g}_{_{F}}(n)\tilde{E\!g}_{_{F^\perp}}(n)
 =\tilde{I\!g}_{_{F^\perp}}(n)\tilde{E\!g}_{_{F}}(n)
\}
\end{equation}
But we are also interested in meromorphic inputs $F$
that may have zeros and poles at 0 or 1. Since the mappings
$F\mapsto \tilde{I\!g}_{_{F}} $ and
$F\mapsto \tilde{E\!g}_{_{F}} $
are clearly multiplicative and since meromorphic functions $F$ possess convergent
Hadamard products\,:
\begin{equation}   \label{a43}
F(x)=c\,x^d\,e^{\sum_{s=1}^{s=\infty} c_s\,x^s}
\prod_i\Big( (1-\frac{x}{a_i})^{k_i}\,
e^{k_i\sum_{s=1}^{s=K_i} \frac{1}{s}\frac{x^s}{a_i^s}}
\Big)\quad(k_i\in\doZ,K_i\in\doN)
\end{equation} 
we require the exact form of the ingress factors for monomial, binomial
and even/odd exponential
factors. Here are the results\,:
\[\begin{array}{lllllllll}
F_{\mr{mon}}(x)\!\!&\!\!=\!\!&\!\!
c
&& \tilde{\mi{Ig}}_{\mr{mon}}(n)
\!\!&\!\!=\!\!&\!\! c^{+\frac{1}{2}}
\\[1.5 ex]
F_{\mr{mon}}(x)\!\!&\!\!=\!\!&\!\!
c\,{x^d}\quad (d\not=0)
&& \tilde{\mi{Ig}}_{\mr{mon}}(n)
\!\!&\!\!=\!\!&\!\! c^{-\frac{1}{2}}\,{(2\pi n)}^{\frac{d}{2}}
\\[1.5 ex]
F_{\mr{bin}}(x)
\!\!&\!\!=\!\!&\!\!
\prod_{i}(1-a_i^{-1}{x})^{s_i} 
&& \tilde{\mi{Ig}}_{\mr{bin}}(n)
\!\!&\!\!=\!\!&\!\!
\prod_{i}\big(\tilde{\Theta}(a_i n)\big)^{s_i}
\\[1.5 ex]
F_{\mr{even}}(x)
\!\!&\!\!=\!\!&\!\!
\exp\!\big(\!-\!\sum_{s\geq 1}f^{\mr{even}}_{2s}\,x^{2s}\big)\,
&& \tilde{\mi{Ig}}_{\mr{even}}(n)
\!\!&\!\!=\!\!& 1
\\[1.5 ex]
F_{\mr{odd}}(x)
\!\!&\!\!=\!\!&\!\!
\exp\!\big(\!-\!\sum_{s\geq 0}f^{\mr{odd}}_{2s+1}\,x^{2s+1}\big)\,
&& \tilde{\mi{Ig}}_{\mr{odd}}(n)
\!\!&\!\!=\!\!&\!\!
\exp\!\big(\!+\!\sum_{s\geq 0}f^{\mr{odd}}_{2s+1}\,
\frac{\frak{b}_{2s+1}^\ast(0)}{n^{2s+1}}\big)
\end{array}\]

\noindent
{\bf Monomial factors.}
\\
The discontinuity between the first two expressions of $\mi{Ig}_{\mr{mon}}(n)$
stems from the fact that for $d=0$ the product in (\ref{a2}) start from $k=0$
as usual, whereas for $d\not=0$ it has to start from $k=1$. The case $d=0$
is trivial, and the case $d\not=0$ by multiplicativity reduces
to the case $d=1$. To calculate the corresponding ingress factor,
we may specialise the identity
\begin{equation}   \label{a44}
\prod_{1\leq k\leq n} F(\frac{k}{n})\sim
\tilde{I\!g}_{F}(n)\;e^{-\nu_\ast\,n}\; \tilde{E\!g}_{F}(n)
\end{equation}
to convenient test functions. Here are the two simplest choices\,:
\[\begin{array}{lllllllll}
\mi{test\;function}\;F_1(x)=x 
&|\!|&
\mi{test\;function}\;F_2(x)=\frac{2}{\pi}\sin(\frac{\pi}{2}x) 
\\
\-----------
&|\!|&
\-------------------
\\
\prod_{k=1}^{k=n} F_1(\frac{k}{n})=\frac{n!}{n^n}
&|\!|&
\prod_{k=1}^{k=n} F_2(\frac{k}{n})=2\,\pi^{-n}\,n^{1/2}
\quad \mi{by\, elem. \, trigon.}
\\[1.3 ex]
\tilde{I\!g}_{F_1}(n) =\mi{unknown}
&|\!|& 
\tilde{I\!g}_{F_2}(n) =\tilde{I\!g}_{F_1}(n)
\hspace{9.ex}\mi{by\,parity\,of}\; \frac{F_2(x)}{F_1(x)}
\\[1.3 ex]
\nu_\ast=1
&|\!|&
\nu_\ast=-\log\pi
\\[1.3 ex]
\tilde{E\!g}_{F_1}(n) =\tilde{\Theta}(n)
&|\!|&
\tilde{E\!g}_{F_2}(n) =1 
\hspace{11.ex}\mi{by\,parity\,of}\; F_2(1\!+\!x)
\end{array}\]

With the choice $F_2$, all we have to do is plug the data in the second column into (\ref{a44})
and we immediately get
$\tilde{I\!g}_{F_2}(n)=(2\pi\,n)^{1/2}$ but before that we have to check the first line's
elementary trigonometric identity. With the choice $F_1$, on the other hand, we need to
check that the egress factor does indeed coincide with $\Theta$. This readily follows from\,:
\begin{eqnarray*}
F_1(1+x)&\!\!=\!\!&\exp\big(\sum_{1\leq s}(-1)^s\frac{x^s}{s}\big) \;\;\;
\Longrightarrow \;\;\;
\tilde{E\!g}_{F_1}(n)=\exp(\tilde{eg}_{F_1}(n))\quad\mi{with}
\\
\tilde{eg}_{F_1}(n)&\!\!=\!\!&\!\!\sum_{1\leq s\,
\mi{odd}}(-1)^{s-1}\frac{n^{-s}}{s}\frak{b}^\ast_s(0) 
=\!\!\sum_{1\leq s\,\mi{odd}}{n^{-s}}{(s-1)!}\,\frak{b}_s \quad \Longrightarrow
\\
\imi{eg}_{F_1}(\nu)&\!\!=\!\!&\!\!\sum_{1\leq s}{\nu^{s-1}}\frak{b}_s 
=\frac{1}{\nu}\,\Big( \frac{e^\nu}{e^\nu-1}-\frac{1}{\nu}-\frac{1}{2}\Big) 
=-\frac{1}{\nu^2}+\frac{1}{2\,\nu}\,\frac{1}{\tanh({\nu}/{2})}
=\;\imi{\theta}(\nu)
\end{eqnarray*} 
We then plug everything into (\ref{a44}) and use formula (\ref{a33}) of \S3.2 to eliminate
both $n!/n^n$ and $\Theta(n)$.
\\

\noindent
{\bf Binomial factors}. \\
By multiplicativity and homogeneity, it is enough to check
the idendity $\tilde{I\!g}_{F_3}(n)=\tilde{\Theta}(n)$ for the test function $F_3(x)=1-x$. 
But since $F_3=1/F_1^\perp$ with the notations of the preceding para,
the parity relations yield\,:
$$
\tilde{I\!g}_{F_3}(n)\,=\,
1/\tilde{E\!g}_{1/F_1^\perp}(n)\,=\,
\tilde{E\!g}_{F_1^\perp}(n)\,=\,
\tilde{\Theta}(n)
\quad\quad (\mi{see\,above)}
$$
which is precisely the required identity. We alse notice that\,:
$$
F(x)=(1-\frac{x}{a})(1+\frac{x}{a})
\quad \Longrightarrow\quad
\tilde{I\!g}_{F}(n)=\tilde{\Theta}(an)\tilde{\Theta}(-an)\equiv 1
$$
which agrees with the trivialness of the ingress factor for an {\it even} input $F$.
\\

\noindent
{\bf Exponential factors.}
\\
For them, the expression of the ingress/egress factors directly follows from (\ref{a32}).
Moreover, since the exponentials occurring in the Hadamard product (\ref{a43})
carry only polynomials or entire functions, the corresponding ingress/egress
factors are actually convergent. 
%%%%%%%%%%%%%%%%%%%%%%%%%%%%%%%%%%%%%%%%%%%%%%%%%%%%%%%%%%%%%%%%%%%%%%%%%%%%%%%%%%%%%%%%%%%%
%%%%%%%%%%%%%%%%%%%%%%%%%%%%%%%%%%%%%%%%%%%%%%%%%%%%%%%%%%%%%%%%%%%%%%%%%%%%%%%%%%%%%%%%%%%%
\subsection{Resummability of the total ingress factor.}
As announced, we shall have to change our SP-series 
$j(\zeta)=\sum J(n)\,\zeta^n$
into $j^\#(\zeta)=\sum J^\#(n)\,\zeta^n$,
which involves dividing the coefficients $J(n)$, not by the asymptotic
series $\tilde{I\!g}_F(n)$, but by its exact resummation  ${I\!g}_F(n)$.
Luckily, this presents no difficulty for meromorphic
\footnote{for simplicity, let us
assume that
$F$ has no purely imaginary poles or zeros.}
 inputs $F$. Indeed, 
the contributions to $\tilde{I\!g}_F$
of the isolated factors in (\ref{a43}) are separetely resummable\,:
\\
\--- for the monomial factors $F_{\mi{mon}}$ this is trivial
\\
\--- for the binary factors  $F_{\mi{bin}}$ this follows from
 $\tilde{\Theta}$'s resummability (see \S3.2)
\\
\--- for the exponential factors $F_{\mi{exp}}$ this follows from
$\log F_{\mi{exp}}$ being either polynomial or entire.

As for the global $\tilde{I\!g}_F$, one easily checks that the Hadamard product (\ref{a43}),
rewritten as
\begin{equation}   \label{a45}
F= F_{\mi{mon}}\;F_{\mi{even}}\;F_{\mi{odd}}\;\prod_i F_{\mi{bin}, i}
\end{equation}
 translates into a product of resurgent functions\,:
\begin{equation}   \label{a46}
\tilde{I\!g}_F= 
\tilde{I\!g}_{F_{\mi{mon}}}\;
\tilde{I\!g}_{F_{\mi{even}}}\;\tilde{I\!g}_{F_{\mi{odd}}}\;
\prod_i 
\tilde{I\!g}_{F_{\mi{bin}, i}}
\end{equation}
which converges in all three models (formal, convolutive, geometric \--- respective
to the corresponding topology) to a limit that doesn't depend on the actual
Hadamard decomposition chosen in (\ref{a43}), i.e. on the actual choice of the truncation-defining
integers $K_i$.

%%%%%%%%%%%%%%%%%%%%%%%%%%%%%%%%%%%%%%%%%%%%%%%%%%%%%%%%%%%%%%%%%%%%%%%%%%%%%%%%%%%%%%%%%%%%
%%%%%%%%%%%%%%%%%%%%%%%%%%%%%%%%%%%%%%%%%%%%%%%%%%%%%%%%%%%%%%%%%%%%%%%%%%%%%%%%%%%%%%%%%%%%
%%%%%%%%%%%%%%%%%%%%%%%%%%%%%%%%%%%%%%%%%%%%%%%%%%%%%%%%%%%%%%%%%%%%%%%%%%%%%%%%%%%%%%%%%%%%
%%%%%%%%%%%%%%%%%%%%%%%%%%%%%%%%%%%%%%%%%%%%%%%%%%%%%%%%%%%%%%%%%%%%%%%%%%%%%%%%%%%%%%%%%%%%
%%%%%%%%%%%%%%%%%%%%%%%%%%%%%%%%%%%%%%%%%%%%%%%%%%%%%%%%%%%%%%%%%%%%%%%%%%%%%%%%%%%%%%%%%%%%
%%%%%%%%%%%%%%%%%%%%%%%%%%%%%%%%%%%%%%%%%%%%%%%%%%%%%%%%%%%%%%%%%%%%%%%%%%%%%%%%%%%%%%%%%%%%
%%%%%%%%%%%%%%%%%%%%%%%%%%%%%%%%%%%%%%%%%%%%%%%%%%%%%%%%%%%%%%%%%%%%%%%%%%%%%%%%%%%%%%%%%%%%
%%%%%%%%%%%%%%%%%%%%%%%%%%%%%%%%%%%%%%%%%%%%%%%%%%%%%%%%%%%%%%%%%%%%%%%%%%%%%%%%%%%%%%%%%%%%
%%%%%%%%%%%%%%%%%%%%%%%%%%%%%%%%%%%%%%%%%%%%%%%%%%%%%%%%%%%%%%%%%%%%%%%%%%%%%%%%%%%%%%%%%%%%
%%%%%%%%%%%%%%%%%%%%%%%%%%%%%%%%%%%%%%%%%%%%%%%%%%%%%%%%%%%%%%%%%%%%%%%%%%%%%%%%%%%%%%%%%%%%
%%%%%%%%%%%%%%%%%%%%%%%%%%%%%%%%%%%%%%%%%%%%%%%%%%%%%%%%%%%%%%%%%%%%%%%%%%%%%%%%%%%%%%%%%%%%
%%%%%%%%%%%%%%%%%%%%%%%%%%%%%%%%%%%%%%%%%%%%%%%%%%%%%%%%%%%%%%%%%%%%%%%%%%%%%%%%%%%%%%%%%%%%
%%%%%%%%%%%%%%%%%%%%%%%%%%%%%%%%%%%%%%%%%%%%%%%%%%%%%%%%%%%%%%%%%%%%%%%%%%%%%%%%%%%%%%%%%%%%
%%%%%%%%%%%%%%%%%%%%%%%%%%%%%%%%%%%%%%%%%%%%%%%%%%%%%%%%%%%%%%%%%%%%%%%%%%%%%%%%%%%%%%%%%%%%
%%%%%%%%%%%%%%%%%%%%%%%%%%%%%%%%%%%%%%%%%%%%%%%%%%%%%%%%%%%%%%%%%%%%%%%%%%%%%%%%%%%%%%%%%%%%
%%%%%%%%%%%%%%%%%%%%%%%%%%%%%%%%%%%%%%%%%%%%%%%%%%%%%%%%%%%%%%%%%%%%%%%%%%%%%%%%%%%%%%%%%%%%
\subsection{Parity relations.}
Starting from the elementary parity relations for the Bernoulli
numbers and polynomials\,:
\begin{eqnarray*}
\frak{b}_{2s}=0\hspace{5.ex}(s\geq 1)\quad &;&\quad \beta_{2s}=0\hspace{5.ex}(s\geq 0)
\\
\frak{b}_s(\tau)\equiv (-1)^{s+1}\,\frak{b}_s(-\tau-1)\quad &;&\quad \beta_s(\tau)\equiv (-1)^{s+1}\,\beta(-\tau)
\hspace{7.ex} (\forall s \geq 0)
\end{eqnarray*}
and setting  
\begin{eqnarray*}
F^{\,\vdash}(x):=1/F(1\!-\!x)\hspace{11.ex}&&
\\
P_F(n):=\prod_{m=0}^{m=n}F(\frac{m}{n})\hspace{11.ex}&&
\\
P_F^\#(n):=\frac{P_F(n)}{I\!g_F(n)\,E\!g_F(n)}=\omega_F^n 
\;\;&\mi{with} &\;\; \omega_F^n=\exp(-\!\int_0^\infty\!\!\!f(x)dx)
\end{eqnarray*}
we easily check that\,:
\begin{eqnarray*} 
\tilde{I\!g}_{_{F}}(n)\,\tilde{E\!g}_{_{F^{\vdash}}}(n)=1\hspace{7. ex}
&\mi{and}&\hspace{3.5 ex}
\tilde{I\!g}_{_{F^{\vdash}}}(n)\,\tilde{E\!g}_{_{F}}(n)=1
\\
J_{F^{\vdash}}(n)=J_F(n)/P_F(n) 
\hspace{4 ex}&\mi{and}&\hspace{3 ex}
J^\#_{F^{\vdash}}(n)=J^\#_F(n)/P^\#_F(n)
\\
j_{F^{\vdash}}(\zeta)\not=j_F(\zeta/\omega_F) 
\hspace{8 ex}&\mi{but}&\hspace{3 ex}
j^\#_{F^{\vdash}}(\zeta)=j^\#_F(\zeta/\omega_F) 
\end{eqnarray*}

%%%%%%%%%%%%%%%%%%%%%%%%%%%%%%%%%%%%%%%%%%%%%%%%%%%%%%%%%%%%%%%%%%%%%%%%%%%%%%%%%%%%%%%%%%%%
%%%%%%%%%%%%%%%%%%%%%%%%%%%%%%%%%%%%%%%%%%%%%%%%%%%%%%%%%%%%%%%%%%%%%%%%%%%%%%%%%%%%%%%%%%%%

% \end{document}
%%%%%%%%%%%%%%%%%%%%%%%%%%%%%%%%%%%%%%%%%%%%%%%%%%%%%%%%%%%%%%%%%%%%%%%%%%%%%%%%%%%%%%%%%%%%
%%%%%%%%%%%%%%%%%%%%%%%%%%%%%%%%%%%%%%%%%%%%%%%%%%%%%%%%%%%%%%%%%%%%%%%%%%%%%%%%%%%%%%%%%%%%
%%%%%%%%%%%%%%%%%%%%%%%%%%%%%%%%%%%%%%%%%%%%%%%%%%%%%%%%%%%%%%%%%%%%%%%%%%%%%%%%%%%%%%%%%%%%
%%%%%%%%%%%%%%%%%%%%%%%%%%%%%%%%%%%%%%%%%%%%%%%%%%%%%%%%%%%%%%%%%%%%%%%%%%%%%%%%%%%%%%%%%%%%
%%%%%%%%%%%%%%%%%%%%%%%%%%%%%%%%%%%%%%%%%%%%%%%%%%%%%%%%%%%%%%%%%%%%%%%%%%%%%%%%%%%%%%%%%%%%
%%%%%%%%%%%%%%%%%%%%%%%%%%%%%%%%%%%%%%%%%%%%%%%%%%%%%%%%%%%%%%%%%%%%%%%%%%%%%%%%%%%%%%%%%%%%
%%%%%%%%%%%%%%%%%%%%%%%%%%%%%%%%%%%%%%%%%%%%%%%%%%%%%%%%%%%%%%%%%%%%%%%%%%%%%%%%%%%%%%%%%%%%
%%%%%%%%%%%%%%%%%%%%%%%%%%%%%%%%%%%%%%%%%%%%%%%%%%%%%%%%%%%%%%%%%%%%%%%%%%%%%%%%%%%%%%%%%%%%
%%%%%%%%%%%%%%%%%%%%%%%%%%%%%%%%%%%%%%%%%%%%%%%%%%%%%%%%%%%%%%%%%%%%%%%%%%%%%%%%%%%%%%%%%%%%

%

%%%%%%%%%%%%%%%%%%%%%%%%%%%%%%%%%%%%%%%%%%%%%%%%%%%%%%%%%%
%% Some resurgence properties of knot-related functions.
%%%%%%%%%%%%%%%%%%%%%%%%%%%%%%%%%%%%%%%%%%%%%%%%%%%%%%%%%%

%\documentclass[12pt,a4paper]{article}\input{SP_commands}\begin{document}

%%%%%%%%%%%%%%%%%%%%%%%%%%%%%%%%%%%%%%%%%%%%%%%%%%%%%%%%%%%%%%%%%%%%%%%%%%%%%%%%%%%%%%%%%%%%
%%%%%%%%%%%%%%%%%%%%%%%%%%%%%%%%%%%%%%%%%%%%%%%%%%%%%%%%%%%%%%%%%%%%%%%%%%%%%%%%%%%%%%%%%%%%
%%%%%%%%%%%%%%%%%%%%%%%%%%%%%%%%%%%%%%%%%%%%%%%%%%%%%%%%%%%%%%%%%%%%%%%%%%%%%%%%%%%%%%%%%%%%

\section{Inner generators.} 
%%%%%%%%%%%%%%%%%%%%%%%%%%%%%%%%%%%%%%%%%%%%%%%%%%%%%%%%%%%%%%%%%%%%%%%%%%%%%%%%%%%%%%%%%%%%
%%%%%%%%%%%%%%%%%%%%%%%%%%%%%%%%%%%%%%%%%%%%%%%%%%%%%%%%%%%%%%%%%%%%%%%%%%%%%%%%%%%%%%%%%%%%

\subsection{Some heuristics.}
Consider a simple, yet typical case. Assume the driving function $f$
 to be {\it entire} (or even think of it as {\it polynomial}, for simplicity),
steadily increasing on
the real interval
$[0,1]$, with a unique zero at
$\bar{x}\in ]0,1[$ on that interval, and no other zeros, real or complex, inside the disk 
$\{|x|\leq |\bar{x}|\}$ \;:
\begin{equation}   \label{a47}
0<\bar{x}<1\;\;,\;\; 
f(0)<0\;,\; f(\bar{x})=0\;,\; f(1)>0 \;,\;\;
f^\prime(x)>0\;\forall x \in [0,1] \quad
\end{equation}
As a consequence, the primitive $f^\ast(x):=\int_0^x f(x^\prime)dx^\prime $
will display a unique minimum at $\bar{x}$ and, for any given large $n$, the
products $\prod_{k=0}^{k=m} F(k/n)=\exp(-\sum_{k=0}^{k=m} f(k/n)) $
will be maximal for $m\sim n\bar{x}$.
It is natural, therefore, to split the Taylor coefficients $J(n)$ of our sum-product
 series (\ref{a2}) into
 a {\it global} but fairly elementary 
factor $J_{1,2,3}(n)$, which subsumes all the {\it pre-critical}
terms $F(k/n)$, and a purely {\it local} but 
analytically more challenging factor $J_4(n)$, which
accounts for the contribution of all {\it near-critical} terms $F(k/n)$. 
Here are the definitions\,:
\begin{eqnarray}   \label{a48}
J(n) &:=& J_{1,2,3}(n)\;J_{4}(n)
\\[1. ex]   \label{a49}
J(n) &:=& \sum_{0\leq m\leq n}\;\prod_{0 \leq k \leq m}F(\frac{k}{n})
\\[1. ex]   \label{a50}
J_{1,2,3}(n)\!\!\! &:=& \prod_{0 \leq k \leq \bar{m}}F(\frac{k}{n})\quad\quad \mi{ with}
\quad\quad \bar{m}:=\mr{ent}(n\,\bar{x})
\\[1. ex]   \label{a51}
J_{4}(n)
&:=& \sum_{0\leq m\leq n}\;\;\Bigg( 
 \prod_{0 \leq k \leq m}F(\frac{k}{n})\;\;
\Big/
 \prod_{0 \leq k \leq \bar{m}}F(\frac{k}{n})\Bigg)
\\ [1.ex]   \label{a52}
& =& 
\left\{\begin{matrix}
\dots+\frac{1}{F}(\frac{\bar{m}-2}{n})\,(\frac{\bar{m}-1}{n})\,
\frac{1}{F}(\frac{\bar{m}\!}{n})
 +\frac{1}{F}(\frac{\bar{m}-1}{n})\,\frac{1}{F}(\frac{\bar{m}\!}{n})
\\[0.5 ex]
+\frac{1}{F}(\frac{\bar{m}\!}{n}) +1 +F(\frac{\bar{m}+1}{n})
\\[0.5 ex]
+F(\frac{\bar{m}\!+\!1}{n})\,F(\frac{\bar{m}\!+\!2}{n})
+F(\frac{\bar{m}\!+\!1}{n})\,F(\frac{\bar{m}\!+\!2}{n})\,F(\frac{\bar{m}\!+\!3}{n})+\dots 
\end{matrix}\right\}
\\[-0.5 ex]\nonumber
\end{eqnarray}
The asympotics of the global factor $J_{123}(n)$ as $n\rightarrow \infty$
easily results from the variant (\ref{a32}) of the Euler-MacLaurin formula and $J_{123}(n)$ 
splits into three subfactors\,:
\\
(i) a factor $J_1(n)$, local at $x=0$, which is none other than the ingress factor $Ig_F(n)$
studied at length in \S3.
\\
(ii) an elementary factor $J_2(n)$, which reduces to an exponential and
carries no divergence.
\\
(iii) a factor $J_3(n)$, local at $x=\bar{x}$ and analogous 
 to the `egress factor' of \S3, but with base point $\bar{x}$ instead of $1$.

That leaves the really sensitive factor $J_4(n)$, which like $J_3(n)$
is local at $x=\bar{x}$, but far more complex.
In view of its expression as the discrete sum
(\ref{a52}), we should expect its asymptotics to be described by
a Laurent series $\sum_{k\geq 0}C_{k/2}\,n^{-k/2} $
involving both integral and semi-integral powers of $1/n$. That turns out
 to be the case indeed, but we shall see that there is a way of jettisoning the integral
powers and retaining 
only the semi-integral ones, i.e. $\sum_{k\geq 0}C_{k+1/2}\,n^{-k-1/2} $.
To do this, we must perform 
a little sleight-of-hand
and attach the egress factor $J_3$ to $J_4$ so as to produce 
the joint factor $J_{3,4}$. In fact, as we shall see, the gains that accrue
from merging $J_3$ and $J_4$ go way beyond the elimination of integral powers.

But rather than rushing ahead, let us describe our four factors $J_i(n)$
and their asymptotic expansions $\tilde{J}_i(n)$\,:
\begin{eqnarray} 
\nonumber
J(n) &:=& J_{1,2,3}(n)\;J_{4}(n)=J_{1}(n)\;J_{2}(n)\,J_{3}(n)\;J_{4}(n)=J_{1,2}(n)\;J_{3,4}(n)
\\[1. ex] \nonumber
J_{1,2,3}(n)\!\!\! &:=&\!\! \prod_{0 \leq k \leq \bar{m}}F(\frac{k}{n})
=\exp(\,-\!\!\!\sum_{0 \leq k \leq \bar{m}}f(\frac{k}{n}))
\quad\quad \mi{ with}
\quad\quad \bar{m}:=\mr{ent}(n\,\bar{x})
\\[1. ex]    \label{a53}
\tilde{J}_{1}(n)\!\!\! &:=&\!\! 
\exp(-\frac{1}{2}\,f_0+\sum_{1\leq s} \frac{\frak{b}^\ast_s(0)\,f_s}{n^s})
\quad\quad \mi{with}\;\quad f_s:=\frac{f^{(s)}(0)}{s!}
\\[1. ex]    \label{a54}
\tilde{J}_{2}(n)\!\!\! &:=&\!\! \exp(-n\,\bar{\nu})\hspace{20.ex} \mi{with}\quad\; 
\bar{\nu}:=\int_{0}^{\bar{x}}f(x)\,dx
\\[1. ex]   \label{a55}
\tilde{J}_{3}(n)\!\!\! &:=&\!\!  
\exp(-\frac{1}{2}\,\bar{f}_0-\sum_{1\leq s} \frac{\frak{b}^\ast_s(0)\,\bar{f}_s}{n^s})
\quad\quad \mi{with}\;\quad \bar{f}_s:=\frac{f^{(s)}(\bar{x})}{s!}
\\[1. ex]\nonumber
\tilde{J}_{4}(n)\!\!\! &:=&\!\! \sum_{0\leq \bar{m}\leq \bar{x}\,n}
\exp{\big(\!\!\!\sum_{0\leq k\leq \bar{m}}f(\bar{x}-\frac{k}{n})\big)}
+\sum_{0 \leq \bar{m}\leq (1-\bar{x})\,n}
\exp{\big(-\!\!\!\sum_{1\leq k\leq \bar{m}}f(\bar{x}+\frac{k}{n})\big)}
\end{eqnarray}
In the last identity, the first exponential inside the second sum, namely $\exp(\sum_{1\leq k\leq
0}(\dots))$, should of course be taken as $\exp(0)=1$. Let us now simplify $\tilde{J}_3$ by using
the fact that
$\bar{f}_0=f(\bar{x})=0$ and let us replace in  $\tilde{J}_4$ the finite
$\bar{m}$-summation (up to $\bar{x}\,n $ or $(1-\bar{x})\,n  $) by an infinite $m$-summation,
up to $+\infty$, which won't change 
 the {\it asymptotics} \footnote{
It will merely change the {\it transasymptotics}  by adding exponentially small summands.
In terms of the Borel transforms $\imi{J}_{3,4}\!(\nu)$ or $\smi{J}_{3,4}\!(\nu)$,
it means that their nearest singularities will remain unchanged.
} in $n$: 
\begin{eqnarray}   \label{a56}
\tilde{J}_{3}(n)\!\!\! &:=& \exp(-\sum_{1\leq s}\frak{b}_s\,\bar{f}_s)
\\[1. ex]   \label{a57}
\tilde{J}_{4}(n)\!\!\! 
&:=&\!\! 2+ \sum_{{1\leq m \atop \epsilon =\pm 1 }}\exp{\Big(-\epsilon\!\!\sum_{1\leq k\leq
m}f(\bar{x}+\epsilon\frac{k}{n})\Big)}  
\\[1. ex]   \label{a58}
\tilde{J}_{4}(n)\!\!\! &:=&\!\! 2+ \sum_{{1\leq m \atop \epsilon =\pm 1 }}
\exp{\Big(-\!\!\sum_{1\leq s}\frac{\epsilon^{s+1}}{n^s}
\big(\frak{b}^\ast_s(m)-\frak{b}^\ast_s(0) \big)
\bar{f}_s
\Big)} 
\end{eqnarray}
We can now regroup the factors $\tilde{J}_3$ and  $\tilde{J}_4$  into  $\tilde{J}_{3,4}$
and switch from the Bernoulli-type polynomials $\frak{b}^\ast_s(m)$ over to
their Euler-Bernoulli counterparts $\beta^\ast_s(m)$. These have the advantage of being odd/even
function of $m$ is $s$ is even/odd, which will enable us to replace
$m$-summation on $\doN$ by $m$-summation on $\frac{1}{2}+\doZ$, eventually
easing the change from  $m$-summation to $\tau$-integration. We successively find\,:
\begin{eqnarray}   \label{a59}
\tilde{J}_{3,4}(n)\!\!\! 
&=&\!\! 2\,\exp{\Big(-\!\!\sum_{1\leq s}\frac{1}{n^s}\frak{b}_s\,\bar{f}_s
\Big)} 
+\sum_{{1\leq m \atop \epsilon =\pm 1 }}
\exp{\Big(-\!\!\sum_{1\leq s}
\frac{\epsilon^{s+1}}{n^s}
\frak{b}^\ast_s(m)\,\bar{f}_s
\Big)} 
\\[1. ex]   \label{a60}
\tilde{J}_{3,4}(n)\!\!\! &=&
 \sum_{{0\leq m \atop \epsilon =\pm 1 }}
\exp{\Big(-\!\!\sum_{1\leq s}
\frac{\epsilon^{s+1}}{n^s}
\frak{b}^\ast_s(m)\,\bar{f}_s
\Big)} 
\\[1. ex]   \label{a61}
\tilde{J}_{3,4}(n)\!\!\! &=&
 \sum_{{0\leq m \atop \epsilon =\pm 1 }}
\exp{\Big(-\!\!\sum_{1\leq s}
\frac{\epsilon^{s+1}}{n^s}
\beta^\ast_s(m+\frac{1}{2})\,\bar{f}_s
\Big)} 
\\[1. ex]   \label{a62}
\tilde{J}_{3,4}(n)\!\!\! &=&
 \sum_{{0\leq m \atop \epsilon =\pm 1 }}
\exp{\Big(-\!\!\sum_{1\leq s}
\frac{1}{n^s}
\beta^\ast_s(\epsilon m+\epsilon \frac{1}{2})\,\bar{f}_s
\Big)}
\hspace{8.ex} \mi{(by \,parity !)}\quad
\\[1. ex]   \label{a63}
\tilde{J}_{3,4}(n)\!\!\! &=&
 \sum_{m\,\in\, \frac{1}{2}+\doZ}
\exp{\Big(-\!\!\sum_{1\leq s}
\frac{1}{n^s}
\beta^\ast_s( m)\,\bar{f}_s
\Big)}  
\end{eqnarray} 
This last identity should actually be construed as\,:
\begin{equation}   \label{a64}
\tilde{J}_{3,4}(n) =
 \sum_{m\,\in\, \frac{1}{2}+\doZ}
\exp{\Big(\!\!-
\frac{1}{n}
\beta^\ast_1( m)\,\bar{f}_1
\Big)}
\exp_\#{\Big(\!\!-\!\!\sum_{1\leq s}
\frac{1}{n^s}
\beta^\ast_s( m)\,\bar{f}_s
\Big)}  
\end{equation} 
Here, the first exponential $\mi{exp}$ decreases fast as $m$ grows, since
$$
\beta^\ast_1(m)\,\bar{f}_1 = \frac{1}{2}\,m^2\,\bar{f}_1 = \frac{1}{2} m^2 f^\prime(\bar{x}) >0
$$
The second  exponential $\mi{exp}_\# $, on the other hand, should be {\it expanded}
as a power series of its argument and each of the resulting terms $m^{s_1}\,n^{-s_2}$ 
should be dealt with separately, leading to a string of clearly {\it convergent} series.
We can now replace the {\it discrete} $m$-summation in (\ref{a63}) by
a {\it continuous} $\tau$-integration\,: here again, that may change the {\it transasymptotics}
in $n$, but not the {\it assymptotics}.\footnote{
Indeed, the summation/integration bounds are $\pm\infty$, with a summand/\!/integrand
that vanishes exponentially fast there.}
We find, using the parity properties of the $\beta^\ast_s$ and maintaining throughout the distinction between $\mi{exp}$ (unexpanded)
and $\mi{exp}_{\#}$ (expanded) \,:
\begin{eqnarray}   \label{a65}
\tilde{J}_{3,4}(n) &\!\!=\!\!&
\int_{-\infty}^{+\infty}
\exp{\Big(\!\!-
\frac{1}{n}
\beta^\ast_1( \tau)\,\bar{f}_1
\Big)}
\exp_\#{\Big(\!\!-\!\!\sum_{2\leq s}
\frac{1}{n^s}
\beta^\ast_s( \tau)\,\bar{f}_s
\Big)} d\tau
\\[1.5 ex]   \label{a66}
\tilde{J}_{3,4}(n) &\!\!=\!\!&
\sum_{\epsilon=\pm1}
\int_{-\infty}^{+\infty}
\exp{\Big(\!\!-
\frac{1}{n}
\beta^\ast_1(\epsilon \tau)\,\bar{f}_1
\Big)}
\exp_\#{\Big(\!\!-\!\!\sum_{2\leq s}
\frac{1}{n^s}
\beta^\ast_s(\epsilon \tau)\,\bar{f}_s
\Big)} d\tau
\\[1.5 ex]   \label{a67}
\tilde{J}_{3,4}(n) &\!\!=\!\!&
\sum_{\epsilon=\pm1}
\int_{-\infty}^{+\infty}
\exp{\Big(\!\!-
\frac{1}{n}
\beta^\ast_1( \tau)\,\bar{f}_1
\Big)}
\exp_\#{\Big(\!\!-\!\!\sum_{2\leq s}
\frac{\epsilon^{s+1}}{n^s}
\beta^\ast_s(\tau)\,\bar{f}_s
\Big)} d\tau
\\[1.5 ex]\nonumber   
\tilde{J}_{3,4}(n) &\!\!=\!\!&
2\,\int_{-\infty}^{+\infty}
\exp{\Big(\!\!-
\frac{1}{n}
\beta^\ast_1( \tau)\,\bar{f}_1
\Big)}
\exp_\#{\Big(\!\!-\!\!\sum_{3\leq s\,\mi{odd}}
\frac{1}{n^s}
\beta^\ast_s(\tau)\,\bar{f}_s
\Big)} 
\\[1.5 ex]   \label{a68}
&&\quad\quad\times\; \;\cosh_\#{\Big(\!\!-\!\!\sum_{2\leq s\,\mi{even}}
\frac{1}{n^s}
\beta^\ast_s(\tau)\,\bar{f}_s\Big)}\;d\tau 
\\[1.5 ex]    \label{a69}
\tilde{J}_{3,4}(n) &\!\!=\!\!&
2\,\Big[
\int_{-\infty}^{+\infty}
\exp{\Big(\!\!-
\frac{1}{n}
\beta^\ast_1( \tau)\,\bar{f}_1
\Big)}
\exp_\#{\Big(\!\!-\!\!\sum_{2\leq s}
\frac{1}{n^s}
\beta^\ast_s(\tau)\,\bar{f}_s
\Big)} d\tau\Big]_{{\mi{} \atop \mi{demi\;}}\atop {\mi{integ.} \atop \mi{part\;\;}}}
\end{eqnarray}
The notation in (\ref{a69}) means that we retain only the demi-integral 
powers of $n^{-1}$ in $\big[\dots \big]$. In view of the results of \S2.3 about the
correspondence between singularities and Taylor coefficient asymptotics,  the Borel transform
$\imi{\ell i}(\nu):=\imi{J}_{3,4}(\nu)$ of $\tilde{J}_{3,4}(n)$,
or rather its counterpart 
$\imi{Li}(\zeta):=\imi{\ell i}(\log(1+\frac{\zeta}{\omega})) $
in the $\zeta$-plane, must in our test-case (\ref{a47}) describe the closest singularities of
the {\it sum-product} function $j(\zeta)$ or rather its `cleansed' variant $j^\#(\zeta)$.

Singularities such as $Li$ shall be referred to as {\it inner generators} of the resurgence
algebra.  They differ from the three other types of generators 
({\it original, exceptional, outer})
 first and foremost by their stability\,: unlike these, they self-reproduce
indefinitely under alien differentiation. Another difference is this\,:
inner generators (minors and majors alike)
tend to carry only {\it demi-integral}\footnote{
in our test-case, i.e. for a driving function $f$ with a simple zero
at $\bar{x}$. For zeros of odd order $\tau>1$ ($\tau$ has to be odd
to produce an extremum in $f^\ast$) we would get ramifications
of  order $(-\tau+2s)/(\tau+1)$ with $s \in \doN$, which again
rules out entire powers. See \S4.2-5 and also \S6.1
} powers of $\zeta$ or $\nu$, as we just saw, whereas the other
types of generators tend to carry only integral powers (in the minors)
and logarithmic terms (in the majors).

So far, so good. But
what about the two omitted factors $J_1(n)$ and $J_2(n)$? The second one, $J_2(n)$,
which is a mere exponential  $\mi{exp}(-n\bar{\nu})$, 
simply accounts for the location $\bar{\zeta}=e^{\bar{\nu}}$ at which $Li$ is seen in the
$\zeta$-plane. 
As for the ingress factor $J_1(n)$,  {\it keeping it} (i.e.
 merging it with $J_{3,4}(n)$) would have rendered $Li$
dependent on the ingress point $x=0$, whereas {\it removing it} ensures that $Li$
(and by extension the whole {\it inner algebra}) is totally independent of the
`accidents' of its construction, such as the choice of ingress point in the $x$-plane. 

As for the move from $\{\frak{b}_s^\ast \}$ to $\{\beta_s^\ast\}$,
apart from easing the change from summation to integration, it brings another,
even greater benefit\,: by removing the crucial coefficient $\beta_0$
in (\ref{a25}) (which vanishes, unlike $\frak{b}_0$ in (\ref{a18})),
 it shall enable us to express 
the future {\it mir}-transform as a purely integro-differential operator.
\footnote{ Indeed,
if it $\beta_0$ didn't vanish, that lone coefficient would suffice
to ruin nearly all the basic
formulae ({\it infra}) about {\it mir} and {\it nir}.}

One last remark, before bringing these heuristics to a close. We have chosen
here the simplest possible way of producing an {\it inner generator}, namely
directly from the {\it original generator} i.e. the sum-product series itself.
To do this, rather stringent assumptions on the driving function $f$
had to be made\footnote{like (\ref{a47})}. However, even when
these assumptions are not met, the {\it original generator} always
produces so-called {\it outer generators} (at least one, but generally two),
which in turn always produce {\it inner generators}. So these two types -- outer and inner --
are a universal feature of sum-product series. 
%%%%%%%%%%%%%%%%%%%%%%%%%%%%%%%%%%%%%%%%%%%%%%%%%%%%%%%%%%%%%%%%%%%%%%%%%%%%%%%%%%%%%%%%%%%%
%%%%%%%%%%%%%%%%%%%%%%%%%%%%%%%%%%%%%%%%%%%%%%%%%%%%%%%%%%%%%%%%%%%%%%%%%%%%%%%%%%%%%%%%%%%%

\subsection{The long chain behind {\it nir/\!/mir}.}
Let us now introduce two non-linear functional transforms central to
this investigation. The {\it nir}-transform is directly inspired
by the above heuristics. It splits into a chain of subtransforms,
all of which are elementary, save for one\,: the {\it mir}-transform.

Both {\it nir} and {\it mir} depend on a coherent choice of scalars $\beta_k$ and
polynomials $\beta^\ast_k(\tau)$. The {\it standard choice}, or Euler-Bernoulli choice,
corresponds to the definitions (\ref{a24}) - (\ref{a30}). It is the one that is relevant in
most applications to analysis and SP-series. However, to gain a better insight into 
the $\beta$-dependence of {\it nir/\!/mir}, it is also useful to consider the
{\it non-standard choice}, with free coefficients $\beta_k$ and
 accordingly redefined
polynomials $\beta_k^\ast(\tau)$\,:
\begin{eqnarray*}
\mi{standard\;choice}\hspace{9.ex}&&\hspace{6.ex}\mi{non{\scriptstyle -}standard\;choice}
\\
\beta(\tau)\!:=\!\frac{1}{e^{\tau/2}-e^{-\tau/2}}=\sum_{-1\leq k}\beta_k\,\tau^k
&&
\beta(\tau):= \sum_{-1\leq k}\beta_k\,\tau^k
\\
\beta^\ast_k(\tau)\!:=\!\beta(\partial_\tau)\,\tau^k=\frac{B_{k+1}(\tau\!+\!\frac{1}{2})}{k+1}
&&
\beta^\ast_k(\tau):=\beta(\partial_\tau)\,\tau^k=\!\!\sum_{s=-1}^{s=k}\!\!\beta_s\,
\tau^{k-s}\frac{\Gamma(k\!+\!1)}{\Gamma(k\!+\!1\!-\!s)}
\end{eqnarray*}
As we shall see, even in the non-standard case it is often necessary
to assume that $\beta_{-1}=1$ and $\beta_0=0$ (like in the standard case) to
get interesting results. The further coefficients, however, can be anything.
\\

\noindent
{\bf The long,  nine-link chain:}
\[\begin{array}{cccccccccccc}
 
\stackrel{\mathbf{nir}}{\rightarrow}   &\rightarrow &\rightarrow
&\rightarrow &\stackrel{\mathbf{nir}}{\rightarrow}\quad   & H 
  &\quad\quad\quad & H&\!=\!&\imi{Li}\;\mi{or} \;\imi{Le}
\\
\uparrow\quad\quad &        &            &           &  \downarrow  
&\;\;\uparrow\stackrel{9}{}
\\
\uparrow\quad\quad & \;\,f^*      &\stackrel{3}{\rightarrow} & \;\,g^*        
&   
\downarrow    &
\;h^{\,\prime}    &\quad\quad\quad & h^\prime&\!=\!&\imi{\ell i}\;\mi{or}\;\imi{\ell e}
\\
 \uparrow\quad\quad&\;\;\uparrow \stackrel{2}{}&  &\;\;\downarrow\stackrel{4}{} &   
\downarrow    &
\;\;\uparrow\stackrel{8}{}\\
\stackrel{\mathbf{nir}} {\longleftarrow }  & f      &            & g  &\quad
\quad\stackrel{\mathbf{nir}}{\longrightarrow}  &
h   &\quad\quad\quad & h&\!=\!&\smi{\ell i}\;\mi{or}\;\smi{\ell e}
\\
              &\;\;\uparrow \stackrel{1}{} & &\;\;\downarrow\stackrel{5}{} &             
&\;\;\uparrow\stackrel{7}{}\\
              & F      &            & {\gbar} &      
\stackrel{6}{\longrightarrow\;\longrightarrow}    &
 {\hbar} \\
              & &  &&\mathbf{mir} & \\

\end{array}\]
{\bf Details of the nine steps:}
\[\begin{array}{cccccccc}
\stackrel{1}{\rightarrow}&:& \mi{precomposition} \; &:&
\; F\rightarrow f\;\;
&\mi{with}& \;\; f(x):=-\log F(x)
\\
\stackrel{2}{\rightarrow}&:& \mi{integration} \; &:&
\; f \rightarrow f^* \;\;
&\mi{with}&\;\; f^*(x):=\int_0^x f(x_0)\,dx_0
\\
\stackrel{3}{\rightarrow}&:& \mi{reciprocation} \;&:&
\; f^* \rightarrow g^* \;\;
&\mi{with}&\;\; f^*\circ g^*=\mi{id}
\\
\stackrel{4}{\rightarrow}&:& \mi{derivation} \;&:&
\; g^*\rightarrow g\;\;
&\mi{with}&\;\; g(y):=\frac{d}{dy}g^*(y)
\\
\stackrel{5}{\rightarrow}&:& \mi{inversion} \;&:&
\; g\rightarrow \gbar \;\;
&\mi{with}&\;\; {-\hspace{-1.5 ex}g}(y)\, := 1/g(y) 
\\
\stackrel{6}{\rightarrow}&:& \mg{mir}\; \mi{functional} \;&:&
\; {-\hspace{-1.5 ex}g} \;
\rightarrow \;
{\hbar}\;\;
&\mi{with}& \mi{int.-}\mi{diff.}\, \mi{expression}\, \mi{below}
\\
\stackrel{7}{\rightarrow}&:& \mi{inversion} \;&:&
\; 
{\hbar}\rightarrow h\;\;
&\mi{with}& \;
\; 
h(\nu):=1/{\hbar}(\nu)
\\
\stackrel{8}{\rightarrow}&:& \mi{derivation} \;&:&
\; h \rightarrow \;h^\prime
&\mi{with}&\;
\; 
h^\prime(\nu):=\frac{d}{d\nu}h(\nu)
\\
\stackrel{9}{\rightarrow}&:& \mi{postcomposition} \;&:&
\; 
h^\prime\rightarrow \;H
&\mi{with}&
\; H(\zeta):=h^\prime(\log(1+\frac{\zeta}{\omega})
\\[0.8 ex]
\stackrel{2...7}{\rightarrow}&:& \mg{nir}\; \mi{functional} \;&:&
\; {g} \;
\rightarrow \;
{ h}\;\;
&\mi{with}&  \mi{see\;\S4.3\; infra}
\end{array}\]

\noindent
{\bf ``Compact" and ``layered" expansions of {\it mir}. }
\\

\noindent
The `sensitive' part of the nine-link chain, namely the {\it mir}-transform, is a non-linear integro-differential functional of infinite
order. Pending its detailed description in \S4-5, let us write down the general shape of its two expansions: the `compact' expansion, which 
merely isolates the $r$-linear parts, and the more precise `layered' expansion, which takes the differential order into account.
We have\,:
\begin{eqnarray*}
\frac{1}{\hbar}
&:=&
\frac{1}{-\hspace{-1.45 ex}g}
+\sum_{1 \leq r\;\in\; \mi{odd}}\! {\doH}_r(\gbar)
\;=\;
\frac{1}{-\hspace{-1.45 ex}g}
+\sum_{1 \leq r\;\in\; \mi{odd}}\; \partial^{1-r}\;{\doD}_r(\gbar)
\hspace{2.ex } (\mi{``compact"})
\\[1.ex]
\frac{1}{\hbar}
&:=&
\frac{1}{\gbar}
+\sum_{{{1 \leq r\;\in\; \mi{odd}}\atop{\frac{1}{2}(r+1) \leq s \leq r}}}\!
{\doH}_{r,s}(\gbar)
\;=\;
\frac{1}{\gbar}
+\sum_{{{1 \leq r\;\in\; \mi{odd}}\atop{\frac{1}{2}(r+1) \leq s \leq r}}}\!
\partial^{-s}\;{\doD}_{r,s}(\gbar)
\hspace{2.ex } (\mi{``layered"})
\end{eqnarray*}
with $r$-linear, purely differential operators $^\ast\!{D}_{r},{D}_{r,s}$ of the form
\begin{eqnarray*}
{\doD}_{r}(\gbar) \!\!&:=&\!\!\!\!\!
\sum_{{\sum n_i=r \atop \sum i\,n_i=r-1 }} {\mr{^\ast Mir}}^{n_0,n_1,\dots,n_{r-1}}
\;\;
 \prod_{0 \leq i \leq r-1} ({\gbar}^{(i)})^{n_i}
\hspace{10.ex } (\mi{``compact"})
\\[1.ex]
{\doD}_{r,s}(\gbar) \!\!&:=&\!\!
\sum_{{\sum n_i=r \atop \sum i\,n_i=s }}\;\; {\mr{ Mir}}^{n_0,n_1,\dots,n_{s}}\;\;
\;\;\;\;
 \prod_{0 \leq i \leq s} ({\gbar}^{(i)})^{n_i}
\hspace{12.ex } (\mi{``layered"})
\end{eqnarray*}
and connected by:
\begin{equation*}
\partial\;\;{\doD}_{r}(\gbar)
=
\sum_{{{\frac{1}{2}(r+1) \leq s \leq r}}}\!
\partial^{r+s}\;\;{\doD}_{r,s}(\gbar)
\hspace{5.ex}(\forall r \in \,\{1,3,5...\})
\end{equation*}
The $\beta$-dependence is of course hidden in the definition of the 
differential operators ${\doD}_{r},{\doD}_{r,s}$\,: cf \S4.5 {\it infra}. All the information about the {\it mir} transform is thus carried by the
two rational-valued, integer-indexed moulds $^\ast\mi{Mir}$ and $\mi{Mir}$.
%%%%%%%%%%%%%%%%%%%%%%%%%%%%%%%%%%%%%%%%%%%%%%%%%%%%%%%%%%%%%%%%%%%%%%%%%%%%%%%%%%%%%%%%%%%%
%%%%%%%%%%%%%%%%%%%%%%%%%%%%%%%%%%%%%%%%%%%%%%%%%%%%%%%%%%%%%%%%%%%%%%%%%%%%%%%%%%%%%%%%%%%%

\subsection{The {\it nir} transform.}
{\bf Integral expression of {\it nir}.}
\\
Starting from $f$ we define  $f^{\Uparrow\beta^\ast}$ and  
$f^{\uparrow\beta^\ast}$ as follows\,:
\begin{eqnarray}
\nonumber
f(x)&\,\,=&\sum_{k\ge \kappa}\;f_k\,x^k\hspace{15. ex} (\kappa \ge 1\;,\;f_\kappa\not=0)
\\   \label{a70}
f^{\Uparrow\beta^\ast}(n,\tau) &:=& \beta(\partial_\tau)\,f(\frac{\tau}{n})
\hspace{13.ex}\mi{with}\;\;\;\partial_\tau:=\frac{\partial}{\partial\tau}
\\   \label{a71}
 &:=& \sum_{k\ge \kappa} 
\,f_k\,\,n^{-k}\,\beta^\ast_{k}(\tau)\hspace{25.ex}
\\   \label{a72}
 &:=& 
\,f_\kappa\,\frac{n^{-\kappa}\,\tau^{\kappa+1} }{\kappa+1}
+f^{\uparrow\beta^\ast}(n,\tau)\hspace{10.ex}
\end{eqnarray}
These definitions apply in the standard and non-standard cases alike.
Recall that in the standard case 
$
\beta^\ast_{k}(\tau)=\frac{B_{k+1}(\tau+\frac{1}{2})}{k+1}
$
is an even/\!/odd function of $\tau$ for $k$ odd/\!/even,
with leading term $\frac{\tau^{k+1}}{k+1}$. 
\\

The {\it nir}-transform $f\mapsto h$ is then defined as follows\,:
\begin{eqnarray}   \label{a73}
h(\nu) \!\!&=&\!\! \frac{1}{2\pi i}\,
\int_{c-i\infty}^{c+i\infty}\exp(n\,\nu)\,\frac{dn}{n}\,
\int_0^{\infty}\, \exp^\#(-f^{\Uparrow\beta^\ast}(n,\tau))\,d\tau
\\ \nonumber
\!&=&\!\! \frac{1}{2\pi i}\,
\int_{c-i\infty}^{c+i\infty}\!\exp(n\,\nu)\,\frac{dn}{n}
\int_0^{\infty}\!
\exp\!\big(\!-\!f_\kappa\, \frac{n^{-\kappa}\,\tau^{\kappa+1}}{\kappa+1}\big)
\exp_{\#}\!\big(-f^{\uparrow\beta^\ast}(n,\tau)\big)\,d\tau
\end{eqnarray}
where $\exp_{\#}(X)$ (resp. $\exp^{\#}(X)$) denotes
 the exponential  expanded as a power series of $X$ (resp. of $X$ minus its leading term). Here, we first perform 
term-by-term, ramified Laplace integration in $\tau$\,:
$$
\int_0^{\infty}\!
\exp\!\big(\!-\!f_\kappa\, \frac{n^{-\kappa}\,\tau^{\kappa+1}}{\kappa+1}\big)\,
\tau^{p}d\tau
=
{(f_\kappa)}^{-\frac{p+1}{\kappa+1}}\;\;
{(\kappa\!+\!1)}^{\frac{p+1}{\kappa+1}-1}\;\;
\Gamma\Big(\frac{p\!+\!1}{\kappa\!+\!1}\Big)\;\;
{n^{\frac{\kappa (p+1)}{\kappa+1}}}
$$
(with the main determination of ${(f_\kappa)}^{-\frac{p}{\kappa+1}}$
when $\Re(f_\kappa)>0$)
and then term-by-term (upper) Borel integration in $n$\,:
$$
\frac{1}{2\pi i}\;\int_{c-i\infty}^{c-i\infty}e^{n\nu}\;n^{-q}\;\;\frac{dn}{n}
\;=\;\frac{\nu^q}{q!}
$$ 
\begin{lemma}[The {\it nir}-transform preserves convergence]{\hspace{3.ex}}\\
Starting from a (locally) convergent $f$, the $\mg{\tau}$-integration in (\ref{a73}) 
usually destroys convergence,
but the subsequent  $\mr{n}$-integration always restores it. This holds not only in the standard
case, but also in the non-standard one, provided $\beta(\tau)$ has positive
convergence radius.
\end{lemma}
%%%%%%%%%%%%%%%%%%%%%%%%%%%%%%%%%%%%%%%%%%%%%%%%%%%%%%%%%%%%%%%%%%%%%%%%%%%%%%%%%%%%%%%%%%%%
%%%%%%%%%%%%%%%%%%%%%%%%%%%%%%%%%%%%%%%%%%%%%%%%%%%%%%%%%%%%%%%%%%%%%%%%%%%%%%%%%%%%%%%%%%%%
 
%%%%%%%%%%%%%%%%%%%%%%%%%%%%%%%%%%%%%%%%%%%%%%%%%%%%%%%%%%%%%%%%%%%%%%%%%%%%%%%%%%%%%%%%%%%%
\subsection{The reciprocation transform.}
Let us first examine what becomes of the nine-link chain in the simplest non-standard case,
i.e. with $\beta(\tau):=\tau^{-1}$.
\begin{lemma}[The simplest instance of {\it nir}-transform]{\hspace{5.ex}}\\
For the choice $\beta(\tau):=\tau^{-1}$, the pair $\{h,\hbar\}$ coincides with
the pair $\{g,\gbar\}$. In other words, $\mr{ mir}$ degenerates into the identity,
and  $\mr{ nir}$ essentially reduces to changing the germ $f^\ast$ into its functional
inverse $g^\ast$ (``reciprocation").
\end{lemma}
Since  $f=g$, the third column in the `long chain' becomes
redundant here, and the focus shifts to the
first two columns, to which we adjoin a new entry $\fbar:=1/f$ for the sake of symmetry.
Lagrange's classical inversion formula fittingly describes the involutions
 $f^\ast \leftrightarrow g^\ast$ and  $f \leftrightarrow g $,
and the simplest way of proving the above lemma is indeed by using Lagrange's formula.
On its own, however, that formula gives no direct information 
about the involution $\fbar
\leftrightarrow
\gbar
$ or the cross-correspondences $\fbar \leftrightarrow g $ and
$f \leftrightarrow \gbar $ which are highly relevant to an understanding
of the nine-link chain, including in the general case, i.e. for an arbitrary
$\beta(\tau)$. So let us first redraw the nine-link chain in 
the ``all-trivial case" $\{\beta(\tau)=\tau^{-1}\,,\,\kappa=0\,,\,f(0)=g(0)=1\}$ 
and then proceed with a description of the three afore-mentioned correspondences.
\[\begin{array}{cccccccccccc} 
\stackrel{\mathbf{nir}}{\rightarrow}   &\rightarrow &\rightarrow
&\rightarrow &\stackrel{\mathbf{nir}}{\rightarrow}\quad   & H
& &\;\;{\bf (Trivial\; case:} &\mg{\beta(\tau)}=\mg{\tau^{-1})}\hspace{5.ex}
 \\
\uparrow\quad\quad &                 &           &  \downarrow  
&\;\;\uparrow\stackrel{9}{}  
\\
\uparrow\quad\quad & \;\,f^*      &\stackrel{3}{\rightarrow} & \;\,g^*        
&   
\downarrow    &
\;g^{\,\prime} && f^\ast(x)=x+a^\ast(x) &g^\ast(y)=y+b^\ast(y) 
\\
 \uparrow\quad\quad&\;\;\uparrow \stackrel{2}{}&  &\;\;\downarrow\stackrel{4}{} &   
\downarrow    &
\;\;\uparrow\stackrel{8}{}\\
\stackrel{\mathbf{nir}} {\longleftarrow }  & f      &            & g  &
\stackrel{\mathbf{nir}}{\longleftarrow}\;\;
\stackrel{\mathbf{nir}}{\longrightarrow}  &
g&& f(x)=1+a(x) &g(y)=1+b(y) 
\\
    \;\;\;\;\;  \nearrow \stackrel{1}{}\! \!\!\!\!      &\;\;\uparrow \; &
&\;\;\downarrow\stackrel{5}{} &             
  &\;\;\uparrow\stackrel{7}{}\\
F \quad         & \fbar      &            & {\gbar} &      
\stackrel{6}{\longrightarrow\;\longrightarrow}    &
 {\gbar} 
&& \fbar(x)=1+\abar(x) &\gbar(y)=1+\bbar(y) 
\\
              & &  &&\stackrel{\mathbf{mir=id}}{} & 

\end{array}\]
\begin{lemma}[Three variants of Lagrange's inversion formula]\hspace{4.ex}\\
The entries $\mr{a,b,\abar,\bbar}$ in the above diagram are connected by\,:
\begin{eqnarray}   \label{a74}
a=\sum_{r\geq1}b_{<r>}&\mi{with}& b_{<r>}=\sum_{r_1+\dots+r_s=r}
M^{n_1,\dots,n_r}\,b^{[n_1,\dots,n_r]}
\\   \label{a75}
a=\sum_{r\geq1}\bbar_{\{r\}}&\mi{with}&
\bbar_{\{r\}}=\sum_{n_1+\dots+n_s=r} P^{n_1,\dots,n_s}\,\bbar^{[n_1,\dots,n_r]}
\\   \label{a76}
\abar=\sum_{n\geq1}\bbar_{[ r]}&\mi{with}& \bbar_{[r]}=\sum_{n_1+\dots+n_r=r}
Q^{n_1,\dots,n_r}\,\bbar^{[n_1,\dots,n_s]}
\end{eqnarray}
with differentially neutral \footnote{indeed, since $n_1+\dots+n_r=r$, we
integrate as many times as we differentiate.} and symmetral \footnote{meaning
that for any two sequences 
$\mathbf{n^\prime}=(n_i^\prime)$ and $\mathbf{n^{\prime\prime}}=(n_i^{\prime\prime})$,
we have the multiplication rule 
$\varphi^{\mathbf{[n^\prime]}}\varphi^{\mathbf{[n^{\prime\prime]}}}
\equiv\sum{\varphi^{\mathbf{[n]}}}$
with a sum running through all  $\mathbf{n}\in \mi{shuffle}
(\mathbf{n^{\prime}},\mathbf{n^{\prime\prime}})$.  }
 integro-differential expressions $ \varphi^{[n_1,\dots,n_s]} $
defined as follows\,:
\begin{equation}   \label{a77}
 \varphi^{[n_1,\dots,n_s]}(t):=\int_{0<t_1<..<t_r<t}
\varphi^{(n_1)}(t_1)\dots\varphi^{(n_r)}(t_r) \;dt_1\dots dt_r\hspace{10.ex}
\end{equation}
and with scalar moulds $M^\bu,P^\bu,Q^\bu$ easily inferred from the relations\,:
\begin{eqnarray}   \label{a78}
 \sum_{\|\bu\|=r} M^{\bu}\,b^{[\bu]}&=& 
\frac{(-1)^r}{r!}\partial^r (I\,b)^r
\\   \label{a79}
 \sum_{\|\bu\|=r} P^{\bu}\,\bbar^{[\bu]}&=& 
\partial_R\bbar\,I_L\dots \partial_R\bbar\,I_L\;\;\quad\quad(r\;\;\mi{times})
\\   \label{a80}
1+\sum_{\bu}Q^\bu\bbar^{[\bu]} &=&\Big(1+\sum_{\bu}P^\bu\bbar^{[\bu]}\Big)^{-1}
\end{eqnarray}
\end{lemma}
{\bf Remark 1\,:} In (\ref{a78}), $\partial$ as usual stands for differentiation  and $I=\partial^{-1}$
for integration from $0$. 
In (\ref{a79}), $\partial_R$ denotes the differentiation operator {\it 
acting on everything
to its right} and
$I_L=\partial_L^{-1}$ denotes the integration operator (with integration starting, again, from $0$) {\it acting on
everything to its left}.\\

Thus we find\,:
\begin{eqnarray*}
a &=& b_{<1>}+b_{<2>}+b_{<3>}+\dots
\\
b_{<1>} &=& -b^{[1]}
\\
b_{<2>} &=& +b^{[0,2]}+b^{[1,1]}
\\
b_{<3>} &=& -b^{[0,0,3]}-4\,b^{[0,1,2]}-4\,b^{[1,0,2]}-3\,b^{[0,2,1]}-15\,b^{[1,1,1]}
\\
\dots &&
\end{eqnarray*}
\begin{eqnarray*}
a &=& \bbar_{\{1\}}+\bbar_{\{2\}}+\bbar_{\{3\}}+\dots
\\
\bbar_{\{1\}} &=& +\bbar^{[1]}
\\
\bbar_{\{2\}} &=& +\bbar^{[0,2]}+\bbar^{[1,1]}
\\
\bbar_{\{3\}} &=& +\bbar^{[0,0,3]}+2\,\bbar^{[0,1,2]}+\bbar^{[1,0,2]}
+\bbar^{[0,2,1]}+\bbar^{[1,1,1]}\hspace{3.6 ex}
\\
\dots &&
\end{eqnarray*}
\begin{eqnarray*}
\abar &=& \bbar_{[1]}+\bbar_{[2]}+\bbar_{[3]}+\dots
\\
\bbar_{[1]} &=&  -\bbar^{[1]}
\\
\bbar_{[2]} &=& +\bbar^{[1,1]}-\bbar^{[0,2]}
\\
\bbar_{[3]} &=& -\bbar^{[1,1,1]}-\bbar^{[0,0,3]}+\bbar^{[0,2,1]}+\bbar^{[1,0,2]}\hspace{7.ex}
\\
\dots &&
\end{eqnarray*}
{\bf Remark 2\,:}
The coefficients $M^\bu,P^\bu,Q^\bu$ verify the following identities,
all of which are elementary, save for the last one (involving $\sum |Q^\bu|$)\,:
\begin{eqnarray}  \label{a81}
\sum_{\|\bu\|=r} (-1)^r\,M^\bu &=&\sum_{\|\bu\|=r} |M^\bu |=r^r
\\  \label{a82}
\sum_{\|\bu\|=r}  P^\bu &=&\sum_{\|\bu\|=r} |P^\bu |= r!
\\  \label{a83}
\sum_{\|\bu\|=r\geq 2}  Q^\bu =0\;&,&\; \sum_{\|\bu\|=r\geq 2} |Q^\bu
|=\frac{(2\,r-1)!}{(r-1)!\,r!}
\end{eqnarray}
\\
{\bf Remark 3\,:} $a$ in terms of $b$ is an elementary consequence of Lagrange's formula
for functional inversion,
but
$a$ in terms of $\bbar$ and $\abar$ in terms of $\bbar$ are not.
\\
\\
{\bf Remark 4\,:}
The formulas (\ref{a74}) through (\ref{a76}) involve only {\it sublinear} sequences
$\mathbf{n}=\{n_1,\dots,n_r\,;\,n_i\geq 0\}$, i.e. sequences verifying\,:
\begin{equation}  \label{a84}
\quad n_1+\dots+n_i\leq i\quad \quad \forall i \quad \quad \mi{and}\quad  
n_1+\dots+n_r=r
\end{equation}
The number of such series is exactly $\frac{(2\,r)!}{r!(r+1)!}$ 
(Catalan number), which puts them in one-to-one correspondence
with $r$-node binary trees.  
Moreover, these sublinear sequences are stable under  {\it
shuffling} and this establishes a link with the `classical product' on
binary trees\footnote{In fact, that product is of recent introduction\,: see Loday, Ronco, Novelli, Thibon, Hivert in [LR] and [HNT].}.
\\  
\\
{\bf Remark 5\,:}
The various $\varphi^{[n_1,\dots,n_r]}$, even for sublinear sequences $[n_1,\dots,n_r]$, are not
linearly independent,
but this does not detract from
the canonicity of the expansions in (\ref{a74}),(\ref{a75}),(\ref{a76}) because the induction
rules (\ref{a78}),(\ref{a79}),(\ref{a80}) behind the definition of $M^\bu,P^\bu,Q^\bu$ unambiguously
define a privileged set of coefficients.

%%%%%%%%%%%%%%%%%%%%%%%%%%%%%%%%%%%%%%%%%%%%%%%%%%%%%%%%%%%%%%%%%%%%%%%%%%%%%%%%%%%%%%%%%%%%
%%%%%%%%%%%%%%%%%%%%%%%%%%%%%%%%%%%%%%%%%%%%%%%%%%%%%%%%%%%%%%%%%%%%%%%%%%%%%%%%%%%%%%%%%%%%

\subsection{The {\it mir} transform.} 
\begin{lemma}[Formula for {\it mir} in the standard case]\hspace{0.ex}\\
The {\it mir} transforms $\gbar\mapsto \hbar $  is explicitely given by\,:
\begin{equation}  \label{a85}
\frac{1}{\hbar(\nu)}=\Big[\frac{1}{\gbar(\nu)}\,
\exp_{\#}\Big(-\sum_{r\geq 1}\beta_r\,I^r\, ({\gbar}(\nu)\partial_{\nu})^r \gbar(\nu)
\Big)\Big]_{I=\partial_{\nu}^{-1}}
\end{equation}
\end{lemma}
Mind the proper sequence of operations\,:
\\ 
\--- {\it first}, we expand the blocks 
$({\gbar}(\nu)\partial_{\nu})^r \gbar(\nu)$.
\\  
\--- {\it second}, we expand  $\exp_{\#}(\dots)$, which involves taking the suitable powers of the formal variable $I$ 
(with $``I"$ standing for {\it ``integration"}).
\\
\--- {\it third }, we divide by $\gbar(\nu)$. 
\\
\--- {\it fourth}, we move each $I^r$ to the left-most
position\footnote{this means that all the powers of $\gbar,\gbar^\prime,\gbar^{\prime\prime}$
etc must be put to the right of $I^r$}.\\
\--- {\it fifth}, we replace each $I^r$ by the operator $\partial_{\nu}^{-r}$ which 
stands for $n$ successive formal integrations from $0$ to $\nu$.  
\\
\--- {\it sixth}, we carry
out these integrations.
\begin{lemma}[The integro-differential components ${I\!D}_{r,s} $ of {\it
mir}]\hspace{0.ex}\\ The $\mr{mir}$ functional admits a canonical expansion\,:
\begin{eqnarray}  \label{a86}
\frac{1}{\hbar}
&:=&
\frac{1}{\gbar}
+\sum_{{{1 \leq r\;\in\; \mi{odd}}\atop{\frac{1}{2}(r-1) \leq s \leq r}}}\!
{\doH}_{r,s}(\gbar)
\;=\;
\frac{1}{\gbar}
+\sum_{{{1 \leq r\;\in\; \mi{odd}}\atop{\frac{1}{2}(r-1) \leq s \leq r}}}\!
\partial^{-s}\;{\doD}_{r,s}(\gbar)
\hspace{2.ex } \hspace{10.ex}
\end{eqnarray}
with $r$-linear differential operators ${\doD}_{r,s}$ of total order $d$\,:
\begin{eqnarray}\label{a87}
{\doD}_{r,s}(\gbar) \!\!&:=&\!\!
\sum_{{\sum n_i=r \atop \sum i\,n_i=s }}\;\; {\mr{ Mir}}^{n_0,n_1,\dots,n_{s}}\;\;
\;\;\;\;
 \prod_{0 \leq i \leq s} ({\gbar}^{(i)})^{n_i}
\hspace{12.ex } \hspace{10.ex}
\end{eqnarray}
and coefficients\,
$
{\mr{ Mir}}^{n_0,n_1,\dots,n_{s}}
\in \frac{1}{s!}\;\doZ[\beta_1,\beta_2,\beta_3,\dots]
$
which are themselves homogeneous of ``degree" $\mr{r+1}$ and ``order" $\mr{s}$ if to each
$\beta_i$ we assign the ``degree" $\mr{i+1}$ and ``order" $\mr{i}$.
\end{lemma} 
For the standard choice $\beta(\tau):=\frac{1}{e^{\tau/2}-e^{-\tau/2}}$, we have
$0=\beta_{2}=\beta_{4}=\dots$, and so we get only integro-differential components
${I\!D}_{r,s}$ which have all {\it odd} degrees $r=1,3,5\dots$. Thus\,:
\begin{eqnarray*}
\doD_{1,1}&=&+\frac{1}{24}\big( \gbar^{\prime} \big)
\\[1.0 ex]
\doD_{3,2}&=&+\frac{1}{1152}\big(\gbar\, \gbar^{\prime\,2} \big)
\\[1.0 ex]
\doD_{3,3}&=&-\frac{7}{5760}\big(
\gbar^{\prime\,3}
+\gbar^{2}\gbar^{\prime\prime\prime}
+4\,\gbar\gbar^\prime\gbar^{\prime\prime} 
\big)
\\[1.0 ex]
\doD_{5,3}&=&+\frac{1}{82944}\big(\gbar^2\gbar^{\prime\,3} \big)
\\[1.0 ex]
\doD_{5,4}&=&-\frac{7}{138240}\big(
\gbar\gbar^{\prime\,4}+\gbar^3\gbar^\prime\gbar^{\prime\prime\prime}
+4\,\gbar^2\gbar^{\prime\,3}
 \big)
\\[1.0 ex]
\doD_{5,5}&=&+\frac{31}{967680}\big(
\gbar^{\prime\, 5}+\gbar^4\gbar^{(5)}+11\,\gbar^3\gbar^\prime\gbar^{(4)}
+32\,\gbar^2\gbar^{\prime\,2}\gbar^{\prime\prime\prime}
\\[1.0 ex]
&&\hspace{10.ex}
+15\,\gbar^3\gbar^{\prime\prime}\gbar^{\prime\prime\prime}
+26\,\gbar\gbar^{\prime\,3}\gbar^{\prime\prime}
+34\,\gbar^2\gbar^{\prime}\gbar^{\prime\prime\,2}
 \big)
\end{eqnarray*}

However the {\it mir} formula has a wider range\,:
\begin{lemma}[Formula for {\it mir} in the non-standard case]\hspace{0.ex}\\
The  formula (\ref{a85}) and (\ref{a86}) for $\mr{mir}$
remains valid if we replace $\beta(\tau):=\frac{1}{e^{\tau/2}-e^{-\tau/2}}$ by any series
of the form
$\beta(\tau):=\sum_{n\geq-1}\;\beta_n\;\tau^n $  subject to
$\beta_{-1}=1,\beta_{0}=0$. The even-indexed coefficients $\beta_{2n}$ need not vanish.
When they don't, the expansion (\ref{a86})
may involves homogeneous components
${\doI\doH}_{r,s}$ of any degree $\mr{r}$, odd or even.
\end{lemma}
Dropping the condition $\beta_{-1}=1$ would bring about only minimal changes, but allowing a
non-vanishing $\beta_0$ would deeply alter and complicate the shape of the {\it mir}
transform\,: it would cease to be
a purely integro-differential functional. We must therefore be thankful for the parity 
phenomenon
(see \S4.1 {\it supra}) responsible for the occurence, in the {\it nir} integral, of the
Bernoulli
polynomials with shift $1/2$ rather than $1$.
\begin{lemma}[Alternative interpretation for the  {\it mir} formula]\hspace{0.ex}\\
The procedure implicit in formula (\ref{a85}) can be rephrased as follows\,:
\\ [1.ex]
\textup{(i)}\;
Form $h(w,y):=\sum_{r\geq 1}\frac{w^r}{r!}\,(\gbar(y)\partial_y)^r.y $
\\[1.ex]
\textup{(ii)}\;
 Form $k(w,y):=\sum_{r\geq 1}\beta_r\frac{w^r}{r!}\,
(\gbar(y)\partial_y)^r.\gbar(y) $
\\[1.ex]
\textup{(iii)}\;
 Interpret $\gbar(y)\partial_y$\; as an infinitesimal generator and\, 
$h^{\circ w}(y)=h(w,y)=g^\ast(w\!+\!f^\ast(y))$ as the corresponding
group of iterates\,: 
$ h^{\circ w_1}\circ h^{\circ w_2}=h^{\circ (w_1+w_2)}$.
\\[1.ex]
\textup{(iv)}\; 
Interpret $k(w,y)$\;\, as the Hadamard product, with repect to the  w variable, of\;
$\beta(w)$
\, and \, $\partial_w h(w,y)$.
\\[1.ex]
\textup{(v)}\; 
Calculate the convolution exponential\;\, $K(w,y):=\exp_{\star}(-k(w,y)) $
relative to the unit-preserving convolution $\star$ acting on the $w$ variable.
\\[1.ex]
\textup{(vi)}\; Integrate $\int_0^\nu K(\nu-\nu_1,\nu_1)\,
(\gbar(\nu_1))^{-1}\,d\nu_1 =:\ell(\nu)$. $ \Box$
\end{lemma}

%%%%%%%%%%%%%%%%%%%%%%%%%%%%%%%%%%%%%%%%%%%%%%%%%%%%%%%%%%%%%%%%%%%%%%%%%%%%%%%%%%%%%%%%%%%%
%%%%%%%%%%%%%%%%%%%%%%%%%%%%%%%%%%%%%%%%%%%%%%%%%%%%%%%%%%%%%%%%%%%%%%%%%%%%%%%%%%%%%%%%%%%%

\subsection{Translocation of the {\it nir} transform.}
If we set $  \eta:=\int_0^\epsilon f(x)\,dx$ and then wish to compare\,:
\\
(i)\; $f(x)$ and its translates 
$^{\epsilon}\!f(x)=e^{\epsilon\partial_x}f(x)=f(x+\epsilon)$
\\
(ii)\; $h(\nu)$ and its translates
${^\eta}h(\nu)=e^{\eta\partial_\nu}h(\nu)=f(\nu+\eta)$\\
there are {\it a priori} four possibilities to choose from\,:
\[\begin{array}{ccccc}
\mi{choice\;1\,:}&
( e^{\eta\partial_\nu}\;\mr{nir}
-\mr{nir}\;e^{\epsilon\partial_x} )\,f 
&\mi{as\; a\; function\; of}&(\epsilon,f)
\\[1.ex]
\mi{choice\;2\,:}&
( e^{\eta\partial_\nu}\;\mr{nir}
-\mr{nir}\;e^{\epsilon\partial_x} )\,f 
&\mi{as\; a\; function\; of}&(\eta,f)
\\[1.ex]
\mi{choice\;3\,:}&
( \mr{nir}
-e^{-\eta\partial_\nu}\,\mr{nir}\;e^{\epsilon\partial_x} )\,f 
&\mi{as\; a\; function\; of}&(\epsilon,f)
\\[1.ex]
\mi{choice\;4\,:}&
(\mr{nir}
-e^{-\eta\partial_\nu}\,\mr{nir}\;e^{\epsilon\partial_x} )\,f 
&\mi{as\; a\; function\; of}&(\eta,f)
\\[1.ex]
\end{array}\]
In the event, however, the best option turns out to be choice 3. So let us define
the finite (resp. infinitesimal) increments $\nabla h $ (resp. $ \delta_m h$)
accordingly\,:
\begin{eqnarray}  \label{a88}
\nabla h(\epsilon,\nu)\!&=&\!\sum (\delta_m h)(\nu)\; {\epsilon^m}
\;:=\;\mr{nir}(f)(\nu)-\mr{nir}(^{\epsilon}\!f)(\nu\!-\!\eta)
\\[1.ex]  \label{a89}
\mi{with}\;\;&& ^{\epsilon}\!f(x):=f(x+\epsilon)\quad\mi{and} \quad
 \eta:=\int_0^\epsilon f(x)\,dx
\end{eqnarray}
Going back to \S4.3, we can calculate
$\mr{nir}(f)(\nu) $
by means of the familiar double integral (\ref{a73}), and then
$\mr{nir}(^{\epsilon}\!f)(\nu\!-\!\eta)$ by using that same double
integral, but after carrying out the substitutions\,:
\begin{eqnarray}  \label{a90}
\nu &\mapsto & \nu -\eta
=\nu-\sum_{k\geq0}f^{(k)}(0)\,\frac{\epsilon^{k+1}}{(k+1)!}
=\nu-\sum_{k\geq0}f_k\,\frac{\epsilon^{k+1}}{k+1}
\\[0.5 ex]  \label{a91}
f(x)&\mapsto &^{\epsilon}\!f(x)=\sum_{k_1,k_2\geq 0}\,
f_{k_1+k_2}\frac{(k_1+k_2)!}{k_1!k_2!}\,\epsilon^{k_1}\,x^{k_2}
\\[0.8 ex]  \label{a92}
f^{\Uparrow\beta}(n,\tau)&\mapsto &^{\epsilon}\!f^{\Uparrow\beta}(n,\tau)
=\sum_{k_1,k_2\geq
0}\,f_{k_1+k_2}\frac{(k_1+k_2)!}{k_1!k_2!}\,
\epsilon^{k_1}\,n^{-k_2}\,\beta^\ast_{k_2}(\tau)
\end{eqnarray}
Singling out the contribution of the various powers of $\epsilon$, we
see that each infinitesimal increment $\delta_m h(\nu)$
is, once again, given by the double {\it nir}-integral, after
multiplication of the integrand by an elementary factor $D_m(n,\tau)$. Thus\,:
\begin{eqnarray*}
D_1(n,\tau)\!\!&=&\!\! f_0\,n+(\partial f)^\dagger(n,\tau)
\\
D_2(n,\tau)\!\!&=&\!\!f_1\,n-f_0^2\,n^2-2\,f_0\,n\,(\partial f)^\dagger(n,\tau)
+ (\partial^2 f)^\dagger(n,\tau) 
-( (\partial f)^\dagger(n,\tau) )^2
\\
\mi{etc}\quad &&
\end{eqnarray*}
Massive cancellations occur, which wouldn't occur under any of the other choices
1, 2 or 4. We can then regroup all
the {\it infinitesimal} increments into one remarkably simple expressions for the
{\it finite} increment\,: 
\begin{lemma}[The finite increment $\nabla h$\;: compact expression.]\hspace{0.ex}\\
Like $\mr{ nir}$ itself, its finite increment is given by a double integral\,:
\begin{equation}  \label{a93}
\nabla h(\epsilon,\nu)=\frac{1}{2\pi i}\int_{c-i\infty}^{c+i\infty}e^{n\nu}\,\frac{dn}{n}
\int_0^{\epsilon\,
n}\exp_{\#}\!\Big(-\beta(\partial_\tau)f(\frac{\tau}{n})
\Big)\,d\tau
\end{equation}
but with truncated Laplace integral and with 
$\exp_{\#}$ instead of $\exp^{\#}$.
\end{lemma}
The presence in (\ref{a93}) of 
$\mi{exp}_{\#}$ instead of $\mi{exp}^{\#}$
means that we must now expand {\it everything} within $\mi{exp}$,
including the leading term (unlike in (\ref{a73})). So we no longer have proper Laplace
integration here. Still, due to the truncation
 $\int_0^{\epsilon n}(\dots)d\tau$ of the integration interval, the integral continues to make sense, 
at least term-by-term. 
Due to the form $\epsilon n$ of the upper bound, it yields 
infinitely many summands $n^{-s}$, with positive {\it and} negative $s$. 
However, the
second integration  $\int_{c-i\infty}^{c+i\infty}(\dots)\frac{dn}{n}$
kill off the  $n^{-s}$ with negative $s$, and turns those with positive $s$
into $\frac{\nu^s}{s!}$\,. If we correctly interpret and carefully execute the
above procedure, we are led to the following analytical expressions for
the increment\,:
\begin{lemma}[The finite increment $\nabla h$\;: analytical expression.]\hspace{0.ex}\\
We have\,:
\begin{eqnarray}  \label{a94}
\nabla h(\epsilon,\nu)&:=&\sum_{s\geq 1}\frac{(-1)^s}{s!}
\sum_{{ p_i\geq 0\,,\, p_i\geq q_i \atop q_i \geq -1\,,\, q_i\not=0}}
\frac{\epsilon^m}{m}\;\frac{\nu^n}{n!}\;
\prod_{i=1}^{i=s}\Big(f_{p_i}\beta_{q_i}\frac{p_i!}{(p_i-q_i)!}\Big)\quad\quad
\\[1.ex]  \label{a95}
\mi{with }&& m:=1+\sum_i p_i -\sum_i q_i\;\;\;,\;\;\; n:=-1+\sum_i q_i
\end{eqnarray}
Equivalently, we may write\,:
\begin{eqnarray}  \label{a96}
\nabla h(\epsilon,\nu)&:=&\sum_{{m\geq 1 \atop n\geq 1 }} \delta_{m,n}(f,\beta)\;
{\epsilon^m}\;{\nu^n}\;
\\
\mi{with}\;\;\;&& \nonumber
\\
\delta_{m,n}(f,\beta)&=&\sum_{s=1}^{m+n}\frac{(-1)^s}{m\,n!\,s!}
\sum_{{m_i\geq 0\,,\, m_i\geq -n_i \atop n_i\not=0\,,\,n_i\geq -1 }}
\sum_{{ \sum m_i=m-1 \atop \sum n_i=n+1}}
\prod_{i=1}^{i=s}\Big(f_{m_i+n_i}\,\beta_{n_i}\,\frac{(m_i+n_i)!}{m_i!} \Big)
\nonumber
\end{eqnarray}
\end{lemma} 
Let us now examine the {\it infinitesimal} increments $\delta_m h$ of (\ref{a88}). Their {\it analytical}
expression clearly follows from (\ref{a96}), but they also admit very useful
{\it compact} expressions. To write these down, we require two sets of power series,
the $f^{\sharp m}$ and their upper Borel transforms $\smi{f^{\sharp m}}$. These series
enter the $\tau$-expansion of $ f^{\Uparrow\beta}$\,:
\begin{equation}  \label{a97}
f^{\Uparrow\beta}(n,\tau)=
f^{\sharp 0}(n)+\tau \,f^{\sharp 1}(n)+\tau^2 \,f^{\sharp 2}(n)+\dots\hspace{19.ex}
\end{equation}
As a consequence, they depend bilinearly on the coefficients of $f$ and $\beta$\,:  
\begin{eqnarray*}
f^{\sharp 0}(n):=\sum_{0\leq p}p!f_p\beta_p \,n^{-p}\hspace{6.ex}&&
\smi{f^{\sharp 0}}(\nu):=\sum_{0\leq p}f_p\beta_p \,\nu^{p}
\hspace{10.ex}(r= 0) 
\\[1.ex]
f^{\sharp m}(n):=\sum_{m-1\leq p}\frac{p!}{m!}f_p\beta_{p-m} \,n^{-p}&&
\smi{f^{\sharp 0}}(\nu):=\sum_{m-1\leq p}\frac{p!}{m!}f_p\beta_{p-m}\, \nu^{p}
\hspace{2.ex}(r\geq 1)
\end{eqnarray*} 
We also require the `upper' variant $\overline{\ast}$ of the finite-path convolution\,:
\begin{eqnarray}  \label{a98}
(A\overline{\ast} B)(t)&:=& \int_0^t A(t-t_1)\,dB(t_1)=\int_0^t B(t-t_1)\,dA(t_1)
\\  \label{a99}
1\overline{\ast}1\equiv 1&,& \frac{(.)^p}{p!}\overline{\ast} \frac{(.)^q}{q!}\equiv
\frac{\;\;(.)^{p+q}}{(p+q)!}
\end{eqnarray}
along with the corresponding convolution
exponential $\mi{exp}_{\,\overline{\star}}$\,:
\begin{equation}  \label{a100}
\exp_{\,\overline{\ast}}{A}:= 1+A+\frac{1}{2} A
\overline{\ast} A+\frac{1}{6} A\overline{\ast} A\overline{\ast} A+\dots
\end{equation}
\begin{lemma}[Infinitesimal increments $\delta_m h$\;: compact expression.]
The infinitesimal increments $\delta_m h$, as defined by the $\epsilon$-expansion\\
$\nabla h(\epsilon,\nu)=\sum_{0\leq m}\epsilon^m\,(\delta_m h)(\nu) $\;, admit the compact
expressions\,:
\begin{eqnarray}  \label{a101}
\delta_1 h&=& \partial_{\nu}\exp_{\overline{\ast}}(-f^{\sharp0})
\\  \label{a102}
\delta_2 h&=&\frac{1}{2}\; \partial_{\nu}^2
\Big(\big(-f^{\sharp1}\big)\star\exp_{\star}(-f^{\sharp0})\Big)
\\  \label{a103}
\delta_3 h&=&\frac{1}{3}\; \partial_{\nu}^3
\Big(\big(-f^{\sharp2}+\frac{1}{2}\,(-f^{\sharp1}){\overline{\ast}}
(-f^{\sharp1})\big){\overline{\ast}}\exp_{\overline{\ast}}(-f^{\sharp0})\Big)
\\
\dots &&\nonumber \dots
\\  \label{a104}
\delta_m h&=&\frac{1}{m}\; \partial_{\nu}^m
\Big(\big(\sum_{{\sum i\,k_i=m-1 \atop i\,\geq 1}}
\prod\frac{(-f^{\sharp\,i})^{{\overline{\ast}} k_i}}{ k_i!}\big)
{\overline{\ast}}\exp_{{\overline{\ast}}}(-f^{\sharp0})\Big)
\quad\quad
\end{eqnarray}
\end{lemma}

\begin{lemma}[The increments in the non-standard case]\hspace{0.ex}\\
The above expressions for $\delta_m h$ and $\nabla h$
remain valid even if we replace $\beta(\tau):=\frac{1}{e^{\tau/2}-e^{-\tau/2}}$ by any series
of the form
$\beta(\tau):=\sum_{n\geq-1}\;\beta_n\;\tau^n $ \quad {subject\; only\; to}\quad
$\beta_{-1}=1,\beta_{0}=0$.
\end{lemma}
\begin{lemma}[Entireness of $\delta_m h$ and $\nabla h$]\hspace{0.ex}\\ 
For any polynomial or entire input $f$, each $\delta_m h(\nu)$ is an entire
function of $\nu$ and $\nabla h(\epsilon,\nu)$ is an entire function of
$(\epsilon,\nu)$. This holds not only for the standard choice
$\beta(\tau):=\frac{1}{e^{\tau/2}-e^{-\tau/2}}$ but also for any series
$\beta(\tau):=\sum_{n\geq-1}\;\beta_n\;\tau^n $ with positive convergence
radius\,\footnote{subject as usual to $\beta_{-1}=1,\beta_0=0$.}.
\end{lemma}
This extremely useful lemma actually results from a sharper
statement\,:
\begin{lemma}[$\nabla h$ bounded in terms of $f$ and $\beta$.]\hspace{0.ex}\\ 
If $f(x)\prec \frac{A}{1-ax}$ and $\beta(\tau)\prec \frac{B}{1-b\tau}$
then $ \nabla h(\epsilon,\nu)\prec 
\frac{Const}{1-2a\epsilon\nu}\exp\Big(\frac{2AB}{b}\frac{\epsilon}{(1-ab\nu)} \Big)$.
\end{lemma}
Here, of course, for any two power series $\{\varphi,\psi\}$, the
notation $\varphi \prec \psi$ is short-hand for ``$\psi$ dominates $\varphi$",
i.e. $|\varphi_n|\leq \psi_n\,\,\forall n$. Under the assumption
$f(x)\prec \frac{A}{1-ax}$ and $\beta(\tau)\prec \frac{B}{1-b\tau}$\, we get\,:
\begin{eqnarray*}
\smi{f^{\# 0}}(\nu)\prec K_0(\nu)&:=&\frac{AB}{1-ab\,\nu}
\\
\smi{f^{\# m}}(\nu)\prec K_m(\nu)&:=&\frac{AB}{ab}\;\frac{a^m}{m!}\;\frac{\nu^{m-1}}{1-ab\,\nu}
\end{eqnarray*}
After some easy majorisations, this leads to\,:
\begin{eqnarray*}
\delta_m h(\nu) &\prec& \partial_\nu^m \sum_{{1\leq s\leq m \atop m_1+\dots+ m_s=m }}
\frac{1}{s!}\big(K_{m_1}\overline{\ast}K_{m_2}\overline{\ast}
\dots K_{m_s} \big)(\nu)
\\[1.ex]
&\prec& \sum_{1\leq s\leq m}\frac{\mi{Const}}{s!}\,\Big(\frac{AB}{ab}\Big)^s\,
\frac{(2a)^r\nu^{r-s}}{(1-ab\nu)^s}
\end{eqnarray*}
and eventually to\,:
\begin{eqnarray*}
\nabla h(\epsilon,\nu)&\prec& \frac{Const}{s!}\,\Big(\frac{AB}{ab}\Big)^s\,
\frac{1}{(1-ab\nu)^s}\,\sum_{s\leq m} (2a\,\epsilon)^m\,\nu^{m-s}
\\[1.ex]
&\prec & 
\frac{Const}{1-2a\epsilon\nu}\exp\Big(\frac{2AB}{b}\frac{\epsilon}{(1-ab\nu)} \Big)
\end{eqnarray*}
In the standard case we may take $B=1,b=\frac{1}{2\pi}$ so that the bound becomes\,:
$$ \nabla h(\epsilon,\nu) \prec \frac{\mi{Cons}}{1-2a\,\epsilon\,\nu} 
\exp\Big(\frac{4\pi A\,\epsilon}{1-\frac{a}{2\pi}\,\nu}\Big)$$
Since {\it Const} is independent of $a,A$, this immediately implies
that $\nabla h(\epsilon,\nu)$ is bi-entire (in $\epsilon$ and $\nu$) if $f(x)$ is
entire in $x$.\footnote{ and {\it only} if $f$ is entire \-- but this part is harder to prove
and not required in practice.}
%%%%%%%%%%%%%%%%%%%%%%%%%%%%%%%%%%%%%%%%%%%%%%%%%%%%%%%%%%%%%%%%%%%%%%%%%%%%%%%%%%%%%%%%%%%%
%%%%%%%%%%%%%%%%%%%%%%%%%%%%%%%%%%%%%%%%%%%%%%%%%%%%%%%%%%%%%%%%%%%%%%%%%%%%%%%%%%%%%%%%%%%%
 
\subsection{Alternative factorisations of {\it nir}. The {\it lir} transform.}
{\bf The {\it nir} transform and its two factorisations.}
\\ \noindent
In some applications, two alternative factorisations of the {\it nir}-transform are preferable to the one corresponding to the
nine-link chain of \S4.2. Graphically\,:
\[\begin{array}{ccccccccccccc}
f&{\longrightarrow}&\!\!\stackrel{\bf{nir}}{\longrightarrow}&{\longrightarrow}\!\!& h
\\[3.ex]
f^\ast&\stackrel{\bf{rec}}{\longrightarrow}&     g^\ast          
 &  \stackrel{\bf{imir}}{\longrightarrow}               & h^\ast
&\quad\quad\quad\quad &(f^\ast:=\partial^{-1}f \!&\!,\!&\! g^\ast:=\partial^{-1}g )
\\
\uparrow&          & \downarrow      &               &  \downarrow 
\\
f&                &      g          &    \dots            & h
\\
   &          & \downarrow      &                     & \uparrow
\\
 &                &       \gbar          &\stackrel{\bf{mir}}{\longrightarrow}& \hbar
&\quad\quad\quad\quad &(\gbar:=1/g \!&\!,\!&\! \hbar:=1/g)
\\[3.ex]
f^\ast&   \stackrel{\bf{ilir}}{\longrightarrow}  &     q^\ast           & 
\stackrel{\bf{rec}}{\longrightarrow}      & h^\ast
&\quad\quad\quad\quad &(q^\ast:=\partial^{-1}q \!&\!,\!&\! h^\ast:=\partial^{-1}h )
\\
\uparrow&          & \uparrow      & & \downarrow
\\
f&    \stackrel{\bf{lir}}{\longrightarrow}&       q          &                & h
\end{array}\]
In the first alternative, we go by {\it imir} (``integral'' {\it mir}) from the indefinite integral $g^\ast$ of $g$ to the indefinite integral 
$h^\ast$ of $h$, rather than from $\gbar$ to $\hbar$. In the second alternative, the middle column $(g,g^\ast)$ gets replaced by
$(q,q^\ast)$ with $q^\ast$ denoting the functional inverse of $h^\ast$. In that last scenario, the non-elementary
 factor-transform becomes {\it lir} (from $f$ to $q$) or {\it ilir} (from $f^\ast$ to $q^\ast$). We get get for these transforms expansions
 similar to, but in some respects simpler than, the expansions for {\it mir}\,:
\begin{eqnarray*}
\mr{mir}\,:&\gbar\rightarrow \hbar\;\;\;\mi{with} & \frac{1}{\hbar}=\frac{1}{\gbar}
+\sum_{1\leq r \,\mi{odd}} \;\; {\doH_r}(\gbar)
\\
\mr{lir}\,:&f\rightarrow q\;\;\;\mi{with} & q = f
+\sum_{3 \leq r \,\mi{odd}} \;\; {{\doQ}_r}(f)
\end{eqnarray*}

\begin{eqnarray*}
\mr{imir}\,:& g^\ast \rightarrow h^\ast \;\;\;\mi{with} & h^\ast =g^\ast
+\sum_{1\leq r \,\mi{odd}} \;\; {\doI\doH_r}(\gbar)
\\
\mr{ilir}\,:&f^\ast \rightarrow q^\ast \;\;\;\mi{with} & q^\ast= f^\ast
+\sum_{3 \leq r \,\mi{odd}} \;\; {{\doI\doQ}_r}(f)
\end{eqnarray*}
Each term on the right-hand sides is a polynomial in the $f^{(i)}$ and the following integro-differential
expressions\,:
\begin{eqnarray}  \label{a105}
f_m^{(d)\{s_1,\dots,s_r\}}&:=&
(I_R \,{\scriptscriptstyle \bu} f)^{s-d}  
\,{\scriptscriptstyle \bu} \,
(f^{-1}\,f^{(s_1)}\dots f^{(s_r)})
\\
&\,=&   \label{a106}
(I_R \,{\scriptscriptstyle \bu} f)^{m-r} 
\,{\scriptscriptstyle \bu} I_R 
\,{\scriptscriptstyle \bu} \,
(f^{(s_1)}\dots f^{(s_r)})
\\
&\,=& I_R \,{\scriptscriptstyle \bu} f 
\,{\scriptscriptstyle \bu}I_R 
\,{\scriptscriptstyle \bu} f \dots I_R \,{\scriptscriptstyle \bu}
 f \,{\scriptscriptstyle \bu} I_R \,{\scriptscriptstyle \bu}\, 
(f^{(s_1)}\dots f^{(s_r)})\nonumber
\end{eqnarray}
with $d\geq -1,m\geq r, s_1,\dots, s_r\geq 1$ and
 $1\!+\!m\!+\!d\!=\!r\!+\!s$.
Here $I_R:=\partial^{-1}=\int_0^{\dots}$ denotes the integration operator that starts from 0
 and
acts on {\it everything standing on the right}. The `monomial' $f_m^{(d)\{s_1,\dots,s_r\}} $
has total degree $m$ (i.e. it is $m$-linear in $f$) and total differential order $d$.
The notation is slightly redundant since $1\!+\!m\!+\!d\!\equiv\!r\!+\!s\!\equiv\!\sum(1+s_i)$
but very convenient, since it makes it easy to check that each summand in the 
expression of $\doH_r(f)$ (resp. $\doI\doH_r(f)$)  has global degree $r$ and global order 
$0$ (resp. -1).
The operators ${\doI\doH}_r$ and ${\doI\doQ}_r$   are simpler and in a sense more basic than the
${\doH}_r$ and ${\doQ}_r$. 
\medskip

\noindent
{\bf Proof\,:} Let us write the two reciprocal (formal) functions $h^\ast$ (known)
and $q^\ast$ (unknown) as sums of a leading term plus a perturbation\,:
\begin{eqnarray*}
h^\ast(x) &=& g^\ast(x)+\doI\doH (x)
\\
q^\ast(x) &=& f^\ast(x)+\doI\doQ (x)
\end{eqnarray*}
The identity $\mi{id}=h^\ast\circ q^\ast$ may be expressed as\,:
\begin{eqnarray*}
\mi{id} &=& (g^\ast+\doI\doH) \circ (f^\ast+\doI\doQ)
\\
 &=& \mi{id}+\doI\doH\circ f^\ast
+\sum_{1\leq r}\;\frac{1}{r!} \;(\doI\doQ)^r\; 
 \Big(\partial^r (g^\ast+\doI\doH)\Big)\circ f^\ast
\end{eqnarray*}
But $h^\ast$ may be written as 
\begin{equation}  \label{a107}
h^\ast(x) = (x+\doJ\doH)\circ g^\ast(x)
\end{equation}
and the identity  $\mi{id}=h^\ast\circ q^\ast$ now becomes\,:
\begin{equation}  \label{a108}
0=\doJ\doH+\sum_{1\leq r}\frac{1}{r!}\;(\doI\doQ)^r\; 
 (f^{-1}\!{\scriptscriptstyle \bu}\;\partial)^r 
 {\scriptscriptstyle \bu}\,(x+\doJ\doH) 
\end{equation}
The benefit from changing $\doI\doH$ into $\doJ\doH$ is that we are now handling direct
functions of $f$. Indeed, in view of the argument in \S4.2  we have\,:
\begin{equation}  \label{a109}
\doJ\doH = \sum_{1\leq r, 1\leq s_i}(-1)^r\Big(\prod_{i=1}^{i=r}\beta_i\Big)
(I_R{\scriptscriptstyle \bu}\,f)^{1+\sum s_i}
{\scriptscriptstyle \bu}\,\Big( f^{-1}\,\prod_{i=1}^{i=r}f^{s_i}\Big)
\end{equation}
The right-hand side turns out to be a linear combination of monomials (\ref{a105}) of order $d=-1$\,:
\begin{equation}  \label{a110}
\doJ\doH = \sum_{1\leq r, 1\leq s_i}(-1)^r\Big(\prod_{i=1}^{i=r}\beta_i\Big)\;
f^{(-1)\{s_1,\dots,s_r\}}_{r+\sum s_i}
\end{equation}
If we now adduce the obvious rules for differentiating these monomials\,:
\begin{eqnarray*}
(f^{-1}{\scriptscriptstyle \bu}\,\partial)^\delta\;
{\scriptscriptstyle \bu}\,f^{(d)\{s_1,\dots,s_r\}}_m
\!&\!=\!&\! 
f^{(d+\delta)\{s_1,\dots,s_r\}}_{m-\delta}  \hspace{20.5 ex}(\mi{if}\; \delta\! \leq\! m\!-\!r)
\\
\!&\!=\!&\! 
f^{-1}\,f^{(s_1)}\dots f^{(s_1)}  \hspace{17.ex}(\mi{if}\; \delta\! =\!1\!+\!m\!-\!r)
\\
\!&\!=\!&\! 
(f^{-1}\!{\scriptscriptstyle \bu}\,\partial)^{\delta+r-m-1}\!\!
{\scriptscriptstyle \bu}\,
f^{-1}\,f^{(s_1)}\dots f^{(s_1)}  
\hspace{1.ex}(\mi{if}\, \delta\! \geq 2\!+\!m\!-\!r)
\end{eqnarray*} 
we see at once that the identity (\ref{a108}) yields an inductive rule for calculating,
for each $m$, the $m$-linear part ${\doI\doQ}_m$ of $\doI\doQ$. At the same time, 
it shows that
any such ${\doI\doQ}_m$  will be exactly of global differential order -1,
and {\it a priori} expressible as a polynomial in $f^{-1},f,f^{(1)},f^{(2)},f^{(3)}\dots $
and finitely many monomials $ f^{(\delta)\{s_1,\dots,s_\rho \}}_\mu $. The only point left to check
is the {\it non-occurence} of negative powers of $f$, which would seem to result from the
above differentiation rules, but actually cancel out in the end result.
%%%%%%%%%%%%%%%%%%%%%%%%%%%%%%%%%%%%%%%%%%%%%%%%%%%%%%%%%%%%%%%%%%%%%%%%%%%%%%%%%%%%%%%%%%%%
%%%%%%%%%%%%%%%%%%%%%%%%%%%%%%%%%%%%%%%%%%%%%%%%%%%%%%%%%%%%%%%%%%%%%%%%%%%%%%%%%%%%%%%%%%%%
 
\subsection{Application\,: kernel of the {\it nir} transform.}
For any input $f$ of the form $p\,\mi{log}(x)+\mi{Reg}(x)$ with $p\in\doZ$ and $\mi{Reg}$ a regular analytic germ, the image $h$ of $f$
under {\it nir} is also a regular analytic germ\,: 
\begin{equation*}
 \mi{nir}  \quad : \;\; f(x)=p\,\log(x)+\mr{Reg_{1}}(x) \; \mapsto \mr{Reg_{2}}(x) 
 \end{equation*} 
The singular part of $h$, which alone has intrinsic significance, is thus 0. In other words, germs $f$ with logarithmic singularities
that are {\it entire} multiples of $\mi{log}(x)$ belong to the kernel of {\it nir} and produce {\it no inner generators}. This important and totally non-trivial fact is essential
when it comes to describing the inner algebra of SP series $j_F$ constructed
 from a meromorphic $F$. It may be proven (see [SS1]) either by using the alternative factorisations of the {\it nir} transform mentioned in the preceding
subsection, or by using an exceptional generator $f(x-x_{0})$ with base-point $x_{0}$ arbitrarily close to 0. An alternative proof, valid in
the special case when $\mi{Reg_{1}}=0$ and relying on the existence in that case of a simple ODE for the {\it nir}-transform, shall be given in 
\S6.6-7 below.

%%%%%%%%%%%%%%%%%%%%%%%%%%%%%%%%%%%%%%%%%%%%%%%%%%%%%%%%%%%%%%%%%%%%%%%%%%%%%%%%%%%%%%%%%%%%
%%%%%%%%%%%%%%%%%%%%%%%%%%%%%%%%%%%%%%%%%%%%%%%%%%%%%%%%%%%%%%%%%%%%%%%%%%%%%%%%%%%%%%%%%%%%
 
\subsection{Comparing/extending/inverting {\it nir} and {\it mir}.}

%%%%%%%%%%%%%%%%%%%%%%%%%%%%%%%%%%%%%%%%%%%%%%%%%%%%%%%%%%%%%%%%%%%%%%%%%%%%%%%%%%%%%%%%%%%%
%%%%%%%%%%%%%%%%%%%%%%%%%%%%%%%%%%%%%%%%%%%%%%%%%%%%%%%%%%%%%%%%%%%%%%%%%%%%%%%%%%%%%%%%%%%%

\begin{lemma}[The case of generalised power-series $f$]\hspace{0.ex}\\
The $\mr{nir}$ transform can be extended to generalised power series
\begin{equation}  \label{a111}
f(x):=\sum_{k_i\geq m} f_{k_i}\,x^{k_i} \quad
 \big(\,k_i\uparrow +\infty\;;\; k_i\in \doR \dot{-}\{-1\} \big)
\end{equation}
in a consistent manner (i.e. one that 
agrees with $\mr{mir}$ and ensures that  $\ell(\nu)$ 
converges whenever $f(x)$ does)  
by replacing in the double {\textup nir}-integral (\ref{a73}) the polynomials $\beta^\ast_k(\tau)$ by
the Laurent-type series\,:
\begin{equation}  \label{a112}
\beta_k^\ast(\tau)
:=\sum_{s=-1}^{+\infty}\,\beta_k\,\tau^{k-s}\frac{\Gamma(k+1)}{\Gamma(k+1-s)}
=\frac{\tau^{k+1}}{k+1}+\sum_{s=1}^{+\infty}\,(\dots)
\end{equation}
As usual, this applies both to the standard and non-standard \footnote{For
 non-standard choices, the series $\beta(\tau):=\sum \beta_s\,\tau^s$
has to be convergent if {\it nir} is to preserve convergence.} choices of
$\beta$.
\end{lemma}
We may also take advantage of the identity  
$f^{\Uparrow\beta}:=\beta(\partial_\tau)f(\frac{\tau}{n}) $
to {\it formally} extend the {\it nir}-transform to functions
$f$ derived from an $F$ with a zero/\!/pole of order $p$ at $x=0$\,:
$$
F(x)=e^{-f}=x^p\,e^{-w(x)}\quad \mi{with}\quad 
p\in\doZ^{\ast}\quad \mi{with}\quad  w(.) \;\mi{regular\; at \;0}
$$
That formal extension would read\,:
\begin{eqnarray*}
h(\nu) \!\!&\stackrel{\mr{formally}}{=}&\!\! \frac{1}{2\pi i}\,
\int_{c-i\infty}^{c+i\infty}\exp(n\,\nu)\,\frac{dn}{n}\,
\int_0^{\infty}\,n^{-p\,\tau}\exp\Big(p\,\lambda(\tau)-\beta(\partial_\tau)\,w(\frac{\tau}{n})\Big)\,d\tau
\\
\tilde{\lambda}(\tau)\!&=&\!\!\beta(\partial_\tau)(\log\tau)
=\tau\log\tau-\tau+\sum_{0\leq s}\beta_{2\,s+1}\;(2\,s)!\;\tau^{-2\,s-1}
\end{eqnarray*}
with $\lambda$ denoting the Borel-Laplace resummation of the divergent series
$\lambda$.
However, in the above formula for $h(\nu)$,
the first integration (in $\tau$) makes no sense at infinity\footnote{even when interpreted
term-by-term,
i.e. after expanding $\mi{exp}\!\big(\!-\!\beta(\partial_\tau)(w(\frac{\tau}{n})\big)$ .}
and one would have to exchange the order of integration (first $n$, then $\tau$),
among other things, to make sense of the formula and arrive at the correct result,
namely that the {\it nir}-transform turns functions
of the form $f(x)=p\log x+\mi{Reg}$ into {\it Reg}\footnote{
as long as $p\in\doZ$.}. In other words, {\it there is no {\it inner generator} attached to the
corresponding base point $x=0$.} But it would be difficult to turn the argument into a rigorous proof, and so
the best approaches remain the ones just outlined in the preceding subsection.

{\it Directly} extending {\it nir} to even more general test functions $f$ would be
possible, but increasingly difficult and of doubtful advantage. Extending {\it mir},
on the other hand, poses no difficulties. 
\begin{lemma}[Extending {\it mir}\,'\,s domain]\hspace{0.ex}\\
The $\mr{mir}$ transform $\gbar\rightarrow \hbar $ extends, formally \textup{and}
analytically, to general transserial inputs $\gbar$ of infinitesimal type, i.e.
$\gbar(y)=o(y)\,,\, y\sim 0$, and even to those with moderate growth
$\gbar(y)=O(y^{-\sigma})\,,\,\sigma>0$.
\end{lemma}
We face a similar situation when investigating the behaviour of $h(\nu)$ {\it over}
$\nu=\infty$\footnote{ 
under the change $\zeta=e^{\nu}$,
this behaviour at $-\!\infty$ in the $\nu$-plane translates into the
behaviour over 0 in the $\zeta$-plane. } for inputs $f$ of the form\,:
$$
f(x)=\mi{polynomial}(x)\quad\mi{or}\quad
\mi{polynomial}(x)+\sum_{i=1}^{i=N}p_i\log(x-x_i)\quad (p_i\in doZ^\ast)
$$
Then the (\ref{a86}) expansion for $h(\nu)$ still converges in some suitable
(ramified) neighbourhood of $\infty$ to some analytic germ, 
but the latter is no longer described by a power series
(or a Laurent series, as we might expect a infinity)
nor even by a (well-ordered) transseries,
but by a complex combination both kinds of infinitesimals\,:
small and large.\footnote{When re-interprented as a germ over $0$ in the
$\zeta$-plane, it typically produces an essential singularity there,
with Stokes phenomena and exponential growth or decrease, depending on the sector.}
\begin{lemma}[Inverting {\it mir}]\hspace{0.ex}\\
The $\mr{mir}$ transform admits a formal inverse $\mr{mir}^{-1}:
\hbar\rightarrow \gbar $\, that acts, not just formally but also 
analytically, on general transserial inputs $\hbar$ of infinitesimal type.
Like $\mr{mir}$, this inverse $\mr{mir}^{-1}$ admits well-defined integro-differential
components ${I\!D}_{r,s}$ of degree r and order s, but these are no longer
of the form $\partial^{-s}{D}_{r,s}$ with a neat separation of the differentiations
(coming first) and integrations (coming last).
\end{lemma}

%%%%%%%%%%%%%%%%%%%%%%%%%%%%%%%%%%%%%%%%%%%%%%%%%%%%%%%%%%%%%%%%%%%%%%%%%%%%%%%%%%%%%%%%%%%%
%%%%%%%%%%%%%%%%%%%%%%%%%%%%%%%%%%%%%%%%%%%%%%%%%%%%%%%%%%%%%%%%%%%%%%%%%%%%%%%%%%%%%%%%%%%%

\subsection{Parity relations.}
With the standard choice for $\beta$, we have the following parity relations
for the {\it nir}-transform\,:
\begin{eqnarray*}
F^{\,\bot}(x):=1/F(-x)\;\;\;,\;\;\;f^{\,\bot}(x):=-f(-x)  \hspace{7.ex}
&&\Longrightarrow 
\\
&& (\mi{tangency}\;\, \kappa=0)
\\
\mr{nir}(f^\bot)(\nu)=-\mr{nir}(f)(\nu)
\hspace{24.ex} 
&&(\mi{tangency}\;\, \kappa=0)
\\
\mr{nir}(f^\bot)\;\; \mi{and}\;\;\mr{nir}(f)\;\;\mi{unrelated}
\hspace{18.ex}
&& (\mi{tangency}\;\, \kappa\;\mi{even}\geq 2)
\\
\mr{nir}(f^\bot)(\nu)=-\mr{nir}(f)(\epsilon_\kappa\nu)\;\;\; \mi{with}\;\;\;
\epsilon_\kappa^{\frac{1}{\kappa+1}}=-1
\hspace{3.3 ex}
&& (\mi{tangency}\;\, \kappa\;\mi{odd}\geq 1)
\\
\Rightarrow h^{\bot}_{\frac{k}{\kappa+1}}=(-1)^{k-1}h_{\frac{k}{\kappa+1}}
\;\;\mi{with}\;\; \;:(f,f^\bot)\stackrel{\mr{nir}}{\mapsto} (h,h^\bot)
\hspace{2.ex}
&& (\mi{tangency}\;\, \kappa\;\mi{odd}\geq 1)
\end{eqnarray*}
For the {\it mir}-transform the parity relation 
doesn't depend on $\kappa$ and assumes the elementary form \,:
$$\mr{mir}(-\gbar)=-\mr{mir}(\gbar)\hspace{26.ex} $$ 
%%%%%%%%%%%%%%%%%%%%%%%%%%%%%%%%%%%%%%%%%%%%%%%%%%%%%%%%%%%%%%%%%%%%%%%%%%%%%%%%%%%%%%%%%%%%
%%%%%%%%%%%%%%%%%%%%%%%%%%%%%%%%%%%%%%%%%%%%%%%%%%%%%%%%%%%%%%%%%%%%%%%%%%%%%%%%%%%%%%%%%%%%

%\end{document}
%%%%%%%%%%%%%%%%%%%%%%%%%%%%%%%%%%%%%%%%%%%%%%%%%%%%%%%%%%%%%%%%%%%%%%%%%%%%%%%%%%%%%%%%%%%%
%%%%%%%%%%%%%%%%%%%%%%%%%%%%%%%%%%%%%%%%%%%%%%%%%%%%%%%%%%%%%%%%%%%%%%%%%%%%%%%%%%%%%%%%%%%%
%%%%%%%%%%%%%%%%%%%%%%%%%%%%%%%%%%%%%%%%%%%%%%%%%%%%%%%%%%%%%%%%%%%%%%%%%%%%%%%%%%%%%%%%%%%%
%%%%%%%%%%%%%%%%%%%%%%%%%%%%%%%%%%%%%%%%%%%%%%%%%%%%%%%%%%%%%%%%%%%%%%%%%%%%%%%%%%%%%%%%%%%%
%%%%%%%%%%%%%%%%%%%%%%%%%%%%%%%%%%%%%%%%%%%%%%%%%%%%%%%%%%%%%%%%%%%%%%%%%%%%%%%%%%%%%%%%%%%%
%%%%%%%%%%%%%%%%%%%%%%%%%%%%%%%%%%%%%%%%%%%%%%%%%%%%%%%%%%%%%%%%%%%%%%%%%%%%%%%%%%%%%%%%%%%%
%%%%%%%%%%%%%%%%%%%%%%%%%%%%%%%%%%%%%%%%%%%%%%%%%%%%%%%%%%%%%%%%%%%%%%%%%%%%%%%%%%%%%%%%%%%%
%%%%%%%%%%%%%%%%%%%%%%%%%%%%%%%%%%%%%%%%%%%%%%%%%%%%%%%%%%%%%%%%%%%%%%%%%%%%%%%%%%%%%%%%%%%%
%%%%%%%%%%%%%%%%%%%%%%%%%%%%%%%%%%%%%%%%%%%%%%%%%%%%%%%%%%%%%%%%%%%%%%%%%%%%%%%%%%%%%%%%%%%%

%

%%%%%%%%%%%%%%%%%%%%%%%%%%%%%%%%%%%%%%%%%%%%%%%%%%%%%%%%%%
%% Some resurgence properties of knot-related functions.
%%%%%%%%%%%%%%%%%%%%%%%%%%%%%%%%%%%%%%%%%%%%%%%%%%%%%%%%%%
 
%\documentclass[12pt,a4paper]{article}\input{SP_commands}\begin{document}

%%%%%%%%%%%%%%%%%%%%%%%%%%%%%%%%%%%%%%%%%%%%%%%%%%%%%%%%%%%%%%%%%%%%%%%%%%%%%%%%%%%%%%%%%%%%
%%%%%%%%%%%%%%%%%%%%%%%%%%%%%%%%%%%%%%%%%%%%%%%%%%%%%%%%%%%%%%%%%%%%%%%%%%%%%%%%%%%%%%%%%%%%
%%%%%%%%%%%%%%%%%%%%%%%%%%%%%%%%%%%%%%%%%%%%%%%%%%%%%%%%%%%%%%%%%%%%%%%%%%%%%%%%%%%%%%%%%%%%

\section{Outer generators.}
%%%%%%%%%%%%%%%%%%%%%%%%%%%%%%%%%%%%%%%%%%%%%%%%%%%%%%%%%%%%%%%%%%%%%%%%%%%%%%%%%%%%%%%%%%%%
%%%%%%%%%%%%%%%%%%%%%%%%%%%%%%%%%%%%%%%%%%%%%%%%%%%%%%%%%%%%%%%%%%%%%%%%%%%%%%%%%%%%%%%%%%%%
\subsection{Some heuristics.}
In the heuristical excursus at the beginning of the preceding section,
we had chosen the driving function $f$ such as to make the nearest
singularity an {\it inner generator}. We must now hone $f$ to ensure
that the nearest singularity be an {\it outer generator}. For maximal
simplicity, let us assume that\,:
\begin{equation}  \label{a113}
0<f(0)\leq +\infty\quad \mi{and}\quad 0<f(x)<+\infty\quad \mi{for}\quad 0<x\leq 1
\end{equation} 
Thus  $f^\ast(x):=\int_0^x f(x^\prime)dx^\prime$ will be $>0$ on
the whole interval $]0,1]$. 
Since we insist, as usual, on $F:=\exp(-f)$ being meromorphic,
(\ref{a113}) leaves but three possibilities\,:
\[\begin{array}{llllllllllll}
\mi{Case \;1} &:& 0=F(0) &;& f(x)=-p\,\log(x)+\sum_{k=0}^{k=\infty} f_k x^k
\\[1.5 ex]
\mi{Case \;2} &:& 0<F(0)<1 &;& f(x)=\sum_{k=0}^{k=\infty} f_k x^k\quad(f_0>0)
\\[1.5 ex]
\mi{Case \;3} &:& \;\;\;\quad F(0)=1 &;& f(x)=\sum_{k=\kappa}^{k=\infty} f_k
x^k\quad(f_\kappa >0\, ,\,\kappa\geq1)
\hspace{10.ex}
\end{array}\]
In all three cases, the nearest singularity of $j(\zeta)$  (cf (\ref{a1})) is located at $\zeta=1$
and reflects the $n$-asymptotics of the Taylor coefficients $J(n)$.
\begin{equation}  \label{a114}
J(n):=\sum_{m=\epsilon}^{m=n}\prod_{k=\epsilon}^{k=m}F(\frac{k}{n})
\quad \big(\epsilon \in\{0,1\}\big)
\quad \Longrightarrow \quad
\tilde{J}(n)\,:=\,
\sum_{k\geq 0}\,j_k\;n^{-k}
\end{equation}

{\it Case 1}\; is simplest.\footnote{Here, we must take $\epsilon=0$
to avoid an all-zero result.} By truncating the $\sum\prod$ expansion
at $m=m_0$, we get the exact values of all coefficients $j_k$ up to
$ k= p\, m_0$.
\\

{\it Case 2}\; corresponds to tangency 0. Here, finite truncations yield
only approximate values. To get the exact coefficients, we must harness the
full $\sum\prod $ expansion but we still end up with closed expressions
for each $j_k$\,.
\\

{\it Case 3}\;  corresponds to tangency $\kappa\geq 1$ 
case. Here, again, the
full $\sum\prod$ expansion must be taken into account, but the difference is that we 
now get coefficients
$j_k$
which, though exact, are no longer neatly expressible in terms of elementary functions.
%%%%%%%%%%%%%%%%%%%%%%%%%%%%%%%%%%%%%%%%%%%%%%%%%%%%%%%%%%%%%%%%%%%%%%%%%%%%%%%%%%%%%%%%%%%%
%%%%%%%%%%%%%%%%%%%%%%%%%%%%%%%%%%%%%%%%%%%%%%%%%%%%%%%%%%%%%%%%%%%%%%%%%%%%%%%%%%%%%%%%%%%%

\subsection{The short and long chains behind {\it nur/mur}.}
Let us now translate the above heuristics into precise (non-linear) functionals.
For case 1, the definition is straightforward\,:
\\[1.5 ex]
\noindent
{\bf The short, four-link chain\,:}
\begin{eqnarray*}
F\;\;\stackrel{1}{\longrightarrow}\;\;
k\;\;\stackrel{2}{\longrightarrow}\;\;
{h}\;\;\stackrel{3}{\longrightarrow}\;\; 
{h}^\prime\;\;\stackrel{3}{\longrightarrow}\;\; 
H\quad \quad \quad
h =\smi{\ell u}\;,\, 
h^\prime =\imi{\ell u}\;,\, 
H=\imi{Lu}
\end{eqnarray*}
{\bf Details of the four steps:}
\[\begin{array}{lllll}
F(x)&\!\!:=\!\!& F_1\,x+F_2\,x^2+F_3\,x^3+\dots \hspace{15.ex} &(\mi{converg^t})
\\
\downarrow \stackrel{1}{} &&
\\
k(n)&\!\!:=\!\!&\Big(
F(\frac{1}{n})+F(\frac{1}{n})F(\frac{2}{n})
+F(\frac{1}{n})F(\frac{2}{n})F(\frac{3}{n})+\dots )/I\!g_{F}(n)
\hspace{1.ex}
& (\mi{diverg^t})
\\
\| \; && 
\\
{k}(n)&\!\!:=\!\!&\sum_{1\le s} \;k_s\; \frac{1}{n^k} \hspace{30.ex}& (\mi{diverg^t})
\\
\downarrow \stackrel{2}{} &&
\\
{h}(\nu)&\!\!:=\!\!&\sum_{1\le s} \;k_s\; \frac{\nu^{s}}{s!} \hspace{30.ex}
&(\mi{converg^t})
\\
\downarrow \stackrel{3}{} &&
\\
{h}^\prime(\nu)&\!\!:=\!\!&\sum_{1\le s} \;k_s\; \frac{\nu^{s-1}}{(s-1)!} \hspace{30.ex}
&(\mi{converg^t})
\\
\downarrow \stackrel{4}{} &&
\\
H(\zeta)&\!\!:=\!\!& {h}^\prime(\log(1+\zeta))
\end{array}\]
Mark the effect of removing the ingress factor $I\!g_F$ after the first step.
 If
$$
F(x)=c_0\,x^d\,F_\ast(x)\quad\mi{with}\quad F_\ast(x)=1+\dots\in\doC\{x^2\}
\;(resp.\; \doC\{x\})
$$
then, according to the results of \S3,  removing  $I\!g_F$
 amounts to dividing $k(n)$ by 
$c_0^{-1/2}(1\pi n)^{d/2} $ and integrating $d/2$ times the functions
\footnote{
or more accurately $c_0^{1/2}(2\pi)^{-d/2}h(\nu) $
and  $c_0^{1/2}(2\pi)^{-d/2}h^\prime(\nu) $.} 
$h(\nu)$ or $h^\prime(\nu)$. The removal of the ingress factor thus has three
main effects\,:
\\

(i) as already pointed out, it makes the outer generators
independent of the ingress point\footnote{just as was the case with the
inner generators.}
\\

(ii) depending on the sign of $d$, it renders the singularities
smoother (for $d>0$) or less smooth (for $d<0$), in the $\nu$- or
$\zeta$-planes alike.
\\

(iii) depending on the parity of $d$, it leads in the Taylor
expansions of the minors $\smi{\ell u}\!\!(\nu):=h(\nu)$ and 
$\imi{\ell u}\!\!(\nu):=h^\prime(\nu)$ 
either to integral powers
of $\nu$ (for $d$ even) or to strictly semi-integral powers (for $d$ odd). 
This means that the corresponding majors 
$\sma{\ell u}$ and 
$\ima{\ell u}$ and, by way of consequence, the inner generators themselves,  will
carry {\it logarithmic} singularities (for $d$ even) or strictly semi-integral
powers (for $d$ odd).\footnote{
Removing the ingress factor has exactly the opposite effect on 
{\it inner generators}\,: these
generically carry semi-integral powers for $d$ even and logarithmic
singularities for $d$ odd.
}.
\\

Time now to deal with the cases 2 and 3 (i.e. $F(0)\not=0 $).
These cases lead to a nine-link 
chain quite similar to that which in \S4.2 did service for the {\it inner generators},
but with the key steps {\it nir} and {\it mir} significantly altered
into {\it nur} and {\it mur}\,:
\\[1.5 ex]
\noindent
{\bf The long,  nine-link chain:}
\[\begin{array}{ccccccccccccc}
\stackrel{\mathbf{nur}}{\rightarrow}   &\rightarrow &\rightarrow
&\rightarrow &\stackrel{\mathbf{nur}}{\rightarrow}\quad   & H 
&\quad\quad\quad& H&\!=\!&\imi{Lu}
\\
\uparrow\quad\quad &        &            &           &  \downarrow  
&\;\;\uparrow\stackrel{9}{}
\\
\uparrow\quad\quad & \;\,f^*      &\stackrel{3}{\rightarrow} & \;\,g^*        
&   
\downarrow    &
\;h^{\,\prime}
&\quad\quad\quad& h^\prime&\!=\!&\imi{lu}
 \\
 \uparrow\quad\quad&\;\;\uparrow \stackrel{2}{}&  &\;\;\downarrow\stackrel{4}{} &   
\downarrow    &
\;\;\uparrow\stackrel{8}{}\\
\stackrel{\mathbf{nur}} {\longleftarrow }  & f      &            & g  &\quad
\quad\stackrel{\mathbf{nur}}{\longrightarrow}  &
h &\quad\quad\quad& h&\!=\!&\smi{lu}
\\
              &\;\;\uparrow \stackrel{1}{} & &\;\;\downarrow\stackrel{5}{} &             
&\;\;\uparrow\stackrel{7}{}\\
              & F      &            & {\gbar} &      
\stackrel{6}{\longrightarrow\;\longrightarrow}    &
 {\hbar} \\
              & &  &&\mathbf{mur} & \\

\end{array}\]
{\bf Details of the nine steps:}
\[\begin{array}{cccccccc}
\stackrel{1}{\rightarrow}&:& \mi{precomposition} \; &:&
\; F\rightarrow f\;\;
&\mi{with}& \;\; f(x):=-\log F(x)
\\
\stackrel{2}{\rightarrow}&:& \mi{integration} \; &:&
\; f \rightarrow f^* \;\;
&\mi{with}&\;\; f^*(x):=\int_0^x f(x_0)\,dx_0
\\
\stackrel{3}{\rightarrow}&:& \mi{reciprocation} \;&:&
\; f^* \rightarrow g^* \;\;
&\mi{with}&\;\; f^*\circ g^*=\mi{id}
\\
\stackrel{4}{\rightarrow}&:& \mi{derivation} \;&:&
\; g^*\rightarrow g\;\;
&\mi{with}&\;\; g(y):=\frac{d}{dy}g^*(y)
\\
\stackrel{5}{\rightarrow}&:& \mi{inversion} \;&:&
\; g\rightarrow \gbar \;\;
&\mi{with}&\;\; {-\hspace{-1.5 ex}g}(y)\, := 1/g(y) 
\\
\stackrel{6}{\rightarrow}&:& \mg{mur}\; \mi{functional} \;&:&
\; {-\hspace{-1.5 ex}g} \;
\rightarrow \;
{ \hbar}\;\;
&&  \mi{see\;\S5.4\; infra}
\\
\stackrel{7}{\rightarrow}&:& \mi{inversion} \;&:&
\; 
{\hbar}\rightarrow h\;\;
&\mi{with}& \;
\; 
h(\nu):=1/{\hbar}(\nu)
\\
\stackrel{8}{\rightarrow}&:& \mi{derivation} \;&:&
\; h \rightarrow \;h^\prime
& &\;
\; 
h^\prime(\nu):=\frac{d}{d\nu}h(\nu)
\\
\stackrel{9}{\rightarrow}&:& \mi{postcomposition} \;&:&
\; 
h^\prime\rightarrow \; H
&\mi{with}&
\; H(\zeta):=h^\prime(\log(1+\zeta)
\\[1.ex]
\stackrel{2...7}{\rightarrow}&:& \mg{nur}\; \mi{functional} \;&:&
\; {g} \;
\rightarrow \;
{ h}\;\;
&\mi{with}&  \mi{see\;\S5.3\; infra}
\end{array}\]
%%%%%%%%%%%%%%%%%%%%%%%%%%%%%%%%%%%%%%%%%%%%%%%%%%%%%%%%%%%%%%%%%%%%%%%%%%%%%%%%%%%%%%%%%%%%
%%%%%%%%%%%%%%%%%%%%%%%%%%%%%%%%%%%%%%%%%%%%%%%%%%%%%%%%%%%%%%%%%%%%%%%%%%%%%%%%%%%%%%%%%%%%

\subsection{The {\it nur} transform.}
\noindent\\
{\bf Integral-serial expression of {\it nur}.}
\\

Starting from $f$ we define  
$f^{\Uparrow\frak{b}^\ast},f^{\uparrow\frak{b}^\ast}$ and 
$f^{\Uparrow\beta^\ast},f^{\uparrow\beta^\ast}$ 
as follows\,:
\begin{eqnarray*}
f(x)&\,\,=&\sum_{k\ge \kappa}\;f_k\,x^k\hspace{5. ex} (\kappa \ge 1\;,\;f_\kappa\not=0)
\\
f^{\Uparrow\frak{b}^\ast}(n,\tau) &:=&\!\! \frak{b}(\partial_\tau)f(\frac{\tau}{n})
=
\sum_{k\ge \kappa} 
\,f_k\,\,n^{-k}\,\frak{b}^\ast_{k}(\tau)\;=:\;
\,f_\kappa\,\frac{n^{-\kappa}\,\tau^{\kappa+1} }{\kappa+1}
+f^{\uparrow\frak{b}^\ast}(n,\tau)\hspace{10.ex}
\\
f^{\Uparrow\beta^\ast}(n,\tau) &:=& \!\! \beta(\partial_\tau)f(\frac{\tau}{n})
=
\sum_{k\ge \kappa} 
\,f_k\,\,n^{-k}\,\beta^\ast_{k}(\tau)\;=:\;
\,f_\kappa\,\frac{n^{-\kappa}\,\tau^{\kappa+1} }{\kappa+1}
+f^{\uparrow\beta^\ast}(n,\tau)\hspace{10.ex}
\\
\mi{with}\!\!\!&\mi{the}&\!\!\!\mi{ usual\; definitions\; in \; the\; standard\;case:} 
\\[0.6 ex]
\frak{b}^\ast_{k}(\tau)&:=&\frac{B_{k+1}(\tau+1)}{k+1}
\,=\,\beta^\ast_{k}(\tau+\frac{1}{2})
\\
\beta^\ast_{k}(\tau)&:=&\frac{B_{k+1}(\tau+\frac{1}{2})}{k+1}
\end{eqnarray*}
where $B_k$ stands for the $k^{th}$ Bernoulli polynomial.  Recall that
$\frak{b}^\ast_k(m)$ is a polynomial in $m$ of degree $k+1$, with leading
term $\frac{m^{k+1}}{k+1} $. 
Then the {\it nur}-transform is defined as follows\,:
\begin{eqnarray}  \label{a115}
\mi{nur}&:& f \mapsto h
\\[0.7 ex]  \label{a116}
h(\nu) \!\!&=&\!\! \frac{1}{2\pi i}\,
\int_{c-i\infty}^{c+i\infty}\exp(n\,\nu)\,\frac{dn}{n}\,
\sum_{m=0}^{\infty}\,\exp^{\#}(-f^{\Uparrow\frak{b}^\ast}(n,m))
\\ \nonumber
\!&=&\!\! \frac{1}{2\pi i}\,
\int_{c-i\infty}^{c+i\infty}\!\exp(n\,\nu)\,\frac{dn}{n}
\sum_{m=0}^{\infty}\!
\exp\!\big(\!\!-\!\!f_\kappa\, \frac{n^{-\kappa}\,m^{\kappa+1}}{\kappa+1}\big)
\exp_{\#}\!\big(\!-\!f^{\uparrow\frak{b}^\ast}(n,m)\big)
\end{eqnarray}
or equivalently (and preferably)\,:
\begin{eqnarray}  \label{a117}
h(\nu) \!\!&=&\!\! \frac{1}{2\pi i}\,
\int_{c-i\infty}^{c+i\infty}\exp(n\,\nu)\,\frac{dn}{n}\,
\sum_{m\,\in\frac{1}{2}+\doN}^{\infty}\,
\exp^{\#}(-f^{\Uparrow\beta^\ast}(n,m))
\\ \nonumber
\!&=&\!\! \frac{1}{2\pi i}\,
\int_{c-i\infty}^{c+i\infty}\!\exp(n\,\nu)\,\frac{dn}{n}
\sum_{m\,\in\frac{1}{2}+\doN}^{\infty}\!
\exp\!\big(\!\!-\!\!f_\kappa\, \frac{n^{-\kappa}\,m^{\kappa+1}}{\kappa+1}\big)
\exp_{\#}\!\big(\!-\!f^{\uparrow\beta^\ast}(n,m)\big)
\end{eqnarray}
where $\mi{exp}_{\#}(X)$
 denotes
 the usual exponential function, but  expanded  as a power series of $X$.
Similarly,
$\mi{exp}^{\#}(X)$  denotes
 the exponential  expanded as a power series
of $X$ {\it minus} the leading term of $X$, which remains within the exponential.
 An unmarked 
 $\mi{exp}(X)$, 
on the other hand, should be construed as the usual exponential function.
\\
The analytical expressions vary depending on the tangency order $\kappa$.
Indeed, after expanding $\mi{exp}_{\#}(...)$, we are left with
the task of calculating individual sums of type\,:
\[\begin{array}{llllllllll}
S_{0,k}(f_0)&\!=\!&\sum_{m\,\in\doN}
\;m^k\;\exp({-f_0\,m}) &\mi{in \; case\;
2}&(\kappa=0,f_0>0)
\\[1.0 ex]
S_{\kappa,k}(f_\kappa)&\!=\!&\sum_{m\,\in\doN}\;m^k\;
\exp\big({-f_\kappa\,\,\frac{n^{-\kappa}\,m^{\kappa+1}}{\kappa+1}}\big)
\quad\quad
 &\mi{in \; case\;3}&(\kappa\geq 1,f_\kappa>0)
\end{array}\]
or of type\,:
\[\begin{array}{llllllllll}
Z_{0,k}(f_0)&\!=\!&\sum_{m\,\in\frac{1}{2}+\doN}\;m^k\;\exp({-f_0\,m}) &\mi{in \; case\;
2}&(\kappa=0,f_0>0)
\\[1.0 ex]
Z_{\kappa,k}(f_\kappa)&\!=\!&\sum_{m\,\in\frac{1}{2}+\doN}\;m^k\;
\exp\big({-f_\kappa\,\,\frac{n^{-\kappa}\,m^{\kappa+1}}{\kappa+1}}\big)
\quad
 &\mi{in \; case\;3}&(\kappa\geq 1,f_\kappa>0)
\end{array}\]
Since we assumed $f_\kappa$ to be positive in all cases,
convergence is immediate and precise bounds are readily found.
However, only for $\kappa=0$ do the sums 
$S_{\kappa,k}\,,\,Z_{\kappa,k}$ admit closed
expressions for all $k$. For the former sums we get\,:
\begin{equation*}
S_{0,k}(\alpha)=\frac{L_k(a)\;\;}{(1-a)^{k+1}}
\;\;\mi{with}\;\; a:=e^{-\alpha}\;\; ;\;\; 
L_k(a):=\mi{tr}_k\Big( (1-a)^{k+1}\!\!\!\!\sum_{0\leq s\leq 2k}s^k\,a^s\Big)
\end{equation*}
\\[-3.ex]
where $\mi{tr}_k$ means that we truncate after the $k^{th}$ power of $a$,
which leads to self-symmetrical polynomials of the form\,:
$$
L_k(a)=a+(1+k+2^k)\,a^2+\dots+(1+k+2^k)\,a^{k-1}+a^k
\;\;\mi{with}\;\;L_k(1)=k!
$$
For the latter sums we get the generating function\,:
\begin{equation}  \label{a118}
\sum_{0\leq k} Z_{0,k}(\alpha)\;\frac{\sigma^k}{k!}=
\frac{1}{e^{\frac{1}{2}(\alpha-\sigma)}-e^{-\frac{1}{2}(\alpha-\sigma)}}
=\frac{1}{\alpha-\sigma}+\mi{Regular}(\alpha-\sigma)\quad
\end{equation}
Hence\,:
\begin{equation}  \label{a119}
Z_{0,k}(\alpha)=\frac{k!}{\alpha^{k+1}}+\mi{Regular}(\alpha)
\end{equation}
Let us now justify the above definition of {\it nur}. For a tangency order $\kappa \geq
0$ and a driving function $f(x):=\sum_{s\geq \kappa} f_s\,x^s $ as in the cases 2 or 3
of \S5.1, our Taylor coefficients $J(n)$ will have the following asymptotic
expansions, {\it before} and {\it after} division by the ingress factor $I\!g_F(n)$\,:
\begin{eqnarray}  \label{a120}
\tilde{J}(n)&:=&\sum_{0\leq m}\exp\Big(\!-\!\sum_{0\leq k\leq m}f(\frac{k}{m})\Big)
\\[1.5 ex]  \label{a121}
&=&\sum_{0\leq m}\exp\Big(\!- (m+1)f_0
-\!\sum_{1\leq s}n^{-s}
\big(\frak{b}_s^\ast(m)-\frak{b}_s^\ast(0)\big)\,f_s
\Big)\quad
\\[1.5 ex]  \label{a122}
\tilde{I\!g}_F(n)&:=&
\exp\big(-\frac{1}{2}f_0+\sum_{1\leq
s}n^{-s}\frak{b}_s^\ast(0)f_s\big)
\\[1.5 ex]  \label{a123}
\tilde{J}(n)/\tilde{I\!g}_F(n)
&=&\sum_{0\leq m}\exp\Big(\!
-(m+\frac{1}{2})\,f_0
-\sum_{1\leq s}n^{-s}
\frak{b}_s^\ast(m)\,f_s
\Big)
\\[1.5 ex]  \label{a124}
&=&\sum_{0\leq m}\exp\Big(\!
-\sum_{0\leq s}n^{-s}
\frak{b}_s^\ast(m)\,f_s
\Big)
\end{eqnarray}
Of course, the summand $\frac{1}{2}f_0$ automatically disappears when
the tangency order $\kappa$ is $>0$. But, whatever the value of
$\kappa$, the hypothesis $f_\kappa>0$ ensures the convergence of 
the $m$-summation\footnote{
after factoring $\exp(-\sum_{\kappa\leq s}(...)) $
into $\exp(-\sum_{\kappa= s}(...))\exp_{\#}(-\sum_{\kappa < s}(...)) $ 
and expanding the second factor as a power series of $ (\sum_{\kappa < s}(...))$.
} 
in (\ref{a124}), which yields, in front of any given power $n^{-s}$, a well-defined, finite
coefficient. If we then suject the right-hand side of (\ref{a124}), term-wise, to
the (upper) Borel transform $n\rightarrow \nu$, we are led straightaway
to the above definition of the {\it nur}-transform $f(x)\mapsto h(\nu)$.
%%%%%%%%%%%%%%%%%%%%%%%%%%%%%%%%%%%%%%%%%%%%%%%%%%%%%%%%%%%%%%%%%%%%%%%%%%%%%%%%%%%%%%%%%%%%
%%%%%%%%%%%%%%%%%%%%%%%%%%%%%%%%%%%%%%%%%%%%%%%%%%%%%%%%%%%%%%%%%%%%%%%%%%%%%%%%%%%%%%%%%%%%
\subsection{Expressing {\it nur} in terms of {\it nir}.}
\begin{lemma}[Decomposition of {\it nur}.] \hspace{3.4 ex}\\
The {\it nur}-transforms reduces to an alternating
sum of  {\it nir}-transforms\,:
\begin{equation}  \label{a125}
\mr{nur}(f)=\sum_{p\in\doZ}(-1)^p\,\mr{nir}(2\pi i\,p +f) 
\end{equation} 
\end{lemma}
It suffices to show that this holds term-by-term, i.e. for the coefficient
of each monomial $\nu^n$ on the left- and right-hand sides of (\ref{a125}). For $\kappa=0$
for instance, this results from the identities\,:
\begin{equation}  \label{a126}
\sum_{m\in \frac{1}{2}+\doZ}m^k\exp(-f_0\,m)=
\sum_{p\in\doZ}(-1)^p\frac{(k+1)!}{(2\pi i\,p+f_0)^{k+1}}
\end{equation}
which are a direct consequence of Poisson's summation formula\footnote{
Decompose the left-hand side of (\ref{a126}) as
$\sum_{m\in\frac{1}{2}+\doZ}=
\sum_{m\in\frac{1}{2}\doZ}
-\sum_{m\in\doZ}
$ and formally apply Poisson's formula separately to each sum.
}. The same argument applies for $\kappa >0$.

As a consequence of the above lemma, we see that whereas the {\it nir}-transform
depends on the exact determination of $\mi{log}F$, 
the {\it nur}-transform
depends only on the determination of $F^{1/2}$. This was quite predictable,
in view of the interpretation of {\it nur}.\footnote{
The square root of $F$ comes from our having replaced
$j_F$ by $j_F^\#$, i.e. from dividing by the ingress factor, which
carries the term $e^{-f_0/2}=F(0)^{1/2}$ }. 
%%%%%%%%%%%%%%%%%%%%%%%%%%%%%%%%%%%%%%%%%%%%%%%%%%%%%%%%%%%%%%%%%%%%%%%%%%%%%%%%%%%%%%%%%%%%
%%%%%%%%%%%%%%%%%%%%%%%%%%%%%%%%%%%%%%%%%%%%%%%%%%%%%%%%%%%%%%%%%%%%%%%%%%%%%%%%%%%%%%%%%%%%

\subsection{The {\it mur} transform.}
Since in this new nine-link chain (of \S5.2) all the steps but {\it mur}
are elementary, and the composite step {\it nur} has just been defined,
that indirectly determines {\it mur} itself, just as knowing {\it nir}
determined {\it mir} in the preceding section. There are, however,
two basic differences between {\it mur} and {\it mir}.
\\
 
(i) {\it Analytic difference}\,: whereas the singularities of a 
{\it mir}-transform were {\it mir}-transforms of singularities 
(reflecting the essential closure of the inner algebra),
the singularities of {\it mur}-transforms are {\it mir}-transforms,
(not {\it mur}-transforms!) of singularities (reflecting the non-recurrence
of outer generators under alien derivation).
\\
 
(ii) {\it Formal difference}\,:
unlike {\it mir}, {\it mur} doesn't reduce to a purely integro-differential
functional. It does admit interesting, if complex, expressions\footnote{
somewhat similar to the expression for the generalised (non-standard) {\it mir}-transform
 when we drop the condition $\beta_0\not= 0$. 
}
but we needn't bother with them, since the whole point of deriving
an exact analytical expression for {\it mir} was to account for the
closure phenomenon just mentioned in (i) but which
no longer applies to {\it mur}.

%%%%%%%%%%%%%%%%%%%%%%%%%%%%%%%%%%%%%%%%%%%%%%%%%%%%%%%%%%%%%%%%%%%%%%%%%%%%%%%%%%%%%%%%%%%%
%%%%%%%%%%%%%%%%%%%%%%%%%%%%%%%%%%%%%%%%%%%%%%%%%%%%%%%%%%%%%%%%%%%%%%%%%%%%%%%%%%%%%%%%%%%%

\subsection{Translocation of the {\it nur} transform.}

Like with {\it nir}, it is natural to `translocate' {\it nur}, i.e.
to measure its failure to commute with translations. To do this,
we have the choice, once again, between four expressions
(where $  \eta:=\int_0^\epsilon f(x)\,dx$)
\[\begin{array}{ccccc}
\mi{choice\;1\,:}&
(\mr{nur}\;e^{\epsilon\partial_x} - e^{\eta\partial_\nu}\;\mr{nur})\,f 
&\mi{as\; a\; function\; of}&(\epsilon,f)
\\[1.ex]
\mi{choice\;2\,:}&
(\mr{nur}\;e^{\epsilon\partial_x} - e^{\eta\partial_\nu}\;\mr{nur})\,f 
&\mi{as\; a\; function\; of}&(\eta,f)
\\[1.ex]
\mi{choice\;3\,:}&
( \mr{nur} -  e^{-\eta\partial_\nu}\,\mr{nur}\;e^{\epsilon\partial_x})\,f 
&\mi{as\; a\; function\; of}&(\epsilon,f)
\\[1.ex]
\mi{choice\;4\,:}&
(  \mr{nur}-e^{-\eta\partial_\nu}\,\mr{nur}\;e^{\epsilon\partial_x})\,f 
&\mi{as\; a\; function\; of}&(\eta,f)
\\[1.ex]
\end{array}\]
but whichever choice we make (let us think of choice 3, for consistency)
two basic differences emerge between {\it nir}'s and {\it nur}'s
translocations\,:
\\

(i) {\it Analytic difference}\,: the finite or infinitesimal increments
$\nabla h(\epsilon,\nu) $ or $\delta h_m(\nu) $ defined
as in \S4.6 but with respect to {\it nur}, are no longer entire
functions of their arguments, even when the driving function $f$
is entire or polynomial. The reason for this is quite simple\,: with
the {\it nir}-transform, to a shift $\epsilon$ in the $x$-plane
there answers a well-defined shift $\eta=\int_0^\infty f(x)dx$
in the $\nu$-plane, calculated from a well-defined
determination of $f=-\log F$, but this no longer holds
with the {\it nur}-transform, whose construction involves {\it all}
determinations of $f$.
\\

(ii) {\it Formal difference}\,: these  increments
still admit exact analytical expansions somewhat similar
to (\ref{a94}) and (\ref{a104}) but the formulas are now more complex\footnote{
with twisted equivalents of the convolution (\ref{a100}), under replacement
of the factorials by $q$-factorials.
}
and above all less useful.
Indeed,  the main point of these formulas in the {\it nir} version 
was to establish that the increments 
$\nabla h(\epsilon,\nu) $ or $\delta h_m(\nu) $ 
were entire functions of $\epsilon$ and $\nu$, but
with {\it nur} this is no longer the case, 
as was just pointed out.

%%%%%%%%%%%%%%%%%%%%%%%%%%%%%%%%%%%%%%%%%%%%%%%%%%%%%%%%%%%%%%%%%%%%%%%%%%%%%%%%%%%%%%%%%%%%
%%%%%%%%%%%%%%%%%%%%%%%%%%%%%%%%%%%%%%%%%%%%%%%%%%%%%%%%%%%%%%%%%%%%%%%%%%%%%%%%%%%%%%%%%%%%

\subsection{Removal of the ingress factor.}
As we saw, changing $j_F$ into $j_F^\#$ brings rather different
changes to the construction of
 the {\it inner} and {\it outer} generators\,:
for the {\it inner} generators it means merging 
the critical stationary factor $J_4$
with the {\it egress} factor $E\!g_F$ ; 
for the {\it outer} generators it means pruning
the critical stationary factor $J$
of the {\it ingress} factor $I\!g_F$. Nonetheless,
the end effect is exactly the same\,:  the parasitical summands $\frak{b}_s^\ast(0)$
vanish from (\ref{a58}) and (\ref{a124}) alike.

%%%%%%%%%%%%%%%%%%%%%%%%%%%%%%%%%%%%%%%%%%%%%%%%%%%%%%%%%%%%%%%%%%%%%%%%%%%%%%%%%%%%%%%%%%%%
%%%%%%%%%%%%%%%%%%%%%%%%%%%%%%%%%%%%%%%%%%%%%%%%%%%%%%%%%%%%%%%%%%%%%%%%%%%%%%%%%%%%%%%%%%%

\subsection{Parity relations.}
\begin{eqnarray*}
F^{\,\bot}(x):=1/F(-x)\;\;\;,\;\;\;f^{\,\bot}(x):=-f(-x)  \hspace{7.ex}
&&\Longrightarrow 
\\[1.5 ex]
\mr{nur}(f^\bot)(\nu)=-\mr{nur}(f)(\nu)
\hspace{23.ex}
&& (\mi{tangency}\;\, \kappa=0)
\end{eqnarray*}

\section{Inner generators and ordinary differential equations.} 
%%%%%%%%%%%%%%%%%%%%%%%%%%%%%%%%%%%%%%%%%%%%%%%%%%%%%%%%%%%%%%%%%%%%%%%%%%%%%%%%%%%%%%%%%%%%
%%%%%%%%%%%%%%%%%%%%%%%%%%%%%%%%%%%%%%%%%%%%%%%%%%%%%%%%%%%%%%%%%%%%%%%%%%%%%%%%%%%%%%%%%%%%
In some important instances, namely for {\it all} polynomial inputs $f$ and {\it some}
rational inputs $F$, 
 the corresponding inner  generators happen to verify ordinary differential
equation of a rather simple type -- {\it linear homogeneous with polynomial coefficients}
-- but often of high degree. These ODEs are interesting on three accounts
\\
(i) they lead to an alternative, more classical derivation of the properties of these inner
generators
\\ 
(ii) they yield a precise description of their behaviour {\it over} $\infty$ in
the $\nu$-plane, i.e. {\it over} $0$ in the $\zeta$-plane.
\\
(iii) they stand out, among similar-looking ODEs, as leading to a {\it rigid}
resurgence pattern, with essentially {\it discrete} Stokes constants, insentitive to the {\it
continuously} varying parameters.

\subsection{``Variable" and ``covariant" differential equations.}
As usual, we consider four types of shift operators $\beta(\partial_{\tau})$, relative to the choices
\begin{eqnarray}\label{diff1}
\mi{trivial\; choice}\quad\quad&\beta(\tau):=&\tau^{-1}
\\   \label{diff2}
\mi{standard\; choice}\quad&\beta(\tau):=&(e^{\tau/2}-e^{-\tau/2}\Big)^{-1}=\tau^{-1}-\frac{1}{24}\tau+\dots
\\   \label{diff3}
\mi{odd\; choice}\quad\quad\quad&\beta(\tau):=&\tau^{-1}+\sum_{s\geq 0}\beta_{2\,s+1}\,\tau^{2\,s+1}
\\  \label{diff4}
\mi{general\; choice}\quad\;&\beta(\tau):=&\tau^{-1}+\sum_{s\geq 0}\beta_{s}\,\tau^{s}
\end{eqnarray}
We the apply the {\it nir}-transform to a driving function $f$ such that $f(0)=0$, with special emphasis on the
case $f'(0)\not=0$\,:
\begin{eqnarray}\label{diff5}
f(x)&:=&\sum_{1\leq s \leq r} f_{s}\,x^s
\\\label{diff6}
\varphi(n,\tau)&:=&\beta(\tau)\,f(\frac{\tau}{n})=\frac{1}{2}\frac{\tau^2}{n}\,f_{1}+\dots \quad\in\doC[n^{-1},\tau]
\\\label{diff7}
\varphi(n,\tau)&:=&\varphi^{+}(n,\tau)+\varphi^{-}(n,\tau)\quad\mi{with}
\quad \varphi^{\pm}(n,\pm\tau) \equiv\pm\varphi^{\pm}(n,\tau)
\end{eqnarray}
\begin{eqnarray}\label{diff8}
k(n)&:=&\Big[\int_{0}^{\infty} \exp^{\#}\!(\varphi(n,\tau))\;d\tau\Big]_{\rm{singular}}
\\   \label{diff9}
&:=&\;\int_{0}^{\infty} \exp^{\#}\!(\varphi(n,\tau))\;\cosh_{\#}\!(\varphi^{-}(n,\tau))\;d\tau\quad (\mi{if}\;\;f_{1}\not=0)
\\    \label{diff10}
\imi{k}(\nu)&:=&\Big[ \frac{1}{2\pi i}\int_{c-i\infty}^{c+i\infty}k(n)\,e^{\nu\,n}\;dn\Big]_{\rm{formal}}\;\;\;=h'(\nu)
\\   \label{diff11}
\smi{k}(\nu)&:=&\Big[ \frac{1}{2\pi i}\int_{c-i\infty}^{c+i\infty}k(n)\,e^{\nu\,n}\;\frac{dn}{n}\Big]_{\rm{formal}}\;\;\;=h(\nu)
\end{eqnarray}
But the case $f(0)\not=0$ also matters, because it corresponds the so-called ``exceptional'' or ``movable'' generators. In that case
the {\it nir}-transform produces no fractional powers. So we set\,:
\begin{equation}
f(x):=\sum_{0\leq s \leq r} f_{s}\,x^s \hspace{20.ex}
\end{equation}
\begin{eqnarray}\label{diff12}
k^{\rm{total}}(n)&:=&\int_{0}^{\infty} \exp^{\#}\!(\varphi(n,\tau))\;d\tau
\\    \label{diff13}
\imi{k}^{\rm{total}}(\nu)&:=&\Big[ \frac{1}{2\pi i}\int_{c-i\infty}^{c+i\infty}k^{\rm{total}}(n)\,e^{\nu\,n}\;dn\Big]_{\rm{formal}}\;
\\   \label{diff14}
\smi{k}^{\rm{total}}(\nu)&:=&\Big[ \frac{1}{2\pi i}\int_{c-i\infty}^{c+i\infty}k^{\rm{total}}(n)\,e^{\nu\,n}\;\frac{dn}{n}\Big]_{\rm{formal}}
\end{eqnarray}
The above definitions also extend to the case $f(0)=0$. The  {\it nir}-transform then produces a mixture of entire
and fractional powers, and the index {\it total} affixed to $k$ signals that we take them all.
\\
\noindent
For polynomial inputs, both $k^{\rm{total}}$ and $k$ along with their Borel transforms verify remarkable linear-homogeneous ODEs. The ones verified by  $k^{\rm{total}}$
are dubbed {\it variable} because there is no simple description of how they change when the base point changes in the $x$-plane
(i.e. when the driving function undergoes a shift from $f$ to $^{\epsilon}\!f$ ). The ODEs verified by $k$, on the other hand, deserve to be called {\it covariant}, for two reasons\,:
\\
(i) when going from a proper base-point $x_{i}$ to another proper base-point $x_{j}$ ({\it proper} means that $f(x_{i})=0, f(x_{j})=0$), these covariant ODEs verified
by $\imi{k}(\nu)$ simply undergo a shift $\nu=\int_{x_{i}}^{x_{j}}f(x)dx$ in the $\nu$-plane.
\\
(ii) there is a unique extension of the covariant ODE even to non-proper base-points $x_{i}$ (i.e. when $f(x_{i})\not=0$), under the same
formal covariance relation as above. That extension, of course, doesn't coincide with the {\it variable} ODE.\footnote{ For a proper base-point, on the other hand, the variable  ODEs, though still distinct from the covariant ones, are {\it also} verified by $k$.}
\\

\noindent
{\bf ``Variable" and ``covariant'' linear-homogeneous polynomial  ODEs:}
\\
They are of the form\,:
\begin{eqnarray} \label{a127}
\mi{variable\; ODE\,:}\quad &P_{v}(n,-\partial_{n})\,k^{\mr{total}}(n)=0&
\Leftrightarrow\;\; P_{v}(\partial_{\nu},\nu)\,\imi{k}^{\mr{total}}\!\!(\nu)=0\quad\quad
\\ \label{a128}
\mi{covariant\; ODE\,:}\quad & P_{c}(n,-\partial_{n})\,k(n)=0&
\Leftrightarrow\;\;  P_{c}(\partial_{\nu},\nu)\,\imi{k}(\nu)=0
\end{eqnarray}
with polynomials
\begin{eqnarray} 
P_{v}(n,\nu)&=&\sum_{0\leq p \leq d}\sum_{0\leq q \leq \delta}\;\rm{dv}_{p,q}\;n^p\,\nu^q
\\
P_{c}(n,\nu)&=&\sum_{0\leq p \leq d}\sum_{0\leq q \leq \delta}\;\rm{dc}_{p,q}\;n^p\,\nu^q
\end{eqnarray}
of degree $d$ and $\delta$ in the non-commuting variables $n$ and $\nu$\,: \;\; $[n,\nu]=1$.
\\ \
The covariance relation reads\,:
\begin{equation}  \label{a129}
P_{c}^{\,^{\epsilon}\!f}(n,\nu\!-\eta)\equiv P_{c}^{f}(n,\nu)\;\forall \epsilon
\quad \mi{with}\quad ^{\epsilon}\!f(x)=f(x+\epsilon)\quad \mi{and}\quad
\eta:=\int_0^\epsilon f(x)dx
\end{equation}
{\bf Existence and calculation of the variable ODEs.}\\
For any $s\in \doN$ let $\varphi_{s},\psi_{s}$ denote the polynomials in $(n^{-1},\tau)$ characterised by the identities\,:
\begin{eqnarray}\label{diff15}
\partial_{n}^s\; k^{\rm{total}}(n)&=&\int_{0}^{0}\varphi_{s}(n,\tau)\,\exp^\#\!(\varphi(n,\tau))\;d\tau
\\
\mi{with}&&\varphi_{s}(n,\tau)\in
\doC[\partial_{n}\varphi,\partial_{n}^2\varphi,\dots,\partial_{n}^s\varphi]\in \doC[n^{-1},\tau]
\\\int_{0}^\infty d_{\tau}^s\big( \tau^k\;\exp^\#\!(\varphi(n,\tau))\big)&=&
\int_{0}^\infty\psi_{s}(n,\tau)\;\exp^\#\!(\varphi(n,\tau))\;d\tau\;=0
\\
\mi{with}&&\psi_{s}(n,\tau)=\tau^s\,\partial_{\tau}\varphi(n,\tau)+s\,\tau^{s-1}\in \doC[n^{-1},\tau]
\end{eqnarray}
For $\delta,\delta'$ large enough, the polynomials $\varphi_{s},\psi_{s}$ become linearly dependent on $\doC[n^{-1}]$
or, what amounts to the same, on $\doC[n]$. So we have relations of the form\,:
\begin{equation}
0=\sum_{0\leq s\leq \delta}A_{s}(n)\,\varphi(n,\tau)+\sum_{0\leq s\leq \delta'}B_{s}(n)\,\varphi(n,\tau)\quad
\mi{with}\quad A(n),B(n)\in \doC[n]
\end{equation}
and to each such relation there corresponds a linear ODE for $k^{\rm{total}}$\,:
\begin{equation}
\Big(\sum_{0\leq s\leq \delta}A_{s}(n)\,\partial_{n}^s\Big)\;k^{\rm{total}}(n)\;=\;0
\end{equation}
\\
{\bf Existence and calculation of the covariant ODEs for $f(0)=0$.}\\
For each $s\in \doN$ let $\varphi_{s}^\pm$ and  $\psi_{s}^{\pm\pm},\psi_{s}^{\pm\mp}$  denote the polynomials in
$(n^{-1},\tau)$ characterised by the identities\,:
\begin{eqnarray*}
\partial_{n}^s\;k(n)\!&\!=\!&\!
\int_{0}^\infty\Big(\varphi^+_{s}\!(n,\tau)\cosh(\varphi^-(n,\tau))+\varphi^-_{s}\!(n,\tau)\sinh(\varphi^-(n,\tau))
 \Big)\,e^{\varphi^+(n,\tau)}\,d\tau\quad\quad
 \\
\mi{with}&&\varphi^\pm_{s}(n,\tau)\in
\doC[\partial_{n}\varphi^+,\dots,\partial_{n}^s\varphi^+,
\partial_{n}\varphi^-,\dots,\partial_{n}^s\varphi^- ]\in \doC[n^{-1},\tau]
\end{eqnarray*}
\begin{eqnarray*}
\int_{0}^\infty d_{\tau}^s\big( \tau^k\;e^{\varphi^+(n,\tau)} \cosh(\varphi^-(n,\tau))\big)\!&\!=\!&\! 
\int_{0}^\infty d_{\tau}^s\big( \tau^k\;e^{\varphi^+(n,\tau)} \sinh(\varphi^-(n,\tau))\big) \,=\, 0
\\ 
\mi{with}\quad\quad&&
\\
 d_{\tau}^s\big( \tau^k\;e^{\varphi^+(n,\tau)} \cosh(\varphi^-(n,\tau))\big)&=& 
 +\varphi_{s}^{++}(n,\tau)\; e^{\varphi^+(n,\tau)}\cosh(\varphi^-(n,\tau) d\tau
 \\
 && + \varphi_{s}^{+-}(n,\tau)\; e^{\varphi^+(n,\tau)}\sinh(\varphi^-(n,\tau) d\tau
 \\ 
 d_{\tau}^s\big( \tau^k\;e^{\varphi^+(n,\tau)} \sinh(\varphi^-(n,\tau))\big)&=& 
 +\varphi_{s}^{--}(n,\tau)\; e^{\varphi^+(n,\tau)}\cosh(\varphi^-(n,\tau) d\tau
 \\
 && + \varphi_{s}^{-+}(n,\tau)\; e^{\varphi^+(n,\tau)}\sinh(\varphi^-(n,\tau) d\tau
\end{eqnarray*}
Here again, for $\delta,\delta'$ large enough, there are going to be dependence relations of the form\,:
\begin{eqnarray}
0&=&\sum_{0\leq s\leq \delta}A_{s}(n)\,\varphi^{+}(n,\tau)
        +\sum_{0\leq s\leq \delta'}B_{s}(n)\,\psi^{++}(n,\tau)
        +\sum_{0\leq s\leq \delta'}C_{s}(n)\,\psi^{++}(n,\tau)\quad\quad
        \\
 0&=&\sum_{0\leq s\leq \delta}A_{s}(n)\,\varphi^{-}(n,\tau)
        +\sum_{0\leq s\leq \delta'}B_{s}(n)\,\psi^{+-}(n,\tau)
        +\sum_{0\leq s\leq \delta'}C_{s}(n)\,\psi^{-+}(n,\tau)\quad\quad
        \\
        \mi{with}&& A(n),B(n),C(n) \in \doC[n]
\end{eqnarray}
and to each such relation there will corresponds a linear ODE for $k$\,:
\begin{equation}
\Big(\sum_{0\leq s\leq \delta}A_{s}(n)\,\partial_{n}^s\Big)\;k(n)\;=\;0
\end{equation}
\\
Remark\,: although the above construction applies, strictly speaking, only to the case of tangency $\kappa=1$, i.e. to
 the case $f_{0}=0,f_{1}\not=0$, it is in fact universal. Indeed, if we set $f_{0}=f_{1}=\dots=f_{\kappa-1},f_{\kappa}\not=0$ in the covariant ODEs
 thus found, we still get the correct covariant ODEs for a general tangency order $\kappa>1$.
 \\
 
 \noindent
 {\bf Existence and calculation of the covariant ODEs for $f(0)\not=0$.}\\
There are five steps to follow\,:
\\
(i) Fix a degree $r$ and calculate $P^f(n,\nu)$ by the above method for
an arbitrary $f$ of degree $r$ such that  $f(0)=0$.
\\
(ii) Drop the assumption $f(0)=0$ but subject $f$ to a shift $\epsilon$
such that $^{\epsilon}\!f(0)=f(\epsilon)=0$ and apply (i) to calculate
$P^{\,^{\epsilon}\!f}(n,\nu)$ without actually solving the equation $f(\epsilon)=0$
(keep $\epsilon$ as a free variable).
\\
(iii) Calculate the $\epsilon$-polynomial $P^{\,^{\epsilon}\!f}(n,\nu\!-\!f^\ast(\epsilon))$
with $f^\ast(x):=\int_0^xf(t)dt $ as usual.
\\
(iv) Divide it by the $\epsilon$-polynomial $f(\epsilon)$ (momentarily assumed to be $\not=0$)
and calculate the remainder $P_0$ and quotient $P_1$ of that division\,:
\begin{equation*}
P^{\,^{\epsilon}\!f}(n,\nu-f^\ast(\epsilon))
=: P_0^{f}(n,\nu,\epsilon)+P_1^{f}(n,\nu,\epsilon)\, f(\epsilon)
\end{equation*}
(v) Use the covariance identity 
$P^{\,^{\epsilon}\!f}(n,\nu\!-\!f^\ast(\epsilon))\equiv P^{f}(n,\nu)\;\forall \epsilon $
to show that the remainder $P_0^{f}(n,\nu,\epsilon)$ is actually constant in
$\epsilon$. Then set  
\begin{equation*}
 P^{f}(n,\nu):= P_0^{f}(n,\nu,0)
\end{equation*}

%%%%%%%%%%%%%%%%%%%%%%%%%%%%%%%%%%%%%%%%%%%%%%%%%%%%%%%%%%%%%%%%%%%%%%%%%%%%%%%%%%%%%%%%%%%%
%%%%%%%%%%%%%%%%%%%%%%%%%%%%%%%%%%%%%%%%%%%%%%%%%%%%%%%%%%%%%%%%%%%%%%%%%%%%%%%%%%%%%%%%%%%%

\subsection{ODEs for polynomial inputs $f$. Main statements.}
{\bf Dimensions of spaces of variable ODEs:}
\\
For $r:=\deg(f)$ and for each pair $(x.y.)$ with \\
$x\in \{v,c\}=\{\mi{variable}, \mi{covariant}\}$ 
\\
$y\in \{t,s,o,g\}=\{\mi{trivial}, \mi{standard}, \mi{odd},\mi{general}\}$ 
\\
the dimension of the corresponding space of ODEs is always of the form:
\begin{eqnarray}  \label{a132}
\mr{dim_{\,x.y.}}(r,d,\delta)&\equiv&
\big(d-A_{\,x.y.}(r)\big)\big(\delta-B_{\,x.y.}(r)\big)-C_{\,x.y.}(r)
\end{eqnarray}
with $\delta$(resp.$d$) denoting the differential order of the the ODEs in the $n$-variable (resp. in
the $\nu$-variable).
Of special interest are the extremal pairs $(\underline{d},\overline{\delta}) $
and  $(\overline{d},\underline{\delta}) $ with
\begin{eqnarray}  \label{a133}
\underline{d}=  1+A_{\,x.y.}(r)\hspace{10.ex}
&\quad& 
\overline{\delta}= 1+B_{\,x.y.}(r) +C_{\,x.y.}(r)
\\ [1.ex]  \label{a134}
\overline{d}= 1+A_{\,x.y.}(r)+C_{\,x.y.}(r) 
&\quad&
\underline{\delta}= 1+B_{\,x.y.}(r)
\end{eqnarray}
( $\underline{d}$ and $\underline{\delta}$ minimal;
 $\overline{d}$ and $\overline{\delta}$ co-minimal)
because the corresponding dimension is exactly\;1.
{\bf Dimensions of spaces of variable ODEs:}
\begin{eqnarray*}
\mr{dim_{\,v.t.}}(r,d,\delta)
\!&\!=\!&\! (d\!-\!r)\,(\delta\!-\!r\!-\!1)
-\frac{1}{2}r^2+\frac{1}{2}r-1
\\
\mr{dim_{\,v.s.}}(r,d,\delta)
\!&\!=\!&\! (d\!-\!r)\,(\delta\!-\!r^2\!-\!2r\!+\!1)
-\frac{1}{2}r^2(r\!+\!1)\hspace{7 ex} (r\;\mi{even})
\\
\!&\!=\!&\! (d\!-\!r)\,(\delta\!-\!r^2\!-\!2r)
-\frac{1}{2}(r^3\!+\!r^2\!-\!5r\!+\!5)\hspace{2 ex} (r\;\mi{odd})
\\
\mr{dim_{\,v.o.}}(r,d,\delta)
\!&\!=\!&\! (d\!-\!r)\,(\delta\!-\!r^2\!-\!2r\!+\!1)
-\frac{1}{2}r^2(r\!+\!1)\hspace{7 ex} (r\;\mi{even})
\\
\!&\!=\!&\! (d\!-\!r)\,(\delta\!-\!r^2\!-\!2r)
-\frac{1}{2}(r^3\!+\!r^2\!-\!3r\!+\!3)\hspace{2 ex} (r\;\mi{odd}\not=3)
\\
\mr{dim_{\,v.g.}}(r,d,\delta)
\!&\!=\!&\! (d\!-\!r)\,(\delta\!-\!r^2\!-\!2r)
-\frac{1}{2}r^2(r\!+\!1)
\end{eqnarray*}
{\bf Dimensions of spaces of covariant ODEs:}
\begin{eqnarray*}
\mr{dim_{\,c.t.}}(r,d,\delta)
\!&\!=\!&\! (d\!-\!r\!+\!1)\,(\delta\!-\!r\!+\!1)
-\frac{1}{2}(r\!-\!1)(r\!-\!2)
\\
\mr{dim_{\,c.s.}}(r,d,\delta)
\!&\!=\!&\! (d\!-\!r\!+\!1)\,(\delta\!-\!r^2\!-\!r\!+\!1)
-\frac{1}{2}r^2(r\!-\!1)\hspace{7 ex} (r\;\mi{even})
\\
\!&\!=\!&\! (d\!-\!r\!+\!1)\,(\delta\!-\!r^2\!-\!r\!+\!1)
-\frac{1}{2}(r^2\!-\!5)(r\!-\!1)\hspace{2 ex} (r\;\mi{odd})
\\
\mr{dim_{\,c.o.}}(r,d,\delta)
\!&\!=\!&\! (d\!-\!r\!+\!1)\,(\delta\!-\!r^2\!-\!r\!+\!1)
-\frac{1}{2}r^2(r\!-\!1)\hspace{7 ex} (r\;\mi{even})
\\
\!&\!=\!&\! (d\!-\!r\!+\!1)\,(\delta\!-\!r^2\!-\!r\!+\!1)
-\frac{1}{2}(r^2\!-\!3)(r-1)\hspace{2 ex} (r\;\mi{odd}\not=3)
\\
\mr{dim_{\,c.g.}}(r,d,\delta)
\!&\!=\!&\! (d\!-\!r\!+\!1)\,(\delta\!-\!r^2\!-\!r\!+\!1)
-\frac{1}{2}r^2(r\!-\!1)
\end{eqnarray*}
{\bf Tables of dimensions for low degrees $r=\deg(f)$\,:}
\[\begin{array}{cccccccccrrrrrrrr}
 \mi{degree} &\quad& \mi{variable} & \mi{variable} 
&\mi{variable}  &  \mi{variable}   &&   
\\[0. ex]
 r &\quad& \mi{trivial} & \mi{standard} &\mi{odd}  &  \mi{general}     
  &&
\\[1.7 ex]
1  && (\underline{d} ,\overline{\delta} ) &   
      (\underline{d} ,\overline{\delta} )&   
      (\underline{d} ,\overline{\delta} )&   
      ( \underline{d} ,\overline{\delta} )&& 
\\[1.7 ex]
1  && ( 2,2 ) &   (2 , 4)&   (2 ,4 )&   (2 ,5 )&&  
\\
2   && (3 ,4 ) &   (3 ,14 )&   (3 ,14 )&   (3 ,15 )&&  
\\
3    && (4 ,7 ) &   (4 ,28 )&   (4 ,28 )&   (4 ,34 )&&  
\\
4    && (5 , 11) &   (5 ,64 )&   (5 ,64 )&   (5 ,65 )&&   
\\
5    && (6 , 16) &   (6 ,100 )&   (6 ,104 )&   (6 ,111 )&&  
\\
6    && (7 ,22 ) &   (7 ,174 )&   (7 ,174 )&   (7 ,175 )&&  
\\
7    && (8 ,29 ) &   (8 ,244 )&   (8 ,250 )&   (8 ,260 )&&   
\\
8    && (9 ,37 ) &   (9 ,368 )&   (9 ,368 )&   (9 ,369 )&&  
\\
9    && (10 ,46 ) &   (10 ,484 )&   (10 ,492 )&   (10 ,505 )&&   
\\
10    && (11 ,56 ) &   (11 ,670 )&   (11 ,670 )&   (11 ,671 )&& 
\\
\dots  && \dots & \dots &\dots  &\dots  && 
\end{array}\]

\[\begin{array}{cccccccccrrrrrrrr}
 \mi{degree} &\quad& \mi{variable} & \mi{variable} 
&\mi{variable}  &  \mi{variable}   &&   
\\[0. ex]
 r &\quad& \mi{trivial} & \mi{standard} &\mi{odd}  &  \mi{general}     
  &&
\\[1.7 ex]
1  && (\overline{d} ,\underline{\delta} ) &   
      (\overline{d} ,\underline{\delta} )&   
      (\overline{d} ,\underline{\delta} )&   
      (\overline{d} ,\underline{\delta} )&& 
\\[1.7 ex]
1       && (3 ,1 ) &   ( 2, 4)&   (2 ,4 )&   ( 3, 4)&&   
\\
2       && ( 5,2 ) &   (9 ,8 )&   (9 ,8 )&   (9 ,9 )&&   
\\
3        && ( 8, 3) &   (16 ,16 )&   (16 ,16 )&   (22 ,16 )&&   
\\
4        && (12 ,4 ) &   (45 ,24 )&   (45 ,24 )&   (45 ,25 )&&    
\\
5       && (17 ,5 ) &   (70 ,36 )&   (74 ,36 )&   (81 ,36 )&& 
\\
6       && (23 ,6 ) &   (133 ,48 )&   (133 ,48 )&   (133 ,49 )&& 
\\
7        && (30 ,7 ) &   (188 ,64 )&   (194 ,64 )&   (204 ,64 )&&   
\\
8        && (38 , 8) &   (297 ,80 )&   (297 ,80 )&   (297 ,81 )&&   
\\
9       && (47 , 9) &   (394 ,100 )&   (402 ,100 )&   (415 ,100 )&&   
\\
10    && (57 ,10 ) &   (561 ,120 )&   (561 ,120 )&   (561 ,121 )&& 
\\
\dots  && \dots & \dots &\dots  &\dots  && 
\end{array}\]

\[\begin{array}{cccccccccrrrrrrrr}
 \mi{degree} &\quad& \mi{covariant} & \mi{covariant} 
&\mi{covariant}  &  \mi{covariant}   &&   
\\[0. ex]
 r &\quad& \mi{trivial} & \mi{standard} &\mi{odd}  &  \mi{general}     
  &&
\\[1.7 ex]
1  && (\underline{d} ,\overline{\delta} ) &   
      (\underline{d} ,\overline{\delta} )&   
      (\underline{d} ,\overline{\delta} )&   
      ( \underline{d} ,\overline{\delta} )&& 
\\[1.7 ex]
1    && ( 1,1 ) &   (1 , 2)&   (1 , 2)&   (1 ,2 )&&  
\\
2    && (2 ,2 ) &   (2 , 8 )&   (2 , 8 )&   (2 ,8 )&&  
\\
3    && (3 ,4 ) &   (3 , 16 )&   (3 , 21 )&   (3 ,21 )&&  
\\
4    && (4 , 7) &   (4 ,44 )&   (4 , 44 )&   (4 ,44 )&&   
\\
5    && (5 , 11) &   (5 ,70 )&   (5 , 80 )&   (5 ,80 )&&  
\\
6    && (6 ,16 ) &   (6 ,132 )&   (6 ,132 )&   (6 ,132 )&&  
\\
7    && (7 ,22 ) &   (7 ,188 )&   (7 ,203 )&   (7 ,203 )&&   
\\
8    && (8 ,29 ) &   (8 ,296 )&   (8 ,296 )&   (8 ,296 )&&  
\\
9    && (9 ,37 ) &   (9 ,394 )&   (9 ,414 )&   (9 ,414 )&&   
\\
10   && (10 ,46 ) &   (10 ,560 )&   (10 ,560 )&   (10 ,560 )&& 
\\
\dots  && \dots & \dots &\dots  &\dots  && 
\end{array}\]

\[\begin{array}{cccccccccrrrrrrrr}
 \mi{degree} &\quad& \mi{covariant} & \mi{covariant} 
&\mi{covariant}  &  \mi{covariant}   &&   
\\[0. ex]
 r &\quad& \mi{trivial} & \mi{standard} &\mi{odd}  &  \mi{general}     
  &&
\\[1.7 ex]
1  && (\overline{d} ,\underline{\delta} ) &   
      (\overline{d} ,\underline{\delta} )&   
      (\overline{d} ,\underline{\delta} )&   
      (\overline{d} ,\underline{\delta} )&& 
\\[1.7 ex]
1       && (1 ,1 ) &   ( 1, 2)&   (1 ,2 )&   (1 ,2 )&&   
\\
2       && ( 2, 2) &   (4 , 6)&   ( 4,6 )&   (4 ,6 )&&   
\\
3        && ( 4, 3) &   (7 ,12 )&   ( 7,12 )&   (12 ,12 )&&   
\\
4        && (7 ,4 ) &   (28 ,20 )&   (28 ,20 )&   (28 ,20 )&&    
\\
5       && (11 ,5 ) &   (45 ,30 )&   (49 ,30 )&   (55 ,30 )&& 
\\
6       && (16 ,6 ) &   (96 ,42 )&   (96 ,42 )&   (96 ,42 )&& 
\\
7        && (22 ,7 ) &   (139 ,56 )&   (145 ,56 )&   (154 ,56 )&&   
\\
8        && (29 ,8 ) &   (232 ,72 )&   (232 ,72 )&   (232 ,72 )&&   
\\
9       && (37 ,9 ) &   (313 ,90 )&   (321 ,90 )&   (333 ,90 )&&   
\\
10    && (46 ,10 ) &   (460 ,110 )&   (460 ,110 )&   (460 ,110 )&& 
\\
\dots  && \dots & \dots &\dots  &\dots  && 
\end{array}\]
{\bf Differential polynomial $P $ in the noncommuting variables $(n,\nu)$. }
\\
Our differential operators will be written as polynomials $P(n,\nu)$ of degree $(d,\delta)$ in the non-commuting variables
$(n,\nu)$, which are capable of two realisations:
$$
(n,\nu) \longrightarrow 
(n,-\partial_n) 
\quad \mi{or} \quad 
(\partial_\nu,\nu) 
$$
Both realisation are of course compatible with $[n,\nu]=1 $ and the ODE interpretation goes like this\,:
\begin{equation}  \label{a135}
P(n,-\partial_n)\; k(n) = 0 
\quad\quad\Longleftrightarrow \quad\quad
P(\partial_\nu,\nu) \imi{k}(\nu) = P(\partial_\nu,\nu)\,\partial_{\nu}\,h(\nu)=0 
\end{equation} 
{\bf Compressing the covariant ODEs.}\\
To get more manageable expressions, we can take advantage of the covariance relation to express everything in terms of shift-invariant
data. This involves three steps\,:
\\
(i) Apply the above the ODE-finding algorithm of \S6.1 to a centered polynomial $f(x)=\sum_{i=0}^{r-2}f_i\,x^i+f_r\,x^r$.
\\
(ii) Replace the coefficients $\{f_0,f_1,\dots,f_{r-2},f_r\}$
by the shift-invariants $\{\mgf_0,\mgf_1,\dots,\mgf_{r-2},\mgf_r\}$ defined in \S6.2 {\it infra}.
\\
(iii) Replace the $\beta$-coefficients by the `centered' $\mg{\beta}$-coefficients defined {\it infra}.
\\

\noindent
{\bf Basic polynomials $f(x)$ and $p(\nu)$.}
\begin{eqnarray}  \label{a136}
f(x)&=&f_0+f_1x+\dots f_rx^r \;=\;(x-x_1)\dots (x-x_r)\,f_r
\\  \label{a137}
p(\nu)&=&p_0+p_1\nu+\dots p_r\nu^r \;=\;(\nu-\nu_1)\dots (\nu-\nu_r)\,p_r
\\  \label{a138}
\mi{with}\;\;\;\; \nu_i &=& f^\ast(x_i)\,=\,\int_0^{x_i}\,f(x)\,dx
=\sum_{0\leq s \leq r} f_s\,\,\frac{x_i^{s+1}}{s+1} 
\end{eqnarray}
{\bf Basic symmetric functions 
$x_s^\ast,x_s^{\ast\ast},\nu_s^\ast,\nu_s^{\ast\ast}  $}.
\begin{eqnarray}  \label{a139}
x_1^\ast:=\sum_{1\leq i \leq r } x_i\;\;\;,\;\;\; 
x_2^\ast:=\sum_{1\leq i < j\leq r } x_i\,x_j &,\dots,& x_r^\ast:=x_1\dots x_r
\quad\quad
\\  \label{a140}
\nu_1^\ast:=\sum_{1\leq i \leq r } \nu_i\;\;\;,\;\;\; 
\nu_2^\ast:=\sum_{1\leq i < j\leq r } \nu_i\,\nu_j &,\dots,& \nu_r^\ast:=\nu_1\dots
\nu_r
\\  \label{a141}
x_s^{\ast\ast}:=\sum_{1\leq i \leq r}x_i^s\hspace{12 ex} (\forall s \in \doN )
\\  \label{a142}
\nu_s^{\ast\ast}:=\sum_{1\leq i \leq r}\nu_i^s\hspace{12 ex} (\forall s \in \doN )
\end{eqnarray}
The change from the $x$-data to the $\nu$-data goes like this\,:
$$
\{f_s\}
\longrightarrow
\{x_s^\ast\}
\stackrel{i}{\longrightarrow}
\{x_s^{\ast\ast}\}
\stackrel{ii}{\longrightarrow}
\{\nu_s^{\ast\ast}\}
\stackrel{iii}{\longrightarrow}
\{\nu_s^{\ast}\}
\longrightarrow
\{p_s\}
 $$
\begin{eqnarray*}
(i)&\quad& \sum_{1\leq s \leq \infty}\frac{1}{s}\frac{x_s^{\ast\ast}}{x^s}
\equiv -\log\Big(1+\sum_{1\leq s \leq r}(-1)^r\frac{x_s^{\ast}}{x^s}\Big)
\\[1.ex]
(ii)&\quad& \nu_s^{\ast\ast}\equiv 
\sum_{s\leq t \leq (r+1)s}\!\!\!f^\ast_{s,t}\, x_t^{\ast\ast}
\;\;\;\;\;\mi{with}\;\;\;\;
\sum_{s\leq t \leq (r+1)s}\!\!\!\!f^\ast_{s,t}\, x^t := (f^\ast(x))^s
\\[1.ex]
(iii)&\quad&
1+\sum_{1\leq s \leq r}(-1)^r\frac{\nu_s^{\ast}}{\nu^s}
\equiv
\exp\Big(- \sum_{1\leq s\leq \infty}\frac{1}{s}\frac{\nu_s^{\ast\ast}}{\nu^s}\Big)
\end{eqnarray*}
{\bf Centered polynomials. Invariants.}
\begin{eqnarray*}
x_0 &\!:=& \frac{1}{r}(x_1+\dots+x_r)\,=\, -\frac{1}{r}\frac{f_{r-1}}{f_r}
\\
\nu_0&\!:=& f^\ast(x_0)=\int_0^{x_0}f(x)\,dx \,=\,  
\sum_{0\leq s \leq r} f_s\,\,\frac{x_0^{s+1}}{s+1}
\\
\underline{\nu}_0 &\!:=& \frac{1}{r}(\nu_1+\dots+\nu_r)\,=\,
-\frac{1}{r}\frac{p_{r-1}}{p_r}
\quad\quad \quad (\nu_0\not=\underline{\nu}_{0}\;\; \mi{in\;general})
\end{eqnarray*}
\begin{eqnarray*}
\mgf(x):= f(x+x_0)=\sum_{0\leq s \leq r}\mgf_s\,x^s \;
& (\mathbf{f}_{r-1}=0)&
\\
\mgp(\nu):= p(\nu+\nu_0)=\sum_{0\leq s \leq r}\mgp_s\,\nu^s
&&\quad\quad \mathbf{P}(\nu):= P(\nu+\nu_0)
\\
\mathbf{\underline{p}}(\nu):= p(\nu+\underline{\nu}_0)=
\sum_{0\leq s \leq r}\mathbf{\underline{p}}_s\,\nu^s 
& (\mathbf{\underline{p}}_{r-1}=0)&\quad\quad
\mathbf{\underline{P}}(\nu):= P(\nu+\underline{\nu}_0)
\end{eqnarray*} 
{\bf Centered $\beta$-coefficients\,:}
\begin{equation}  \label{a130}
\beta(\tau)=\tau^{-1}+\sum_{0\leq k}\beta_k\,\tau^k 
=\tau^{-1}\,\big(1+\sum_{1\leq k}\frac{b_k}{k!}\tau^k\big)
\end{equation}
\begin{equation}  \label{a131}
1+\sum_{2\leq k}\frac{\mathbf{b}_k}{k!}\tau^k 
=
\big(1+\sum_{1\leq k}\frac{b_k}{k!}\tau^k\big)\,
\big(1+\sum_{1\leq k}\frac{(-b_1)^k}{k!}\tau^k\big)
\end{equation}
\begin{eqnarray*}
\BEE{1} &\!=\!& 0\, =\, 0
\\
\BEE{2} &\!=\!&
b_2-b_1^2 
\,=\, 
2\,\beta_1-\beta_0^2
\\
\BEE{3}  &\!=\!& 
b_3-3\,b_1 b_2+2\,b_1^3 
\,=\, 
6\,\beta_2-6\,\beta_0\beta_1+2\,\beta_0^3
\\
\BEE{4}  &\!=\!&
b_4
-4\,b_1 b_3
+6\,b_1^2 b_2
-3\,b_1^4
\,=\, 
24\,\beta_3
-24\,\beta_0 \beta_2
+12\,\beta_0^2\beta_1
-3\,\beta_0^4
\\
\BEE{5}  &\!=\!& 
b_5
-5\,b_1 b_4
+10\,b_1^2 b_3
-10\,b_1^3 b_2
+4\,b_1^5
\\
 &\!=\!&  
120\,\beta_4
-120\,\beta_0 \beta_3
+60\,\beta_0^2 \beta_2
-20\,\beta_0^3 \beta_1
+4\,\beta_0^5
\end{eqnarray*}
{\bf Invariance and homogeneousness 
under $f(\bu)\mapsto \lambda f(\gamma \bu+ \epsilon)$}.
\\
{ Invariance 
under $f(\bu)\mapsto f( \bu+ \epsilon)$}.
\begin{eqnarray*}
(x,n,\nu)\quad\stackrel{\partial_\epsilon}{\mapsto}\quad(1,0,-f_0)
&&
\\
\partial_\epsilon x_i=-1 \hspace{9 ex}(1\leq i \leq r) 
&&\quad\partial_\epsilon x_0=-1
\\
\partial_\epsilon \nu_i=-f_0\hspace{9 ex} (1\leq i \leq r) 
&&\quad
\partial_\epsilon \nu_0=\partial_\epsilon \underline{\nu}_0=-f_0 
\\
\partial_\epsilon f_s=(1+s)f_{1+s}\;\;\;\; (0\leq s<r)
& \partial_\epsilon f_r=0
&\quad\partial_\epsilon \mathbf{f}_s=0\;\; (0\leq s\leq r)
\\
\partial_\epsilon p_s=(1+s)p_{1+s}f_0\;\; (0\leq s<r)
& \partial_\epsilon p_r=0
&\quad\partial_\epsilon \mathbf{p}_s=0\;\; (0\leq s\leq r)
\end{eqnarray*}

{ Homogeneousness under $f(\bu)\mapsto f(\gamma \bu)$}.
\begin{eqnarray*}
(x,n,\nu) &\mapsto& (\gamma^{-1}x,\gamma n,\gamma^{-1}\nu)
\\
(f_s,\mathbf{f}_s)&\mapsto& (\gamma^s f_s,\gamma^s \mathbf{f}_s)
\\
(p_s,\mathbf{p}_s)&\mapsto& (\gamma^{s-r} p_s,\gamma^{s-r} \mathbf{p}_s)
\end{eqnarray*}
{ Homogeneousness 
under $f(\bu)\mapsto \lambda f(\bu)$}.
\begin{eqnarray*}
(x,n,\nu) &\mapsto& (x,\lambda^{-1} n,\lambda\nu)
\\
(f_s,\mathbf{f}_s)&\mapsto& (\lambda f_s,\lambda \mathbf{f}_s)
\\
(p_s,\mathbf{p}_s)&\mapsto& (\lambda^{r-s} p_s,\lambda^{r-s} \mathbf{p}_s)
\\
\beta_{s-1}&\mapsto& \lambda^{-s}\beta_{s-1}
\end{eqnarray*}
%%%%%%%%%%%%%%%%%%%%%%%%%%%%%%%%%%%%%%%%%%%%%%%%%%%%%%%%%%%%%%%%%%%%%%%%%%%%%%%%%%%%%%%%%%%%
%%%%%%%%%%%%%%%%%%%%%%%%%%%%%%%%%%%%%%%%%%%%%%%%%%%%%%%%%%%%%%%%%%%%%%%%%%%%%%%%%%%%%%%%%%%%
\subsection{Explicit ODEs for low-degree polynomial inputs $f$.}
To avoid glutting this section, we shall restrict ourselves to the standard choice for $\beta$ and mention only the covariant ODEs.
\footnote{But we keep extensive tables for all 8 cases $(v,c)\times(t,s,o,g)$ at the disposal of the interested reader. } Concretely, for all 
values of the $f$-dregree $r$ up to 4 we shall write down a complete set of {\it minimal} polynomials $P_{(d_{i},\delta_{i})}(n,\nu)$, of degrees 
$(d_{i},\delta_{i})$ in $(n,\nu)$, that generate all the other convariant  polynomials by non-commutative pre-multiplication by
covariant  polynomials in $(n,\nu)$.\footnote{In fact, all {\it variable} ODEs can also be expressed as suitable combinations of
the minimal {\it covariant} ODEs.}. For each $r$, the sequence
$$ (\underline{d},\overline{\delta})^1,\dots, (d_{i},\delta_{i})^{m_{i}},\dots, (\overline{d},\underline{\delta})^1$$
indicates the degrees $(d_{i},\delta_{i}) $ of all minimal spaces with their dimensions $m_{i}$, i.e. the number of polynomials in
them. For the extreme cases, right and left, that dimension is always 1.
\\

\noindent
{\bf Input $f$ of degree 1. }
\\
Invariant coefficients\,:
$\mgf_1:=f_1 $
\\
Covariant shift\,:
$\nu_0:=-\frac{1}{2}\,\frac{f_0^2}{f_1} $
\\
First leading polynomial (shifted)\,:\quad $\mgp(\nu)=p(\nu+\nu_0)= \nu$
\\
Second leading polynomial\,:\quad $\mgq(n)= n^2$
\\
Covariant differential equations\,:$ (1,2) $
$$\mathbf{P}_{(1 ,2 ) }(n,\nu)={P}_{(1 ,2 ) }(n,\nu\!+\!\nu_0)=
 n^2\,\nu+\frac{1}{2}\,n -\frac{1}{24}\,\mgf_1
$$
Variable differential equations\,:$(2,4)$
%%%%%%%%%%%%%%%%%%%%%%%%%%%%%%%%%%%%%%%%%%%%%%%%%%%%%%%%%%%%%%%%%%%%%%%%%%%%%%%
%%%%%%%%%%%%%%%%%%%%%%%%%%%%%%%%%%%%%%%%%%%%%%%%%%%%%%%%%%%%%%%%%%%%%%%%%%%%%%%
\\

\noindent
{\bf Input $f$ of degree 2. }
\\
Invariant coefficients\,:
$$
\mgf_0=f_0-\frac{1}{4}\frac{f_1^2}{f_2} \;=\; -\frac{1}{4}\,(x_1-x_2)^2\,f_2
\quad,\quad\mgf_2=f_2 
$$
Covariant shift\,:
$$
\nu_0=  -\frac{1}{2}\frac{f_0 f_1}{f_2}+\frac{1}{12}\frac{f_1^3}{f_2^2}
\,=\,-\frac{1}{12}\,(x_1+x_2)(x_1^2-4 x_1 x_2+x_2^2)\,f_2
$$
Leading scalar factor\,:
$$
\mgf_0=-\frac{1}{4}\,(x_1-x_2)^2\,f_2;
$$
First leading polynomial (shifted) 
$$
\mgp(\nu)=p(\nu+\nu_0)= \frac{4}{9}\frac{\mgf_0^3}{\mgf_2}+\nu^2 
$$
Second leading polynomial\,:
$$
 \mgq(n)= \frac{1}{6}\,\frac{\mgf_2}{\mgf_0}\,n^6+n^8
$$
Covariant differential equations\,:\;$ (2,8),(3,7)^2,(4,6) $
\begin{eqnarray*}
\mathbf{P}_{(2 ,8 ) }(n,\nu)={P}_{(2 ,8 ) }(n,\nu\!+\!\nu_0)=
n^8\,\mgf_0\,\mgp(\nu)
+n^7\,\mgf_0\,\nu
+n^6\,(\frac{1}{6}\mgf_2\,\nu^2+\frac{5}{27}\mgf_0^3+\frac{8}{9}\mgf_0)\hspace{20 ex}
\\
-n^5\,(\frac{1}{6}\mgf_2\,\nu)
+n^4\,(\frac{1}{54}\mgf_0^2\,\mgf_2-\frac{2}{27}\mgf_2)
-n^2\,(\frac{1}{972}\mgf_0\,\mgf_2^2)
-\frac{1}{583}\,\mgf_2^3
\hspace{40 ex}
\end{eqnarray*}
\begin{eqnarray*}
\mathbf{P}_{( 3, 7) }(n,\nu)={P}_{(3 ,7 ) }(n,\nu\!+\!\nu_0)=
n^7\,\mgf_0\,\mgp(\nu)
+n^6(-\frac{1}{32}\mgf_2\,\nu^3-\frac{1}{72}\mgf_0^3\,\nu+\mgf_0\,\nu)
\hspace{40 ex}
\\
+n^5\,(\frac{7}{32}\mgf_2\,\nu^2+\frac{8}{9}\mgf_0\!+\!\frac{2}{9}\mgf_0^3)
+n^4\,(-\frac{1}{288}\mgf_2\,\mgf_0^2\nu-\frac{41}{288}\mgf_2\,\nu)
\hspace{40 ex}
\\
+n^3\,(\frac{7}{432}\mgf_0^2\,\mgf_2-\frac{1}{18}\mgf_2)
-n\,(\frac{11}{7776}\mgf_0\,\mgf_2^2)
+\frac{1}{31104}\,\mgf_2^3\,\nu
\hspace{40 ex}
\end{eqnarray*}
\begin{eqnarray*}
\mathbf{P}^\dag_{( 3,7 ) }(n,\nu)={P}^\dag_{(3 , 7) }(n,\nu\!+\!\nu_0)=
n^7\,\mgf_0\,\nu\,\mgp(\nu)
+n^6\,(-\frac{5}{3}\mgf_0\,\nu^2-\frac{32}{27}\frac{\mgf_0^4}{\mgf_2})
\hspace{77 ex}
\\
+n^5\,(-\frac{7}{9}\mgf_0\,\nu+\frac{1}{9}\mgf_0^3\,\nu)
+n^4\,(\frac{2}{27}\mgf_0^3-\frac{16}{27}\mgf_0)
+n^2\,(\frac{1}{81}\mgf_0^2\,\mgf_2)
-n\,(\frac{1}{972}\mgf_0\,\mgf_2^2\nu)
-\frac{4}{729}\,\mgf_0\,\mgf_2^2
\hspace{60 ex}
\end{eqnarray*}
\begin{eqnarray*}
\mathbf{P}_{(4 ,6 ) }(n,\nu)={P}_{( 4, 6) }(n,\nu\!+\!\nu_0)=
n^6\,\mgf_0\,(\nu^2+\frac{416}{3}\frac{\mgf_0}{\mgf_2})\,\mgp(\nu)
\hspace{80. ex} &&
\\
+n^5\,(-13\,\mgf_0\,\nu^3+\frac{416}{3}\frac{\mgf_0^2}{\mgf_2}\,\nu
-\frac{56}{9}\frac{\mgf_0^4}{\mgf_2}\,\nu)
+n^4\,(\frac{356}{9}\mgf_0\,\nu^2+\frac{1}{9}\mgf_0^3\,\nu^2
         +\frac{3328}{27}\frac{\mgf_0^2}{\mgf_2}+\frac{1024}{27}\frac{\mgf_0^4}{\mgf_2})
\hspace{55 ex}&&
\\
+n^3\,(-\frac{148}{9}\mgf_0\,\nu-\frac{26}{27}\mgf_0^3\,\nu)
+n^2\,(\frac{158}{81}\mgf_0^3-\frac{16}{3}\mgf_0)
+n\,(\frac{1}{81}\mgf_0^2\,\mgf_2\nu)
-\frac{1}{972}\,\mgf_0\,\mgf_2^2\,\nu^2-\frac{161}{729}\mgf_0^2\,\mgf_2
\hspace{51 ex}&&
\end{eqnarray*}
\\
Variable differential equations\,:\;$(3,14),(4,11)^2,(5,10)^3,(6,9)^2,(9,8)$
\\
%%%%%%%%%%%%%%%%%%%%%%%%%%%%%%%%%%%%%%%%%%%%%%%%%%%%%%%%%%%%%%%%%%%%%%%%%%%%%%%
%%%%%%%%%%%%%%%%%%%%%%%%%%%%%%%%%%%%%%%%%%%%%%%%%%%%%%%%%%%%%%%%%%%%%%%%%%%%%%%

\noindent
{\bf Input $f$ of degree 3. }
\\
Invariant coefficients\,:
\begin{eqnarray*}
\mgf_0 &\!\!=\!\!& f_0-\frac{1}{3}\frac{f_1\,f_2}{f_3}+\frac{2}{27}\frac{f_2^3}{f_3^2} 
=
\frac{1}{27}(x_1+x_2-2\,x_3)(x_2+x_3-2\,x_1)(x_3+x_1-2\,x_2)f_3
\\
 \mgf_1 &\!\!=\!\!& f_1-\frac{1}{3}\frac{f_2^2}{f_3}
= -\frac{1}{3}(x_1^2+x_2^2+x_3^2-x_1x_2-x_2x_3-x_3x_1)f_3
\\
\mgf_3 &\!\!=\!\!& f_3
\end{eqnarray*}
Covariant shift\,:
\begin{eqnarray*}
\nu_0 &\!\!=\!\!&  -\frac{1}{3}\frac{f_0f_2}{f_3}/
+\frac{1}{18}\frac{f_1f_2^2}{f_3^2}
-\frac{1}{108}\frac{f_2^4}{f_3^3}
=  -\frac{1}{108}
(x_1+x_2+x_3)(x_1^3+x_2^3+x_3^3+
\\ &&
\;\;\;24\,x_1x_2x_3-3\,x_1x_2^2-3\,x_1^2x_2-3\,x_3x_1^2-3\,x_3^2x_1-3\,x_3x_2^2-3\,x_3^2x_2)
f_3
\end{eqnarray*}
Leading scalar factor\,:
\begin{eqnarray*}
 \mga:=4\,\mgf_1^3+27\,\mgf_0^2\,\mgf_3=-(x_1-x_2)^2\,(x_2-x_3)^2\,(x_3-x_2)^2\,f_3^3
\end{eqnarray*}
First leading polynomial (shifted) 
$$
\mgp(\nu)=p(\nu+\nu_0)=
\frac{1}{32}\frac{\mgf_0^2\mgf_1^3}{\mgf_3^2}
+\frac{27}{64}\frac{\mgf_0^4}{\mgf_3}
+\Big(\frac{9}{8}\frac{\mgf_0^2\mgf_1}{\mgf_3}+\frac{1}{16}\frac{\mgf_1^4}{\mgf_3^2}\Big)\nu
+\frac{1}{2}\frac{\mgf_1^2}{\mgf_3}\nu^2+\nu^3
$$
Second leading polynomial\,:
$$
\mgq(n)= 
\frac{9}{4}\,\frac{\mgf_1\,\mgf_3^2}{\mga}\,n^{12}
+\frac{81}{4}\,\frac{\mgf_3^2}{\mga}\,n^{13}
+6\,\frac{\mgf_1^2\,\mgf_3}{\mga}\,n^{14}
+3\,\frac{\mgf_1\,\mgf_3}{\mga}\,n^{15}
+n^{16}
$$
Covariant differential equations\,:\;$ (3,16),(4,14)^2,(5,13)^2,(7,12)$
$$
\mathbf{P}_{(3,16)}(n,\nu)={P}_{(3,16)}(n,\nu+\nu_0)=  \,\mga\,\mgp(\nu)+O(n^{15})\,O(\nu^3)
\hspace{25 ex}
$$
\begin{eqnarray*}
\mathbf{P}_{(4,14)}(n,\nu) &={P}_{(4,14)}(n,\nu+\nu_0) &=
 \mgf_1\,\mga\,\mgb\,n^{14}\,   p(\nu)+O(n^{13})\,O(\nu^4)
\\
\mathbf{P}^\dag_{(4,14)}(n,\nu)   &={P}^\dag_{(4,14)}(n,\nu+\nu_0) &=
       \mga\,\mgb\,n^{14}\nu\,p(\nu)+O(n^{13})\,O(\nu^4)\hspace{15 ex}
\end{eqnarray*}
with the following invariant coefficient $\mgb$\,:
\begin{eqnarray*}
\mgb&\!\!:=\!&
2097152\,\mgf_1^{12}
-766779696\,\mgf_1^3\,\mgf_3^3
-520497152\,\mgf_1^9\,\mgf_3
-36074005128\,\mgf_0^4\,\mgf_1^3\,\mgf_3^3
\\&+\!\!\!\!&  
1428879744\,\mgf_1^6\,\mgf_3^2\,
-1314579456\,\mgf_0^6\mgf_1^3\,\,\mgf_3^3\,
+1099865088\,\mgf_0^4\,\mgf_1^6\,\mgf_3^2
\\&+\!\!\!\!&
205963264\,\mgf_1^9\,\mgf_0^2\,\mgf_3
-8872609536\,\mgf_0^2\,\mgf_1^6\,\mgf_3^2
+73222472421\,\mgf_0^4\,\mgf_3^4
\\&+\!\!\!\!&
20602694736\,\mgf_0^2\,\mgf_1^3\,\mgf_3^3
+5971968\,\mgf_0^6\,\mgf_1^6\,\mgf_3^2
+884736\,\mgf_0^4\,\mgf_1^9\,\mgf_3
-5165606520\,\mgf_0^2\,\mgf_3^4
\end{eqnarray*}
\noindent
\begin{eqnarray*}
\mathbf{P}_{(5,13)}(n,\nu)    
      \!\!&\!\!={P}_{(5,13)}(n,\nu+\nu_0)    
\!\!&= n^{13}\,(\mgf_1^2\,\mgc_1-180\,\mgf_3\,\mgc_2\,\nu)\,\mgp(\nu)+O(n^{12})\,O(\nu^5)
\\
\mathbf{P}^\dag_{(5,13)}(n,\nu)
      \!\!&\!\!={P}^\dag_{(5,13)}(n,\nu+\nu_0)
\!\!&= n^{13}\,(\mgc_3\,\nu+180\,\mgf_1\,\mgf_3\,\mgc_1\,\nu^2)\,\mgp(\nu)+O(n^{12})\,O(\nu^5)
\end{eqnarray*}
with the following invariant coefficients $\mgc_1,\mgc_2,\mgc_3$\,:
\begin{eqnarray*}
\mgc_1 &\!\!:=\!\!&
917290620205793280\,\mgf_0^2\,\mgf_1^{12}\,\mgf_3^2
+78717609050112\,\mgf_0^6\,\mgf_1^{12}\,\mgf_3^2 \hspace{20 ex}
\\&+\!\!\!\!&
4163751641088\,\mgf_0^4\,\mgf_1^{15}\,\mgf_3
+50281437903388672\,\mgf_1^{15}\,\mgf_3
\\&+\!\!\!\!&
1581069280739328\,\mgf_0^2\,\mgf_1^{15}\,\mgf_3
+17755411807125504\,\mgf_0^4\,\mgf_1^{12}\,\mgf_3^2
\\&+\!\!\!\!&
99407759207731200\,\mgf_0^6\,\mgf_1^9\,\mgf_3^3
+5640800181652267776\,\mgf_0^4\,\mgf_1^9\,\mgf_3^3
\\&+\!\!\!\!&
344140580192256\,\mgf_0^8\,\mgf_1^9\,\mgf_3^3
+11726669550606570432\,\mgf_0^6\,\mgf_1^6\,\mgf_3^4
\\&+\!\!\!\!&
326589781381042176\,\mgf_0^8\,\mgf_1^6\,\mgf_3^4
-498496347843530688\,\mgf_1^6\,\mgf_3^4
\\&+\!\!\!\!&
16926659444736\,\mgf_0^{10}\,\mgf_1^6\,\mgf_3^4
-85405328111733120\,\mgf_0^8\,\mgf_1^3\,\mgf_3^5
\\&-\!\!\!\!&
1691608028258304\,\mgf_0^{10}\,\mgf_1^3\,\mgf_3^5
-15390509185018432260\,\mgf_0^6\,\mgf_1^3\,\mgf_3^5
\\&+\!\!\!\!&
98766738625551624\,\mgf_0^8\,\mgf_3^6
-7432537028329878624\,\mgf_0^2\,\mgf_1^3\,\mgf_3^5
\\&+\!\!\!\!&
10331678048256\,\mgf_1^{18}
-27319961213550950355\,\mgf_0^4\,\mgf_3^6
\\&+\!\!\!\!&
1500717585045441600\,\mgf_0^2\,\mgf_3^6
+226960375516131600\,\mgf_1^3\,\mgf_3^5
\\&-\!\!\!\!&
88258622384581632\,\mgf_1^{12}\,\mgf_3^2
-8253051882421560660\,\mgf_0^4\,\mgf_1^6\,\mgf_3^4
\\&-\!\!\!\!&
1478991931831367424\,\mgf_0^2\,\mgf_1^9\,\mgf_3^3
-450793967617700928\,\mgf_3^3\,\mgf_1^9
\\&+\!\!\!\!&
3821964710670454374\,\mgf_0^6\,\mgf_3^6
-17792355610879876332\,\mgf_0^4\,\mgf_1^3\,\mgf_3^5
\\&-\!\!\!\!&
5761805211034236864\,\mgf_0^2\,\mgf_1^6\,\mgf_3^4
\end{eqnarray*}
\begin{eqnarray*}
\mgc_2 &\!\!:=\!\!& 
346056266653347168\,\mgf_0^2\,\mgf_1^3\,\mgf_3^5
-561701191680\,\mgf_1^{18}
\\&-\!\!\!\!&
12492140160024576\,\mgf_0^2\,\mgf_1^{12}\,\mgf_3^2
-1568573227008\,\mgf_0^6\,\mgf_1^{12}\,\mgf_3^2 \hspace{20 ex}
\\&-\!\!\!\!&
225501511680\,\mgf_0^4\,\mgf_1^{15}\,\mgf_3
-812740325605376\,\mgf_1^{15}\,\mgf_3
\\&-\!\!\!\!&
83614219370496\,\mgf_0^2\,\mgf_1^{15}\,\mgf_3
-1675724924436480\,\mgf_0^4\,\mgf_1^{12}\,\mgf_3^2
\\&-\!\!\!\!&
7670187447717888\,\mgf_0^6\,\mgf_1^9\,\mgf_3^3
-50653665706906368\,\mgf_0^4\,\mgf_1^9\,\mgf_3^3
\\&-\!\!\!\!&
313456656384\,\mgf_0^8\,\mgf_1^9\,\mgf_3^3
-23099325980303808\,\mgf_0^6\,\mgf_1^6\,\mgf_3^4
\\&+\!\!\!\!&
24015789981646272\,\mgf_1^6\,\mgf_3^4
+28152325951488\,\mgf_0^8\,\mgf_1^6\,\mgf_3^4
\\&-\!\!\!\!&
64268410079232\,\mgf_0^8\,\mgf_1^3\,\mgf_3^5
-24631997881011588\,\mgf_0^6\,\mgf_1^3\,\mgf_3^5
\\&+\!\!\!\!&
1250788627474992675\,\mgf_0^4\,\mgf_3^6
+308789626552560\,\mgf_1^3\,\mgf_3^5
\\&+\!\!\!\!&
1527911696305152\,\mgf_1^{12}\,\mgf_3^2
+60223210699403700\,\mgf_0^4\,\mgf_1^6\,\mgf_3^4
\\&+\!\!\!\!&
19863849419539200\,\mgf_0^2\,\mgf_1^9\,\mgf_3^3
+18320806414937664\,\mgf_1^9\,\mgf_3^3
\\&+\!\!\!\!&
3706040377703682\,\mgf_0^6\,\mgf_3^6
+945719068633781580\,\mgf_0^4\,\mgf_1^3\,\mgf_3^5
\\&+\!\!\!\!&
263746805911956288\,\mgf_0^2\,\mgf_1^6\,\mgf_3^4
\end{eqnarray*}

\begin{eqnarray*}
\mgc_3 &\!\!:=\!\!&
188945409245184\,\mgf_1^{21}
-4265434334643431940864\,\mgf_0^8\,\mgf_1^6\,\mgf_3^5 
\hspace{25 ex} 
\\&+\!\!\!\!&
10475616970801152\,\mgf_0^8\,\mgf_1^{12}\,\mgf_3^3
-6409779863684795640576\,\mgf_0^2\,\mgf_1^{12}\,\mgf_3^3
\\&+\!\!\!\!&
1643585979933664752384\,\mgf_0^4\,\mgf_1^{12}\,\mgf_3^3
+31374503650787328\,\mgf_0^{10}\,\mgf_1^9\,\mgf_3^4
\\&+\!\!\!\!&
3977644977060256658880\,\mgf_0^6\,\mgf_1^9\,\mgf_3^4
+1290160568497752489024\,\mgf_1^9\,\mgf_3^4
\\&+\!\!\!\!&
1205510242413389496768\,\mgf_0^2\,\mgf_1^9\,\mgf_3^4
-48405699843949469920500\,\mgf_0^4\,\mgf_1^9\,\mgf_3^4
\\&+\!\!\!\!&
129087554262282601920\,\mgf_1^{12}\,\mgf_3^3
-280615966839399140352\,\mgf_1^{15}\,\mgf_3^2
\\&+\!\!\!\!&
9973443990092156928\,\mgf_0^6\,\mgf_1^{12}\,\mgf_3^3 
-1492256344300529883948\,\mgf_0^4\,\mgf_1^6\,\mgf_3^5
\\&+\!\!\!\!&
914039610015744\,\mgf_0^{12}\,\mgf_1^6\,\mgf_3^5
-122919033279447568214604\,\mgf_0^6\,\mgf_1^6\,\mgf_3^5
\\&+\!\!\!\!&
29574529753446346752\,\mgf_0^8\,\mgf_1^9\,\mgf_3^4
+13297895549157703680\,\mgf_0^{10}\,\mgf_1^6\,\mgf_3^5
\\&+\!\!\!\!&
30311992402755395296032\,\mgf_0^2\,\mgf_1^6\,\mgf_3^5
-24763502547539307564426\,\mgf_0^6\,\mgf_1^3\,\mgf_3^6
\\&+\!\!\!\!&
169418141509536645120\,\mgf_0^2\,\mgf_1^{15}\,\mgf_3^2 
-908625541799649020400\,\mgf_1^6\,\mgf_3^5
\\&+\!\!\!\!&
5603533051087967184\,\mgf_0^{10}\,\mgf_3^7
-67837564266533409024\,\mgf_0^{10}\,\mgf_1^3\,\mgf_3^6
\\&+\!\!\!\!&
3613883506978117107600\,\mgf_0^8\,\mgf_3^7
-41025866981179759740000\,\mgf_0^4\,\mgf_3^7
\\&+\!\!\!\!&
581286688880237080992900\,\mgf_0^6\,\mgf_3^7
-12210659652336342667200\,\mgf_0^2\,\mgf_1^3\,\mgf_3^6
\\&+\!\!\!\!&
231937500459010111367085\,\mgf_0^4\,\mgf_1^3\,\mgf_3^6
-12977326722621245045184\,\mgf_0^8\,\mgf_1^3\,\mgf_3^6
\\&+\!\!\!\!&
1049410426382106624\,\mgf_0^4\,\mgf_1^{15}\,\mgf_3^2
-100601484577357824\,\mgf_0^{12}\,\mgf_1^3\,\mgf_3^6 
\\&+\!\!\!\!&
27598162056708096\,\mgf_0^2\,\mgf_1^{18}\,\mgf_3
+76397618921472\,\mgf_0^4\,\mgf_1^{18}\,\mgf_3
\\&+\!\!\!\!&
4008794200000102400\,\mgf_1^{18}\,\mgf_3
+1381995569479680\,\mgf_0^6\,\mgf_1^{15}\,\mgf_3^2
\end{eqnarray*}
\begin{eqnarray*}
\mathbf{P}_{(7,12)}(n,\nu)\!=\!P_{(7,12)}(n,\nu\!+\!\nu_0)\!=\!
n^{12}p(\nu)\,(\mgf_1^4\,\nu^4\!+\!\mgd_3\,\nu^3\!+\!\mgf_1^2\,\mgd_2\,
\nu^2\!+\!\mgf_1\,\mgd_1\,\nu\!+\!\mgd_0)\!
\\
+ O(n^{11})\,O(\nu^7) \hspace{0 ex}
\end{eqnarray*}
with the following invariant coefficients $\mgd_0,\mgd_1,\mgd_2,\mgd_3$\,:
\begin{eqnarray*}
\mgd_0 &\!\!:=\!\!&
+\frac{29859111}{128}\,\mgf_0^2
+\frac{3664683}{1600}\,\mgf_0^4
+\frac{29889}{4000}\,\mgf_0^6
+\frac{81}{10000}\,\mgf_0^8
+\frac{22240737}{640}\frac{\mgf_1^3}{\mgf_3}\hspace{20 ex}
\\ &\!\!\!\!&
+\frac{336626989}{3200}\frac{\mgf_0^2\,\mgf_1^3}{\mgf_3}
+\frac{3493333}{24000}\frac{\mgf_0^4\,\mgf_1^3}{\mgf_f3}
+\frac{159}{5000}\frac{\mgf_1^3\,\mgf_0^6}{\mgf_3}
+\frac{1969}{60000}\frac{\mgf_0^4\,\mgf_1^6}{\mgf_3^2}
\\ &\!\!\!\!&
+\frac{242977752829}{15552000}\frac{\mgf_1^6}{\mgf_3^2}
+\frac{40541647}{1296000}\frac{\mgf_0^2\,\mgf_1^6}{\mgf_3^2}
+\frac{15317}{4860000}\frac{\mgf_0^2\,\mgf_1^9}{\mgf_3^3}
\\ &\!\!\!\!&
+\frac{203363491}{69984000}\frac{\mgf_1^9}{\mgf_3^3}
+\frac{83521}{1049760000}\frac{\mgf_1^{12}}{\mgf_3^4}
\end{eqnarray*}
\begin{eqnarray*}
\mgd_1 &\!\!:=\!\!&
-\frac{368631}{160}
+\frac{3305043}{800}\,\mgf_0^2
+\frac{93339}{2000}\,\mgf_0^4
+\frac{27}{250}\,\mgf_0^6
-\frac{642277459}{1296000}\,\frac{\mgf_1^3}{\mgf_3}
+\frac{123}{500}\,\frac{\mgf_0^4\,\mgf_1^3}{\mgf_3}\hspace{20 ex}
\\ &\!\!\!\!&
+\frac{10657943}{18000}\,\frac{\mgf_0^2\,\mgf_1^3}{\mgf_3}
+\frac{697}{9000}\,\frac{\mgf_0^2\,\mgf_1^6}{\mgf_3^2}
+\frac{101072021}{1944000}\,\frac{\mgf_1^6}{\mgf_3^2}
+\frac{4913}{1458000}\,\frac{\mgf_1^9}{\mgf_3^3}
\end{eqnarray*}
\begin{eqnarray*}
\mgd_2 &\!\!:=\!\!&
+\frac{361809}{800}
+\frac{10467}{200}\,\mgf_0^2
+\frac{27}{50}\,\mgf_0^4
+\frac{479929}{2160}\,\frac{\mgf_1^3}{\mgf_3}
+\frac{29}{50}\,\frac{\mgf_0^2\mgf_1^3}{\mgf_3}
+\frac{289}{5400}\,\frac{\mgf_1^6}{\mgf_3^2}
\hspace{19 ex}
\end{eqnarray*}
\begin{eqnarray*}
\mgd_3 &\!\!:=\!\!&
-\frac{6561}{40}\,\mgf_3
-\frac{1347}{20}\,\mgf_1^3
+\frac{6}{5}\,\mgf_0^2\,\mgf_1^3
+\frac{17}{45}\,\frac{\mgf_1^6}{\mgf_3}
\hspace{59 ex}
\end{eqnarray*}
Variable differential equations\,:
$$(4,28),(5,22)^2,(6,20)^3,(7,19)^4,(8,18)^3,(10,17)^2,(16,16)$$
%%%%%%%%%%%%%%%%%%%%%%%%%%%%%%%%%%%%%%%%%%%%%%%%%%%%%%%%%%%%%%%%%%%%%%%%%%%%%%%
%%%%%%%%%%%%%%%%%%%%%%%%%%%%%%%%%%%%%%%%%%%%%%%%%%%%%%%%%%%%%%%%%%%%%%%%%%%%%%%
\\

\noindent
{\bf Input $f$ of degree 4. }
\\
Invariant coefficients\,:
\begin{eqnarray*}
\mgf_0 &\!\!=\!\!&
f0-\frac{1}{4}\frac{f_1f_3}{f4}+\frac{1}{16}\frac{f_2\,f_3^2}{f_4^2}
-\frac{3}{256}\frac{f_3^4}{f_4^3}
=\frac{1}{256}\prod_{i=1}^{i=4}(4x_i\!-\!x_1\!-\!x_2\!-\!x_3\!-\!x_4)
\\
\mgf_1 &\!\!=\!\!& f_1-\frac{1}{2}\frac{f_2f_3}{f_4}+\frac{1}{8}\frac{f_3^3}{f_4^2}
 =
-\frac{1}{8}(x_1\!+\!x_2\!-\!x_3\!-\!x_4)(x_2\!+\!x_3\!-\!x_1\!-\!x_4)
(x_2\!+\!x_4\!-\!x_1\!-\!x_3)\,f_4
\\
\mgf_2 &\!\!=\!\!& f_2-\frac{3}{8}\frac{f_3^2}{f_4}
= -\frac{1}{8}\,\Big(4(x_1^2+x_2^2+x_3^2+x_4^2) -(x_1+x_2+x_3+x_4)^2\Big)\,f_4
\\
\mgf_4 &\!\!=\!\!& f_4
\end{eqnarray*}
Covariant shift\,:
$$
\nu_0= -\frac{1}{4}\frac{f_0f_3}{f_4}
+\frac{1}{32}\frac{f_1f_3^2}{f_4^2}
-\frac{1}{192}\frac{f_2f_3^3}{f_4^3}
+\frac{1}{1280}\frac{f_3^5}{f_4^4}
$$
Leading scalar factors:
\begin{eqnarray*}
\mga &\!\!=\!\!& 256\,\mgf_0^3\,\mgf_4^2 -128\,\mgf_0^2\,\mgf_2^2\,\mgf_4
+16\,\mgf_0\,\mgf_2^4+144\,\mgf_0\,\mgf_2\,\mgf_1^2\,\mgf_4
-4\,\mgf_1^2\,\mgf_2^3-27\,\mgf_1^4\,\mgf_4
\\
&\!\!=\!\!& \prod_{1\leq i< j\leq 4}(x_i-x_j)^2\,f_4^5
\\
\mgb &\!\!=\!\!& -1280\,\mgf_2^6
+32256\,\mgf_0\,\mgf_2^4\,\mgf_4
-269568\,\mgf_0^2\,\mgf_2^2\,\mgf_4^2
+746496\,\mgf_0^3\,\mgf_4^3
+69984\,\mgf_0\,\mgf_1^2\,\mgf_2\,\mgf_4^2
\\&&
-9504\,\mgf_1^2\,\mgf_2^3\,\mgf_4
+19683\,\mgf_1^4\,\mgf_4^2
\\
&\!\!=\!\!&\prod_{{1\leq i < j \leq 4}\atop {1\leq k < l \leq 4}}^{(i,j)\not=(k,l)}
\frac{1}{128}\;\Big(5\,(x_i+x_j-x_k-x_l)^2+(x_i-x_j)^2-5\,(x_k-x_l)^2 \Big)\;f_4^6
\end{eqnarray*}
First leading polynomial (shifted) 
\begin{eqnarray*}
\mgp(\nu)&\!\!=\!\!&
\frac{12}{125}\frac{\mgf_0^3\mgf_1^2\mgf_2}{\mgf_4^2}
-\frac{27}{2000}\frac{\mgf_0^2\mgf_1^4}{\mgf_4^2}
+\frac{256}{625}\frac{\mgf_0^5}{\mgf_4}
+\frac{16}{2025}\frac{\mgf_0^3\mgf_2^4}{\mgf_4^3}
-\frac{128}{1125}\frac{\mgf_0^4\mgf_2^2}{\mgf_4^2}
-\frac{1}{675}\frac{\mgf_0^2\mgf_1^2\mgf_2^3}{\mgf_4^3}
\\&\!\!+\!\!&\!\!
\Big(
\frac{32}{25}\frac{\mgf_0^3\mgf_1}{\mgf_4}
-\frac{56}{225}\frac{\mgf_0^2\mgf_1\mgf_2^2}{\mgf_4^2}
+\frac{21}{100}\frac{\mgf_0\mgf_1^3\mgf_2}{\mgf_4^2}
-\frac{27}{1000}\frac{\mgf_1^5}{\mgf_4^2}
+\frac{4}{225}\frac{\mgf_0\mgf_1\mgf_2^4}{\mgf_4^3}
-\frac{2}{675}\frac{\mgf_1^3\mgf_2^3}{\mgf_4^3}
\Big)\,\nu
\\&\!\!+\!\!&\!\!
\Big(
\frac{16}{15}\frac{\mgf_0^2\mgf_2}{\mgf_4}
+\frac{9}{10}\frac{\mgf_0\mgf_1^2}{\mgf_4}
+\frac{11}{60}\frac{\mgf_1^2\mgf_2^2}{\mgf_4^2}
-\frac{4}{15}\frac{\mgf_0\mgf_2^3}{\mgf_4^2}
+\frac{4}{225}\frac{\mgf_2^5}{\mgf_4^3}
\Big)\,\nu^2
+\frac{\mgf_1\mgf_2}{\mgf_4}\,\nu^3
+\nu^4
\end{eqnarray*} 
Second leading polynomial\,:
$$
\mgq(n)= 
-\frac{2^{14}\,7}{3^3\,5^6}\frac{\mgf_4^{11}}{\mga\mgb}\,n^{20}
-\frac{2^{13}\,11}{3^2\,5^5}\,\frac{\mgf_2\,\mgf_4^{10}}{\mga\,\mgb}\,n^{22}
+\dots
+8\,\Big(
\frac{\mgf_1\,\mgf_2^2\,\mgf_4}{\mga}
+12\,\frac{\mgf_0\,\mgf_1\,\mgf_4^2}{\mga}
\Big)\,n^{43}
+n^{44}
$$
Covariant differential equations\,:
$$ (4,44),
(5,32)^2,
(6,28)^3,
(7,26)^4,
(8,24),
(10,23)^4,
(11,22)^3,
(16,21)^2,
(28,20) $$
$$
\mathbf{P}_{(4,44)}(n,\nu)=\mathbf{P}_{(4,44)}(n,\nu+\nu_0)=
\mga\,\mgb\,n^{44}\,\mgp(\nu)+O(n^{43})\,O(\nu^4)\hspace{40 ex}
$$
Variable differential equations\,:
$$
(5,64),
(6,44)^2,
(7,37)^2,
(8,34)^4,
(9,32)^5,
(10,30)^2,
$$
$$
(11,29)^2,
(13,28)^5,
(15,27)^4,
(18,26)^2,
(25,25)^2,
(45,24)
$$
%%%%%%%%%%%%%%%%%%%%%%%%%%%%%%%%%%%%%%%%%%%%%%%%%%%%%%%%%%%%%%%%%%%%%%%%%%%%%%%%%%%%%%%%%%%%
%%%%%%%%%%%%%%%%%%%%%%%%%%%%%%%%%%%%%%%%%%%%%%%%%%%%%%%%%%%%%%%%%%%%%%%%%%%%%%%%%%%%%%%%%%%%
\subsection{The global resurgence picture for polynomial inputs $f$.}
The covariant ODEs enable us to describe the exact singular behaviour of $\imi{k}\!\!(\nu)=h(\nu)$ at infinity in the $\nu$-plane, and by way of consequence  all singularities {\it over} 0 in the $\zeta$-plane. In the $\nu$-plane, the singularities in question consist of linear combinations of rather elementary exponential factors multiplied by series in negative powers of $\nu$. These are always divergent, resurgent,
and resummable. The case of radial inputs $f$ (i.e. $f(x)=f_{r}\,x^r$) is predictably much simpler and deserves special mention. We find\,:
\begin{eqnarray}  \label{a143}
\Big(\sum_{r+1\leq k} c_{s}(\omega)\,\nu^{-\frac{s}{r+1}}\Big)\!&\!\!&\!
\exp\Big(\omega\,\nu^{\frac{r}{r+1}}\Big)
\hspace{20 ex} (\mi{for\; radial} \; f)\quad
\\  \label{a144}
\Big(\sum_{r+1\leq k} c_{s}(\omega)\,\nu^{-\frac{s}{r+1}}\Big)\!&\!\!&\!
\exp\Big(\omega\,\nu^{\frac{r}{r+1}}+\sum_{s=1}^{r-2}\omega_s\,\nu^{\frac{s}{r+1}}\Big)
\quad \quad\quad (\mi{for\, general} \; f)\quad
\end{eqnarray}
The ``leading'' frequencies $\omega$ featuring in the exponential factors depend only on the leading coefficient $f_{r}$ of $f$.
Via the variable $\theta$ thus defined\,:
\begin{equation}  \label{a145}
\theta:=\Big(\frac{r+1}{r}\Big)^{\!r}\,\;\frac{\mathbf{f}_r}{\omega^{r+1}}
\;=\;
\Big(\frac{r+1}{r}\Big)^{\!r}\,\frac{f_r}{\omega^{r+1}}\;
\end{equation}
the leading frequencies $\omega$ correspond, for each degree $r$, to the roots of the following polynomials $\mg{\pi}_{r}(\theta)$ of degree
$r$\,:
\begin{eqnarray*}
\mg{\pi}_1(\theta) &\!=\!& -12+\theta 
\\
\mg{\pi}_2(\theta) &\!=\!& -432+\theta^2= -2^4\,3^3+\theta^2 
\\
\mg{\pi}_3(\theta) &\!=\!& (240+7\,\theta)\,(-30+\theta)^2=
(2^4\,3\times5+7\,\theta)\,(-2\times3\times5+\theta)^2
\\
\mg{\pi}_4(\theta) &\!=\!& (1749600000-1620000\,\theta^2+343\,\theta^4)
\\
                     &\!=\!& (2^8\,3^7\,5^5-2^5\,3^4\,5^4\,\theta^2+7^3\,\theta^4)
\\
\mg{\pi}_5(\theta) &\!=\!& (-1344+31\,\theta)\,(189+\theta)^2\,(42+\theta)^2
\\
    &\!=\!& (-2^6\,3\times7+31\,\theta)\,(3^3\,7+\theta)^2\,(2\times3\times7+\theta)^2 
\\
\mg{\pi}_6(\theta) &\!=\!& (-66395327975424+152320630896\,\theta^2-116688600\,\theta^4+29791\,\theta^6)
\\
       &\!=\!& (-2^{12}\,3^9\,7^7+2^4\,3^7\,7^6\,37\,\theta^2-
2^3\,3^5\,5^2\,7^4\,\theta^4+31^3\,\theta^6) 
\\
\mg{\pi}_7(\theta) &\!=\!& (3840+127\,\theta)\,(-30+\theta)^2\,(24300+1080\,\theta+37\,\theta^2)^2
\\
   &\!=\!&
(2^8\,3\times
5+127\,\theta)\,(-2\times3\times5+\theta)^2\,(2^2\,3^5\,5^2+2^3\,3^3\,5\,\theta+37\,\theta^2)^2
\end{eqnarray*}
For a non-standard choice of $\beta$ and with the `centered' coefficients $\BEE{i}$ introduced at the end of \S6.2,
these polynomials $\mg{\pi}_{r}(\theta)$ become\,:
\begin{eqnarray*}
\mg{\pi}_1(\theta) &\!=\!& 
1
+\BEE{2}\,\theta
\\[1.5 ex]
\mg{\pi}_2\,(\theta) &\!=\!& 
1
+2\,\BEE{3}\theta
+(\BEE{3}^2
+4\,\BEE{2}^3)\,\theta^2 
\\[1.5 ex]
\mg{\pi}_3(\theta) &\!=\!&  
1
+3\,(\BEE{4} 
-6\,\BEE{2}^2)\,\theta
+3\,(\BEE{4}^2
+ 18\,\BEE{2}\BEE{3}^2
-12\,\BEE{2}^2\BEE{4}
+ 27\,\BEE{2}^4)\,\theta^2
\\&\!\!&+\,
( \BEE{4}^3 
-27\,\BEE{3}^4
+54\,\BEE{2}\BEE{3}^2\BEE{4}
-18\,\BEE{2}^2\BEE{4}^2 
-54\,\BEE{2}^3\BEE{3}^2
+81\,\BEE{2}^4\BEE{4})\,\theta^3 
\\[1.5 ex]
\mg{\pi}_4(\theta) &\!=\!&  
1
+4\,(\BEE{5}
-30\,\BEE{2}\BEE{3})\,\theta
+2\,(3\,\BEE{5}^2
+80\,\BEE{2}\BEE{4}^2
+180\,\BEE{3}^2\,\BEE{4}
-180\,\BEE{2}\BEE{3}\BEE{5}
\\&&
+1320\,\BEE{2}^2\BEE{3}^2
-720\,\BEE{2}^3\BEE{4}
+1728\,\BEE{2}^5)\,\theta^2
+4\,(\BEE{5}^3 
-160\,\BEE{3}\BEE{4}^3
\\&&
-90\,\BEE{2}\BEE{3}\BEE{5}^2
+180\,\BEE{3}^2\BEE{4}\BEE{5}
+80\,\BEE{2}\BEE{4}^2\BEE{5}
+1120\,\BEE{2}^2\BEE{3}\BEE{4}^2
+864\,\BEE{3}^5
\\&&
-2520\,\BEE{2}\BEE{3}^3\BEE{4}
+1320\,\BEE{2}^2\BEE{3}^2\BEE{5}
-720\,\BEE{2}^3\BEE{4}\BEE{5}
+ 1728\,\BEE{2}^5\BEE{5}
+1280\,\BEE{2}^3\BEE{3}^3
\\&&
-2880\,\BEE{2}^4\BEE{3}\BEE{4})\,\theta^3
+(\BEE{5}^4 
-120\,\BEE{2}\BEE{3}\BEE{5}^3
+160\,\BEE{2}\BEE{4}^2\BEE{5}^2
+360\,\BEE{3}^2\BEE{4}\BEE{5}^2
\\&&
-640\,\BEE{3}\BEE{4}^3\BEE{5}
+256\,\BEE{4}^5 
-2560\,\BEE{2}^2\BEE{4}^4
+3456\,\BEE{3}^5\BEE{5}
+5760\,\BEE{2}\BEE{3}^2\BEE{4}^3
\\&&
+2640\,\BEE{2}^2\BEE{3}^2\BEE{5}^2
-1440\,\BEE{2}^3\BEE{4}\BEE{5}^2
-2160\,\BEE{3}^4\BEE{4}^2
+4480\,\BEE{2}^2\BEE{3}\BEE{4}^2\BEE{5}
\\&&
-10080\,\BEE{2}\BEE{3}^3\BEE{4}\BEE{5}
+ 3456\,\BEE{2}^5\BEE{5}^2
-3200\,\BEE{2}^3\BEE{3}^2\BEE{4}^2
+5120\,\BEE{2}^3\BEE{3}^3\BEE{5}
\\&&
+6400\,\BEE{2}^4\BEE{4}^3
-11520\,\BEE{2}^4\BEE{3}\BEE{4}\BEE{5}
)\,\theta^4 
\end{eqnarray*}
Thus, for $r=1$ we get two basic singular summands\,:
\begin{eqnarray}  \label{a146}
&&\Big(\sum_{2\leq k} c_{s}(\omega)\,\nu^{-\frac{s}{2}}\Big)
\exp\Big(\omega\,\nu^{\frac{1}{2}}\Big)\hspace{14 ex} (\mi{for\,all\; f\; of\; degree\;1})
\quad
\\  \label{a147}
&&\Big( {\nu}^{-1}-\frac{1}{\omega}\,{\nu^{-\frac{1}{2}}}\Big)
\exp\Big(\omega\,\nu^{\frac{1}{2}}\Big)\hspace{14 ex} (\mi{if}\; f(x)=f_1\,x)
\end{eqnarray}
with   frequencies $\omega$ corresponding to the solutions of
$\mg{\pi}_1(\theta)=0 $ i.e. 
$\mg{\pi}_1(2\,\frac{\mathbf{f}_1}{\omega^2})=0 $ i.e. $\omega=(-2\,\BEE{2}f_1)^{\frac{1}{2}}$.
\\

\noindent
For $r=2$ we have $6=2\times 3$ basic summands
\begin{eqnarray}  \label{a148}
&&\Big(\sum_{3\leq k} c_{s}(\omega)\,\nu^{-\frac{s}{3}}\Big)
\exp\Big(\omega\,\nu^{\frac{2}{3}}\Big)\hspace{14 ex} (\mi{for\,all\; f\; of\; degree\;2})
\quad
\end{eqnarray}
with  frequencies $\omega$ corresponding to the solutions of
$\mg{\pi}_2(\theta)=0 $ i.e. 
$\mg{\pi}_2(\frac{9}{4}\,\frac{\mathbf{f}_2}{\omega^3})=0 $.
\\

\noindent
For $r=3$ we have $12=3\times 4$ basic summands
\begin{eqnarray}  \label{a149}
&&\Big(\sum_{4\leq k} c_{s}(\omega)\,\nu^{-\frac{s}{4}}\Big)
\exp\Big(\omega\,\nu^{\frac{3}{4}}\Big)\hspace{12 ex} (\mi{for\;all\;radial\; f\; of\;
degree\;3})
\quad\quad \\  \label{a150}
&&\Big(\sum_{4\leq k} c_{s}(\omega)\,\nu^{-\frac{s}{4}}\Big)
\exp\Big(\omega\,\nu^{\frac{3}{4}}+\omega_1\,\nu^{\frac{1}{4}}\Big)
\hspace{8 ex} (\mi{for\;all\; f\; of\; degree\;3})\quad\quad
\end{eqnarray}
with  main frequencies $\omega$ solution of
$\mg{\pi}_3(\theta)=0 $ i.e. 
$\mg{\pi}_3(\frac{64}{27}\,\frac{\mathbf{f}_3}{\omega^4})=0 $,
and  with secondary frequencies $\omega_1$ dependent on the main ones and given by\,:
\begin{eqnarray}  \label{a151}
\omega_1 &=&\frac{2}{3}\,\frac{\mathbf{f}_1}{\omega}\,
\frac{\big( 1+(\BEE{4}+3\,\BEE{2}^2)\,\theta \big)\,
\big(\BEE{2}+(\BEE{2}\BEE{4}-3\,\BEE{3}^2-9\,\BEE{2}^3)\,\theta \big)}
{\big(1+(\BEE{4}-6\,\BEE{2}^2)\,\theta \big)^2+9\,\BEE{2}\,\big(2\,\BEE{3}^2-\BEE{2}^3
\big)\,\theta^2}
\\  \label{a152}
\mi{with}&& \mathbf{f}_1=f_1-\frac{1}{3}\,\frac{f_2^2}{f_3}\;\; ,\;\;
\mathbf{f}_3=f_3\;\;\mi{and} \;\;
\theta =\Big(\frac{4}{3}\Big)^{\!3}\,\frac{\mathbf{f}_3}{\omega^4}
\end{eqnarray}
\\

\noindent
Lastly, for $r=4$ we have $20=4\times 5$ basic summands
\begin{eqnarray}  \label{a153}
&&\Big(\sum_{5\leq k} c_{s}(\omega)\,\nu^{-\frac{s}{5}}\Big)
\exp\Big(\omega\,\nu^{\frac{4}{5}}\Big)\hspace{14 ex} (\mi{for\;all\;radial\; f\; of\;
degree\;4})
\quad\quad \\  \label{a154}
&&\Big(\sum_{4\leq k} c_{s}(\omega)\,\nu^{-\frac{s}{5}}\Big)
\exp\Big(
\omega\,\nu^{\frac{4}{5}}
+\omega_2\,\nu^{\frac{2}{4}}
+\omega_1\,\nu^{\frac{1}{4}} \Big)
\hspace{4 ex} (\mi{for\;all\; f\; of\; degree\;4})\quad\quad\quad
\end{eqnarray}
with  main frequencies $\omega$ solution of
$\mg{\pi}_4(\theta)=0 $ i.e. 
$\mg{\pi}_4(\frac{625}{256}\,\frac{\mathbf{f}_4}{\omega^5})=0 $,
and  with secondary frequencies $\omega_1,\omega_2$ that
depend on the main ones and vanish {\it iff}  the shift-invariants
$\mathbf{f}_1$ resp.  $\mathbf{f}_2$ vanish.
%%%%%%%%%%%%%%%%%%%%%%%%%%%%%%%%%%%%%%%%%%%%%%%%%%%%%%%%%%%%%%%%%%%%%%%%%%%%%%%%%%%%%%%%%%%%
%%%%%%%%%%%%%%%%%%%%%%%%%%%%%%%%%%%%%%%%%%%%%%%%%%%%%%%%%%%%%%%%%%%%%%%%%%%%%%%%%%%%%%%%%%%%
\subsection{The antipodal exchange for polynomial inputs $f$.}
As noted in the preceding subsection, the behaviour of our {\it nir}-transforms  $h(\nu)$ at infinity in the $\nu$-plane involves
 elementary exponential factors multiplied by divergent-resurgent power series 
$\varphi_{\omega}(\nu)=\sum_{r+1\leq k} c_{s}(\omega)\,\nu^{-\frac{s}{r+1}}$, which verify simple linear ODEs
easily deducible from the frequencies $\omega$ and the original ODE verified by $h(\nu)$. Therefore, to resum the 
$\varphi_{\omega}(\nu) $, which are local data at infinity, we must subject them to a formal Borel transform, 
which takes us back to the origin, with a new set of linear ODEs. This kicks off a resurgence ping-pong between 0 and $\infty$.
\footnote{which is quite distinct from the ping-pong between two inner generators  associated with two proper base points
$x_{i},x_{j}$ in the $x$-plane.} Before taking a closer look at it, let us state a useful lemma\,:
\begin{lemma}[Deramification of linear homogeneous ODEs].\\
 Let $\rho$ be a positive integer and
  $ \Phi(t)$ any power series in $\doC\{t^{\frac{1}{\rho}}\}$ or
$\doC\{t^{-\frac{1}{\rho}}\}$ that verifies a linear homogeneous differential equation
$P^\ast(t,\partial_t)\,\Phi(t)=0$ of order $\delta^\ast$ and
 with coefficients polynomial in $t^{\frac{1}{\rho}}$ of degree $d^\ast$. Then $\Phi$
automatically
verifies a new linear homogeneous differential equation
$P(t,\partial_t)\,\Phi(t)=0$ of order $\delta $ and
 with coefficients polynomial in $t$ of degree $d$
such that
$$\delta \leq \delta^\ast\,\rho\;\;,\;\;d\leq (1+d^\ast)\\(1+\delta^\ast\,(\rho-1))^2  $$
\end{lemma}
Proof\,: The initial, ramified differential equation, after division by
the leading coefficient and deramification of the denominators,
can be written uniquely in the form
\begin{equation}  \label{a155}
\Phi^{(\delta^\ast)}=\sum_{0\leq j < \rho}\;\sum_{0\leq s <\delta^\ast}
a_{\delta^\ast\!,j,s}\, t^{\frac{j}{\rho}}\,\Phi^{(s)}
\end{equation}
with unramified coefficients $a_{\delta^\ast\!,j,s}$ that are rational in $t$.
Under successive differentiations and eliminations of the derivatives of order larger than
$\delta^\ast$ but $\not= i$, we then get a sequence of similar-looking equations\,:
\begin{equation}  \label{a156}
\Phi^{(i)}=\sum_{0\leq j < \rho}\;\sum_{0\leq s <\delta^\ast}
a_{i,j,s}\, t^{\frac{j}{\rho}}\,\Phi^{(s)}\hspace{4.ex} 
(\forall i, \delta^\ast\leq i \leq \delta^\ast\rho)
\end{equation}
again with  unramified coefficients $a_{\delta^\ast\!,j,s}$ rational in $t$.
One then checks that there always exists a linear combination of the $(\rho\!-\!1)\,\delta^\ast$
equations (\ref{a156}) with coefficients $L_i(t)$ polynomial in the 
$a_{i^\prime,j^\prime,s^\prime}(t)$ and therefore rational in $t$, that eliminates
the (at most) $(\rho\!-\!1)\,\delta^\ast$ terms of the form $t^{\frac{j}{\rho}}$
with $1\leq j<\rho$ and $0\leq s<\delta^\ast$. After multiplication by a suitable
$t$-polynomial, this yields the required unramified equation $P(t,\partial_t)\Phi(t)=0$.
A closer examination of the process shows that the coefficients $a$ are of the form\,:
$$ a_{i,j,s}(t)=\frac{b_{i,j,s}(t)}{t^{i-\delta^\ast}\,c(t)^{1+i-\delta^\ast}}
\;\mi{with}\; \deg_t(c)\leq d^\ast\;,\; \deg_t(b_{i,j,s}) \leq (1+i-\delta^\ast)\,d^\ast
$$
Plugging this into the elimination algorithm, we get the bound
\begin{eqnarray*}
d &\leq &
(1+(\rho-1)\delta^\ast)\,(d^\ast+(\rho-1)(d^\ast+1)\delta^\ast)\quad \Longleftrightarrow
\\
(1+d) &\leq &
(1+d^\ast)\,(1+\delta^\ast\,\rho^\ast)^2  \;\;\;\quad\;\mi{with}\;\;\;\quad\; \rho^\ast:=\rho-1
\end{eqnarray*}
which, barring unlikely simplifications, is probably near-optimal.\;\;$\Box$.
\\
Let us now return to the resurgence ping-pong $0\leftrightarrow \infty$. Graphically, we get the following sequence
of transforms\,:
\\
\[\begin{array}{ccccccccccccccc}
 \rightarrow\!\!&\!\! \rightarrow\!\!&\!\! \rightarrow\!\!&\!\! \rightarrow\!\!&
\!\! \rightarrow\!\!&\!\! \rightarrow\!\!&\!\! 
\rightarrow\!\!&\!\! \rightarrow\!\!& \!\! \rightarrow
\\
\uparrow\;\;\; &&&&&&&&\;\;\; \downarrow
\\
(P_1,k_1,n_1) & & 
(P_3^\ast,k_3,n_3) &\!\! \stackrel{33}{\rightarrow}\!\! & 
(P_3,k_3,n_3) &  & 
(P_5^\ast,k_5,n_5) &\!\! \stackrel{55}\rightarrow\!\! & 
(P_5,k_5,n_5) &
\\
 &\!\!\stackrel{\;\;12}{ \searrow} \!\! & \uparrow\scriptstyle{23}
 & & 
 &\stackrel{\;\;34}{\searrow} & \uparrow\scriptstyle{45}
 & & 
 &
\\
              & & 
(P_2,k_2,\nu_2) & & 
              & & 
(P_4,k_4,\nu_4) & & 
 &
\end{array}\]

\noindent
{\bf Step 1}:
we have the polynomial-coefficient linear ODE\; 
$P_1(n_1,\partial_{n_1})\,k_1(n_1) =0 $ with 
$$ n_1\equiv n \sim \infty , k_1(n_1)\equiv k(n),P_1(n_1,\partial_{n_1})\equiv
P(n_1,-\partial_{n_1})
$$ 
\medskip

\noindent
{\bf Arrow 12}: we perform the Borel transform from the variable $n_1=n$ to the conjugate
variable $\nu_2=\nu$. Thus\;:\;
$ n_1^{-s}\mapsto \frac{\nu_2^{s-1}}{\Gamma(s)}, 
n_1\mapsto \partial_{\nu_2},  
\partial_{n_1}\mapsto -\nu_2
$.
\\
\medskip

\noindent
{\bf Step 2}: 
we have the  polynomial-coefficient linear  ODE\;
 $ P_2(\nu_2,\partial_{\nu_2})\,k_2(\nu_2) =0 $ with 
$$ \nu_2\equiv \nu\sim 0, k_2(\nu_2)\equiv \hat{k}(\nu),
P_2(\nu_2,\partial_{\nu_2})\equiv
P_1(\partial_{\nu_2},{-\nu_2})
$$ {\bf Arrow 23}:
we go from 0 to $\infty $ and from increasing power series of the variable $\nu_2 $ to 
decreasing power series of the variable $n_3$.
For an input $f$ of degree $r$, we set $n_3:=\nu_2^{\frac{r}{r+1}}=:
\nu_2^{\frac{1}{\kappa_3}}$, the new variable $n_3$ being the ``critical resurgrence variable"
at $\infty$.
\medskip

\noindent
{\bf Step $\mathbf{3}^\ast$}: 
we have the  ramified-coefficient linear  ODE\;
$ P^\ast_3(n_3,\partial_{n_3})\,k_3(n_3) =0 $ with 
$$ n_3^{\kappa_3}\equiv \nu_2\;\;\mi{but}\; n_3\sim \infty\,,\, k_3(n_3)\equiv k_2(\nu_2),\,
P_3^\ast(n_3,\partial_{n_3})\equiv
P_2(n_3^{\kappa_3},\frac{n_3^{\kappa_3-1}}{\kappa_3}\partial_{n_3})
$$ {\bf Arrow 33}: since for an input $f$ of degree $r$, we must take $\kappa_3=\frac{r+1}{r}$,
this leads to a ramification of order $r$ in the coefficients of $P^\ast_3$. We then
apply the above Lemma ... with $\rho=r$ to deramify $P^\ast_3$ to $P_3$.
\medskip

\noindent
{\bf Step 3}:  
we have the  polynomial-coefficient linear  ODE\;
$ P_3(n_3,\partial_{n_3})\,k_3(n_3) =0 $ with 
$n_3$ and $k_3 $ as in step $3^\ast$ but with a linear homogeneous differential
operator $P_3$ which,
unlike
$P_3^\ast$, is polynomial in $n_3$.
\medskip

\noindent
{\bf Arrow 34}:
we perform the Borel transform from the variable $n_3$ to the conjugate
variable $\nu_4$. Thus\;:\;
$ n_3^{-s}\mapsto \frac{\nu_4^{s-1}}{\Gamma(s)}, 
n_3\mapsto \partial_{\nu_4},  
\partial_{n_3}\mapsto -\nu_4
$.
\medskip

\noindent
{\bf Step 4}: 
we have the  polynomial-coefficient linear  ODE\;
 $ P_4(\nu_4,\partial_{\nu_4})\,k_4(\nu_4) =0 $ with 
$ \nu_4$ conjugate to $n_3$ and
$ P_4(\nu_4,\partial_{\nu_4})\equiv P_3(\partial_{\nu_4},-{\nu_4}) $
\medskip

\noindent
{\bf Arrow 45}:
we go from 0 to $\infty $ and from increasing power series of the variable $\nu_4 $ to 
decreasing power series of the variable $n_5$.
For an input $f$ of degree $r$, we set $n_5:=\nu_4^{\frac{r+1}{r}}=:
\nu_4^{\frac{1}{\kappa_5}}=\nu_4^{\kappa_3}$, the new variable $n_5$ being the ``critical
resurgrence variable" at $\infty$.
\medskip

\noindent
{\bf Step $\mathbf{5}^\ast$}:  
we have the  ramified-coefficient linear  ODE\;
$ P^\ast_5(n_5,\partial_{n_5})\,k_5(n_5) =0 $ with 
$$ n_5^{\kappa_5}\equiv \nu_4\;\;\mi{but}\; n_5\sim \infty\,,\, k_5(n_5)\equiv k_4(\nu_4),\,
P_5^\ast(n_5,\partial_{n_5})\equiv
P_4(n_5^{\kappa_5},\frac{n_5^{\kappa_5-1}}{\kappa_5}\partial_{n_5})
$$
\medskip

\noindent
{\bf Arrow 55}: since for an input $f$ of degree $r$, we must take 
$\kappa_5=\frac{r}{r+1}=\frac{1}{\kappa_3}$,
this leads to a ramification of order $r+1$ in the coefficients of $P^\ast_5$. We then
apply once again the above Lemma 6.1 with $\rho=r+1$ to deramify $P^\ast_5$ to $P_5$.
\medskip

\noindent
{\bf Step 5}:  
we have the  polynomial-coefficient linear  ODE\;
$ P_5(n_5,\partial_{n_5})\,k_5(n_5) =0 $ with 
$n_5$ and $k_5 $ as in step $5^\ast$ but with a linear homogeneous differential
operator $P_5$ which,
unlike
$P_5^\ast$, is polynomial in $n_5$.
%%%%%%%%%%%%%%%%%%%%%%%%%%%%%%%%%%%%%%%%%%%%%%%%%%%%%%%%%%%%%%%%%%%%%%%%%%%%%%%%%%%%%%%%%%%%
%%%%%%%%%%%%%%%%%%%%%%%%%%%%%%%%%%%%%%%%%%%%%%%%%%%%%%%%%%%%%%%%%%%%%%%%%%%%%%%%%%%%%%%%%%%%
\subsection{ODEs for monomial inputs $F$.}
General meromorphic inputs $F$, with more than one zero or pole, shall be investigated in \S7.2 and \S8.3-4 with
the usual {\it nir-mir} approach. Here, we shall restrict ourselves to strictly {\it monomial} $F$, i.e. with only one
zero {\it or} pole (but of abitrary order $p$), for these monomial inputs, and only they, give rise to {\it nir} transforms that
verify linear ODEs with polynomial coefficients. So for now our inputs shall be\,:
\begin{eqnarray}  \label{a157}
f(x):=+p\,\log(1+p\,x) &,& F(x):=(1+p\,x)^{-p}\quad\quad(p\in \doN^\ast)
\\  \label{a158}
f(x):=-p\,\log(1-p\,x) &,& F(x):=(1-p\,x)^{+p}\quad\quad\;(p\in \doN^\ast)
\end{eqnarray}
and we shall set as usual\,:
\begin{eqnarray*}
k(n) &:=& \mr{singular}\Big( \int_0^\infty
e_{\#}^{-\beta(\partial_\tau)\,f(\frac{\tau}{n})\ }\;d\tau
\Big)\hspace{7.ex} \in \Gamma(1/2)\, n^{1/2}\,\doQ[[n^{-1}]]
\\
\smi{h}(\nu) &:=& \mr{formal}
\Big( \frac{1}{2\pi i}\int_{c-i\infty}^{c+i\infty}\,k(n)\,e^{\nu\,n}\frac{dn}{n} \Big)=h(\nu)
\hspace{6.ex} \in \nu^{-1/2}\,\doQ\{\nu\}
\\
\imi{k}(\nu) &:=& \mr{formal}
\Big( \frac{1}{2\pi i}\int_{c-i\infty}^{c+i\infty}\,k(n)\,e^{\nu\,n}{dn}\Big)=h^\prime(\nu)
\quad\quad\;\; \in \nu^{-3/2}\,\doQ\{\nu\}
\end{eqnarray*}
Unlike  with the polynomial inputs $f$ of \S6.2-5, the global {\it nir} transforms now verify no
(variable) polynomial linear-homogeneous ODEs. Only their singular parts,
which in the present case ($\forall p$) always consist of semi-entire powers of the variable,
do verify (covariant) linear ODEs with polynomial coefficients. These equations
depend only on the absolute value $|p|$ and read\,:
\begin{eqnarray*}
\big(n+n\partial_n-\frac{|p|}{2}\big)^{|p|}\, k(n)  
&=& n^{|p|}\, k(n)
\\
\big(\partial_\nu-\nu\,\partial_\nu-\frac{|p|}{2}\big)^{|p|}\, h(\nu)  
&=& (\partial_\nu)^{|p|}\,h(\nu)
\\
\big(\partial_\nu-\nu\,\partial_\nu-1-\frac{|p|}{2}\big)^{|p|}\, \hat{k}(\nu)  
&=& (\partial_\nu)^{|p|}\,\hat{k}(\nu)
\end{eqnarray*}
If we regard $n$ and $\nu$ no longer as commutative variables (as in \S4 and \S5),
but as non-commutative ones bound by $[n,\nu]=1$ (as in the preceding sections),
our covariant ODEs read\,:
\begin{eqnarray*}
P(n,-\partial_n)\,k(n)=0\;\;,\;\;
\partial_\nu^{-1}\,P(\partial_\nu,\nu)\,\partial_\nu\,h(\nu)=0\;\;,\;\;
\,P(\partial_\nu,\nu)\,\hat{k}(\nu)=0\;\;\;\;
\\
\mi{with} \;\;\; \quad
P(n,\nu):= \big(n-n\,\nu-\frac{|p|}{2}\big)^{|p|} - n^{|p|}
= (n-\nu\,n-1-\frac{|p|}{2})^{|p|} - n^{|p|}\hspace{4 ex}
\end{eqnarray*} 
If we now apply the covariance relation (\ref{a129}) to the shifts $(\epsilon,\eta)$\,:
$$\epsilon:=-1/|p|\;,\; \eta:=\int_0^\epsilon f(x)\,dx=1\,,\,
 ^{\epsilon\!}f(x)=|p|\,\log(|p|x)  $$
we find a centered polynomial $P_\ast$ predictably simpler than $P$\,:
$$ P_\ast(n,\nu)=P(n,\nu+\eta)=\big(-n\,\nu-\frac{|p|}{2}\big)^{|p|} -n^{|p|} $$
Although our covariant operators $P(n,\nu)$ are now much simpler,
and of far lower degree in $n$, than was the case for polynomial inputs $f$,
their form is actually harder to derive. As for their dependence on $|p|$
rather than $p$, it follows from the general parity relation for the $nir$
transform (cf \S4.10), but here it also makes direct formal sense. Indeed, in view of $[n,\nu]=1$, we have the
chain of formal equivalences:
\begin{eqnarray*}
\{\big(n-n\,\nu-\frac{p}{2}\big)^{p}\,k(n)=n^{p}\,k(n)\} &\Longleftrightarrow &\hspace{10. ex}
\\
\{k(n)=\big(n-n\,\nu-\frac{p}{2}\big)^{-p}\,n^{p}\,k(n)\} &\Longleftrightarrow &
\\
\{k(n)=n^{p}\,\big(n-n\,\nu+\frac{p}{2}\big)^{-p}\,k(n)\} &\Longleftrightarrow &
\\
\{ n^{-p}\,k(n)= \big(n-n\,\nu +\frac{p}{2}  \big)^{-p}\,k(n) \} &&
\end{eqnarray*}
which reflects the invariance of $P(n,\nu)$ \,\footnote{or more accurately\,:
the invariance of the relation $P(n,\nu)\,k(n)=0 $.} under the change
$p\mapsto -p$.

\medskip
From the form of the centered differential operator, it is clear that $h(1-\nu)$ has all its irregular singular points
over the unit roots, plus a regular singular point at infinity.

\medskip
Remark\,:  Although both inputs $f_{1}(x)=\frac{1}{p}\,x^p -1$ and $f_{2}(x)=\pm p\log(1\pm p\,x)$ lead to {\it nir}-transforms
$h_{1}(1-\nu)$ and $h_{2}(1-\nu)$ with radial symmetry and singular points over the unit roots of order $p$, there are far-going differences\,:
\\
(i) $h_{2}$ verifies much simpler ODEs than $h_{1}$
\\
(ii) conversely,  $h_{1}$ verifies much simpler resurgence equations than $h_{2}$ (see {\it infra})
\\
(iii)
the singularities of $h_{1}$ over $\infty$ are of divergent-resurgent type (see \S6.4-5) whereas those of $h_{2}$
are merely ramified-convergent (see \S6.7).

\medskip
Let us now revert to our input (\ref{a157}) or (\ref{a158}) with the corresponding {\it nir} transform $h(\nu)$ and its linear
ODE. That ODE always has very explicit power series solutions at $\nu=0$ and $\nu=\infty$ and, as we shall see, this is what really matters.
At $\nu=0$ the solutions are of the form\,:
\begin{eqnarray*}
k(n)=\sum_{s\in-\frac{1}{2}+\doN} k_{s}\,n^{-s} \quad &,&
\quad
h(\nu)=\sum_{s\in-\frac{1}{2}+\doN} k_{s}\,\nu^{s}
\quad\quad {\it (relevant)}
\\
k^{\rm{en}}\!(n)=\sum_{s \in \doN} k_{s}\,n^{-s} \quad &,&
\quad 
h^{\rm{en}}\!(\nu)=\sum_{s \in \doN} k_{s}\,\nu^{s}
\hspace{6. ex} {\it (irrelevant)}
\end{eqnarray*}
but only for $ p \in \{ \pm1,\pm 2,\pm 3 \}$ are the coefficients explicitable.

\medskip
\noindent
{\bf The case $p=\pm 1$.}
\begin{eqnarray*}
k_{-\frac{1}{2}+r}&=& 0 \quad\mi{if}\quad r\geq 1\quad\mi{and}\quad k_{-\frac{1}{2}}=\Big(\frac{\pi}{2}\Big)^{\frac{1}{2}}
\\
h_{-\frac{1}{2}+r}&=& 0 \quad\mi{if}\quad r\geq 1\quad\mi{and}\quad k_{-\frac{1}{2}}=\Big(\frac{1}{2}\Big)^{\frac{1}{2}}
\end{eqnarray*}
\\
{\bf The case $p=\pm 2$.}
\begin{eqnarray*}
k_{-\frac{1}{2}+r}&=& 2^{-5\,r}\frac{(2\,r)!(2\,r)!}{r!\,r!\,r!}k_{-\frac{1}{2}}\quad\mi{with}\quad k_{-\frac{1}{2}}=\Big(\frac{\pi}{8}\Big)^{\frac{1}{2}}
\\
h_{-\frac{1}{2}+r}&=&  2^{-3\,r}\frac{(2\,r)!}{r!\,r!}h_{-\frac{1}{2}}\quad\quad\quad\mi{with}\quad h_{-\frac{1}{2}}=\Big(\frac{1}{8}\Big)^{\frac{1}{2}}
\end{eqnarray*}
\\
{\bf The case $p=\pm 3$.}
 The coefficients of $k,h$ have no simple multiplicative structure, but the entire analogues $k^{\mr{en}},h^{\mr{en}}$
are simple superpositions of hypergeometric series.
%%%%%%%%%%%%%%%%%%%%%%%%%%%%%%%%%%%%%%%%%%%%%%%%%%%%%%%%%%%%%%%%%%%%%%%%%%%
%%%%%%%%%%%%%%%%%%%%%%%%%%%%%%%%%%%%%%%%%%%%%%%%%%%%%%%%%%%%%%%%%%%%%%%%%%%
\subsection{Monomial inputs $F$\,: global resurgence.}
Let us replace the pair $(h,P)$ by $(\mi{h\!o},\mi{P\!o})$ with
\begin{equation}
\mi{h\!o}(\nu):=h(1-\nu)\quad ;\quad \mi{P\!o}(n,\nu):=(-1)^p\,P(-n,1-\nu)=(\nu\,n+\frac{p}{2})^p-n^p\quad \quad\quad
\end{equation}
so as to respect the radial symmetry and deal with a function {\it ho} having all its singular points over the unit roots 
$e_{j}=\exp(2\pi i j/p)$. At the crucial points $\nu_{0}\in \{0,\infty,e_{0},\dots,e_{p-1}\}$ the $p$-dimensional kernel of the operator
\begin{equation}
\mi{P\!o}(\partial,\nu+\nu_{0})=\big((\nu+\nu_{0})\partial_{\nu}+\frac{p}{2}\big)^p\,-\,\big(\partial_{\nu})^p
\end{equation}
is spanned by the following systems of fundamental solutions
\begin{eqnarray*}
{\mi at}\; \nu_{0}=0:&
\mi{h\!o}_{s}(\nu) &\in \nu^s\; \doC\{\nu^p\} \hspace{35. ex}   (0\leq s \leq p-1)
\\
{\mi at}\; \nu_{0}=\infty:&
\mi{hi}_{s}(\nu) &\in\bigoplus_{0\leq \sigma \leq s-1}
\Big(  \nu^{-p/2}\,  \doC\{ \nu^{-1} \}\, \frac{(\log\nu)^\sigma}{\sigma!}\Big) \hspace{12. ex}   (0\leq s \leq p-1)
\\
{\mi at}\; \nu_{0}=e_{j}:&
\mi{h\!a}_{j}(\nu) &\in \nu^{-1/2}\; \doC\{\nu\} \hspace{30. ex}   
\\
&\mi{h\!a}_{j,s}(\nu) &\in  \doC\{\nu\} \hspace{38. ex}   (1\leq s \leq p-1)
\end{eqnarray*}
The singular solutions $\mi{h\!a}_{j}$ (normalised in a manner consistent with the radial symmetry) are, up to sign, none other
than the
{\it inner generators} whose resurgence properties we want to describe. For $p\geq 3$, their coefficients have
no transparent expression, but the coefficients of the $\mi{hi}_{s}$ and, even more so, those of the $\mi{h\!o}_{s}$
do possess a very simple multiplicative structure, which allows us to apply the method of coefficient asymptotics
in \S2.3 to derive the resurgence properties of the $\mi{h\!a}_{j}$, and that too from `both sides' --- from 0 and $\infty$.
A complete treatment shall be given in [ES] but here we shall only state the result and describe the closed resurgence system
governibg the behaviour of the $\mi{h\!a}_{j}$. To that end, we consider their Laplace integrals along any given
axis $\arg\nu=\theta$, with the ``location factor'' $e^{ -e_{j}\,n }$\,:
\begin{equation}
\mi{h\!a\!a}_{j}^\theta(n):=e^{ -e_{j}\,n } \int_{0}^{ e^{i\theta}\infty } \mi{h\!a}_{j}(\nu)\,d\nu 
\quad \quad\quad \quad (\theta\in\doR\,,\, j\in\doZ/p\doZ)
\end{equation}
Everything boils down to describing the effect on the system 
$\{  \mi{h\!a\!a}_{1}^{\theta},\dots,\mi{h\!a\!a}_{p}^{\theta}  \}$ of crossing a singular axis 
$\theta_{0}=\arg(e_{j_{2}}-e_{j_{1}})$, i.e. of going from $\theta_{0}-\epsilon$ to  $\theta_{0}+\epsilon$.
The underlying ODE being linear, such a crossing will simply subject 
$\{  \mi{h\!a\!a}_{1}^{\theta},\dots,\mi{h\!a\!a}_{p}^{\theta}  \}$
to a linear transformation with constant coefficients. Moreover, since all $\mi{h\!a}_{j}(\nu)$ are in $\nu^{-1/2}\,\doC\{\nu\}$,
two full turns (i.e. changing $\theta$ to $\theta+4\pi$) ought to leave $\{  \mi{h\!a\!a}_{1}^{\theta},\dots,\mi{h\!a\!a}_{p}^{\theta}  \}$ unchanged. All the above facts can be derived in a rather straightforward manner by resurgence analysis (see [ES]) but, when
translated into matrix algebra, they lead to rather complex matrices and to remarkable, highly non-trivial relations between
these. Of course, the relations in question  also admit {\it `direct' algebraic} proofs, but these are rather difficult -- and in any case much longer
than their {\it `indirect' analytic} derivation. The long subsection which follows is entirely devoted to this `algebraic'
description of the resurgence properties of the $\mi{h\!a}_{j}$.

%%%%%%%%%%%%%%%%%%%%%%%%%%%%%%%%%%%%%%%%%%%%%%%%%%%%%%%%%%%%%%%%%%%%%%%%%%%
%%%%%%%%%%%%%%%%%%%%%%%%%%%%%%%%%%%%%%%%%%%%%%%%%%%%%%%%%%%%%%%%%%%%%%%%%%%
\subsection{Monomial inputs $F$\,: algebraic aspects.}
{\bf Some elementary matrices.}
\\
Eventually, $\epsilon$ will stand for -1 and $\epsilon^{q}$ for $e^{\pi i q}\,,\,\forall q\in\doQ$, but for greater clarity $\epsilon$
shall be  kept free (unassigned) for a while. We shall encounter both $\epsilon$-carrying matrices, which we shall
underline, and $\epsilon$-free matrices. For each $p$, we shall also require the following elementary square matrices
$(p\times p)$\,:
\\
$\caI$ \;\;:\;\; {identity}
\\
$\underline{\caI} $ \;\;:\;\; $\epsilon$-{carrying diagonal}
\\
$\caJ$ \;\;:\;\; {Jordan correction}
\\
$\caP$ \;\;:\;\; {unit shift}
\\
$\caQ$  \;\;:\;\; {twisted unit shift}
\\
These are hollow matrices, whose only nonzero entries are\,:
\[\begin{array}{llrllllrlllll}
\underline{\caI}[i,j]& =& \epsilon^{j/p}&\mi{if}& j=i 
\\
{\caJ}[i,j]&\!=\!& 1&\mi{if}& j=i+1
 \\
{\caP}[i,j]&\!=\!& 1&\mi{if}& j=i+1 \mod p &{\caP}^{k}[i,j]&\!=\!& 1&\mi{if}& j=i+k \mod p
 \\
{\caQ}[i,j]&\!=\!& 1&\mi{if}& j=i+1  &{\caQ}^{k}[i,j]&\!=\!& 1&\mi{if}& j=i+k 
 \\
{\caQ}[i,j]&\!=\!& -1&\mi{if}& j=i+1-p  &{\caQ}^{k}[i,j]&\!=\!& -1&\mi{if}& j=i+k - p
\end{array}\]
{\bf The simple-crossing matrices $\underline{\caM}_{k},{\caM}_{k}$.}
\\
Let $\mi{fr}(x)$ resp.  $\mi{en}(x)$  denote the {\it fractional} resp. {\it entire} part of any real  $x $\,:
$$x\equiv\mr{fr}(x)+\mr{en}(x)\quad \mi{with}\quad x\in \doR\;\;,\;\; \mr{fr}(x)\in [0,1[\;\;,\;\; \mr{en}(x)\in \doZ  $$
Fix $p\in\doN^{\ast}$ and set $e_{j}:=\exp(2\pi i j/p), \forall j\in \doZ$. For any $k\in \frac{1}{2}\doZ$, it is convenient
to denote $\theta_{k}$ the axis of direction $2\pi (\frac{k}{p}+\frac{3}{4})$, i.e. the axis from $e_{j_{1}}$
to  $e_{j_{2}}$ for any pair $j_{1},j_{2}\in \doZ$ such that $j_{1}+j_{2}=2 k \mod p$ and $(k<j_{1}<j_{2})^{\mr{circ}}_{p}$.
The matrix $\underline{\caM}_{k} $ corresponding to the (counterclockwise) crossing of the axis $\theta_{k}$ has the
following elementary entries\,:
\[\begin{array}{lllll}
\underline{\caM}_{k}[i,j]&=&1  &\mi{if}& i=j
\\
\underline{\caM}_{k}[i,j]&=&-\epsilon^{\mr{fr}(\frac{j-k}{p})-\mr{fr}(\frac{i-k}{p})}\frac{p!}{(|i-j|)!(p-|i-j|)!}  &\mi{if}& \mr{fr}(\frac{i+j-2k}{p})=0 \;\;
\\&&&\mi{and}& \mr{fr}(\frac{i-k}{p})>\mr{fr}(\frac{j-k}{p})
\\
\underline{\caM}_{k}[i,j]&=&0  & \mi{otherwise}&
\end{array}\]
Alternatively, we may start from the simpler matrix $\underline{\caM}_{0}$\,:
\[\begin{array}{lllll}
\underline{\caM}_{0}[i,j]&=&1  &\mi{if}& i=j
\\
\underline{\caM}_{0}[i,j]&=&-\epsilon^{\frac{j-i}{p}}\frac{p!}{(i-j)!(p-i+j)!}  & \mi{if} & i>j \;\mi{and}\; i+j=2k\mod p
\\
\underline{\caM}_{0}[i,j]&=&0  & \mi{otherwise}&
\end{array}\]
and deduce the general  $\underline{\caM}_{k}$ under the rules\,:
$$
\underline{\caM}_{k}[i,j]= \underline{\caM}_{0}[[i-k]_{p},[j-k]_{p}]
\quad \quad
\mi{with}\;\;\;\; [x]_{p}:=p\,.\,\mr{en}(\frac{x}{p})
$$
$\underline{\caM}_{k} $ carries unit roots of order $2 p$ (hence the underlining) but 
can be turned into a unit root-free matrix  ${\caM}_{k} $ under a $k$-independent conjugation\,:
\begin{equation}  \label{b1}
\caM_{k}= \underline{\caI}\;\; \underline{\caM}\;\;\underline{\caI}^{-1}
\end{equation}
with the elementary diagonal matrix $\underline{\caI}$ defined above. We may therefore work with the simpler matrices $\caM_{k}$ whose
entries are\,:
\[\begin{array}{lllll}
\caM_{k}[i,j]&=&1  &\mi{if}& i=j
\\
\caM_{k}[i,j]&=&-\epsilon^{\mr{en}(\frac{j-k}{p})-\mr{en}(\frac{i-k}{p})}\frac{p!}{(|i-j|)!(p-|i-j|)!}  &\mi{if}& \mr{fr}(\frac{i+j-2k}{p})=0 \;\;
\\&&&\mi{and}&\mr{fr}(\frac{i-k}{p})>\mr{fr}(\frac{j-k}{p})
\\
\caM_{k}[i,j]&=&0  & \mi{otherwise}&
\end{array}\]
However, since $\underline{\caI}$ and $\caP$ do not commute, we go 
from $\underline{\caM}_{k}$ to  $\underline{\caM}_{k+1}$ under  the {\it regular} shift $\caP$ but
from ${\caM}_{k}$ to  ${\caM}_{k+1}$ under  the {\it twisted} shift $\caQ$\,:
\begin{equation}  \label{b2}
\underline{\caM}_{k+1}=\caP^{-1}\;\;\underline{\caM}_{k}\;\;\caP
\quad,\quad
{\caM}_{k+1}=\caQ^{-1}\;\;{\caM}_{k}\;\;\caQ
\end{equation}
\\
{\bf The multiple-crossing matrices $\underline{\caM}_{k_{2},k_{1}},{\caM}_{k_{2},k_{1}}$.}
\\ \\
For any $k_{1},k_{2} \in \frac{1}{2}\doZ$ such that $k_{2} > k_{1} $ we set \,:
\begin{eqnarray}  \label{b3}
\underline{\caM}_{k_{2},k_{1}}:=
\underline{\caM}_{k_{2}}\;\;\underline{\caM}_{k_{2}-\frac{1}{2}}\;\;\underline{\caM}_{k_{2}-\frac{2}{2}}\;\dots\;
\underline{\caM}_{k_{1}+\frac{3}{2}}\;\;\underline{\caM}_{k_{1}+\frac{2}{2}}\;\;\underline{\caM}_{k_{1}+\frac{1}{2}}
\\  \label{b4}
\caM_{k_{2},k_{1}}:=
\caM_{k_{2}}\;\;\caM_{k_{2}-\frac{1}{2}}\;\;\caM_{k_{2}-\frac{2}{2}}\;\dots\;
\caM_{k_{1}+\frac{3}{2}}\;\;\caM_{k_{1}+\frac{2}{2}}\;\;\caM_{k_{1}+\frac{1}{2}}
\end{eqnarray}
 For $k_{2} < k_{1}$ or $k_{2} = k_{1}$ we set of course\,:
 $$
 \underline{\caM}_{k_{2},k_{1}}:= \underline{\caM}_{k_{1},k_{2}}^{-1}
  \quad,\quad
 {\caM}_{k_{2},k_{1}}:= {\caM}_{k_{1},k_{2}}^{-1}
 \quad,\quad
\underline{\caM}_{k,k}:= {\caM}_{k,k}:=\caI$$
thus ensuring the composition rule\,:
  $$\underline{\caM}_{k_{3},k_{2}}\;\;\underline{\caM}_{k_{2},k_{1}}=\underline{\caM}_{k_{3},k_{1}}
  \quad,\quad
  {\caM}_{k_{3},k_{2}}\;\;{\caM}_{k_{2},k_{1}}={\caM}_{k_{3},k_{1}}\quad \big(\forall k_{i}\in\frac{1}{2}\doZ\big)
  $$
  Since $\underline{\caM}_{p+k}\equiv \underline{\caM}_{k}$ and $\caM_{p+k}\equiv \caM_{k}$
  for all $k$ ($p$-periodicity), each full-turn matrix    $\underline{\caM}_{p+k,k}$ or ${\caM}_{p+k,k}$
  is conjugate to any other. It turns out, however, that just two of them (corresponding to $ k \in \{0,1\} $ if $p=0$ or $1$\!$\mod 4$,
  and  to $k \in \{\pm\frac{1}{2} \} $ if $p=2$ or $3$\!$\mod 4$) admit a simple or at least tolerably
   explicit normalisation (i.e. a conjugation to the canonical Jordan form, or in this case, a more convenient
   variant thereof). That normalisation involves remarkable lower diagonal matrices $\caL$ and $\caR$. 
  To construct  $\caL$ and $\caR$, however, we require a set of rather intricate polynomials $H^{\delta}_{d}$.
  \\
\\ \noindent
{\bf The auxiliary polynomials $H^{\delta}_{d}(x,y)$.}
\\
These polynomials, of global degree $d$ in each of their two variables  $x,y$,  also depend on an 
integer-valued parameter $\delta \in \doZ$. They are $d$-inductively determined by the following system of
difference equations in $y$, along with the initial conditions for $y=0$\,:
\begin{eqnarray} \label{asa1}
H^{\delta}_{d}(x,y)&=&H^{\delta}_{d}(x,y-1)+(x-d)\,H^{\delta}_{d-1}(x,y-1)
\\  \label{asa2}
H^{\delta}_{d}(x,0)&=& \frac{(x+\delta+d)!}{(x+\delta)!}\,=\,\prod_{0<d_{1}\leq d} (x+\delta+d_{1})
\end{eqnarray}
One readily sees that this induction leads to the direct expression\,:
\begin{eqnarray}  \label{asa3}
H^{\delta}_{d}(x,y)&=&\sum_{d_{1}=0}^{d}\frac{(x-1-d+d_{1})!}{(x-1-d)!}\frac{(x+\delta+d-d_{1})!}{(x+\delta)!}\frac{y!}{d_{1}!(y-d_{1})!}
\\ \nonumber
&=& \sum_{d_{1}=0}^{d}\frac{1}{d_{1}!}
\prod_{0 \leq k_{1} < d_{1}} \!(x-d+k_{1})\!
\prod_{1 \leq k_{2} \leq d-d_{1}} \!(x+\delta+k_{2})\!
\prod_{0 \leq k_{3} < d_{1}}\! (y-k_{3})
\end{eqnarray}
which is turn can be shown to be equivalent to\,:
\begin{eqnarray} \label{asa4}
H^{\delta}_{d}(x,y)&\!\!=\!\!&\sum_{d_{1}=0}^{d}\;
 \Big[\!\Big[\frac{\delta\!+\!2d_{1}}{\delta\!+\!d_{1}}\Big]\!\Big]!!\;\;\;
 \frac{(2d\!+\!d_{1}\!-\!y)!}{(d\!+\!2d_{1}\!-\!y)!}\;\;\;
\prod_{0 \leq d_{2}\leq d}^{d_{2}\not=d_{1}}\frac{(x-d_{2})}{(d_{1}-d_{2})}\;
\\  \label{asa5}
&\!\!=\!\!&\sum_{d_{1}=0}^{d}\;  \Big[\!\Big[\frac{\delta\!+\!2d_{1}}{\delta\!+\!d_{1}}\Big]\!\Big]!!
\prod_{d_{1} < d_{3} \leq d}\!\!\! (d\!+\!d_{1}\!+\!d_{3}\!-\!y) 
\prod_{0 \leq d_{2}\leq d}^{d_{2}\not=d_{1}}\frac{(x-d_{2})}{(d_{1}-d_{2})}\;
\end{eqnarray}
with\,:
\begin{eqnarray} \nonumber
 \Big[\!\Big[\frac{a}{b}\Big]\!\Big]!!&:=& \frac{a!}{b!}\hspace{21.ex} \mi{if}\quad a,b\in \doN
 \\  \nonumber
 &:=& (-1)^{a-b}\frac{(-1-b)!}{(-1-a)!}\quad \quad \mi{if}\quad a,b\in -\doN^{\ast}
 \\
 &:=& 0 \quad\hspace{21.ex} \mi{otherwise}  \label{abc3}
\end{eqnarray}
\\
\noindent
{\bf The left normalising matrix $\caL$.}
\[\begin{array}{lllllllll}
\mi{if} \;\; i<j    &     (\forall p)                            &:   \caL[i,j]=& 0
\\[2.ex]
\mi{if} \;\; p=0  &\!\!\!\mod 4\;\;\;\; \mi{and}      & \dots
\\[1.ex]
 2\,j \leq  p     ,& i\!+\!j \leq p\!+\!2                   &:  \caL[i,j]=& (-1)^{i}\frac{(i-1)!}{(j-1)!(i-j)!}
\\[1.ex]
 2\,j \leq  p     ,& i\!+\!j > p\!+\!2                     &:    \caL[i,j]=& (-1)^{j}\frac{(i-1)!}{(2j-3)!(p-2j+2)!}\,H^{1}_{p-i}(j\!-\!2,p)
\\[1.ex]
 2\,j >    p        &                                            &:    \caL[i,j]=& (-1)^{j}\frac{(p-j)!}{(p-i)!(i-j)!}
\\[2.ex]
\mi{if} \;\; p=1  &\!\!\!\mod 4\;\;\;\; \mi{and}      & \dots
\\[1.ex]
 2\,j \leq  p\!+\!1     ,& i\!+\!j \leq p\!+\!2                   &:  \caL[i,j]=&         (-1)^{i}\frac{(i-1)!}{(j-1)!(i-j)!}
\\[1.ex]
 2\,j \leq  p\!+\!1     ,& i\!+\! j > p\!+\!2                     &:    \caL[i,j]=&\!\!\!\!-(-1)^{j}\frac{(i-1)!}{(2j-3)!(p-2j+2)!}\,H^{1}_{p-i}(j\!-\!2,p)
\\[1.ex]
 2\,j >    p\!+\!1       &                                            &:    \caL[i,j]=&        (-1)^{j}\frac{(p-j)!}{(p-i)!(i-j)!}
\\[2.ex]
\mi{if} \;\; p=2  &\!\!\!\mod 4\;\;\;\; \mi{and}      & \dots
\\[1.ex]
 2\,j \leq  p     ,& i\!+\!j \leq p\!+\!1                   &:  \caL[i,j]=&        (-1)^{i}\frac{(i-1)!}{(j-1)!(i-j)!}
\\[1.ex]
 2\,j \leq  p     ,& i\!+\!j > p\!+\!1                     &:    \caL[i,j]=&\!\!\!\!-(-1)^{j}\frac{(i-1)!}{(2j-2)!(p-2j+1)!}\,H^{0}_{p-i}(j\!-\!1,p)
\\[1.ex]
 2\,j >    p        &                                            &:    \caL[i,j]=&         (-1)^{j}\frac{(p-j)!}{(p-i)!(i-j)!}
\\[2.ex]
\mi{if} \;\; p=3  &\!\!\!\mod 4\;\;\;\; \mi{and}      & \dots
\\[1.ex]
 2\,j \leq  p\!-\!1     ,& i\!+\!j \leq p\!+\!1                  &:  \caL[i,j]=& (-1)^{i-1}\frac{(i-1)!}{(j-1)!(i-j)!}
\\[1.ex]
 2\,j \leq  p\!-\!1     ,& i\!+\!j > p\!+\!1                     &:    \caL[i,j]=& (-1)^{i-1}\frac{(i-1)!}{(2j-2)!(p-2j+1)!}\,H^{0}_{p-i}(j\!-\!1,p)
\\[1.ex]
 2\,j >    p\!-\!1        &                                            &:    \caL[i,j]=& (-1)^{p-j}\frac{(p-j)!}{(p-i)!(i-j)!}
\end{array}\]
\noindent
{\bf The right normalising matrix $\caL$.}
\[\begin{array}{lllllllll}
\mi{if} \;\; i<j    &     (\forall p)                            &:   \caR[i,j]=& 0
\\[2.ex]
\mi{if} \;\; p=0  &\!\!\!\mod 4\;\;\;\; \mi{and}      & \dots
\\[1.ex]
 2\,j \leq  p\!-\!2     ,& i\!+\!j \leq p                  &:  \caR[i,j]=& (-1)^{i}\frac{(i-1)!}{(j-1)!(i-j)!}
\\[1.ex]
 2\,j \leq  p\!-\!2     ,& i\!+\!j > p                     &:    \caR[i,j]=& (-1)^{j}\frac{(i-1)!}{(2j-1)!(p-2j)!}\,H^{-1}_{p-i}(j,p)
\\[1.ex]
 2\,j >    p\!-\!2        &                                            &:    \caR[i,j]=& (-1)^{j}\frac{(p-j)!}{(p-i)!(i-j)!}
\\[2.ex]
\mi{if} \;\; p=1  &\!\!\!\mod 4\;\;\;\; \mi{and}      &\dots 
\\[1.ex]
 2\,j \leq  p\!-\!1     ,& i\!+\!j \leq p                  &:  \caR[i,j]=&         (-1)^{i}\frac{(i-1)!}{(j-1)!(i-j)!}
\\[1.ex]
 2\,j \leq  p\!-\!1     ,& i\!+\!j > p                     &:    \caR[i,j]=&\!\!\!\!-(-1)^{j}\frac{(i-1)!}{(2j-1)!(p-2j)!}\,H^{-1}_{p-i}(j,p)
\\[1.ex]
 2\,j >    p\!-\!1       &                                            &:    \caR[i,j]=&        (-1)^{j}\frac{(p-j)!}{(p-i)!(i-j)!}
\\[2.ex]
\mi{if} \;\; p=2  &\!\!\!\mod 4\;\;\;\; \mi{and}      & \dots
\\[1.ex]
 2\,j \leq  p\!-\!2     ,& i\!+\!j \leq p\!-\!1                   &:  \caR[i,j]=&        (-1)^{i}\frac{(i-1)!}{(j-1)!(i-j)!}
\\[1.ex]
 2\,j \leq  p\!-\!2     ,& i\!+\!j > p\!-\!1                     &:    \caR[i,j]=&\!\!\!\!-(-1)^{j}\frac{(i-1)!}{(2j)!(p-2j-1)!}\,H^{-2}_{p-i}(j\!+\!1,p)
\\[1.ex]
 2\,j >    p\!-\!2        &                                            &:    \caR[i,j]=&         (-1)^{j}\frac{(p-j)!}{(p-i)!(i-j)!}
\\[2.ex]
\mi{if} \;\; p=3  &\!\!\!\mod 4\;\;\;\; \mi{and}      & \dots
\\[1.ex]
 2\,j \leq  p\!-\!3     ,& i\!+\!j \leq p\!-\!1                  &:  \caR[i,j]=& (-1)^{i-1}\frac{(i-1)!}{(j-1)!(i-j)!}
\\[1.ex]
 2\,j \leq  p\!-\!3     ,& i\!+\!j > p\!-\!1                     &:    \caR[i,j]=& (-1)^{i-1}\frac{(i-1)!}{(2j)!(p-2j-1)!}\,H^{-2}_{p-i}(j\!+\!1,p)
\\[1.ex]
 2\,j >    p\!-\!3        &                                            &:    \caR[i,j]=& (-1)^{p-j}\frac{(p-j)!}{(p-i)!(i-j)!}
\end{array}\]
{\bf Normalisation identities for the full-turn matrices $\caM_{p+k,k}$\,: }
\begin{eqnarray*}
\caL\;\;\caM_{p+1,\,1}\;\;\caL^{-1}&=&(-1)^{p-1}(\caI+\caJ)^{p} \quad\quad \mi{if}  \;\; p=0 \; \mi{or} \;1 \;\; \mod 4
\\
\caR\;\;\caM_{p+0,\,0}\;\;\caR^{-1}&=&(-1)^{p-1}(\caI+\caJ)^{p} \quad\quad \mi{if}  \;\; p=0 \; \mi{or} \;1 \;\; \mod 4
\\[1.5 ex]
\caL\;\;\caM_{p+\frac{1}{2},+\frac{1}{2}}\;\;\caL^{-1}&=&(-1)^{p-1}(\caI+\caJ)^{p} \quad\quad \mi{if}  \;\; p=2 \; \mi{or} \;3 \;\; \mod 4
\\
\caR\;\;\caM_{p-\frac{1}{2},-\frac{1}{2}}\;\;\caR^{-1}&=&(-1)^{p-1}(\caI+\caJ)^{p} \quad\quad \mi{if}  \;\; p=2 \; \mi{or} \;3 \;\; \mod 4
\end{eqnarray*}
with $\caI$ denoting the identity matrix and $\caJ$ the matrix carrying a maximal upper-Jordan side-diagonal:
$$ \caJ[i,j]=1\;\; \;\mi{if}\;\;\; j=1\!+\!i \;\;\;\; \mi{and} \;\;\;  \caJ[i,j]=0 \;\;\; \mi{otherwise} \ $$
This result obviously implies that {\it all} full-rotation matrices $\caM_{p+k,\,k} $ are also conjugate to $(-1)^{p-1}(\caI+\caJ)^{p}$ but the point,
as already mentioned, is that only for $k\in\{0,1\}$ or $\{\pm\frac{1}{2}\}$ do we get an explicit conjugation with simple, lower-diagonal
matrices like $\caL,\caR$. As for the choice of  $(-1)^{p-1}(\caI+\caJ)^{p}$ rather than  $(-1)^{p-1}\caI+\caJ $ as normal form,
it is simply a matter of convenience, and a further, quite elementary conjugation, immediately takes us from the one to the other.
\medskip

\noindent
{\bf Defining identities for the normalising matrices $\caL,\caR$\,. }
\begin{eqnarray}  \label{b7}
\caR&=& (\caI+\caJ)\;\caL\;\caQ^{-1}
\\  \label{bb7}
\caR &=& \caL\; \caW
\end{eqnarray}
with the twisted shift matrix $\caQ$ defined right at the beginning of \S6.8 and with
\begin{eqnarray*} 
\caW&=\;\caM_{1\,,\,0}  &\quad \mi{if}\;\;p=0\;\;\mi{or}\;\;1\; \mod 4
\\
\caW&=\;\caM_{\frac{1}{2},-\frac{1}{2}} &\quad \mi{if}\;\;p=2\;\;\mi{or}\;\;3\; \mod 4
\end{eqnarray*}
The matrix entries of $\caW$ are elementary binomial coefficients\,:
\[\begin{array}{lllllll}
\mi{if} \;\; i<j   &:\caW[i,j]=0
\\
\mi{if} \;\; i=j   &:\caW[i,j]=1
\\
\mi{if}\;\;  i>j \;\; \mi{and} \dots&
\\
p\in\{0,1\}\!\!\!\!\mod 4\;\;\mi{and}\;\; p\!-\!i\!-\!j\in\{1,2\} &: \caW[i,j]= -  \frac{p!}{(i-j)!(p-i+j)!}
\\
p\in\{2,3\}\!\!\!\!\mod 4\;\;\mi{and}\;\;p\!-\!i\!-\!j\in\{0,1\} &: \caW[i,j]= -  \frac{p!}{(i-j)!(p-i+j)!}
\\
\mi{otherwise} &:\caW[i,j]=0
\end{array}\]
If we now eliminate either $\caR$ (resp. $\caL$) from the system (\ref{b7}),(\ref{bb7}) and express the remaining
matrix as a sum of an {\it elementary part} (which corresponds to the two extreme subdiagonal zones and
carries only binomial entries) and a {\it complex part} (which corresponds to the middle subdiagonal zone
and involves the intricate polynomials $H^{\delta}_{d}$), we get a linear system which, as it turns out,
completely determines  $\caL^{\mi{comp.}}$ or  $\caR^{\mi{comp.}}$ (viewed as unknown) in terms 
$\caL^{\mi{elem.}} $ or $\caR^{\mi{elem.}} $
(viewed as known). Thus\,:
\begin{equation}  \label{b8}
(\caI+\caJ)\; (\caL^{\mi{elem.}}+\caL^{\mi{comp.}})=
(\caL^{\mi{elem.}}+\caL^{\mi{comp.}})\; \caW\; \caQ
\end{equation}
To understand just how special the value $\epsilon=-1$ and the case of full-rotation matrices are, let us briefly examine,
first, the case of full rotations with unassigned $\epsilon$, then the case of partial rotations with $\epsilon=-1$.
\medskip

\noindent
{\bf Complement: full rotations with $\epsilon \not=-1 $\,.}
\\
Keeping $\epsilon$ free and setting $V_{p}(t,\epsilon):=\det(t\,\caI-\caM_{p,0})$ we get\,:
\begin{eqnarray*}
V_{2}(t,\epsilon)&=&(t+1)^2\,-2^{2} \,\,(1+\epsilon)\, t
\\
V_{3}(t,\epsilon)&=&(t-1)^3\,+3^{3}\,\, (1+\epsilon)\, t\, \epsilon
\\
V_{4}(t,\epsilon)&=&(t+1)^4\,-2^{3}\,\,(1+\epsilon)\,t\,( 1+16\,\epsilon+32\,\epsilon^2 +14\epsilon\, t+t^2 ) 
\\
V_{5}(t,\epsilon)&=&(t-1)^5\,+5^4\,\,(1+\epsilon)\, t\, \epsilon\, (1+5\,\epsilon+5\,\epsilon^2+3\,t+t^2)
\end{eqnarray*}
Etc\dots.The only conspicuous properties of the $V_{p}$ polynomials seem to be\,:
\begin{eqnarray}  \label{b11}
V_{p}(t,-1)&=& \big(t+(-1)^{p}\big)^{p}
\\  \label{b12}
V_{p}(1,\epsilon)&=& \epsilon^{p}\;V_{p}(1,\epsilon^{-1})
\end{eqnarray}
(\ref{b11}) follows from the short analysis argument given in \S6.7, and we have devoted the bulk of the present
section (\S6.8) to checking it algebraically. As for the self-inversion property (\ref{b12}), it directly follows from the way the
simple-crossing matrices $\caM_{k}$ are constructed. As far as we can see, the  $V_{p}$ polynomials appear to possess
only one additional property, albeit a curious one (we noticed it empirically and didn't attempt a proof). It is this\,: 
for $p$ prime $\geq 5$ and $t=1$ we have (at least up to $p=59$)\,:
\begin{equation}  \label{b13}
V_{p}(1,\epsilon)=\det(\caI-\caM_{p,0})=p^{p}\,\epsilon\,(1+\epsilon)\,(1+\epsilon+\epsilon^{2})^{\kappa(p)}\,W_{p}(\epsilon)
\end{equation}
with $\kappa(p)\!=\!1$ (resp.$2$) if  $p=2$ (resp.$1$) \!$\mod 4$ and some $\doQ$-irreducible polynomial
$W(\epsilon)\in \doZ(\epsilon)$. However, $V_{p}(t,\epsilon)\not=0 \mod 1+\epsilon+\epsilon^{2 }$, which
reduces the above relation (\ref{b13}) to a mere oddity. 
 \footnote{\;
True, we have $V_{p}(1,\epsilon)=\mi{Const} \mod 1+\epsilon+\epsilon^{2 }$, but this is a trivial consequence
of  $V_{p}(1,\epsilon)$ being self-inverse in $\epsilon$.}
More generally, the ``semi-periodicity''
in $k$ of $\caM_{k,0}$ that we noticed for $\epsilon=-1$  has no counterpart for any other value of $\epsilon$, not
even for $\epsilon^{3}=1$ or, for that matter, $\epsilon=1$.
\medskip

 \noindent
{\bf Complement: partial rotations with $\epsilon=-1$\,.}\\
The partial-rotation matrices $\caM_{k_{2},k_{1}}$ with $ |k_{2}-k_{1}| \leq \frac{p}{2}$ all share the same trivial
characteristic polynomial $(t-1)^{p}$, but possess increasingly numerous and increasingly large Jordan blocks as 
$ |k_{2}-k_{1}|$ goes from $0$ to $\frac{p}{2}$. For $\frac{p+1}{2} < |k_{2}-k_{1}| < p $, the Jordan blocks disappear and 
the characteristic polynomials become thoroughly unremarkable, apart from being self-inverse 
({\it always} so if $p$ is even, {\it only when} $k_{2}-k_{1}\in \doZ $ if $p$ is odd). 
For $ |k_{2}-k_{1}| = p $, as we saw earlier in this section, we have
one single Jordan block of maximal size, with eigenvalue $\mp1$ depending on the parity of $p$. That leaves only
the border-line case $ |k_{2}-k_{1}|=\frac{p+1}{2}$. We have no Jordan blocks then, yet the characteristic polynomials
possess a remarkable factorisation on $\doZ$\,:
\begin{eqnarray*}
\mi{If}\;\; k_{2}-k_{1}=\pm\frac{p+1}{2}\;\; \mi{then}\;:
\\
(\mi{for} \;\;p\;\;\mi{odd})\quad\quad\quad\quad\det(t\,\caI-\caM_{k_{2},k_{1}})\!\!&\!\!=\!\!&\!\!(t-1)\,\prod_{s=1}^{\frac{p-1}{2}} P_{s}(p,t)
\\
(\mi{for} \;\;p=0\! \mod 4)\quad\quad\det(t\,\caI-\caM_{k_{2},k_{1}})\!\!&\!\!=\!\!&\quad\quad\quad\prod_{s=1}^{\frac{p}{4}}\big(P_{2s-1}(p,t)\big)^{2}
\\
(\mi{for} \;\;p=2\! \mod 4)\quad\quad\det(t\,\caI-\caM_{k_{2},k_{1}})\!\!&\!\!=\!\!&\!\! P_{\frac{p}{2}}(p,t)\prod_{s=1}^{\frac{p-2}{4}}\big(P_{2s-1}(p,t)\big)^{2}
\end{eqnarray*}
with polynomials $P_{s}(p,t)\in \doQ[p,t]$ quadratic and self-inverse in $t$, 
of degree $2s$ in $p$, and assuming values in $\doZ[t]$ for $p\in\doZ $\,:
\begin{equation*}
P_{s}(p,t) = \big(1-t\big)^{\!2}+\Big(\prod_{i=0}^{s-1}\frac{p-i}{1+i}\Big)^{\!\!2}\,t
\end{equation*}
\\ \noindent
{\bf Complement: some properties of the polynomials $H_{d}^{\delta}$\,.}
\\
For any fixed $n,d\in \doN$ with $n\leq d$, the $H^{\delta}_{d}(x,n)$ and $H^{\delta}_{d}(n,y)$, as polynomials in $x$ or $y$,
factor into a string of fully explicitable one-degree factors. This immediately follows from the expansions (\ref{asa3}),(\ref{asa4}),(\ref{asa5}). Conversely,
the factorisations may be directly derived from the induction (\ref{asa1}),(\ref{asa2}) and then serve to establish the remaining properties. Most zeros $(x,y)$
in $\doZ^{2}$ or $(\frac{1}{2}\doZ)^{2}$ can also be read off the factorisation. All the above properties suggest a measure
of symmetry between the two variables, under the simple exchange $x\leftrightarrow y$. But there also exists a more 
recondite symmetry, which is best expressed in terms of the polynomials 
\begin{equation}  \label{b20}
K^{\delta}_{d}(x,y):=H^{\delta}_{d}(d-x,\frac{1}{2}+2d+\delta-\frac{3}{2}x+\frac{1}{2}\,y)
\end{equation}
under the exchange $y\leftrightarrow -y$. It reads, for $x=n\in\doN\cup[0,d]$\,:
\begin{eqnarray*}
K^{\delta}_{d}(n,y)+K^{\delta}_{d}(n,-y)&=& 
2^{n-1}\Big[\!\Big[\frac{d\!+\!\delta\!-\!n }{2d\!+\!\delta\!-\!2n} \Big]\!\Big]!!
\prod^{i\,\mi{odd}}_{0<i<n}(y^{2}-i^{2})
 \hspace{1.3 ex} (n\;\mi{even})
\\                                                        
                                       &=& \;\;0 \hspace{32.ex} (n\;\mi{odd})
\end{eqnarray*}
with the factorial ratio $ [\![\dots]\!]!!$ defined as in (\ref{abc3})

%%%%%%%%%%%%%%%%%%%%%%%%%%%%%%%%%%%%%%%%%%%%%%%%%%%%%%%%%%%%%%%%%%%%%%%%%%%
%%%%%%%%%%%%%%%%%%%%%%%%%%%%%%%%%%%%%%%%%%%%%%%%%%%%%%%%%%%%%%%%%%%%%%%%%%%
\subsection{Ramified monomial inputs $F$\,: infinite order ODEs.}
If we now let $p$ assume arbitrary complex values $\alpha$, our {\it nir}-transform $h(\nu)$ and its centered
variant  $h_{\ast}(\nu)=h(\nu+\nu_{\ast})=h(\nu+1)$ ought to verify the following ODEs of infinite order
\begin{eqnarray} \label{uhu1}
Q(\partial_{\nu},\nu)\,h(\nu):=
\Big( (\partial_\nu -\nu\,\partial_\nu-\frac{\alpha}{2})^\alpha-\partial_\nu^\alpha\Big) h(\nu)=0
\quad \quad\quad\alpha \in \doC
\\ \label{uhu2}
Q_{\ast}(\partial_{\nu},\nu)\,h_{\ast}(\nu):=
\Big( ( -\nu\,\partial_\nu-\frac{\alpha}{2})^\alpha-\partial_\nu^\alpha \Big) h_{\ast}(\nu)=0
\quad \quad\quad\alpha \in \doC
\end{eqnarray}
to which a proper meaning must now be attached. This is more readily done with the first, non-centered variant, since
\\
\begin{proposition}.
\\ The $\mr{nir}$-transform $h_{\alpha}(\nu)$ of $f_{\alpha}(x):=\alpha\,\log(1+\alpha\,x)$ is of the form
\begin{equation} \label{b30}
h_{\alpha}(\nu)=-h_{-\alpha}(\nu)=\frac{1}{\sqrt{2}\,\alpha}\sum_{n\in \doN}\gamma_{-\frac{1}{2}+n}(\alpha^{2})\,\nu^{-\frac{1}{2}+n}
\end{equation}
with $ \gamma_{-\frac{1}{2}+n}(\alpha^{2})$ polynomial of degree $n$ in $\alpha^{2}$ and it verifies (mark the sign change) an infinite integro-differential
equation of the form
\begin{equation} \label{b31}
\Big(\sum_{1\leq k} \partial_{\nu}^{-k}\;\,S_{k}(\nu\partial_{\nu}+\frac{k}{2},\alpha-k)\Big)\,h_{\alpha}(-\nu)=0
\end{equation}
with  integrations $\partial_{\nu}^{-k} $ from $\nu=0$ and with elementary differential operators $\doS(.,.) $ which, being
polynomial in their two arguments, merely multiply each monomial $\nu^{n}$ by a scalar factor polynomial in $(n,k,\alpha)$,
effectively yielding an infinite induction for the calculation of the coefficients $\gamma_{-\frac{1}{2}+n}(\alpha^{2})$.
\end{proposition}
Thus the first three coefficients are 
$$\gamma_{-\frac{1}{2}}(\alpha^{2})= 1
\;\;,\;\;
\gamma_{\frac{1}{2}}(\alpha^{2})=\frac{1}{12}(\alpha^{2}-1)
\;\;,\;\;
\gamma_{\frac{3}{2}}(\alpha^{2})=\frac{1}{864}(\alpha^{2}-1)\,(\alpha^{2}+23)
 $$
 For $n\geq 1$ all polynomials $\gamma_{-\frac{1}{2}+n}(\alpha^{2})$ are divisible by $(\alpha^{2}-1)$ but this is their only common factor.
 
Remark\,: the regular part of the {\it nir}-transform $h_{\alpha}$ of $f_{\alpha}$ has the same shape 
$ \sum_{n\in\doN}\alpha^{-1}\,\gamma_{n}(\alpha^{2}) $ as the singular part, also with $ \gamma_{n}(\alpha^{2})$
polynomial of degree $n$ in $\alpha^{2}$, but it doesn't verify the integro-differential equation (\ref{uhu1}). We'll need
the following identies\,:
\\
\begin{equation} \label{b35}
[\mgd,\mgD]=\mgd \quad \quad (\mi{here}\quad \quad \mgd=\partial_{\nu}\;\;,\;\; \mgD=\nu\partial_{\nu}+\frac{\alpha}{2})
\end{equation}
\begin{eqnarray*}
(\mgd+\mgD)^{\alpha}&=& \sum_{0\leq k}\hspace{5.5 ex} S_{k}(\mgD+\frac{\alpha-k}{2},\alpha-k)\;\mgd^{\alpha-k}
\\
&=& \sum_{0\leq k} \mgd^{\frac{\alpha-k}{2}}\;S_{k}(\mgD\hspace{8.5 ex},\alpha-k)\;\mgd^{\frac{\alpha-k}{2}}
\\
&=& \sum_{0\leq k} \mgd^{\alpha-k}\;S_{k}(\mgD-\frac{\alpha-k}{2},\alpha-k)
\end{eqnarray*}
The non-commutativity relation $[\mgd,\mgD]=1 $, combined with the above expansions, yields for the polynomials
$S_{k}$ the following addition equation\,:
\begin{equation*}
S_{k}(\mgD,\beta_{1}+\beta_{2})=\sum_{k_{1}+k_{2}=k}
S_{k_{1}}(\mgD-\frac{\beta_{2}-k_{2}}{2},\beta_{1}-k_{1})\,
S_{k_{2}}(\mgD+\frac{\beta_{1}-k_{1}}{2},\beta_{2}-k_{2})\,
\end{equation*}
and the difference equation\,:
\begin{equation} \label{b36}
S_{k}(\frac{\beta}{2},\beta)\equiv S_{k}(\frac{\beta+1}{2},\beta-1)
\end{equation}
That relation, in turn, has two consequences\,: on the one hand, it leads to a finite expansion (\ref{b37}) of $S_{k}(\mgD,\beta)$ in powers
of $\mgD$ with coefficients $T_{2k_{\ast}}(\beta)$ that are polynomials in $\beta$ of degree exactly $k_{\ast}$ with $2k_{\ast}\leq k$.
On the other, it can be partially reversed, leading, for entire values of $b$, to a finite expansion (\ref{b38}) 
of $T_{2k}(b)$ in terms of some special values of $S_{2k-1}(\,.\,,b)$.
\begin{eqnarray} \label{b37}
S_{k}(\mgD,\beta)&=&\Big( \prod_{i=1}^{k}(\beta+i)\Big)
\sum_{k_{1}+2\,k_{2}=k }^{ k_{1},k_{2}\geq 0}\;
\frac{\mgD^{k_{1}}}{k_{1}!}\;\frac{T_{2 k_{2}}(\beta)}{(2\,k_{2})!}\quad\quad\  \forall \beta\in \doC
\\ \label{b38}
T_{2 k}(b) &=& 
 \frac{(2 k)! \;b!}{(2 k+b)!}\,\sum_{0 \leq c \leq b} \; (c-\frac{b}{2}) \; S_{2k-1}( \frac{c-b}{2}, c )
 \quad\;  \forall b \in \doN
\end{eqnarray} 
Together, (\ref{b37}) and its reverse (\ref{b38}) yield an explicit inductive scheme for constructing the polynomials $T_{2k}$\,. We
first calculate $T_{2k}(b)$ for $b$ whole, via the identity (\ref{b39}) whose terms $S_{2k-1}(\,.\,,b)$ 
involve only the earlier polynomials $T_{2h}(c)$, with indices $ h < k$ and $c \leq b $. The identity reads\,:
\begin{equation} \label{b39}
\frac{T_{2k}(b)}{(2k)!b!}=\sum_{0\leq c \leq b}^{0 \leq h < k} \frac{T_{2h}(c)}{(2h)!c!}
\, \frac{(c/2-b/2)^{2k-2h-1}}{(2\,k-2\,h-1)!} \,\frac{(2k+c)!}{(2k+b)!}\,\frac{(c-b/2)}{(c+2\,k)}
\end{equation}
 Then we use Lagrange interpolation (\ref{b40})-(\ref{b41}) to calculate
$T_{2k}(\beta)$ for general complex arguments $\beta$ \,:
\begin{eqnarray} \label{b40}
T_{2 k}(\beta) &=& 
\sum_{1 \leq b \leq k} \;\Lambda_{k}(\beta,b)\; T_{2 k}(b)
 \quad\quad\quad\;  \forall \beta \in \doC
 \\ \label{b41}
 \Lambda_{k}(\beta,b)&:=& \frac{\beta}{b}\prod_{1\leq i \leq k}^{i\not=b}\frac{i-\beta}{i-b}
\end{eqnarray} 
\\
{\bf First values of the $T_{2\,k}$-polynomials\,:}
\begin{eqnarray*}
T_{0}(\beta)&=& 1
\\
T_{2}(\beta)&=& \frac{1}{12}\,\beta
\\
T_{4}(\beta)&=& \frac{1}{240}\,\beta\,(-2+5\,\beta)
\\
T_{6}(\beta)&=& \frac{1}{4032}\,\beta\,(16+42\,\beta+35\,\beta^{2})
\\
T_{8}(\beta)&=& \frac{1}{34560}\,\beta\,(-4+5\,\beta)(36-56\,\beta+35\,\beta^{2})
\\
T_{10}(\beta)&=&\frac{1}{101376}\,\beta\, (768-2288\,\beta+2684\,\beta^{2}-1540\,\beta^{3}+385\,\beta^{4})
\end{eqnarray*}
{\bf Special values of the $T_{2\,k}$-polynomials\,:}
\begin{eqnarray*}
T_{2\,k}(2)&=& \frac{2}{(2\,k+1)(2\,k+2)}
\\
T_{2\,k}(1)&=& \frac{1}{(2\,k+1)\,2^{{2\,k}}}
\\
T_{2\,k}(0)&=&0 \quad \mi{if}\quad  k\not=0 \quad \mi{and} \quad T^{\,\prime}_{2\,k}(0)=\frac{B_{2\,k}}{2\,k}
\\
T_{2\,k}(-1)&=&B_{2\,k}(\frac{1}{2})
\\[1. ex]
T_{2\,k}(-2)&=&-(2\,k-1)\,B_{2\,k}
\\[0.7 ex]
T_{2\,k}(-1-2\,k)&=&(-1)^{k}\,\frac{(2\,k)!}{4^{k}\,k!}
\end{eqnarray*}
with $B_{n}$ and $B_{n}(.)$ denoting the Bernoulli numbers and polynomials.
\\
\\
{\bf Special values of the $S_{k}$-polynomials\,:}
\begin{eqnarray*}
\mi{For}\, k\; \mi{odd}\;:\quad S_{k}(\mgD,-1-k)&=&\prod_{-\frac{k}{2}<s<\frac{k}{2} }^{ k\in \doZ}\;(\mgD+s)
\\
\mi{For}\, k\; \mi{even}\;:\quad S_{k}(\mgD,-1-k)&=&\prod_{-\frac{k}{2}<s<\frac{k}{2} }^{ k\in \doZ-\frac{1}{2}\doZ}\;(\mgD+s)
\end{eqnarray*}
Note that in neither case are the bounds $\pm k/2$ reached by $s$, since the pair $\{k/2,s\}$ always consists of an integer and a half-integer.
$S_{k}(\mgD,b)$ appears to have no simple factorisation structure except (trivially) for $b=1$ and $b=2$ when in view of (\ref{b37}),(\ref{b38}) we have\,:
\begin{eqnarray*}
S_{k}(\mgD,1) &=& 2^{-k-1}\Big( (2\,\mgD+1)^{k+1} -(2\,\mgD-1)^{k+1}\Big)
\\
S_{k}(\mgD,2) &=& 2^{-1}\Big( (\mgD+1)^{k+2} +(\mgD-1)^{k+2}-2\,\mgD^{k+2}\Big)
\end{eqnarray*}
Since $S_{1}(n,\alpha-1)= n\,\alpha $, the induction rule for the $\gamma$-coefficients reads
\begin{eqnarray*}
\gamma_{-\frac{1}{2}}(\alpha^2) \!\!&\!=\!\!& 1
\\
\gamma_{-\frac{1}{2}+n}(\alpha^2) \!\!&\!=\!\!&\!\! 
\sum_{1\leq k\leq n}{(-1)^{k+1}}\,\frac{\Gamma(\frac{1}{2}+n-k)}{\Gamma(\frac{1}{2}+n)}\,
\frac{S_{k+1}(n-\frac{k}{2},a-k-1)}{S_{1}(n,\alpha-1)}\,
\gamma_{-\frac{1}{2}+n-k}(\alpha^2)
\\
 \!\!&\!=\!\!&\!\! 
\sum_{1\leq k\leq n}{(-1)^{k+1}}\,\frac{\Gamma(\frac{1}{2}+n-k)}{\Gamma(\frac{1}{2}+n)}\,
\frac{S_{k+1}(n-\frac{k}{2},a-k-1)}{n\,\alpha}\,
\gamma_{-\frac{1}{2}+n-k}(\alpha^2)
\end{eqnarray*}
Moreover, since $ S(k)(\mgD,-1)=S(k)(\mgD,-2)=\dots =S(k)(\mgD,-k)=0$, when $\alpha$ is a positive integer, the above
induction involves a constant, finite number of terms, with a sum $\sum$ over $1\leq k \leq \alpha-1$ instead of
 $1\leq k \leq n$, which is consistent which the finite differential equations of \S6.6.

%%%%%%%%%%%%%%%%%%%%%%%%%%%%%%%%%%%%%%%%%%%%%%%%%%%%%%%%%%%%%%%%%%%%%%%%%%%%%%%%%%%%%%%%%%%%
%%%%%%%%%%%%%%%%%%%%%%%%%%%%%%%%%%%%%%%%%%%%%%%%%%%%%%%%%%%%%%%%%%%%%%%%%%%%%%%%%%%%%%%%%%%%

\subsection{Ramified  monomial inputs $F$\,: arithmetical aspects.}
In this last subsection, we revert to the case of polynomial inputs $f$ and replace the high-order ODEs verified by $k$ by
a first-order order differential system, so as to pave the way for a future paper [SS2] devoted to understanding, from
a pure ODE point of view, the reasons for the rigidity of the {\it inner algebra's} resurgence, i.e. its surprising insentivity
to the numerous parameters inside $f$.

The normalised coefficients $\gamma_{r},\delta_{r},\delta^{\mr{ev}}_{r}$ of
the series $h_{\alpha},k_{\alpha},k^{\mr{ev}}_{\alpha}$, whose definitions we recall\,:
\begin{eqnarray*}
h_{\alpha}(\nu)&\!=\!& \frac{1}{\sqrt{2}}\frac{1}{\alpha}\sum_{r\in -\frac{1}{2}+\doN}\gamma_{r}(\alpha^{2})\;\nu^{r}
\quad\;\; \mi{with}\quad \gamma_{-\frac{1}{2}}(\alpha^{2})\equiv 1
\\
k_{\alpha}(n)&\!=\!& \frac{\sqrt{\pi}}{\sqrt{2}}\frac{1}{\alpha}\sum_{r\in -\frac{1}{2}+\doN}\delta_{r}(\alpha^{2})\;n^{-r}
\quad \mi{with}\quad \delta_{-\frac{1}{2}}(\alpha^{2})\equiv 1
\\
k_{\alpha}(n)k_{\alpha}(-n) 
=:k^{\mr{ev}}_{\alpha}(n)
&\!=\!& \frac{\pi i}{2\alpha^{2}}\;\;\sum_{r\in -1+2\,\doN}\delta^{\mr{ev}}_{r}(\alpha^{2})\;n^{-r}
\;\;\; \mi{with}\quad \delta^{\mr{ev}}_{-1}(\alpha^{2})\equiv 1
\end{eqnarray*} 
seem to possess remarkable arithmetical properties, whether we view them
\\
(i) as polynomials in $\alpha$
\\
(ii) as polynomials in $r$
\\
(iii) as rational numbers, for $\alpha$ fixed in $\doZ$.
\\
These arithmetical properties, at least some of them, do not obviously
follow from the shape of the {\it nir} transform nor indeed from the above induction. 
Thus, as {\it polynomials}, they $\gamma$ coefficients appear to be
exactly of the form\,:
\begin{equation} \label{b50}
\gamma_{-\frac{1}{2}+r}(\alpha^{2})=
\frac{6^{-r}}{(2r)!}
\frac{(\alpha^{2}-1)\,\gamma_{-\frac{1}{2}+r}^{\ast}(\alpha^{2})}
{\prod_{5 \leq p \leq r+2 }p^{\mu_{r,p}}}
\,=\, \,\frac{1}{(2r)!}\,\frac{a^{r}}{6^{r}}\,\gamma_{-\frac{1}{2}+r}^{\ast\ast}(\alpha^{2}-1)
\end{equation}
with the $\gamma^{\ast}_{-\frac{1}{2}+r}(\alpha^{2}) $ irreducible in $\doZ[\alpha^{2}]$ and 
with on the denominator a product $\prod$ 
involving only prime numbers between 5 and $r\!+\!2$. 
This at any rate holds for all values of $n$ up to 130. The surprising
thing is not the presence of  these $p$ in $[5,r\!+\!2]$ but rather the fact that their powers $\mu_{r,p}$ seem to obey
no exact laws (though they are easily majorised), unlike the powers of 2 and 3 that are {\it exactly} accounted for by the factor $6^{-r}$.
But this 2- and 3-adic regularity seems to go much further. It becomes especially striking if we consider
the polynomials $\gamma_{-\frac{1}{2}+r}^{\ast\ast} $ in the rightmost term of (\ref{b50}) after changing to the variable $a:=\alpha^{2}-1$.
Indeed\,: 
\begin{conjecture}[2- or 3-adic expansions for the $\gamma$ as $\alpha$-polynomials].\\
The polynomials $\gamma^{\ast\ast}$ defined by
\begin{equation} \label{b51}
\gamma_{-\frac{1}{2}+r}(a+1)\,=: \,\frac{1}{(2r)!}\,\frac{a^{r}}{6^{r}}\,\gamma_{-\frac{1}{2}+r}^{\ast\ast}(a)
\quad\quad,\quad\quad
\gamma_{-\frac{1}{2}+r}^{\ast\ast} \in \doQ[a^{-1}]
\end{equation}
possess 2- and 3-adic expansions to all orders\,:
\begin{eqnarray*}
\gamma_{-\frac{1}{2}+r}^{\ast\ast}(a) &=&\sum_{0\leq j}\lambda_{2,j}(r,a)\,2^{j}
\quad\quad\quad\quad \Big(\lambda_{2,j}(r,a)\in \{0,1\}\Big)
\\
\gamma_{-\frac{1}{2}+r}^{\ast\ast}(a) &=&\sum_{0\leq j}\lambda_{3,j}(r,a)\,3^{j}
\quad\quad\quad\quad \Big(\lambda_{3,j}(r,a)\in \{0,1,2\}\Big)
\end{eqnarray*}
with coefficients $\lambda_{p,j}$ that in turn depend only on the first ${j}$ terms of the p-adic expansion of $\mr{r}$. In other words\,:
\begin{eqnarray*}
\lambda_{2,j}(r,a)=\lambda_{2,j}([r_{0},r_{1},\dots,r_{j-1}],a)\quad\quad \mi{with}\quad r=\sum_{0\leq i < j}r_{i}\,2^{i}\;\;\mod 2^{j}
\\
\lambda_{3,j}(r,a)=\lambda_{3,j}([r_{0},r_{1},\dots,r_{j-1}],a)\quad\quad \mi{with}\quad r=\sum_{0\leq i < j}r_{i}\,3^{i}\;\;\mod 3^{j}
\end{eqnarray*}
Moreover, as a polynomial in ${a^{-1}}$, each $\lambda_{2,j}(r,a)$ is of degree $j$ at most.
\end{conjecture}
These facts have been checked up to the p-adic order $j=25$ and for all $r$ up to 130. Moreover, no such
regularity seems to obtain for the other p-adic expansions, at any rate not for $p=5,7,11,13$.
\begin{conjecture}[p-adic expansions for the $\gamma$ as $r$-polynomials].\\
The  $\gamma^{\ast\ast}$ defined as above verify
\begin{equation} \label{b52}
\gamma^{\ast\ast}_{-\frac{1}{2}+r}(a)\,= 1+\sum_{1\leq d \leq r} a^{-d}\,Q_{d}(r)
\end{equation}
with universal polunomials $Q_{d}(r)$ of degree $3d$ in $r$ and of the exact arithmetical form\,:
\begin{eqnarray*}
Q_{d}(r)&=&\Big(6^{d}\prod_{p\,\mi{prime}\geq 2}p^{-\mu_{p}(d)}\Big)\;\Big(Q^{\ast}_{d}(r)\prod_{1\leq i \leq d}(r-i)\Big)
\\
Q^{\ast}_{d}(r)&=&\sum_{1\leq i \leq 2d}c_{d,i}\,r^{i}\hspace{4.ex}\mi{with}\hspace{4.ex} (c_{d,1},\dots,c_{d,2d})\;\;\;\;\mi{coprime}
\\ 
\mu_{p}(d) &=& \sum_{0\leq s}\mr{en}\Big(\frac{2\,d}{(p\!-\!3)\,p^{s}}\Big)\hspace{3.ex} \mi{if}\;\;\; p\geq 5
\\ 
\mu_{3}(d) &=& \sum_{0\leq s}\mr{en}\Big(\frac{d}{2\,p^{s}}\Big)
\\ 
 \mu_{2}(d) &=& -\sum_{0 \leq s} d_{s} \hspace{7.ex}\mi{if} \;\;\; d=\sum_{0 \leq s} d_{s} \,2^{s} \;\;\;(d_{s}\in \{0,1\} )
 \\
 c_{d,2d} &\in& (-1)^{d}\doN^{+}
\end{eqnarray*}
\end{conjecture}
with $\mr{en}(x)$ denoting as usual the entire part of $x$.
\begin{conjecture}[Special values of the $Q_{d}$ ].\\
If for any $ q\in \doQ $ we set 
\begin{eqnarray*}
\mr{pri}(q)&:=& q \quad \mi{if} \quad q\in \doN \quad \mi{and}\quad q\;\;\;\mi{prime}
\\  
 \mr{pri}(q)           &:=& 1 \hspace{7.ex} \mi{otherwise}
\end{eqnarray*}
then for any $d,s\in\doN^{\ast}$ we have 
\begin{equation} \label{b53}
Q_{d}(d\!+\!s)\in \frac{1}{Q_{d,s}}\doZ \quad \mi{with}\quad Q_{d,s}:=\prod_{{1\leq j \leq s \atop s <k \leq s+2\,j}} \mr{pri} 
\Big(\frac{d\!+\!k}{j}\Big)
\end{equation}
Moreover, for $s$ fixed and $d$ large enough, the denominator of $Q_{d}(d\!+\!s)$ is exactly $Q_{d,s}$. Note that by construction
 $Q_{d}(d\!+\!s)$ is automatically quadratfrei as soon as $d>2\,s\,(s\!-\!1)$.
\end{conjecture}
Together with the trivial identities $ Q_{d}(s)=0$ for $s\in[0,d]\cup \doN$, this majorises  $\mi{denom}(Q_{d}(s)$ for all $s\in \doN$.
We have no such simple estimates for negative values of $s$. 
\begin{conjecture}[ Coefficients $\delta^{\mr{ev}}$ ].\\
The normalised coefficients $\delta_{r}$ of $k_{\alpha}$ (with $r \in -1+2\,\doN$) are of the form\,:
\begin{equation}
\delta^{\mr{ev}}_{r}(\alpha^{2})
= \frac{1}{B_{r}}\,R_{r}(\alpha^{2})
= \frac{A_{r}}{B_{r}}\,R_{r}^{\ast}(\alpha^{2})\,\prod_{d | r}(\alpha^{2}-d^{2})
\end{equation}
where
\\
(i) $A_{r}$ is of the form $\prod_{p | r }^{p\, \mr{prime}} p^{\sigma_{p,r}} $ with $\sigma_{r,p}\in \doN$
\\
(ii) $B_{r}$ is of the form $\prod_{(p,r)=1 }^{p\, \mr{prime}\leq r+2} \,p^{\tau_{p,r}}$ with $\tau_{r,p}\in \doN$
\\
(iii) $R^{\ast}_{r}(\alpha^{2})$ is an irreducible polynomial in $\doZ[\alpha^{2}]$
\\ However, when $\alpha$ takes entire values $q$, the arithmetical properties of  $\delta^{\mr{ev}}_{r}(q^{2})$
become more dependent  on q  than r. In particular\,;
\\
(iv)
$ \mr{denom}(\delta_{r}(q^{2}))=\prod^{p \,\mr{prime}}_{p|q} d^{\kappa_{r,q,p}}$  with  $\kappa_{r,q,p} \in \doN$
and $\kappa_{r,p,p} \leq 3$ 
\\
(v) $ \mr{denom}(\delta_{r}(p^{2}))=p^{\kappa_{r,p,p}}$  with 
 $\kappa_{r,p,p} \leq 3$\;\;(  $p$ prime). 
\\
 This suggests a high degree of divisibility for $R_{r}^{\ast}(q^{2})$ and above all $R_{r}(q^{2})$, 
specially for q prime. In particular we surmise that\,:
\\
(vi) $R_{r}(p^{2})\in (p\!-\!1)!\;\doZ$ for $p$ prime.
\end{conjecture}
%%%%%%%%%%%%%%%%%%%%%%%%%%%%%%%%%%%%%%%%%%%%%%%%%%%%%%%%%%%%%%%%%%%%%%%%%%%%%%%%%%%%%%%%%%%%
%%%%%%%%%%%%%%%%%%%%%%%%%%%%%%%%%%%%%%%%%%%%%%%%%%%%%%%%%%%%%%%%%%%%%%%%%%%%%%%%%%%%%%%%%%%%

\subsection{From flexible to rigid resurgence.}
%%%%%%%%%%%%%%%%%%%%%%%%%%%%%%%%%%%%%%%%%%%%%%%%%%%%%%%%%%%%%%%%%%%%%%%%%%%%%%%%%%%%%%%%%%%%
%%%%%%%%%%%%%%%%%%%%%%%%%%%%%%%%%%%%%%%%%%%%%%%%%%%%%%%%%%%%%%%%%%%%%%%%%%%%%%%%%%%%%%%%%%%%
%%%%%%%%%%%%%%%%%%%%%%%%%%%%%%%%%%%%%%%%%%%%%%%%%%%%%%%%%%%%%%%%%%%%%%%%%%%%%%%%%%%%%%%%%%%%
We assume the tangency order to be 1. To get rid of the demi-entire powers, we go from $k$ to a
new unknown 
$K$ such that 
\begin{equation}   \label{a159}
k(n)=n^{\frac{1}{2}}\;K(n)\quad
\quad\quad (k(n)\in n^{\frac{1}{2}}\doC[[n^{-1}]]\;,\; K(n)\in \doC[[n^{-1}]])
\end{equation}
The new differential equation in the $n$-plane reads $ P(n,-\partial_n)K(n)=0$ and may be written
in the form\,:
\begin{equation}   \label{a160}
\Big(\prod_{i=1}^r(\partial_n+\nu_i)\Big)\,K(n)
+\sum_{i=0}^{r-1} \theta^\ast_i(n)\,\partial_n^iK(n)=0
\;\;\;;\;\;\;\theta^\ast_i(n)=O(n^{-1})
\end{equation}
This ODE is equivalent to the following first order differential system with $r$ unknowns
$K_i^\ast =\partial_n^iK\; ( 0\leq i \leq r-1)$\,:
\begin{eqnarray}
\partial_n K^\ast_0 -K^\ast_1 &=& 0 \nonumber
\\
\partial_n K^\ast_1 -K^\ast_2  &=& 0 \nonumber
\\
\dots && \nonumber
\\
\partial_n K^\ast_{r-2} -K^\ast_{r-1}  &=& 0 \nonumber
\\   \label{a161}
\partial_n K^\ast_{r-1} +\sum_{i=0}^{r-1}\big(\nu_{r-i}^\ast+\theta_i^\ast(n)\big)\,K_i  &=& 0
\end{eqnarray}
with $\nu_l$ standing for the symmetric sum of order $l$ of $\nu_1,\dots,\nu_r$\,.
Changing from the unknowns $K_i^\ast \;(0\leq i\leq r-1)$ to the unknowns
 $K_i\; (1\leq i\leq r)$ under the vandermonde transformation
$$
K_i^\ast=\sum_{0\leq j\leq r-1}(-\nu_i)^j\,K_j =\sum_{0\leq j\leq r-1}(-\nu_i)^j\,\partial^j_n K
$$
we arrive at a new differential system in normal form\,:
\begin{equation}   \label{a162}
\partial_n K_i+\nu_i\,K_i+\sum_{i=1}^{r}\theta_{i,j}(n)\,K_j =0
\quad\quad (1\leq i \leq r)
\end{equation}
with constants $\nu_i$ and rational coefficients $\tau_i(n)$ which,
unlike the earlier $\theta_i(n)$, are not merely $O(n^{-1})$ but also, crucially,  $O(n^{-2})$.
Concretely, the rank-1 matrix $\Theta=[\theta_{i,j}]$
is conjugate to a rank-1 matrix $\Theta^\ast=[\theta^\ast_{i,j}]$ with
only one non-vanishing (bottom) line, under the vandermonde 
matrix $V=[v_{i,j}]$ \,:
\begin{equation}   \label{a163}
\Theta=V^{-1}\,\Theta^\ast \,V \quad \mi{with}\;\;\; 
v_{i,j}=(-\nu_i)^{j-1}\;,\;
\theta^\ast_{i,j}=0\;\;\mi{if}\;i<r\;,
\;\theta^\ast_{r,j}=\theta^\ast_{j-1}
\end{equation}
The coefficients $\theta_{i,j}$ have a remarkable structure. They admit a unique
factorisation of the form\,:
\begin{equation}   \label{a164}
\theta_{i,j}(n) =  \frac{1}{\delta_i}\,\frac{\alpha_j(n)}{\gamma(n)}
\end{equation}
with factors $\alpha_i,\beta_j,\gamma $ derived from symmetric
polynomials $\alpha,\beta,\gamma $
of $r\!-\!1$ \;$f$-related variables (not counting the additional $n$-variable)\,:
\begin{eqnarray}   \label{a165}
\delta_i\;\;\; &=&\delta(x_1-x_i,\dots,\widehat{x_i-x_i},\dots,x_r-x_i)
\\   \label{a166}
\alpha_j(n)&=&\alpha(x_1-x_j,\dots,\widehat{x_j-x_j},\dots,x_r-x_j)(n)
\\   \label{a167}
\gamma(n) &=&\gamma(x_1-x_i,\dots,\widehat{x_i-x_i},\dots,x_r-x_i)(n)
\\   \label{a168}
&=&\gamma(x_1-x_j,\dots,\widehat{x_j-x_j},\dots,x_r-x_j)(n)
\end{eqnarray}
Let us take a closer look at all three factors\,:
\\
\\
(i) The $\gamma$ factor is simply the ``second leading polynomial" of \S6.2 and \S6.3
after division by its leading term $n^{\overline{\delta}}$. Being a direct
shift-invariant of $f$, if may also be viewed as a polynomial in $\mgf_0,\mgf_1\dots$.
\\
\\
(ii) The $\delta$ factor comes from the inverse vandermonde matrix 
$V^{-1}=[u_{i,j}]$. Indeed\,:
\begin{equation}   \label{a169}
u_{i,j}=\sigma_{r-j}(\nu_1,\dots,\nu_r)\,\delta_i
\quad\quad\quad\mi{with}\quad\quad\quad
\delta_i=\prod_{1\leq s\leq r}^{s\not= i}\frac{1}{\nu_s-\nu_i}=u_{i,r}
\quad
\end{equation}
Moreover\,:
\begin{eqnarray*}
\nu_i&=&f^\ast(x_i)=f_r\;\int_0^{x_i}\prod_{1\leq j \leq r}^{j\not= i} (x-x_j)\,dx
\\
\nu_i-\nu_j &=& \nu(x_i,x_j\,;\,x_1\,\dots,\hat{x}_i,\dots,\hat{x}_j,\dots,x_r)
\end{eqnarray*}
with a function $\nu$ antisymmetric (resp. symmetric) in its first two (resp. last $r\!-\!2$)
variables, and completely determined by the following two identities\,:
\begin{eqnarray*}
\nu(x_1,x_2\,;\, y_1,\dots,y_{r-2})&\!\equiv\!&
\nu(x_1-t,x_2-t\,;\, y_1-t,\dots,y_{r-2}-t)\quad\quad
\forall t
\\[ 1 ex]
\nu(x,-x\,;\, y_1,\dots,y_{r-2})&\!=\!&-\nu(-x,x\,;\, y_1,\dots,y_{r-2})
\\ 
&\!=\!& (-1)^{r-1}\,4\,f_r\!\sum_{1\leq s\leq [\frac{r-2}{2}]}
\frac{\sigma_{r-2s}(y_1,\dots,y_{r-2})}{4\,s^2-1}\;t^{2s+1}
\\
&\!=\!& \!-4\,f_r\int_0^x\!t_2\,dt_2\!\int_0^{t_2}
\!\Big( \prod_{i=1}^{r-2}(y_i+t_1)+\prod_{i=1}^{r-2}(y_i-t_1)\Big)\,dt_1
\end{eqnarray*}
 where $\sigma_l(y_1,y_2,\dots) $ denotes the  symmetric sum of order $l$ of $y_1,y_2\dots$.
\\
\\
(iii) The $\alpha$ factor stems from the coefficients $\theta_s$ in the
differential equation (\ref{a162}). Indeed, in view of (i) and (ii) and with the
$n$-variable implicit\,:
\begin{equation*}
\theta_{i,j}=\frac{1}{\delta_i}\,\frac{\alpha_j(n)}{\gamma(n)}
=\sum_{1\leq t\leq r \atop 1\leq s\leq r }u_{i,t}\,\theta^\ast_{t,s}\,v_{s,j}
=\sum_{1\leq s\leq r}u_{i,r}\,\theta^\ast_{r,s}\,v_{s,j}
=\sum_{1\leq s\leq r}\frac{1}{\delta_i}\,\theta^\ast_{s-1}\,(-\nu_j)^{s-1}
\end{equation*}
Hence 
\begin{equation} \label{ehe1}
\alpha_j(n) =\gamma(n)\sum_{0\leq l \leq r-1}(-\nu_j)^l\,\theta_l^\ast(n)
\end{equation}
We may also insert the covariant shift $\nu_0$ or $\underline{\nu}_0$
and rewrite the differential equation (\ref{a162}) as\,:
\begin{eqnarray}   \label{a170}
\Big(\prod_{i=1}^r(\partial_n+\nu_i)\Big)\,K(n)
+\sum_{i=0}^{r-1} \theta^\#_i(n)\,(\partial_n+\nu_0)^iK(n)=0
\hspace{9 ex}
\\  \label{a171}
\Big(\prod_{i=1}^r(\partial_n+\nu_i)\Big)\,K(n)
+\sum_{i=0}^{r-1} \underline{\theta}^\#_i(n)\,(\partial_n+\underline{\nu}_0)^iK(n)=0
\hspace{9 ex}
\end{eqnarray}
with new coefficients $\theta^\#_i,\underline{\theta}^\#_i$ 
that are not only shift-invariant but also root-symmetric.\footnote{i.e. symmetric
with respect to the roots $x_i$ of $f$.}
This leads  for the $\alpha$ factors to  expressions\,:
\begin{equation*}
\alpha_j(n) 
=\gamma(n)\sum_{0\leq l \leq r-1}(\nu_0-\nu_j)^l\,\theta_l^{\#}(n)
=\gamma(n)\sum_{0\leq l \leq
r-1}(\underline{\nu}_0-\nu_j)^l\,\underline{\theta}_l^{\#}(n)
\end{equation*} 
which have over (\ref{a166}) the advantage of involving only shift-invariants,
namely the $\theta^\#_l$ or $\underline{\theta}^\#_j$ (root-symmetric)
and the $\nu_0\!-\!\nu_j$ or $\underline{\nu}_0\!-\!{\nu}_j$ (not root-symmetric).
To show the whole extent of the {\it rigidity}, we may even introduce new parameters
by taking a non-standard shift operator $\beta(\partial_{\tau})$,
$$
\beta(t)=t^{-1}+\sum_{0\leq k}\beta_k\,t^k
 =t^{-1}+\sum_{1\leq k}\bee{k}\,t^{k-1}
\quad\quad (\bee{k}\equiv \beta_{k-1})
$$
but with a re-indexation $b_{k}=\beta_{k}$ to do justice to the underlying homogeneity.\footnote{N.B. the present $b_{k}$
differ from those in (\ref{a131}).} The case $r=1$ is uninteresting (no ping-pong, there being only one inner generator), 
and here are the results for $r=2$ and $3$.
\medskip

\noindent
{\bf Input $f $ of degree 2\,: } $f(x)=(x-x_1)\,(x-x_2)\,f_2 $
\begin{eqnarray*}
\delta(y_1) &=& -\frac{1}{6}\,f_2\,y_1^3
\\
\gamma(y_1) &=& 
+36\,f_2^4\,y_1^2
+288\,f_2^4\,b_2\,n^{-2}
\\
\alpha(y_1) &=& 
+(
5\,f_2^4\,y_1^2
-3\,f_2^6\,y_1^6\,b_2
)\,n^{-2}
+8\,(
f_2^6\,y_1^5\,b_3
-12\,f_2^5\,b_2\,y_1^3
)\,n^{-3}
\\&&
+8\,(
9\,f_2^5\,y_1^2\,b_3
-7\,f_2^4\,b_2
-3\,f_2^6\,y_1^4\,b_2^2)\,n^{-4}
\\&&
+64\,f_2^6\,y_1^3\,b_2\,b_3\,n^{-5}
+24\,\big(
12\,f_2^5\,b_2\,b_3
+f_2^6\,y_1^2\,(4\,b_2^3+b_3^2)
\big)\,n^{-6}\\&&
+128\,f_2^6\,b_2\,(4\,b_2^3+b_3^2)\,n^{-8}
\end{eqnarray*}
{\bf Input $f $ of degree 3\,: } $f(x)=(x-x_1)\,(x-x_2)\,(x-x_3)\,f_3 $
\\
$ \sy:=y_1+y_2\,,\,\py:=y_1\,y_2$
\begin{eqnarray*}
\delta(y_1,y_2) &\!\!=\!\!& \frac{1}{144}\,f_3^2\,\py^3\,(9\,\py-2\,\sy^2)
\\[1.5 ex]
\gamma(y_1,y_2) &\!\!=\!\!& 
-2^6\,3^3\,f_3^{10}\,\sy\,\py^2\,(9\,\py-2\,\sy^2)\,(4\,\py-\sy^2)
\\&&
-2^6\,3^3\,f_3^9\,\sy\,(9\,\py-2\,\sy^2)\,(3\,\py-\sy^2)\,n^{-1}
\\&&
+ 2^9\,3^3\,f_3^{10}\,b_2\,\sy\,(9\,\py-2\,\sy^2)\,(3\,\py-\sy^2)^2
\,n^{-2 }
\\&&
+2^6\,3^3\,\big(
243\,f_3^9\,\sy\,b_2\,(9\,\py-2\,\sy^2)
\\&&
+f_3^{10}\,b_3\,(513\,\sy^2\,\py^2+28\,\sy^6-252\,\sy^4\,\py+216\,\py^3)
\big)\,n^{- 3}
\\&&
- 2^6\,3^6\,(3\,\py-\sy^2)\,
\big(
33\,f_3^9\,b_3           
+4\,f_3^{10}\,b_2^2\,\sy  
\,(9\,\py-2\,\sy^2)
\big)\,n^{- 4}
\\&&
-  2^8\,3^7\,f_3^{10}\,b_2\,b_3\,(3\,\py-\sy^2)^2
\,n^{- 5} 
\\&&
+ 2^6\,3^9\,\big(
9\,f_3^9\,b_2\,b_3
-f_3^{10}\,b_3^2\,\sy\,(9\,\py-2\,\sy^2)
\big)\,n^{- 6}
\\&&
+ 2^6\,3^{12}\,f_3^{10}\,b_3\,(4\,b_2^3+b_3^2)
\,n^{- 9}
\end{eqnarray*}
For the $\alpha$ factor, we mention only the two lowest and highest powers of $n^{-1}$\,:
\begin{eqnarray*}
\alpha(y_1,y_2)\!\! &\!\!=\!\!&\! 
\Big(
\frac{1}{27}\,f_3^{11}\,\sy\,(9\,\py\!-\!2\,\sy^2)\,
(8748\,\py^5\!-\!13851\,\sy^2\,\py^4\!+\!378\,\sy^4\,\py^3\!+\!2403\,\sy^6\,\py^2
\\&&
-600\,\sy^8\,\py+40\,\sy^{10})
+\frac{2}{27}\,f_3^{13}\,b_2\,\sy\,\py^2\,(9\,\py-2\,\sy^2)\,(4\,\py-\sy^2)\times
\\&&
     \; \;\;(9\,\py^2+6\,\sy^2\,\py-\sy^4)\,(81\,\py^3-36\,\sy^2\,\py^2+9\,\sy^4\,\py-\sy^6)
    \Big)\,n^{- 2}
\\&&
+\Big(
f_3^{10}\,\sy\,(9\,\py\!-\!2\,\sy^2)\,(2835\,\py^3-675\,\sy^2\,\py^2-9\,\sy^4\,\py+\sy^6)
\\&&
+\frac{2}{27}\,f_3^{12}\,b_2\,\sy\,(9\,\py\!-\!2\,\sy^2)\,
(66339\,\py^6\!-\!129033\,\sy^2\,\py^5\!+\!175770\,\sy^4\,\py^4
\\&&
  -119475\,\sy^6\,\py^3+38520\,\sy^8\,\py^2-5742\,\sy^{10}\,\py+319\,\sy^{12})
\\&&
-\frac{8}{9}\,f_3^{13}\,b_3\,\sy^2\,\py^2\,(4\,\py\!-\!\sy^2)\,
(9\,\py^3\!-\!18\,\sy^2\,\py^2+9\,\sy^4\,\py\!-\!\sy^6)\times
\\&&
\,(9\,\py\!-\!2\,\sy^2)^2
\Big)\,n^{- 3}
+\sum_{4 \leq s \leq 18 }(\dots)\,n^{-s}
\\&&
+2\times3^{14}\,f_3^{13}\,b_3\,
(4\,b_2^3+b_3^2)\,
(b_2\,b_4-12\,b_2^3-3\,b_3^2)\,
b_4\,
(3\,\py-\sy^2)
\,n^{- 19}
\\&&
+\frac{3^{15}}{2}\,f_3^{13}\,b_3\,
(4\,b_2^3+b_3^2)\,
\big(
27\,(4\,b_2^3+b_3^2-b_2\,b_4)^2
-b_4^3
\big)\,n^{- 21}  
\end{eqnarray*}

\noindent
For a direct, ODE-theoretical derivation of the {\it rigidity} phenomenon, see [SS2].  General
criteria will also be given there for deciding which parameters inside an ODE contribute 
to the resurgence constants (or Stokes constants) and which don't.
% \end{document}
%%%%%%%%%%%%%%%%%%%%%%%%%%%%%%%%%%%%%%%%%%%%%%%%%%%%%%%%%%%%%%%%%%%%%%%%%%%%%%%%%%%%%%%%%%%%
%%%%%%%%%%%%%%%%%%%%%%%%%%%%%%%%%%%%%%%%%%%%%%%%%%%%%%%%%%%%%%%%%%%%%%%%%%%%%%%%%%%%%%%%%%%%
%%%%%%%%%%%%%%%%%%%%%%%%%%%%%%%%%%%%%%%%%%%%%%%%%%%%%%%%%%%%%%%%%%%%%%%%%%%%%%%%%%%%%%%%%%%%
%%%%%%%%%%%%%%%%%%%%%%%%%%%%%%%%%%%%%%%%%%%%%%%%%%%%%%%%%%%%%%%%%%%%%%%%%%%%%%%%%%%%%%%%%%%%
%%%%%%%%%%%%%%%%%%%%%%%%%%%%%%%%%%%%%%%%%%%%%%%%%%%%%%%%%%%%%%%%%%%%%%%%%%%%%%%%%%%%%%%%%%%%
%%%%%%%%%%%%%%%%%%%%%%%%%%%%%%%%%%%%%%%%%%%%%%%%%%%%%%%%%%%%%%%%%%%%%%%%%%%%%%%%%%%%%%%%%%%%
%%%%%%%%%%%%%%%%%%%%%%%%%%%%%%%%%%%%%%%%%%%%%%%%%%%%%%%%%%%%%%%%%%%%%%%%%%%%%%%%%%%%%%%%%%%%
%%%%%%%%%%%%%%%%%%%%%%%%%%%%%%%%%%%%%%%%%%%%%%%%%%%%%%%%%%%%%%%%%%%%%%%%%%%%%%%%%%%%%%%%%%%%
%%%%%%%%%%%%%%%%%%%%%%%%%%%%%%%%%%%%%%%%%%%%%%%%%%%%%%%%%%%%%%%%%%%%%%%%%%%%%%%%%%%%%%%%%%%%

%

%%%%%%%%%%%%%%%%%%%%%%%%%%%%%%%%%%%%%%%%%%%%%%%%%%%%%%%%%%
%% Some resurgence properties of knot-related functions.
%%%%%%%%%%%%%%%%%%%%%%%%%%%%%%%%%%%%%%%%%%%%%%%%%%%%%%%%%%

%\documentclass[12pt,a4paper]{article}\input{SP_commands}\begin{document}

%%%%%%%%%%%%%%%%%%%%%%%%%%%%%%%%%%%%%%%%%%%%%%%%%%%%%%%%%%%%%%%%%%%%%%%%%%%%%%%%%%%%%%%%%%%%
%%%%%%%%%%%%%%%%%%%%%%%%%%%%%%%%%%%%%%%%%%%%%%%%%%%%%%%%%%%%%%%%%%%%%%%%%%%%%%%%%%%%%%%%%%%%
%%%%%%%%%%%%%%%%%%%%%%%%%%%%%%%%%%%%%%%%%%%%%%%%%%%%%%%%%%%%%%%%%%%%%%%%%%%%%%%%%%%%%%%%%%%%

%\addtocounter{section}{}  

\section{The general resurgence algebra for SP series.} 
We recall the definition of the {\it raw} and {\it cleansed} SP series\,:
\begin{eqnarray}   \label{a172}
j_F(\zeta):=\sum_{0\leq n}\;J_F(n)\;\zeta^n
\quad &\mi{with}&\quad
J_F(n):=\sum_{0\leq m<n}\prod_{0\leq k\leq m}F(\frac{k}{n})\quad\quad
\\[1.ex]  \label{a173}
j^\#_F(\zeta):=\sum_{0\leq n}\;J^\#_F(n)\;\zeta^n
\quad &\mi{with}&\quad
J^\#_F(n):= J_F(n)/{I\!g}_{_{F}}(n)
\end{eqnarray}
We also recall that the $\perp$ transform turns the set $\{F,f,f^\ast,j_{F}^\# \}$
into the set $\{F^\perp,f^\perp,f^{\perp\,\ast},j_{F^\perp}^\#\}$ with\,:
\begin{eqnarray}  \label{a174}
 F^{\vdash}(x)=1/F(1-x) & ; &
 f^{\vdash}(x)=-f(1-x)\quad ; \quad
 f^{\vdash\,\ast}(x)=f^\ast(1-x)\quad  \quad
\\  \label{a175}
j_{F^\vdash}^\#(\zeta)=j_{F}^\#(\frac{\zeta}{\omega_F})
&& \mi{with} \;\;\; \omega_F:=e^{-\eta_F}\;\;\;\mi{and}\;\;\; \eta_F:=\int_0^1f(x)dx
\end{eqnarray}
We shall now (pending a more detailed investigation in [SS1]) sketch how 
the various generators arise
and how they reproduce under alien differentiation. Piecing all this information
together, we shall then get a global description of the Riemann surfaces
of our SP functions. 
\\

\noindent
For convenience, let us distinguish two degrees of difficulty\,:
\\ 
\--- first, the case of holomorphic inputs $f$
\\
\--- second, the case of meromorphic inputs $F$
\\
and split the investigation into two phases\,:
\\
\--- first, focusing on the auxiliary $\nu$-plane
\\
\--- second, reverting to the original $\zeta$-plane.

%%%%%%%%%%%%%%%%%%%%%%%%%%%%%%%%%%%%%%%%%%%%%%%%%%%%%%%%%%%%%%%%%%%%%%%%%%%%%%%%%%%%%%%%%%%%
\subsection{Holomorphic input $f$. The five arrows.}
%%%%%%%%%%%%%%%%%%%%%%%%%%%%%%%%%%%%%%%%%%%%%%%%%%%%%%%%%%%%%%%%%%%%%%%%%%%%%%%%%%%%%%%%%%%%
\subsubsection{From {\it original} to {\it outer}.}
Let us check, in the four simplest instances, that SP series 
(our so-called {\it original} generators)
with  an holomorphic input $f$ always give rise to two {\it outer} generators\footnote{
which exceptionally coalesce into one when $\eta_F=0,\omega_f=1$, which
may occur only in the cases 3 or 4 {\it infra}.}
$$ \{\imi{lo}_{\mr{in}}\!(\nu),
\imi{Lo}_{\mr{in}}\!(\zeta)=\,\imi{lo}_{\mr{in}}\!(1+\zeta) \} 
\quad,\quad
 \{\imi{lo}_{\mr{out}}\!(\nu),
\imi{Lo}_{\mr{out}}\!(\zeta)=\,\imi{lo}_{\mr{out}}\!(1+\omega_F\zeta) \} $$
located respectively over 
$\{\nu=0,\zeta=1 \} $ or $\{\nu=\eta_F,\zeta=1/\omega_F \}$ 
and produced under the {\it nur}-tranform, i.e. by inputting respectively
$f$ or $f^\vdash$  into the
long chain of \S5.2

\noindent
{\bf Case 1:}\;\;{\it $-f^\ast$ decreases on $[0,1]$ .}\\
 To explain the occurence $\mi{Lo}_{\mr{in}}$, apply the argument
at the beginning of \S5.1. To explain  the occurence $\mi{Lo}_{\mr{out}}$,
the shortest way is to pick $\epsilon >0$ small enough for
$-f^\ast$ to be decreasing on the whole of $[0,1+\epsilon]$,
and then to form the SP series $j\!j^\#_F(\zeta)$ defined
exactly as $j^\#_F(\zeta)$ but with a summation ranging over
$ 0\leq m<(1+\epsilon)\,n$ instead of $0\leq m <n$. Then $j\!j^\#_F$
clearly has no singularity at $\zeta=1/\omega_F$. On the other hand,
applying once again the argument of \S.. to the difference
 $j\!j^\#_F-j^\#_F$ we see that it has at $\zeta=1/\omega_F$
a singularity which, up to the dilation factor $\omega_F$, is given
by the {\it nur}-transform of $^1\!f$ with $^1\!f(x):=f(1+x)$.
In view of the parity relation of \S5.8 it is also equal to {\it minus}
the {\it nur}-transform of $(^1\!f)^\perp$. But 
$ (^1\!f)^\perp=f^\vdash$. Hence the result\footnote{Recall that $f^{\perp}(x):=-f(-x)$ and $f^{\vdash}(x):=-f(1-x)$    }.

A trivial \-- but telling \-- example corresponds to the choice
of a constant input $F(x)\equiv \alpha$ with $ 0<\alpha <1$. We then get\,:
\begin{equation}\label{last1}
j_{F}(\zeta)=\frac{\alpha}{1-\alpha}\,\big(\frac{1}{1-\zeta}-\frac{1}{1-\alpha\zeta}\big)\quad ; \quad
j^{\#}_{F}(\zeta)=\frac{1}{\alpha^{-1/2}-\alpha^{1/2}}\,\big(\frac{1}{1-\zeta}-\frac{1}{1-\alpha\zeta}\big)\quad  \quad
\end{equation}
\\
\medskip
\noindent
{\bf Case 2:}\;\;{\it $-f^\ast$ increases on $[0,1]$ .}\\
The $\vdash$ transform turns case 2 into case 1, with $f$ anf $f^\vdash$
exchanged. Hence the result. 
Again, we have the trivial example 
of a constant input $F(x)\equiv \alpha$ but now with $ 1<\alpha $. We then get the same power series as in (\ref{last1}) but
with $\alpha$ changed into $1/\alpha$, which of course agrees with the relation (\ref{a175}) 
between $j^{\#}_{F}$ and $j^{\#}_{F^{\vdash}}$. 

\medskip
\noindent
{\bf Case 3:}\;\;{\it $-f^\ast$ decreases on $[0,x_0]$, then increases on $[x_0,1]$.}\\
Here again, the argument at the beginning of \S5.1 takes care of $\mi{Lo}_{\mr{in}}$.
To justify  $\mi{Lo}_{\mr{out}}$, all we have to do is observe that the $\vdash$
transform turns case 3 into another instance of that same case 3, 
while exchanging the
roles of $\mi{Lo}_{\mr{in}}$ and $\mi{Lo}_{\mr{out}}$.
 
\medskip
\noindent
{\bf Case 4:}\;\;{\it $-f^\ast$ increases on $[0,x_0]$, then decreases on $[x_0,1]$.}\\
 Case 4 is exactly the reverse of case 3. The argument about $j^\#_F$ and
$j\!j^\#_F$ (see case 1) takes care of $\mi{Lo}_{\mr{out}}$ and then the
fact that $\vdash$ turns case 4 into another case 4 justifies the occurence of
$\mi{Lo}_{\mr{in}}$. Case 4, however, presents us with a novel difficulty\,:
the presence for $-f^\ast$ of a maximum at $x=x_0$ gives rise
(see \S7.1.2 infra) to an inner generator $Li$ located at a
point $\omega^\prime_F=e^{\eta^\prime_F}$ 
(with  $\eta^\prime_F=\int_0^{x_0}f(x)dx$)
that is closer to the origin than both 1 (location of $Lo_{\mr{in}}$)
and $\omega_F$ (location of  $Lo_{\mr{out}}$). So the
method of \S5.1 for translating coefficient asymptotics into {\it nearest}
singularity description no longer applies. One must then resort to a suitable
{\it deformation} argument. We won't go into the details, but just mention a
simplifying circumstance\,: from the fact that {\it inner generators}
never produce {\it outer generators} (under alien differentiation),
it follows that the actual manner of pushing $Li$  beyond 
$Lo_{\mr{in}}$ and $Lo_{\mr{out}}$ (i.e. under right or left circumvention)
doesn't matter.

%%%%%%%%%%%%%%%%%%%%%%%%%%%%%%%%%%%%%%%%%%%%%%%%%%%%%%%%%%%%%%%%%%%%%%%%%%%%%%%%%%%%%%%%%%%%
\subsubsection{From {\it original} to {\it inner}.}
\noindent
{\bf Case 4:}\;\;{\it $-f^\ast$ decreases on $[0,x_0]$, then increases on
$[x_0,1]$.}\\ 
When $f$ has a simple zero at $x_0$, i.e.  when the ``tangency order''is $\kappa=1$,
we are back to the heuristics of \S4.1. When $f$ has a multiple zero (necessarily
of odd order, if $f^\ast$ is to have an extremum there), we have a tangency order
$\kappa\in\{3,5,7\dots\}$ and the same argument as in \S4.1 points to the existence of a
singularity over $\eta_F^\prime $ in the $\nu$-plane 
or  $\omega_F^\prime $ in the $\zeta$-plane, with
$\eta_F^\prime,\omega_F^\prime $ as above. In the $\nu$-plane, this
singularity is characterised by an upper-minor $\smi{li}$ given by\,:
\begin{equation}  \label{a176}
\smi{li}:=\mr{nir}(^0f)+\mr{nir}(^0f^\perp)\quad\mi{with}\quad
^0f(x):=f(x_0+x)
\end{equation}
In view of the parity relation (cf \S4.10) this implies\,:
\begin{equation}  \label{a177}
\smi{li}(\nu)=\sum_{0\leq
k}\,h_{\frac{-\kappa+2k}{\kappa+1}}\,\nu^{\frac{-\kappa+2k}{\kappa+1}}
\quad\;\mi{with}\quad
h=\mr{nir}(f)=\sum_{0\leq
k}\,h_{\frac{-\kappa+k}{\kappa+1}}\,\nu^{\frac{-\kappa+k}{\kappa+1}}
\quad
\end{equation}
Thus, only every second coefficient of $h$ goes into the making of $ \smi{li}$.
Moreover, since $\kappa$ here is necessarily odd, the ratio 
$\frac{-\kappa+2k}{\kappa+1} $ can never be an integer. This means that the
corresponding {\it majors}\footnote{whether $\sma{li},\ima{li}$ or 
$\ima{Li}$} never carry any logarithms, but only fractional powers.
%%%%%%%%%%%%%%%%%%%%%%%%%%%%%%%%%%%%%%%%%%%%%%%%%%%%%%%%%%%%%%%%%%%%%%%%%%%%%%%%%%%%%%%%%%%%
\subsubsection{From {\it outer} to {\it inner}.}
The relevant functional transform here is {\it nur}, which according to (\ref{a125})
is an infinite superposition of {\it nir} transforms applied separately to
all determinations of $\log f(.) $. To calculate the alien derivatives of 
$\smi{lo}_{\mr{in}}$ or $\smi{lo}_{\mr{out}}$, we must therefore 
apply
the recipe of the next para (\S7.1.4) to the 
various 
$\mi{nir}(2\pi i\,k+\log f(0)+\dots)$  
or $\mi{nir}(2\pi i\,k+\log(f(1)+\dots)$. Exceptionnally,
if $2\pi i\,k+\log f(0)$ or  $2\pi i\,k+\log f(1)$ 
vanishes for some $k$, we must also deal with tangency orders $\kappa >0$
and apply the recipe of the para after next (\S7.1.5). But in this as in that case,
the result will always be {\it some} inner generator $\smi{li}$, and never an outer one.
%%%%%%%%%%%%%%%%%%%%%%%%%%%%%%%%%%%%%%%%%%%%%%%%%%%%%%%%%%%%%%%%%%%%%%%%%%%%%%%%%%%%%%%%%%%%
\subsubsection{From {\it exceptional} to {\it inner}.}
Let $\smi{le}$ an exceptional generator with base point $x_1$. Assume, in other words,
that 
$f(x_1)\not=0$ and\,:
\begin{equation}  \label{a178}
\smi{le}\,=\,^{\nu_1}\!h=\mr{nir}\,({ ^{x_1}\!f})\;\;\;\;\mi{with}\;\;\;\;
^{x_1}\!f(x):=f(x_1+x)
\;\;,\;\; \nu_1:=\int_0^{x_1}f(x)dx
\end{equation}
To calculate the alien derivatives of $\smi{le}$,
we go back to the long chain \S4.2 and decompose the {\it nir}-transform
into  elementary steps from 1 to 7. The elementary steps 1,2,4,5,7
neither produce nor destroy singularities. The steps that matter are the
reciprocation (step 3) and the {\it mir}-transform (step 6).
The singularities produced by reciprocation are easy to predict.
As for the {\it mir}-transform, its integro-differential expression (\ref{a85})
and the properties of the Euler-Bernoulli numbers\footnote{ more
exactly, the fact that the singularities of $\beta$ are all on $2\pi i\doZ^\ast$.}
 show that the closest singularity or singularities of $\smi{le}$
\footnote{i.e. those lying on the boundary of the disk of convergence.
Recall that for an exceptional generator we have a tangency order $\kappa=0$ and so 
$\smi{le}$ is a regular, unramified germ at the origin.} necessarily correspond to
closest singularity/\!/ies of $\gbar$ (see Lemma 4.7). Now comes the crucial,
non-trivial  fact\,: this one-to-one correspondance between singularities of $\gbar$ and
$\smi{le}$ holds also in the large, at least when the initial input $f$
is holomorphic. This is by no means obvious, since the singularities of
$\gbar$ might combine with those of $\beta$ to produce infinitely
many new ones, farther away, under the Hadamard product mechanism . 
To show that this
{\it doesn't occur}, assume the existence of a point $\nu_2$ in the $\nu$-plane where
$^{\nu_1}h=\smi{le}$ is singular but $g$ is regular. We can then write
$\nu_2=\int_0^{x_2} f(x)dx$ for some $x_2$ and then choose $x_3$ close enough to
$x_2$ to ensure that the exceptional generator $ ^{\nu_3}h$ of base point $x_3$ is
regular at $\nu_2$
\footnote{ by ensuring that  $ ^{\nu_3}h $ has 
$\nu_2$ within its convergence disk.}.
  We then use the bi-entireness of the finite {\it nir}-increment
$\nabla h(\epsilon,\nu)$ with $\epsilon=x_3-x_1\, ,\, \nu=\nu_3-\nu_1$
to conclude that $^{\nu_1}h:=\mi{nir}(^{x_1}h)$, just like 
 $^{\nu_3}h:=\mi{nir}(^{x_3}h)$, is regular at $\nu_2$.
%%%%%%%%%%%%%%%%%%%%%%%%%%%%%%%%%%%%%%%%%%%%%%%%%%%%%%%%%%%%%%%%%%%%%%%%%%%%%%%%%%%%%%%%%%%%
\subsubsection{From {\it inner} to {\it inner}.  Ping-pong resurgence.}
Let $\smi{li}_1$ be an inner generator with base point $x_1$. This means that
$f(x_1)=0$ and\,: 
\begin{equation}  \label{a179}
\smi{li}_1=\,^{\nu_1}\!h=\mr{nir}\,({ ^{x_1}\!f})\;\;\;\;\mi{with}\;\;\;\;
^{x_1}\!f(x):=f(x_1\!+\!x)
\;\;,\;\; \nu_1:=\int_0^{x_1}f(x)dx
\end{equation}
Assuming once again $f$ to be holomorphic, the same argument as above shows that all
singularities of  $\smi{li}_1$, not just the closest ones, correspond to zeros $x_i$
of $f$. They are therefore inner generators 
$\smi{li}_2, \smi{li}_3, \smi{li}_4 \dots $
with base points $x_2,x_3,x_4\dots$ and the resurgence equations
between them\,
\begin{equation}  \label{a180}
\Delta_{\nu_q-\nu_p}\,\smi{li}_p \;= \;\smi{li}_q
\end{equation}
will exactly mirror the resurgence equations between the singularities of $g$.
The only difference is that if $\smi{li}_p$  {\it``sees"} $\smi{li}_q$,
i.e. if (\ref{a180}) holds, then the converse is automatically true\,:
 $\smi{li}_q$  {\it sees}  $\smi{li}_p$. Exceptional generators, on
the other hand, {\it see } but are not {\it seen}. \footnote{
Regarding the inner generators, one may note that what matters is the
 geometry in the $\nu$-plane, not in the $x$-plane. Consider for
instance\,:
$$
f(x):=(x-x_1)(x-x_2)(x-x_3)\;\;\;\mi{with}\;\;\; x_1=1,x_2=2,x_3=2+\epsilon+
\epsilon^2\,i\;\;(0<\epsilon<<1)
$$
Then a simple calculation shows that the inner generator  $\smi{li}_1$ {\it sees}
$\smi{li}_2$ but not $\smi{li}_3$, although $x_1$ {\it sees}  $x_2$ and
$x_3$. (On the other hand, $\smi{li}_2$ {\it sees} both $\smi{li}_1$ and
 $\smi{li}_3$.)
}
%%%%%%%%%%%%%%%%%%%%%%%%%%%%%%%%%%%%%%%%%%%%%%%%%%%%%%%%%%%%%%%%%%%%%%%%%%%%%%%%%%%%%%%%%%%%
\subsubsection{Recapitulation. One-way arrows, two-way arrows.}
Let us sum up pictorially our findings for a holomorphic input $f$ \,:
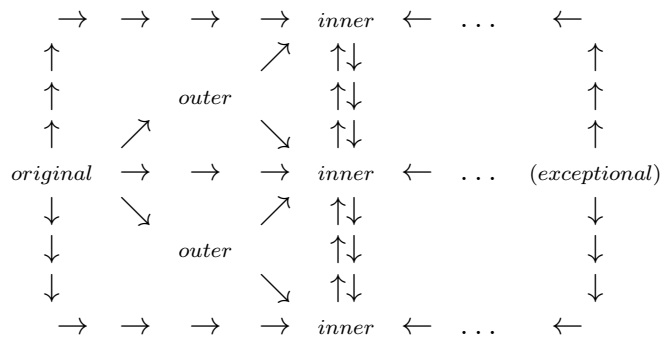
\[\begin{array}{cccccccccc} 
\hspace{3.ex}\rightarrow&\rightarrow&\rightarrow&\rightarrow
&\scriptstyle{inner} &\leftarrow&\dots&\leftarrow\hspace{4.ex}
\\
\uparrow&&&\nearrow&\uparrow\downarrow&&&\uparrow
\\
\uparrow&&\scriptstyle{outer}&&\uparrow\downarrow&&&\uparrow
\\
\uparrow&\nearrow& &\searrow&\uparrow\downarrow&&&\uparrow
\\
\scriptstyle{original}&\rightarrow&\rightarrow
&\rightarrow&\scriptstyle{inner}&\leftarrow&\dots&\scriptstyle{(exceptional)}
\\
\downarrow &\searrow&&\nearrow&\uparrow\downarrow&&&\downarrow
\\
\downarrow &&\scriptstyle{outer}&&\uparrow\downarrow&&&\downarrow
\\
\downarrow &&&\searrow&\uparrow\downarrow&&&\downarrow
\\

\hspace{3.ex}\rightarrow&\rightarrow&\rightarrow
&\rightarrow&\scriptstyle{inner}&\leftarrow&\dots&\leftarrow\hspace{4.ex}
\end{array}\]
The above picture displays four types of generators\,:
\\
\-- one {\it original} generator, which is none other than the `cleansed' SP series 
\\
\-- two {\it outer} generators ({\it in} and {\it out}) which may occasionally coalesce
\\
\-- a countable number of {\it inner} generators\,:  as many as $f$ has zeros
\\
\-- a continous infinity of {\it exceptional} generators\,: any $x_i$ where $f$ doesn't 
vanish can serve as base point.

The picture also shows five types of arrows linking these generators\footnote{
meaning in each case that the {\it target} is generated by the {\it source}
under alien differentiation.}. All these arrows are one-way, except for those
linking pairs of inner generators.

As this ``one-way/two-way traffic" suggests, the various generators differ widely
as to origin, shape, and function.

The original generator clearly stands apart, not just because it kicks off
the whole generation process, but also because it makes (immediate) sense
only in the $\zeta$-plane\,: in the $\nu$-plane it is relegated to infinity.

Directly proceeding from it under the {\it nur}-transform, we have
 two outer generators, which in turn generate the potentially more numerous
inner generators, this time under the {\it nir}-transform, relatively
in each case to a given determination of $f=-\log F$. To each such determination
(corresponding to an additive term\, $2\pi i k$) there answers a
distinct inner algebra $\mi{Inner}_f$ spanned by $K$ inner generators,
with $K:=\mr{card}\{f^{-1}(0)\}$.

Another way of entering the inner algebras is via exceptional generators, but
these are ``artificial" in the sense that they never occur naturally, i.e.
under analytic continuation of the original generator. They are more in the nature of
auxiliary tools\footnote{as components of the {\it nur}-transform
under the Poisson formula (see \ref{a125}) and also, as we just saw,
as mobile tools for sifting out true singularities from illusory ones (see \S7.1.4).
}. Also, since each exceptional generator results from applying the {\it nir}-transform
to a {\it given} determination of $f=\log F$, it gives acces to {\it one}
inner algebra $\mi{Inner}_f$, unlike the outer generators, which give access to 
{\it them all}. 
  
These inner algebras  $\mi{Inner}_f$ are in one-to-one correspondance with $\doZ$.
Though  distinct (and usually  disjoint) from each other, they are essentially
isomorphic. Each of them is also ``of one piece" in the sense that for any
pair $\smi{li}_{p},\smi{li}_{q}\in \mi{Inner}_f$, there is always
a connecting chain $\mi{li}_{n_i}$ starting at $\smi{li}_{p}$, ending at
$\smi{li}_{q}$, and such that any two neighbours $\mi{li}_{n_i}$
and $\mi{li}_{n_{i+1}}$ {\it see} each other.
\\

The emphasis so far has been on the singularities in the $\nu$-plane. Those
in the $\zeta$-plane follow, except  {\it over} the origin
$\zeta=0$, where quite specific and severe singularities
may also occur ({\it at} the origin itself, i.e. on the 
main Riemann leaf, the SP function is of course regular).
For a brief discussion of these 
$0$-based singularities and their resurgence properties, see \S7.2.1 below. 
%%%%%%%%%%%%%%%%%%%%%%%%%%%%%%%%%%%%%%%%%%%%%%%%%%%%%%%%%%%%%%%%%%%%%%%%%%%%%%%%%%%%%%%%%%%%
%%%%%%%%%%%%%%%%%%%%%%%%%%%%%%%%%%%%%%%%%%%%%%%%%%%%%%%%%%%%%%%%%%%%%%%%%%%%%%%%%%%%%%%%%%%%

\subsection{Meromorphic input $F$\,: the general picture.}
Let us briefly review the main changes which take place when we relax the hypothesis
about $f:=-\log(F) $ being holomorphic and simply demand that $F$ be
{\it meromorphic}.\footnote{ Since $F$ and $F^\vdash$ (recall that
$F^\vdash(x):=1/F(1-x)$) are essentially
on the same footing, it would make little sense to assume one to be {\it holomorphic} rather than
the other. So we must assume {\it meromorphy}, even {\it strict} meromorphy, with at least
one zero or pole. 
} 
%%%%%%%%%%%%%%%%%%%%%%%%%%%%%%%%%%%%%%%%%%%%%%%%%%%%%%%%%%%%%%%%%%%%%%%%%%%%%%%%%%%%%%%%%%%%
%%%%%%%%%%%%%%%%%%%%%%%%%%%%%%%%%%%%%%%%%%%%%%%%%%%%%%%%%%%%%%%%%%%%%%%%%%%%%%%%%%%%%%%%%%%%
\subsubsection{Logarithmic/non-logarithmic singularities.}
If $F$ has at $x=0$ a zero or pole of order $d\in \doZ^\ast$, we must
replace the $\sum_{0\leq m <n}$ summation in (\ref{a2}) by  $\sum_{0 < m <n}$ 
for the definition of the SP coefficients $J_F(n)$ to make sense. More significantly,
depending of the parity of $d$, the outer and inner singularities may
exchange their logarithmic/non-logarithmic nature. Recall that for the {\it cleansed}
SP function and $d=0$, the outer generators have purely logarithmic
singularities\footnote{i.e. with majors of type $\mi{Reg}_1(\zeta)+\mi{Reg}_2(\zeta)\log(\zeta) $.}
while the inner generators have power-type singularities, with strictly
rational (non-entire) powers. That doesn't
change when $d$ is $\not=0$, at least where the {\it cleansed} SP series
are concerned. However, when
we revert to the {\it raw} SP series, i.e. to the position prior to coefficient division
by the ingress factor
$\mi{Ig}_F(n) \sim n^{-d/2}\;(c_0+O(n^{-1})) $,
we are faced with a neat dichotomy\,:
\\
(i) \;$d$ even\,: nothing changes.
\\
(ii) $d$  odd\,: everything gets reversed, with the outer singularities
becoming strict rational (semi-integral) powers and the inner singularities becoming
purely logarithmic\footnote{ At least in the generic case, i.e. for a tangency
order
$\kappa=1$. For $\kappa>1$, the inner singularities involve a mixture of 
rational powers and logarithms.}

%%%%%%%%%%%%%%%%%%%%%%%%%%%%%%%%%%%%%%%%%%%%%%%%%%%%%%%%%%%%%%%%%%%%%%%%%%%%%%%%%%%%%%%%%%%%
%%%%%%%%%%%%%%%%%%%%%%%%%%%%%%%%%%%%%%%%%%%%%%%%%%%%%%%%%%%%%%%%%%%%%%%%%%%%%%%%%%%%%%%%%%%%
\subsubsection{Welding the inner algebras into one. }
The presence of even a single zero or pole in $F$, no matter where
 \--- whether at x=0 or x=1 or elsewhere \---
suffices to abolish the distinction between the various 
inner algebras $\mi{Inner}_f$ attached to the various 
determinations of $f:=-\log(F)$, since $f$ itself now becomes
multivalued and assumes the form\,:
\begin{equation}  \label{a181}
f(x)=\sum d_i\,\log(x-x_i)+\mi{holomorphic}(x)\;\; 
\end{equation}
Everything hinges on $ d:=g.c.d.(d_1,d_2,\dots)$. If $d=1$, then all inner algebras
merge into one. If $d>1$, they merge into $d$ distinct but ``isomorphic" copies. 

Notice that no such change affects the outer generators, because
these are constructed, not from $f$, but directly from $F$ (in the case
of the {\it raw} SP function ) or $F^{1/2}$ (in the case of the {\it cleansed}
SP function).
%%%%%%%%%%%%%%%%%%%%%%%%%%%%%%%%%%%%%%%%%%%%%%%%%%%%%%%%%%%%%%%%%%%%%%%%%%%%%%%%%%%%%%%%%%%%
%%%%%%%%%%%%%%%%%%%%%%%%%%%%%%%%%%%%%%%%%%%%%%%%%%%%%%%%%%%%%%%%%%%%%%%%%%%%%%%%%%%%%%%%%%%%
\subsection{The $\zeta$-plane and its violent $0$-based singularities. }
\noindent
{\bf Converting $\nu$-singularities into $\zeta$-singularities.}
\\
So far, we have been describing the outer/inner singularities in the
auxiliary  $\nu$-plane (more exactly, the $\nu$-Riemann surface) which
is naturally adapted to  Taylor coefficient asymptotics. To revert to
the original $\zeta$-plane, we merely apply the formulas for step 9 in the long chain of \S4.2 which
convert $\nu$-singularities into $\zeta$-singularities, for majors as well minors.
The resurgence equations, too, carry over almost unchanged, with the
additive indices $\nu_i$ simply turning into  multiplicative indices $\zeta_i$.
But there is one exception, namely the origin $\zeta=0$. Under the correspondence
$\nu \mapsto \zeta=e^\nu$, the SP function's behaviour {\it over} will
reflect its behaviour {\it over} the ``point" $\Re(\nu)\!=\!-\infty$
on various Riemann leaves. This is the tricky matter we must now look into.
\\

\noindent
{\bf Description/expansion of the $0$-based singularities.}\\
The SP function itself is regular {\it at} $\zeta=0$, i.e.
on the main Riemann leaf,\footnote{or, if $F$ has a zero/pole of odd order $d$,
it is of the form $\zeta^{d/2}\varphi(\zeta)$, but again with a regular $\varphi$.
} but usually not {\it over} $\zeta=0$, i.e. on the other leaves. Studying these
$0$-singularities entirely reduces to studying the $0$-singularities of the outer/inner
generators
$\imi{Lo}/\imi{Li}$, which in turn reduces to investigating the $\infty$-behaviour of
$\imi{lo}/\imi{li}$. This can be done in the standard manner, by going to the long chain
of \S4.2 and applying the {\it mir}-transform to $\gbar$, {\it but locally at
$-\infty$}\footnote{ each time on the suitable leaf, of course.}. The integro-differential expansion
for {\it mir} still converges in this case, but no longer formally so 
(i.e. coefficient-wise),
and it still yields inner generators, but of a very special, quite irregular sort.
Pulled back into the $\zeta$-plane, they produce violent singularities
{\it over} $\zeta= 0$, usually with exponentially 
explosive/implosive radial behaviour,
depending on the sectorial neighbourhood of 0.\\
 
 \noindent
{\bf Resurgence properties of the $0$-based singularities.}\\
 Fortunately, no detailed {\it local} description of the 0-based singularities
is required to calculate their alien derivatives and, therefore, to obtain
a complete system of resurgence equations for our original SP function.
Indeed, turning $k$ times around $\zeta=0$ on {\it some} leaf amounts to making a
$2\pi i k$-shift in the $\nu$-plane, again on {\it some} leaf. But the effect
of that is easy to figure out, especially for an holomorphic input $f$ (in that case, it simply takes us from one
inner algebra $\mi{Inner}_{f}$ to the next)
but also for a general meromorphic input $F$ (for illustrations, see the examples of \S8.3, especially examples 8.7 and 8.8. See also \S6.6-8.)

%%%%%%%%%%%%%%%%%%%%%%%%%%%%%%%%%%%%%%%%%%%%%%%%%%%%%%%%%%%%%%%%%%%%%%%%%%%
%%%%%%%%%%%%%%%%%%%%%%%%%%%%%%%%%%%%%%%%%%%%%%%%%%%%%%%%%%%%%%%%%%%%%%%%%%%
\subsection{Rational inputs $F$\,: the inner algebra.}
Let $F$ be a rational function of degree $ d $ \,:
\begin{equation}\label{s1}
 F(x)=\prod_{1\leq j \leq r} \big(1-\frac{x}{\alpha_{j}} \big)^{d_{j}}\quad \mi{with}  \quad d_{j} \in \doZ^{\ast}
 \quad, \quad \delta:= \mr{g.c.d.}(|d_{1}|,\dots, |d_{r}|)  
 \end{equation}
and let $x_{0},\dots, x_{d-1} $ be the zeros (counted with multiplicities) of the equation $F(x)=1$.
We then fix a determination of the the corresponding $f$\,:
\begin{equation}\label{s2}
f(x) = -\log(F) = -\sum_{1\leq j\leq r}d_{j}\log\big(1-\frac{x}{\alpha_{j}} \big)
\end{equation}
with its Riemann surface $S_{f}$ . We denote $X_{j}^{n}\subset S_{f}$ the set of all $x^{\star}\in S_{f}$ lying
over $x_{j}\in \doC$ and such that $f(x^{\star})=2\pi i\delta$ and select some point $x_{0}\in X_{0}^{0}$ as
base point of $S_{f}$. The internal generators will  then correspond one-to-one to the points of $ \cup_{1\leq j< d} X^{0}_{j}$
and be located at points $\nu_{j}$ of the ramified $\nu$-plane, with projections $\dot{\nu}_{j}$ such that\,:
\begin{equation}\label{s3}
\dot{\nu}_{j}-\dot{\nu}_{i}=\int_{x_{i}}^{x_{j}}f(x)dx\quad\quad \quad(x_{i}\in X_{i}^{0}\;,\;x_{j}\in X_{j}^{0}\;)
\end{equation}
For two distinct points $x_{j}^{\prime},x_{j}^{\prime\prime}$ in the same $X_{j}^{0}$, the above integral is obviously
a multiple of $2\pi i\delta$. Therefore, three cases have to be distinguished.
\\
{\bf Case 1}. $F$ has only one single zero $\alpha_{1}$ (of any multiplicity $d_{1}=p$)
or again one single pole $\alpha_{1}$ (of any multiplicity $d_{1}=p$). In that case, we have exactly $p$ sets $X_{j}^{0}$,
but each one reduces to a single point, since $f$ has only one single logaritmic singularity. That case 
(``monomial input F'') was investigated in detail in \S6.6, \S6.7,\S6.8.
\\
{\bf Case 2}: $F$ has one simple  zero $\alpha_{1}$ and one simple pole $\alpha_{2}$. 
The position is now the reverse\,: we then have only one set $X_{0}^{0}$,but with a countable infinity of
points in it, since $f$ has now two logarithmic singularities, thus allowing integrals (\ref{s3}) with distinct end points
$x_{0}^{\prime} , x_{0}^{\prime\prime}$ both in $X_{0}^{0}$. (See \S8.11 below).
\\
{\bf Case 3}: $F$ has either more than two distinct zeros, or more than two distinct poles, or both. We then have
$p+q$ distinct sets $X_{j}^{0}$ , with $p$ (resp. $q$) the number of distinct zeros (resp. poles). Each such $X_{j}^{0}$
contains a countable number of points $x_{j}$, to which there answer, in the ramified $\nu$-plane,
distinct singular points $\nu_{j}$ that generate a set $\mathcal{N}_{j}$ whose projection 
$\dot{\mathcal{N}}_{j}$ on $\doC$ is of the form $\nu_{j}+\Omega$, with
\begin{eqnarray}\label{u71}
\Omega &=& \Big\{ \omega\;,\; \omega= \!\!\sum_{\sum n_{j}d_{j}=0} n_{j}\,d_{j} \,\eta_{j} = -2\pi i \!\!\!
\sum_{\sum n_{j}d_{j}=0} n_{j}\,d_{j} \,\alpha_{j} \Big\}
\quad\quad 
\\[1.ex]
 \eta_{j}&=&-\int_{\mathcal{I}_{j}}\log(1-\frac{x}{\alpha_{j}})\,dx\,=\,2\pi i\, (x_{0}-\alpha_{j})
\end{eqnarray}
and with integration loops $\mathcal{I}_{j}$ so chosen as to generate the fundamental homotopy group of 
$ \doC \setminus \{\alpha_{1},\dots,\alpha_{r}\}$. Each  $\mathcal{I}_{j}$ describes a positive turn round $\alpha_{j}$
and the choice of the loops' common end-point is immaterial, since changing the end-point merely adds adds a common
constant  to each $\eta_{j}$, which constant cancels out from the sums $\omega$ due to the condition $\sum n_{j}\,d_{j}=0$.

Though $\Omega$ usually fails to be discrete as soon as $  \geq 3$, the sets $\mathcal{N}_{j} $ are of course
always discrete in the ramified $\nu$-plane.
In particular, from any given singular point $\nu_{i}\in \mathcal{N}_{i} $ only finitely
many $\nu_{j} \in \mathcal{N}_{j} $ can be {\it seen} \--- those namely that correspond to {\it `simple'} integration paths in (\ref{s3}),
i.e. typically paths whose projection on $\doC$ is {\it short} and doesn't {\it self-intersect}.\footnote{ Various examples of such situations
shall be given in \S8.3, with {\it simple/complicated} integration paths corresponding to {\it visible/invisible} singularities. The
general situation, with the exact criteria for visibility/invisibility, shall be investigated in [S.S.1]} The other points 
$\nu_{j} \in \mathcal{N}_{j} $ are located on more removed Riemann leaves and therefore {\it hidden from view} (from $\nu_{i} $).

\section{The inner resurgence algebra for SP series.}
%%%%%%%%%%%%%%%%%%%%%%%%%%%%%%%%%%%%%%%%%%%%%%%%%%%%%%%%%%%%%%%%%%%%%%%%%%%%%%%%%%%%%%%%%%%%
%%%%%%%%%%%%%%%%%%%%%%%%%%%%%%%%%%%%%%%%%%%%%%%%%%%%%%%%%%%%%%%%%%%%%%%%%%%%%%%%%%%%%%%%%%%%
%%%%%%%%%%%%%%%%%%%%%%%%%%%%%%%%%%%%%%%%%%%%%%%%%%%%%%%%%%%%%%%%%%%%%%%%%%%%%%%%%%%%%%%%%%%%
%%%%%%%%%%%%%%%%%%%%%%%%%%%%%%%%%%%%%%%%%%%%%%%%%%%%%%%%%%%%%%%%%%%%%%%%%%%%%%%%%%%%%%%%%%%%
\subsection{Polynomial inputs $f$. Examples.} 
\begin{example}:\,$f(x)=x^r$
\textup{
There is only one inner generator $\hat{h}(\nu)$ which, up to to the factor $\nu^{-1/2}$, is an entire
function of $\nu$.
}
\end{example}
\begin{example}:\,$f(x)=x^r-1$
\textup{
There are $r$ inner generators. We have exact radial symmetry, of radius $1$ in the $x$-plane
and radius $\eta=1/(r+1)$ in the $\nu$-plane. Every singular point there sees all the others\,:
we have `multiple ping-pong', governed by a very simple resurgence system (see [SS1]).
}
\end{example}
\begin{example}:\,$f(x)=\prod_{j=1}^{j=r}(x-x_{j})$
\textup{
Every such configuration, including the case of multiple roots, can be realised
by continuous deformations of the radial-symmetric configuration of Example 8.2, and the thing
is to keep track of the $\nu_{j}$-pattern as the $x_{j}$-pattern changes. When $\nu_{j}-\nu_{i}$ becomes
small while $x_{j}-x_{i}$ remains large, that usually reflects mutual invisibility of $\nu_{i}$ and $\nu_{j}$.
Thus, if $r=3$ and $x_{1}=0$, $x_{2}=1$, $x_{3}=1+\epsilon\,e^{i\theta}$ with $0\!<\!\epsilon\!<\!\!<\!1$,
the case $\theta=\pi/2$ with its approximate symmetry $x_{2}\leftrightarrow x_{3}$ corresponds to three mutually
visible singularities $\nu_{1},\nu_{2},\nu_{3}$, but when $\theta$ decreases to $0$, causing $x_{3}$
to make a $-\pi/2$ rotation around $x_{2}$, the point $\nu_{3}$ makes a $-3\pi/2$ rotation around $\nu_{2}$,
so that the projection $\dot{\nu}_{3}$ actually lands on the real interval $[\dot{\nu}_{1},\dot{\nu}_{2}]$.
But the new $\nu_{3}$ has actually moved to an adjacent Riemann leaf and is no longer {\it visible}
from $\nu_{1}$.
}
\end{example}
%%%%%%%%%%%%%%%%%%%%%%%%%%%%%%%%%%%%%%%%%%%%%%%%%%%%%%%%%%%%%%%%%%%%%%%%%%%%%%%%%%%%%%%%%%%%
%%%%%%%%%%%%%%%%%%%%%%%%%%%%%%%%%%%%%%%%%%%%%%%%%%%%%%%%%%%%%%%%%%%%%%%%%%%%%%%%%%%%%%%%%%%%
\subsection{Holomorphic inputs $f$. Examples.} 
\begin{example}:\,$f(x)=\exp(x)$
\\
\textup{
To the unique `zero' $x_{0}=-\infty$ of $f(x)$ there answers a unique inner generator $\hat{h(\nu)}$.
It is of rather exceptional type, in as far as its local behaviour is described by a {\it transseries} rather
than a series, but the said transseries is still produced by the usual mechanism of the nine-link chain.
}
\end{example}
\begin{example}:\,$f(x)=\exp(x)-1$ or $f(x)=\sin^{2}(x)$
\\
\textup{
All zeros $x_{i}$ of $f(x)$ contribute distinct inner generators, identical up to shifts but positioned at different
locations $\nu_{j}$.
}
\end{example}
\begin{example}:\,$f(x)=\sin(x)$\\
\textup{ Here, the periodic $f(x)$ still has infinitely many zeros but is constant-free (i.e. is itself the derivative
of a periodic function). As a consequence, we have just two inner generators, at two distinct locations,
like in the case $f(x)=1-x^{2}$ but of course with a more complex resurgence pattern.
}
\end{example}
%%%%%%%%%%%%%%%%%%%%%%%%%%%%%%%%%%%%%%%%%%%%%%%%%%%%%%%%%%%%%%%%%%%%%%%%%%%%%%%%%%%%%%%%%%%%
%%%%%%%%%%%%%%%%%%%%%%%%%%%%%%%%%%%%%%%%%%%%%%%%%%%%%%%%%%%%%%%%%%%%%%%%%%%%%%%%%%%%%%%%%%%%
\subsection{Rational inputs $F$. Examples.} 
\begin{example}\,: $F(x)=(1-x)$
\\ \textup{
The inner algebra here reduces to one generator $h(\nu)$ and a fairly trivial one at that, since
$ h(\nu)=\mi{const}\,\nu^{-1/2}$, as given by the {\it semi-entire} part of the {\it nir}-transform. In contrast, the {\it entire}
part of the {\it nir}-transform (which lacks intrinsic significance) is, even in this simplest of cases,
a highly transcendental function\,: in particular, it verifies no linear ODE with polynomial coefficients.
}
\end{example}
\begin{example}:\,$F(x)=(1-x)^p$\\
 \textup{ 
Under the change $x\rightarrow p\,x$, this reduces the case of ``monomial F'', which
was extensively investigated in \S6.6, \S6.7, \S6.8. We have now exactly p internal generators $h_{j}(\nu)$ located at the unit roots
$\nu_{j}=-e^{2\pi i\,j/p}$ and verifying a simple ODE of order p, with polynomial coefficients. Each singular
point $\nu_{j}$ `sees' all the others, and the resurgence regimen is completely encapsulated in the matrices $\mathcal{M}_{p,q}$
of \S6.7, which account for the basic closure phenomenon: a $4\pi i$-rotation (around any base point) leaves the whole picture
unchanged.}
\end{example}
\begin{example}:\,$F(x)=\frac{x^{2}-\alpha^{2}}{1-\alpha^{2}}=\frac{x^{2}+\beta^{2}}{1+\beta^{2}}\;\;,\;\; \alpha=i\,\beta$
\\
\textup{
The general results of \S7.4 apply here, with $x_{0}=0,x_{1}=1$ and the lattice $\Omega= 4\pi i \alpha\,\doZ=-4\pi\beta\,\doZ $.
We have therefore two infinite series of internal generators in the $\nu$-plane, located over 
$\dot{\nu}_{0}+\Omega$ and $\dot{\nu}_{1}+\Omega$ respectively, where the difference $\dot{\nu}_{1}-\dot{\nu}_{0}$ may be
taken equal to any determination of
$ -\int_{0}^{1}\log(F(x))dx$. However, depending on the value of the parameters $\alpha,\beta$,
each singular point $\nu_{j}$ `sees' {\it one, two} or {\it three} singular points of the `opposite' series.
 Let us illustrate this on the three `real' cases\,:
\\
{\bf Case 1:} $0<\beta $.\\
The only singularity {\it seen} (resp. {\it half-seen}) from $\nu_{0}$ is $\nu_{1}$ (resp. $\nu_{1}^{\ast}$) with 
\begin{eqnarray*}
\dot{\nu}_{1}&=&\dot{\nu}_{0}+2\,\eta
\\
\dot{\nu}_{1}^{\ast}&=&\dot{\nu}_{0}+2\,\eta+4\pi\beta
\\
\mi{with}\quad\quad\eta &=& 2-2\,\beta \arctan({1}/{\beta})>0
\end{eqnarray*}
All other singularities above $\dot{\nu}_{1}+\Omega$ lie are on further Riemann leaves. The singularity $\nu_{1}$ corresponds
to the straight integration path $\mathcal{I}_{1}$  whereas $\nu_{1}^{\ast}$ corresponds to either of the equivalent
paths $\mathcal{I}_{1}^{\ast}$ and $\mathcal{I}_{1}^{\ast\ast}$.
\\
{\bf Case 2:} $0<\alpha <1$.\\
Only two singularities are {\it seen} from $\nu_{0}$, namely $\nu_{1}^{\ast}$ and $\nu_{1}^{\ast\ast}$ of projections :
\\
\begin{eqnarray*}
 \dot{\nu}_{1}^{\ast}\, &=&\dot{\nu}_{0}+2\,\eta+2\,\pi i\alpha 
\\
 \dot{\nu}_{1}^{\ast\ast} &=&\dot{\nu}_{0}+2\,\eta-2\,\pi i\alpha 
\\
\mi{with}\quad\quad \eta&=&2-\alpha\log\big(\frac{1+\alpha}{1-\alpha}\big)>0
\end{eqnarray*}
They correspond to the integration paths $\mathcal{I}^{\ast}$ and $\mathcal{I}^{\ast\ast}$ .
\\
{\bf Case 3:} $1<\alpha $.
\\
Three singularities are {\it seen} from $\nu_{0}$, namely $\nu_{1}, \nu_{1}^{\ast},\nu_{1}^{\ast\ast}$ of projections :
\begin{eqnarray*}
 \dot{\nu}_{1}\; &=&\dot{\nu}_{0}+2\,\eta
\\
 \dot{\nu}_{1}^{\ast}\, &=&\dot{\nu}_{0}+2\,\eta+2\,\pi i\alpha 
\\
 \dot{\nu}_{1}^{\ast\ast}\! &=&\dot{\nu}_{0}+2\,\eta-2\,\pi i\alpha 
\\
\mi{with}\quad\quad \eta&=&2-\alpha\log\big(\frac{\alpha+1}{\alpha-1}\big)< 0
\end{eqnarray*}
They correspond to the integration paths $\mathcal{I},\mathcal{I}^{\ast},\mathcal{I}^{\ast\ast}$ .
}
\end{example} 

%%%%%%%%%%%%%%%%%%%%%%%%%%%%%%%%%%%%%%
 \begin{tikzpicture}
\draw[dotted,black] (0,-2)--(15,-2);
\draw[dotted,black] (0,-3.5)--(15,-3.5);
\draw[dotted,black] (0,-7)--(15,-7);
\draw[dotted,black] (0,-11)--(15,-11);
%\draw[dotted,black] (0,-15)--(15,-15);
\draw[dotted,black] (0,-18.5)--(15,-18.5);
\draw[dotted,black] (5.5,-2)--(5.5,-23.5);
\draw[dotted,black] (10.5,-2)--(10.5,-23.5);
\draw[dotted,black] (0,-23.5)--(15,-23.5);

\begin{scope} [xshift=3cm,yshift=-5cm]

\node [rectangle,draw] {$ F(x)= \frac{x^2 -\alpha^2}{1- \alpha^2} ; \alpha > 1$};
\end{scope}

%%%%%%%%%%%%%%%%%  A1 %%%%%%%%%%%%%%%%%
\begin{scope}[xshift=4cm]
\node at (3,-3){$x$-plane};
\filldraw [black] (2,-5) circle  (3pt) node [anchor= north] {$x_0$};
\filldraw [black] (5,-5) circle  (3pt) node [anchor= north] {$x_1$};
\draw [black,dotted](2,-5)--(5,-5);
\filldraw [black] (2.5,-5) circle  (2pt) node [anchor= north] {$\alpha_{1}$};
\filldraw [black] (4.5,-5) circle  (2pt) node [anchor= north] {$\alpha_{2}$};
\draw [->,black,out=30,in=150](2,-5) to (3.5,-5);
\draw [->,black,out=-30,in=-150](3.5,-5) to (5,-5);
%\draw [->,black,out=30,in=-165](2,-5) to (3.5,-4.2);
%\draw [->,black,out=30,in=-165](2,-5) to (3.5,-4.2);
\node at (6,-5) {$I_{1,1}$};
\end{scope}

%%%%%%%%%%%%%%%%%%%% A2  %%%%%%%%%%%%%%%%%%%%%%%%
\begin{scope}[xshift=4cm,yshift=-1cm]
\filldraw [black] (2,-5) circle  (3pt) node [anchor= north] {$x_0$};
\filldraw [black] (5,-5) circle  (3pt) node [anchor= north] {$x_1$};
\draw [black,dotted](2,-5)--(5,-5);
\filldraw [black] (2.5,-5) circle  (2pt) node [anchor= north] {$\alpha_{1}$};
\filldraw [black] (4.5,-5) circle  (2pt) node [anchor= north] {$\alpha_{2}$};
\draw [->,black,out=-30,in=-150](2,-5) to (3.5,-5);
\draw [->,black,out=30,in=150](3.5,-5) to (5,-5);
\node at (6,-5) {$I^{}_{1,2}$};
\end{scope}

\begin{scope}[xshift=10cm,yshift=-.5cm]
\node at (3,-2.5){$\nu$-plane};
\filldraw [black] (2,-5) circle  (4pt) node [anchor= north] {$\nu_{0}$};
\filldraw [black] (4,-5.5) circle  (3pt) node [anchor= north] {$\nu_{1,2}$};
\filldraw [black] (4,-4.5) circle  (3pt) node [anchor= north] {$\nu_{1,1}$};
\end{scope}
\begin{scope}[xshift=4cm, yshift= -4cm]
\filldraw [black] (2,-5) circle  (3pt) node [anchor= north] {$x_0$};
\filldraw [black] (5,-5) circle  (3pt) node [anchor= north] {$x_1$};
\draw [black](2,-5)--(5,-5);
\end{scope}
\begin{scope}[xshift=4cm, yshift= -4cm]
\filldraw [black] (2.5,-4.5) circle  (2pt) node [anchor= north] {$\alpha_{1}$};
\filldraw [black] (4.5,-5.5) circle  (2pt) node [anchor= north] {$\alpha_{2}$};
\draw[black](2.5,-4.5)--(4.5,-5.5);
\end{scope}

\begin{scope} [xshift=3cm,yshift=-12cm]

\node [rectangle,draw] {$ F(x)= \frac{x^2 +\beta^2}{1+ \beta^2} ; \beta > 0$};
\end{scope}

%%%%%%%%%%%% B... I11%%%%%%%%%%%
\begin{scope}[xshift=4cm, yshift= -7cm]
\filldraw [black] (2,-5) circle  (3pt) node [anchor= north] {$x_0$};
\filldraw [black] (5,-5) circle  (3pt) node [anchor= north] {$x_1$};
\draw [->,black](2.3,-5)--(4.7,-5);
\node at (6,-5) {$I_{1}$};
\filldraw [black] (3.5,-4.5) circle  (2pt) node [anchor= north] {$\alpha_{1}$};
\filldraw [black] (3.5,-5.5) circle  (2pt) node [anchor= north] {$\alpha_{2}$};
%\draw[black](3.5,-4)--(3.5,-6);
\end{scope}

%%%%%%%%%%%%%%%%% B....I1%%%%%%%%%%%%%%%%%%%%%%%%
\begin{scope}[xshift=4cm, yshift= -9.5cm]
\filldraw [black] (2,-5) circle  (3pt) node [anchor= north] {$x_0$};
\filldraw [black] (5,-5) circle  (3pt) node [anchor= north] {$x_1$};
%\draw [->,black](2.3,-5)--(4.7,-5);
\node at (6,-5) {$I_{1,1}$};
\draw [->,black,out=30,in=-165](2,-5) to (3.5,-4.2);
\draw [->,black,out=-30,in=30](3.5,-4.2) to (3.5,-5);
\draw [->,black,out=-150,in=150](3.5,-5) to (3.5,-5.8);
\draw [->,black,out=15,in=-140](3.5,-5.8) to (5,-5);
\filldraw [black] (3.5,-4.5) circle  (2pt) node [anchor= north] {$\alpha_{1}$};
\filldraw [black] (3.5,-5.5) circle  (2pt) node [anchor= north] {$\alpha_{2}$};
%\draw[black](3.5,-4)--(3.5,-6);
\end{scope}

%%%%%%%%%%%%%%%% B........ I1prime %%%%%%%%%%%%%%%%
\begin{scope}[xshift=4cm, yshift= -12cm]
\filldraw [black] (2,-5) circle  (3pt) node [anchor= north] {$x_0$};
\filldraw [black] (5,-5) circle  (3pt) node [anchor= north] {$x_1$};
%\draw [->,black](2.3,-5)--(4.7,-5);
\filldraw [black] (3.5,-4.5) circle  (2pt) node [anchor= north] {$\alpha_{1}$};
\filldraw [black] (3.5,-5.5) circle  (2pt) node [anchor= north] {$\alpha_{2}$};
%\draw[black](3.5,-4)--(3.5,-6);
\draw [->,black,out=-30,in=165](2,-5) to (3.5,-5.8);
\draw [->,black,out=30,in=-30](3.5,-5.8) to (3.5,-5);
\draw [->,black,out=150,in=-150](3.5,-5) to (3.5,-4.2);
\draw [->,black,out=-15,in=140](3.5,-4.2) to (5,-5);
\node at (6,-5) {$I^{'}_{1,1}$};
\end{scope}

\begin{scope}[xshift=10cm, yshift= -10cm]
\filldraw [black] (2,-5) circle  (2pt) node [anchor= north] {$\nu_{0}$};
\filldraw [black] (5,-5) circle  (3pt) node [anchor= north] {$\nu_{1,1}$};
\draw [->,black](2,-5)--(3.8,-5);
\filldraw [black] (4,-5) circle  (3pt) node [anchor= north] {$\nu_{1}$};
\draw [->,dotted,out=30,in=150] (4,-5) to (5,-5);
\end{scope}

%%%%%%%%%%%%%%%%slanted alpha%%%%%%%%%%%%%%%%%%%%%%%%%%%%%%%%
%\begin{scope}[xshift=4cm, yshift= -12cm]
%\filldraw [black] (2,-5) circle  (3pt) node [anchor= north] {$x_0$};
%\filldraw [black] (5,-5) circle  (3pt) node [anchor= north] {$x_1$};
%\draw [black](2,-5)--(5,-5);
%\filldraw [black] (2.5,-4.5) circle  (2pt) node [anchor= north] {$\alpha_{1}$};
%\filldraw [black] (4.5,-5.5) circle  (2pt) node [anchor= north] {$\alpha_{2}$};
%\draw[black](2.5,-4.5)--(4.5,-5.5);
%\draw[black](3.5,-5)--(4.5,-5.5);
%\end{scope}

%%%%%%%%%%%%%%%C....I1%%%%%%%%%%%%%%%%%

\begin{scope}[xshift=4cm, yshift= -14.5cm]
\filldraw [black] (2,-5) circle  (2pt) node [anchor= north] {$\alpha_{1}$};
\filldraw [black] (5,-5) circle  (2pt) node [anchor= north] {$\alpha_{2}$};
%\draw [black](2,-5)--(5,-5);
\draw [->,black](3,-5)--(3.8,-5);
\filldraw [black] (3,-5) circle  (3pt) node [anchor= north] {$x_{0}$};
\filldraw [black] (4,-5) circle  (3pt) node [anchor= north] {$x_{1}$};
%\draw[black](3.5,-4)--(3.5,-6);
\node at (6,-5) {$I_{1,1}$};
\end{scope}

%%%%%%%%%%%%%%%%%%%%%%%C.....I12%%%%%%%%%%%%%%%%%%%%%%%%%%%%%
\begin{scope}[xshift=4cm, yshift= -16cm]
\filldraw [black] (2,-5) circle  (2pt) node [anchor= north] {$\alpha_{1}$};
\filldraw [black] (5,-5) circle  (2pt) node [anchor= north] {$\alpha_{2}$};
%\draw [black,](2,-5)--(5,-5);
%\draw [->,black](3,-5)--(3.8,-5);
\draw [->,out=150,in=30] (3,-5) to (1.6,-5);
\draw [->,out=-30,in=-150] (1.6,-5) to (3.5,-5);
\draw [->,out=30,in=150] (3.5,-5) to (5.4,-5);
\draw [->,out=-150,in=-30] (5.4,-5) to (4,-5);
\filldraw [black] (3,-5) circle  (3pt) node [anchor= north] {$x_{0}$};
\filldraw [black] (4,-5) circle  (3pt) node [anchor= north] {$x_{1}$};
\node at (6,-5) {$I_{1,3}$};
%\draw[black](3.5,-4)--(3.5,-6);
\end{scope}

%%%%%%%%%%%%%%%%%%%%%C....I13%%%%%%%%%%%%%%%%%%%%%%%%%
\begin{scope}[xshift=4cm, yshift= -17.5cm]
\filldraw [black] (2,-5) circle  (2pt) node [anchor= north] {$\alpha_{1}$};
\filldraw [black] (5,-5) circle  (2pt) node [anchor= north] {$\alpha_{2}$};%%%%%%%%%%%%%
%\draw [black,](2,-5)--(5,-5);
%\draw [->,black](3,-5)--(3.8,-5);
\draw [->,out=-150,in=-30] (3,-5) to (1.6,-5);
\draw [->,out=30,in=150] (1.6,-5) to (3.5,-5);
\draw [->,out=-30,in=-150] (3.5,-5) to (5.4,-5);
\draw [->,out=150,in=30] (5.4,-5) to (4,-5);
\filldraw [black] (3,-5) circle  (3pt) node [anchor= north] {$x_{0}$};
\filldraw [black] (4,-5) circle  (3pt) node [anchor= north] {$x_{1}$};
%\draw[black](3.5,-4)--(3.5,-6);
\node at (6,-5) {$I_{1,2}$};
\end{scope}

\begin{scope} [xshift=3cm,yshift=-20cm]

\node [rectangle,draw] {$ F(x)= \frac{x^2 -\alpha^2}{1- \alpha^2} ; 0<\alpha < 1$};
\end{scope}

\begin{scope}[xshift=10cm, yshift= -16cm]
\filldraw [black] (2,-5) circle  (2pt) node [anchor= north] {};
\filldraw [black] (4,-5) circle  (2pt) node [anchor= north] {$\nu_{1,1}$};
\filldraw [black] (4,-6) circle  (2pt) node [anchor= north] {$\nu_{1,3}$};
\filldraw [black] (4,-4) circle  (2pt) node [anchor= north] {$\nu_{1,2}$};
\draw[->,black](2,-5)--(3.8,-5);
\draw[->,black](4,-5)--(4,-5.9);
\draw[->,black](4,-5)--(4,-4.1);
\end{scope}
 \end{tikzpicture}
%%%%%%%%%%%%%%%%%%%%%%%%%%%%%%%%%%%%%%

\begin{example}:\,$F(x)=\frac{x^{p}-\alpha^{p}}{1-\alpha^{p}}=\frac{x^{p}+\beta^{p}}{1+\beta^{p}}\;\;,\;\; \epsilon=e^{\pi i/p}\beta$.
\\
\textup{
Here $\Omega$ is generated by the unit roots of order $p$. More precisely, due to the condition $\sum n_{j}\,d_{j}=0$
in (219) (with $d_{j}=1$ here) we have
$$
\Omega=2\pi i\alpha\Big( (\epsilon-1)\doZ+(\epsilon^{2}-1)\doZ+\dots (\epsilon^{p-1}-1)\doZ\Big)\quad\,\quad
\mi{with}\quad  \epsilon:=e^{2\pi i/p}
$$
Thus, except for $p\in \{2,3,4,6 \}$ the point set $\Omega$ is never discrete, but this doesn't prevent
there being, from any point of the ramified $\nu$-plane, only a finite number of {\it visible} singularities. 
\\
{\bf Case 1:} $0<\beta $.
\\
\begin{eqnarray}\label{xx12}
\Omega_{j}&=& (\eta+\Omega)\,\epsilon^{j}=\eta\,\epsilon^{j}+\Omega \quad\quad\quad (1\leq j \leq p\;\;,\;\; \epsilon=e^{\pi i/p}\beta)
\\ [1.6 ex]
\label{yy12}
\eta &=& -p \sum_{1\leq k}(-1)^{k}\frac{\beta^{-kp}}{kp+1}>0   \hspace{16 ex} (\mi{if}\;\; 1\leq \beta)
\\
\label{zz12}
\eta &=& p- p\, b_{p}\, \beta -p\sum_{1\leq k}(-1)^{k}\frac{\beta^{kp}}{kp-1}>0   \quad \quad (\mi{if}\;\; 0< \beta \leq 1)
\\
\label{ww12}
b_{p}&=&\int_{0}^{1}\frac{1+t^{p-2}}{1+t^{p}\;\;}dt\,=\,1-2\sum_{1\leq k}\frac{(-1)^{k}}{k^{2}p^{2}-1}
\end{eqnarray}
\\
{\bf Case 2:} $0<\alpha <1$.
\\
\begin{eqnarray}\label{xx13}
\Omega_{j}&=&(\eta+\pi i \alpha+\Omega)\,\epsilon^{j}\;=\;
(\eta-\pi i \alpha+\Omega)\,\epsilon^{j}
\\ \nonumber
&=&(\eta+\pi i \alpha)\,\epsilon^{j}+\Omega\;=\;
(\eta-\pi i \alpha)\,\epsilon^{j}+\Omega \quad\quad (1\leq j \leq p\;\;,\;\; \epsilon=e^{\pi i/p}\beta)
\\[1.ex] \label{yy13}
\eta &=& p - p\, a_{p}\,\alpha -p\sum_{1\leq k}\frac{\alpha^{kp}}{kp-1}>0
\\
\label{ww13}
a_{p}&=&\int_{0}^{1}\frac{1-t^{p-2}}{1-t^{p}\;\;}dt\,=\,1-2\sum_{1\leq k}\frac{1}{k^{2}p^{2}-1}
\end{eqnarray}
\\
{\bf Case 3:} $1<\alpha $.
\\
\begin{eqnarray}\label{xx14}
\Omega_{j}& =& (\eta+\Omega)\,\epsilon^{j}=\eta\,\epsilon^{j}+\Omega \quad\quad\quad (1\leq j \leq p\;\;,\;\; \epsilon=e^{\pi i/p}\beta)
\\[1. ex] \label{yy14}
\eta &=&-p\sum_{1\leq k}\frac{\alpha^{-kp}}{kp+1} < 0
\end{eqnarray} 
Remark: the expressions (\ref{ww12}) for $b_{p}$ are obtained by identifying the two distinct expressions (\ref{yy12}),(\ref{zz12}) for $\eta$ which are
are equally valid when $\beta=1$. The expressions (\ref{ww13}) for $a_{p}$ are {\it formally} obtained in the same way, i.e. by 
equating the expressions (\ref{yy13}),(\ref{yy14}) when $\alpha=1$, but since both diverge in that case, the derivation is illegitimate,
and the proper way to proceed is by rotating $\alpha$ by $e^{\pi i/p}$ so as to fall back on the situation of case 1.
Here are the $\doZ$-irreducible equations verified by the
first algebraic numbers $\alpha_{p}:=\frac{p}{\pi}\,a_{p}$ and $\beta_{p}:=\frac{p}{\pi}\,b_{p}$\,:
\[
\begin{array}{llllllllll}
0&=& \alpha_{2}   &&&&0&=& \beta_{2}-1
\\
0&=&3\,\alpha_{3}^{2}-1 &&&&0&=& 3\,\beta_{3}^{2}-4 
\\
0&=&\alpha_{4}-1 &&&&0&=& \beta_{4}^{2}-2 
\\
0&=&5\,\alpha_{5}^{4}-10\,\alpha_{5}^{2}+1 &&&&0&=& 5\,\beta_{5}^{4}-20\,\beta_{5}^{2}+16
\\
0&=&\alpha_{6}^{2}-3  &&&&0&=& \beta_{6}-2
\\
0&=&\alpha_{8}^{4}-6\,\alpha_{8}^{2}+1 &&&&0&=& \beta_{8}^{4}-8\,\beta_{8}^{2}+8
\\
0&=&\alpha_{10}^{4}-10\,\alpha_{10}^{2}+5 &&&&0&=& \beta_{10}^{2}-2\,\beta_{10}-4
\\
0&=&\alpha_{12}^{2}-4\,\alpha_{12}+1 &&&&0&=& \beta_{12}^{4}-16 \,\beta_{12}^{2}+16
\end{array}\]
Let us illustrate the situation for $p=4$ in all three 'real' cases. We choose one singularity $\nu_{0}$,
corresponding to $x_{0}=-1$ (resp. $x_{0}=1$) in case 1 or 2 (resp. 3), as base point of the $\nu$-plane, and plot as bold (resp. faint) points
all singularities visible or semi-visible from $\nu_{0}$ (resp. the closest invisible ones). For clarity, the scale (i.e. the relative values of 
$\alpha,\eta,\pi$) has not been strictly respected -- only the points' relatives position has.
}
\end{example}  

%
\\
\\

$\nu$\,-plane viewed from $\nu_0$
\newpage
%%%%%%%%%%%%%%%%%%%%%%%%%%%%%%%%%%%%%%%%%%%%%%%%%%%%%%%%%%
%%%%%%%%%%%%%%%%%%%%%%%%%%%%%%%%%%%%%%%%%%%%%%%%%%%%%%%%%%
%%%%%%%The singularity contours for B%%%%%%%%%%%%%%%%%%%%%
$x$ \,-plane viewed from \,\,$x_0$
\\
%

$\nu$\,-plane viewed from $\nu_0$

\newpage
$x$\, -plane viewed from \,$x_0$
\\

%

$\nu$ \, -plane viewed from $\nu_0$.
\newpage
$x$\, -plane viewed from $x_0$

%
%%%%%%%%%%%%%%%%%%%%%%%%%%%%%%%%%%%%%%%%%%%%%%%%%%%%%%%%%%%%%%%%%%%%%%%%%
\begin{example}:\,$F(x)=\frac{1-x}{1+x}$ or $F(x)=\frac{1-x/\alpha}{1+x/\beta}$
\\
\textup{
This interesting case is the only one where, despite the equation $F(x)\!=\!1$ having only one solution
$x_{0}\!=\!0$, the function $f$ has two logarithmic singularities, so that we get a non-trivial set  $\Omega=2\pi i\doZ$
and infinitely many copies of one and the same inner generator. From any given singularity $\nu_{0}$
there are two {\it visible} neighbouring singularities over $\dot{\nu}_{0}\pm2\pi i$ and infinitely many {\it semi-visible}
ones over $\dot{\nu}_{0}\pm2\pi i\,k\,(k\geq2)$.
}
\end{example}
\begin{example}:\,$F(x)=\frac{(1-{x}/{\alpha_{1}})(1-{x}/{\alpha_{2}})}{(1-{x}/{\alpha_{3}})}$
or $F(x)=\frac{(1-{x}/{\alpha_{1}})(1-{x}/{\alpha_{2}})}{(1-{x}/{\alpha_{3}})(1-{x}/{\alpha_{4}})}$
\\
\textup{
Here the equation $F(x)\,=\,0$ has two distinct solutions, so we have two distinct families of `parallel' inner generators,
and a set $\Omega$ which is generically discrete in the first sub-case (no $\alpha_{4}$) and
generically non-discrete in the second sub-case.\footnote{As usual, this discrete/non-discrete
dichotomy applies only to the projection on $\doC$ of the ramified $\nu$-plane, which is itself  always 
a {\it discrete} Riemann surface, with only a {\it discrete} configuration of singular points visible from any
given base point.} 
}
\end{example}
%%%%%%%%%%%%%%%%%%%%%%%%%%%%%%%%%%%%%%%%%%%%%%%%%%%%%%%%%%%%%%%%%%%%%%%%%%%%%%%%%%%%%%%%%%%%
%%%%%%%%%%%%%%%%%%%%%%%%%%%%%%%%%%%%%%%%%%%%%%%%%%%%%%%%%%%%%%%%%%%%%%%%%%%%%%%%%%%%%%%%%%%%
\subsection{Holomorphic/meromorphic inputs $F$. Examples.} 
\begin{example}: \; $F = \prod_{j=1}^{\infty} \big(1-\frac{x}{\alpha_{j}}\big)\, e^{A_{j}(x)} $
or $F =\frac{  \prod_{j=1}^{\infty} (1-{x}/{\alpha_{j}}) e^{A_{j}(x)}}{ \prod_{j=1}^{\infty} (1-{x}/{\beta_{j}}) e^{B_{j}(x)}} $.
\\
\textup{
Predictably enough, we inherit here features from the case of polynomial inputs $f$ and from that of rational
inputs $F$, but three points need to be stressed\,:
\\
(i) the presence of even a single zero $\alpha_{j}$ or of a single pole $\beta_{j}$ is enough to weld all inner algebras 
$\mi{Inner}_{f}$ into one (see \S7.2.2 {\it supra}).
\\
(ii) though, for a given $x_{0}$, the numbers $\eta_{0,j}:=-\int_{x_{0}}^{x_{j}}f(x)dx$ may accumulate 0, the corresponding
singularities $\nu_{j}$ never accumulate $\nu_{0}$ in the ramified $\nu$-plane.
\\
(iii) the question of deciding which integration paths (in the $x$-plane) lead to visible singularities (in the $\nu$-plane)
is harder to decide than for purely polynomial inputs $f$ or purely rational inputs $F$ because, unlike in these two earlier situations,
we don't always have the option of {\it deforming} a  configuration with full radial symmetry. The precise criteria for
visibility/invisibility shall be given in [SS 2].
}
\end{example} 
\begin{example}:\; $F$ = trigonometric polynomial.
\\
\textup{The series associated with knots tend to fall into this class.
Significant simplifications occur, especially when $f=-\log(F)$ is itself the derivative of a periodic function,
because the number of singularities $\nu_{j}$ then becomes finite up to $\Omega$-translations. The
special case $F(x)=4\sin^{2}(\pi x)$, which is relevant to the knot $4_{1}$, is investigated at length in the next section.
}
\end{example}  
%%%%%%%%%%%%%%%%%%%%%%%%%%%%%%%%%%%%%%%%%%%%%%%%%%%%%%%%%%%%%%%%%%%%%%%%%%%%%%%%%%%%%%%%%%%%
%%%%%%%%%%%%%%%%%%%%%%%%%%%%%%%%%%%%%%%%%%%%%%%%%%%%%%%%%%%%%%%%%%%%%%%%%%%%%%%%%%%%%%%%%%%%

%\end{document}
%%%%%%%%%%%%%%%%%%%%%%%%%%%%%%%%%%%%%%%%%%%%%%%%%%%%%%%%%%%%%%%%%%%%%%%%%%%%%%%%%%%%%%%%%%%%
%%%%%%%%%%%%%%%%%%%%%%%%%%%%%%%%%%%%%%%%%%%%%%%%%%%%%%%%%%%%%%%%%%%%%%%%%%%%%%%%%%%%%%%%%%%%
%%%%%%%%%%%%%%%%%%%%%%%%%%%%%%%%%%%%%%%%%%%%%%%%%%%%%%%%%%%%%%%%%%%%%%%%%%%%%%%%%%%%%%%%%%%%
%%%%%%%%%%%%%%%%%%%%%%%%%%%%%%%%%%%%%%%%%%%%%%%%%%%%%%%%%%%%%%%%%%%%%%%%%%%%%%%%%%%%%%%%%%%%
%%%%%%%%%%%%%%%%%%%%%%%%%%%%%%%%%%%%%%%%%%%%%%%%%%%%%%%%%%%%%%%%%%%%%%%%%%%%%%%%%%%%%%%%%%%%
%%%%%%%%%%%%%%%%%%%%%%%%%%%%%%%%%%%%%%%%%%%%%%%%%%%%%%%%%%%%%%%%%%%%%%%%%%%%%%%%%%%%%%%%%%%%
%%%%%%%%%%%%%%%%%%%%%%%%%%%%%%%%%%%%%%%%%%%%%%%%%%%%%%%%%%%%%%%%%%%%%%%%%%%%%%%%%%%%%%%%%%%%
%%%%%%%%%%%%%%%%%%%%%%%%%%%%%%%%%%%%%%%%%%%%%%%%%%%%%%%%%%%%%%%%%%%%%%%%%%%%%%%%%%%%%%%%%%%%
%%%%%%%%%%%%%%%%%%%%%%%%%%%%%%%%%%%%%%%%%%%%%%%%%%%%%%%%%%%%%%%%%%%%%%%%%%%%%%%%%%%%%%%%%%%%

%

%%%%%%%%%%%%%%%%%%%%%%%%%%%%%%%%%%%%%%%%%%%%%%%%%%%%%%%%%%
%% Some resurgence properties of knot-related functions.
%%%%%%%%%%%%%%%%%%%%%%%%%%%%%%%%%%%%%%%%%%%%%%%%%%%%%%%%%%

%\documentclass[12pt,a4paper]{article}\input{SP_commands}\begin{document}

%%%%%%%%%%%%%%%%%%%%%%%%%%%%%%%%%%%%%%%%%%%%%%%%%%%%%%%%%%%%%%%%%%%%%%%%%%%%%%%%%%%%%%%%%%%%
%%%%%%%%%%%%%%%%%%%%%%%%%%%%%%%%%%%%%%%%%%%%%%%%%%%%%%%%%%%%%%%%%%%%%%%%%%%%%%%%%%%%%%%%%%%%
%%%%%%%%%%%%%%%%%%%%%%%%%%%%%%%%%%%%%%%%%%%%%%%%%%%%%%%%%%%%%%%%%%%%%%%%%%%%%%%%%%%%%%%%%%%%
%\addtocounter{section}{6}
\section{Application to some knot-related power series.}
%%%%%%%%%%%%%%%%%%%%%%%%%%%%%%%%%%%%%%%%%%%%%%%%%%%%%%%%%%%%%%%%%%%%%%%%%%%%%%%%%%%%%%%%%%%%
%%%%%%%%%%%%%%%%%%%%%%%%%%%%%%%%%%%%%%%%%%%%%%%%%%%%%%%%%%%%%%%%%%%%%%%%%%%%%%%%%%%%%%%%%%%%

\subsection{The knot $4_1$ and the attached power series $G^{NP},G^P$.}
Knot theory attaches to each knot $\caK$ two types of power series\,:
the so-called {\it non-perturbative} series  $G_{\caK}^{NP}$
and their {\it perturbative} companions $G_{\caK}^{P}$. Both encode
the bulk of the invariant information about $\caK$ and both are largely
equivalent, though non-trivially so\,: each one can be retieved
from the other, either by non-trivial arithmetic manipulations (the Habiro
approach) or under a non-trivial process of analytic continuation (the approach
favoured in this section).

The main ingredient in the construction of  $G_{\caK}^{NP}$ and  $G_{\caK}^{P}$
is the so-called {\it quantum factorial}, classically denoted $(q)_m$\,:
\begin{equation}  \label{a182}
(q)_m:=\prod_{k=1}^{k=m}(1-q^k)
\end{equation}
For the simplest knots, namely $\caK=3_1$ or $4_1$ in standard notation, the
general definitions yield\,:
\[\begin{array}{lllllll}
\Phi_{3_1}(q)&\!:=\!&\sum_{m\geq 1}(q)_m &&
\Phi_{4_1}(q)&\!:=\!&\sum_{m\geq 1}(q)_m (q^{-1})_m 
\\[1.5 ex]
{\imi{G}}^{NP}_{3_1}(\zeta)&\!:=\!&\sum_{n\geq 0}\Phi_{3_1}(e^{2\pi i/n})\,\zeta^n &&
{\imi{G}}^{NP}_{4_1}(\zeta)&\!:=\!&\sum_{n\geq 0}\Phi_{4_1}(e^{2\pi i/n})\,\zeta^n
\\[1.5 ex]
{\tilde{G}}^{P}_{3_1}(n)&\!:=\!&\Phi_{3_1}(e^{-1/n}) = \sum{c_k}\,{n^{-k}}
&&
{\tilde{G}}^{P}_{4_1}(n)&\!:=\!&\Phi_{4_1}(e^{-1/n}) = \sum {c^\ast_k}\,{n^{-k}}
\\[1.5 ex]
{\imi{G}}^{P}_{3_1}(\nu)&\!:=\!& \sum {c_k}\;\frac{\nu^{k-1}}{(k-1)!}
&&
{\imi{G}}^{P}_{4_1}(\nu)&\!:=\!&\sum {c^\ast_k}\;\frac{\nu^{k-1}}{(k-1)!}
\end{array}\]
A few words of explanation are in order here.
\\

{\it First}, when we plug unit roots $q=e^{2\pi i/n}$ into the infinite
series $\Phi_{3_1}(q)$ or $\Phi_{4_1}(q)$, these reduce to finite sums.
\\

{\it Second}, the coefficients $\Phi_{3_1}(e^{2\pi i/n})$ or $\Phi_{4_1}(e^{2\pi i/n})$
thus defined are syntactically of {\it sum-product} type, relative
to the driving functions\,:
\begin{equation}  \label{a183}
F_{3_1}(x):=1-e^{2\pi i x}\;;\;
F_{4_1}(x):=(1-e^{2\pi i x})(1-e^{-2\pi i x})\,=\,4\sin^2(\pi x)\quad
\end{equation}

{\it Third}, whereas the non-perturbative series $\imi{G}^{NP}$ clearly
possess positive radii of convergence, their perturbative counterparts 
 $\tilde{G}^{P}$ are divergent power series of $1/n$, of Gevrey type 1, i.e.
with coefficients bounded by $$ |c_k|< \mi{Const}\; k! \quad,\quad
 |c^\ast_k|< \mi{Const}^\ast\, k!$$

{\it Fourth}, the perturbative series $\tilde{G}^P(n)$ being Gevrey-divergent,
we have to take their Borel transforms $\imi{G}^{P}\!\!(\nu)$ to restore convergence.
\\

Here, we won't discuss the series attached to knot $3_1$, because that
case has already been thoroughly investigated by Costin-Garoufalidis [CG1],[CG2]
and also because it is rather atypical, with an uncharacteristically poor
resurgence structure\,: indeed,  ${G}^{NP}_{3_1}$ and  ${G}^{P}_{3_1}$
give rise to only {\it one} inner generator $Li$, whereas it takes at least two
of them for the phenomenon of ping-pong resurgence to manifest.

So we shall concentrate on the next knot, to wit $4_1$, with its driving function 
$F(x):= 4\sin^2(\pi x)$. That case was/is also being investigated by Costin-Garoufalidis
but with  methods quite different from ours\,: see \S 12.2 below for a comparison. Here,
we shall approach the problem as a special case of {\it sum-product} series,
unravel the underlying resurgence structure, and highlight the typical
interplay between the four types of generators\,: {\it original, exceptional,
outer, inner}.

Our main {\it original} generator $Lo$ and main {\it outer} generator $Lu$,
both corresponding to the same base point $x=0$, shall turn out to be essentially equivalent,
respectively, to the {\it non-perturbative} and {\it perturbative} 
series of the classical theory, with only minor differences
stemming from the {\it ingress factor} (see below) and a trivial $2\pi i$
rotation. The exact correspondence goes like this\,:
\begin{eqnarray}  \label{a184}
\imi{G}_{4_1}^{NP}(\zeta)&\equiv& \zeta\,\partial_\zeta \imi{Lo}(\zeta)
\\[1.5 ex]  \label{a185}
\imi{G}_{4_1}^{P}(\nu)&\equiv& \frac{1}{2\pi i}\,\partial_\zeta \imi{lu}\!(2\pi i \nu)
\;=\;\frac{1}{2\pi i}\,\partial_\zeta \imi{Lu}\!(e^{2\pi i \nu}-1)
\end{eqnarray}
But we shall also introduce other generators, absent from the classical picture\,:
namely an {\it exceptional} generator $Le$, relative to the base-point $x=1/2$,
as well as a new pair consisting of a secondary {\it original} generator $Loo$ and 
a secondary {\it outer} generator $Luu$, also relative to the base-point $x=1/2$.

We shall show that these generators don't self-reproduce under alien differentiation,
but vanish without trace\,: they are mere {\it gates} for entering the true
core of the resurgence algebra, namely the {\it inner algebra}, which in the present
instance will be spanned  by just two {\it inner} generators, $Li$ and $Lii$.
%%%%%%%%%%%%%%%%%%%%%%%%%%%%%%%%%%%%%%%%%%%%%%%%%%%%%%%%%%%%%%%%%%%%%%%%%%%%%%%%%%%%%%%%%%%%
%%%%%%%%%%%%%%%%%%%%%%%%%%%%%%%%%%%%%%%%%%%%%%%%%%%%%%%%%%%%%%%%%%%%%%%%%%%%%%%%%%%%%%%%%%%%

\subsection{Two contingent ingress factors.}
Applying the rules of \S3 we find that to the driving function $F\!o$
and its translate $F\!oo$:
\begin{equation}  \label{a186}
F\!o(x)=F(x)=4\sin^2(\pi\,x)\;\; ;\; \;
F\!oo(x)=F(x\!+\!\frac{1}{2})=4\cos^2(\pi\,x)
\quad
\end{equation}
there correspond the following ingress factors\,:
\begin{equation}   \label{a187}
\mi{I\!g}_{\mi{Fo}}(n):=(4\,pi^2)^{-1/2}\;(2\pi\,n)^{2/2}=n\quad\quad ;\quad\quad
\mi{I\!g}_{\mi{Foo}}(n):= 4^{1/2}=2 \quad
\end{equation}
Their elementary character stems from the fact the only contributing
factors in $F\!o(x)$ and $F\!oo(x)$ are  $4\pi^2 x^2$ and $4$ respectively.
All other {\it binomial} or {\it exponential} factors inside $F\!o(x)$ and $F\!oo(x)$
contribute nothing, since they are {\it even} functions of $x$.

Leaving aside the totally trivial ingress factor $\mi{I\!g}_{\mi{Foo}}(n)=2$, we 
can predict
what the effect will be of removing $\mi{I\!g}_{\mi{Fo}}(n)=n$\; from 
$\imi{G}^{{}_{NP}}\!\!\!\!(\zeta)$ and all its alien derivatives\,: it will {\it smoothen}
all singularities under what shall amount to one $\zeta$-integration. In particular,
it shall replace the leading terms 
$ \mi{C_1} (\zeta-\zeta_1)^{-5/2}$ and   $ \mi{C_3} (\zeta-\zeta_3)^{-5/2}$ 
in the singularities of
$\imi{G}^{{}_{NP}}\!\!(\zeta)$ at $\zeta_1$ and $\zeta_3$ by the leading terms
$ \mi{C^\prime_1} (\zeta-\zeta_1)^{-3/2}$ and   
$ \mi{C^\prime_3}(\zeta-\zeta_3)^{-3/2}$ typical
of {\it inner generators} produced by driving functions $f(x)$ of tangency order $m=1$
(see \S4).
\\

Remark\,: an alternative, more direct but less conceptual way of deriving the form
of the ingress factor $\mi{I\!g}_{\mi{Fo}}(n)=n$ would be to use the following
trigonometric identities\,:
\begin{equation}  \label{a188}
K_{n,n-1}\equiv n^2
\;,\;
K_{2\,n,n-1}\equiv n
\;,\;
K_{2\,n,n}\equiv 4\,n
\;,\;
K_{2\,n+1,n}\equiv 2\,n+1
\end{equation}
with
\begin{equation}  \label{a189}
K_{n,m}:= \prod_{1\leq k\leq m}\;F(\frac{k}{m})
= 4^m \prod_{1\leq k\leq m}\;\sin^2(\pi\frac{k}{m})
\end{equation}
%%%%%%%%%%%%%%%%%%%%%%%%%%%%%%%%%%%%%%%%%%%%%%%%%%%%%%%%%%%%%%%%%%%%%%%%%%%%%%%%%%%%%%%%%%%%
%%%%%%%%%%%%%%%%%%%%%%%%%%%%%%%%%%%%%%%%%%%%%%%%%%%%%%%%%%%%%%%%%%%%%%%%%%%%%%%%%%%%%%%%%%%%

\subsection{Two  original generators ${Lo}$ and ${Loo}$.}
Here are the power series of {\it sum-product} type corresponding to the driving functions
$F\!o$ and $F\!oo$ (mark the lower summation bounds\,: first 1, then 0):
\begin{eqnarray}  \label{a190}
\imi{Jo}(\zeta)&\!\!:=\!\!&\sum_{1\leq n}\mi{Jo}_n\,\zeta^n
\hspace{4.ex} \mi{with}\quad \mi{Jo}_n:=
\sum^{m = n}_{m=1}\prod^{ k = m}_{k=1}\mi{Fo}(\frac{k}{n})\quad\quad
\\[1.4 ex]  \label{a191}
\imi{Joo}(\zeta)&\!\!:=\!\!&\sum_{1\leq n}\mi{Joo}_n\,\zeta^n
\hspace{3.ex} \mi{with}\quad \mi{Joo}_n:=
\sum^{ m= n}_{m=0}\prod^{ k = m}_{k=0}\mi{Foo}(\frac{k}{n})\quad\quad
\end{eqnarray}
After removal of the respective ingress factors $\mi{I\!g}_{\mi{Fo}}(n)$
and $\mi{I\!g}_{\mi{Foo}}(n)$ these become our two {\it `original generators'}\,:
\begin{eqnarray}  \label{a192}
\imi{Lo}(\zeta)&\!\!:=\!\!&
\sum_{1\leq n}\mi{Lo}_n\,\zeta^n
=\sum_{1\leq n}\frac{1}{n}{\mi{Jo}_n}\,\zeta^n
\;\; \Longrightarrow\;\;
\imi{Lo}(\zeta)=\int_0^\zeta{\imi{Jo}}(\zeta^\prime)\frac{d\zeta^\prime}{\zeta^\prime}
\quad\quad
\\[1.4 ex]  \label{a193}
\imi{Loo}(\zeta)&\!\!:=\!\!&\sum_{1\leq n}\mi{Loo}_n\,\zeta^n
=\sum_{1\leq n}\frac{1}{2}{\mi{Joo}_n}\,\zeta^n
\;\; \Longrightarrow\;\;
\imi{Loo}(\zeta)=\frac{1}{2}\imi{Joo}(\zeta)\quad\quad
\end{eqnarray} 
%%%%%%%%%%%%%%%%%%%%%%%%%%%%%%%%%%%%%%%%%%%%%%%%%%%%%%%%%%%%%%%%%%%%%%%%%%%%%%%%%%%%%%%%%%%%
%%%%%%%%%%%%%%%%%%%%%%%%%%%%%%%%%%%%%%%%%%%%%%%%%%%%%%%%%%%%%%%%%%%%%%%%%%%%%%%%%%%%%%%%%%%%

\subsection{Two  outer generators ${Lu}$ and ${Luu}$.}
The two outer generators $\imi{Lu}$ and $\imi{Luu}$ (resp. their
variants $\imi{\ell u}$ and $\imi{\ell uu}$) are produced
as outputs ${H}$ (resp. $\imi{k}$) by
inputting $F=F\!o$ or $F=1/F\!oo$ into the short chain \S5.2 and duly removing
the ingress factor 
$\mi{I\!g}_{\mi{Fo}}$
or $\mi{I\!g}_{\mi{Foo}}$. Since both $F\!o(x)$ and $F\!oo(x)$ are even functions of $x$,
we find\,:
\[\begin{array}{lllllllll}
\tilde{\ell u}(n)&\!\!:=\!\!& \frac{1}{n}\sum_{1\leq m}\prod_{1\leq k\leq m} F\!o(\frac{k}{n})
&&
\tilde{\ell uu}(n)&\!\!:=\!\!& \frac{1}{2}\sum_{1\leq m}\prod_{1\leq k\leq m}
(1/F\!oo)(\frac{k}{n})
\\
\;\;\;\Downarrow &&&&\;\;\;\Downarrow &
\\[0.3 ex]
\tilde{\ell u}(n)&\!\!:=\!\!& 
\sum_{1\leq k} c_{2k+1}\,n^{-2k-1}
&&
\tilde{\ell uu}(n)&\!\!:=\!\!& 
\sum_{0\leq k} c^\ast_{2k}\,n^{-2k}
\\
\;\;\;\Downarrow &&&&\;\;\;\Downarrow &
\\[0.3 ex]
\imi{\ell u}(\nu)&\!\!:=\!\!& 
\sum_{1\leq k} c_{2k+1}\,\frac{\nu^{2k}}{(2k)!}
&&
\imi{\ell uu}(\nu)&\!\!:=\!\!& 
\sum_{0\leq k} c^\ast_{2k}\,\frac{\nu^{2k-1}}{(2k-1)!}
\\
\;\;\;\Downarrow &&&&\;\;\;\Downarrow &
\\[0.3 ex]
\imi{Lu}(\zeta)&\!\!:=\!\!& 
\imi{\ell u}(\log(1+\zeta))
&&
\imi{Luu}(\zeta)&\!\!:=\!\!& 
\imi{\ell uu}(\log(1+\zeta))
\end{array}\]
There is a subtle difference between the two columns, though. Whereas in the left column,
the sum product $\sum\prod$ truncated at order $m$ yields the {\it exact} values
of all coefficients $c_{2k+1}$ up to order $m$, the same doesn't hold true for the
right column\,: here, the truncation of $\sum\prod$ at order $m$ yields only approximate
values of the coefficients $c_{2k}$ (of course, the larger $m$, the better the approximation).
This is because $F\!o(0)=0$ but $1/F\!oo(0)\not=0$. Therefore, whereas the short, four-link
chain of \S5.2 suffices to give the exact coefficients $c_{2k+1}$, 
one must resort to the more complex
$\mi{nur}$-transform, as articulated in the long, nine-link chain of \S5.2, to get the exact value of any given
coefficient $c^\ast_{2k}$. 

%%%%%%%%%%%%%%%%%%%%%%%%%%%%%%%%%%%%%%%%%%%%%%%%%%%%%%%%%%%%%%%%%%%%%%%%%%%%%%%%%%%%%%%%%%%%
%%%%%%%%%%%%%%%%%%%%%%%%%%%%%%%%%%%%%%%%%%%%%%%%%%%%%%%%%%%%%%%%%%%%%%%%%%%%%%%%%%%%%%%%%%%%

\subsection{Two  inner generators ${Li}$ and ${Lii}$.}
The two outer generators $\imi{Li}$ and $\imi{Lii}$ (resp. their
variants $\imi{\ell i}$ and $\imi{\ell ii}$) are produced
as outputs ${h}$ (resp. ${H}$) by
inputting 
\begin{eqnarray}  \label{a194}
f(x)\!=\!fi(x)\!=\!-\log\big(4\sin^2(\pi(x\!+\!\frac{5}{6})))= +2\sqrt{3}\pi\,x+4\pi^2\,x^2+O(x^3)
\\  \label{a195}
f(x)\!=\!fii(x)\!=\!-\log\big(4\sin^2(\pi(x\!+\!\frac{1}{6})))=-2\sqrt{3}\pi\,x+4\pi^2\,x^2+O(x^3)
\end{eqnarray} 
into the long chain of \S4.2 expressive of the {\it nir}-transform.
However, due to an obvious symmetry, it is enough to calculate $\imi{\ell i}(\nu)$
and then deduce  $\imi{\ell ii}(\nu)$ under (essentially) the chance $\nu\rightarrow -\nu$.
Notice that the tangency order here is $m=1$, leading to semi-integral powers of $\nu$\,:
\begin{eqnarray}  \label{a196}
\imi{\ell i}(\nu):=\sum_{0\leq n}\; d_{-\frac{3}{2}+n}\; \nu^{-\frac{3}{2}+n}
\quad ;\quad
\imi{\ell ii}(\nu):=\sum_{0\leq n}\; (-1)^n\,d_{-\frac{3}{2}+n}\; \nu^{-\frac{3}{2}+n}
\quad
\end{eqnarray}
Notice, too, that there is no need to bother about the ingress factors here\,:
the very definition of the {\it nir}-transform automatically provides for their
removal. 
%%%%%%%%%%%%%%%%%%%%%%%%%%%%%%%%%%%%%%%%%%%%%%%%%%%%%%%%%%%%%%%%%%%%%%%%%%%%%%%%%%%%%%%%%%%%
%%%%%%%%%%%%%%%%%%%%%%%%%%%%%%%%%%%%%%%%%%%%%%%%%%%%%%%%%%%%%%%%%%%%%%%%%%%%%%%%%%%%%%%%%%%%

\subsection{One  exceptional generator ${Le}$.}
The exceptional generators $\imi{Le}$ (resp. their
variant $\imi{\ell e}$) is produced
as  output ${h}$ (resp. ${H}$) by
inputting 
\begin{equation}  \label{a197}
f(x)\!=\!f\!o(x)\!=\!-\log\Big(4\sin^2\big(\pi(x\!+\!\frac{1}{2})\big)\Big)= -2\log2+\pi^2\,x^2+O(x^4)
\end{equation} 
into the long chain of \S4.2 expressive of the {\it nir}-transform.
The tangency order here being $m=0$
and $f\!o(x)$ being an even function of $x$, the series $\smi{\ell e}$
(resp. $\imi{\ell e}$)  carries
only integral-even (resp. integral-odd) powers of 
$\nu$\,:
\begin{eqnarray}  \label{a198}
\smi{\ell e}(\nu):=\sum_{0\leq n}\; c^{\ast\ast}_{2\,n}\; \nu^{\,2\,n}
\quad
\imi{\ell e}(\nu):=\sum_{0\leq n}\; 2n\,c^{\ast\ast}_{2\,n}\; \nu^{\,2\,n-1}
\end{eqnarray}
As with the {\it inner generators}, the {\it nir}-transform automatically 
takes care of removing  the ingress factor. 
%%%%%%%%%%%%%%%%%%%%%%%%%%%%%%%%%%%%%%%%%%%%%%%%%%%%%%%%%%%%%%%%%%%%%%%%%%%%%%%%%%%%%%%%%%%%
%%%%%%%%%%%%%%%%%%%%%%%%%%%%%%%%%%%%%%%%%%%%%%%%%%%%%%%%%%%%%%%%%%%%%%%%%%%%%%%%%%%%%%%%%%%%

\subsection{A complete system of resurgence equations.}
Before writing down the exact resurgence equations, let us depict them
graphically, in the two pictures below, where each arrow connecting
two generators signals that the {\it target generator} can be obtained
as an alien derivative of the {\it source generator}.

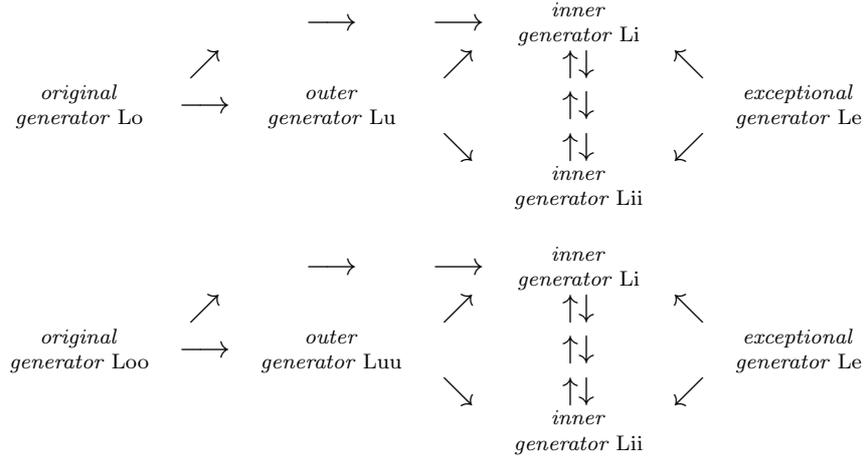
\[\begin{array}{ccccccc}
& &\longrightarrow &\longrightarrow &\mi{{inner \atop generator\;\mr{Li} }} & & 
\\
&\nearrow & &\nearrow  & \uparrow\downarrow &\nwarrow  & 
\\
\mi{{original \atop generator \;\mr{Lo} }} &\longrightarrow &\mi{{outer \atop generator
\;\mr{Lu} }} &
  &\uparrow\downarrow & & \mi{{exceptional \atop generator \;\mr{Le} }} 
\\
& & &\searrow &\uparrow\downarrow  &\swarrow & 
\\
& & & &\mi{{inner \atop generator \;\mr{Lii} }} & &
\\[3.ex]
& &\longrightarrow &\longrightarrow &\mi{{inner \atop generator\;\mr{Li} }} & & 
\\
&\nearrow & &\nearrow  & \uparrow\downarrow &\nwarrow  & 
\\
\mi{{original \atop generator \;\mr{Loo} }} &\longrightarrow &\mi{{outer \atop generator
\;\mr{Luu} }} &
  &\uparrow\downarrow & & \mi{{exceptional \atop generator \;\mr{Le} }} 
\\
& & &\searrow  &\uparrow\downarrow  &\swarrow & 
\\
& & & &\mi{{inner \atop generator \;\mr{Lii} }} & &
\end{array}\]
We observe that whereas each inner generator is both {\it source} and {\it target}, the other
generators (\--- original, outer, exceptional \---) are {\it sources} only. Moreover, although
there is perfect symmetry between $Li$ and $Lii$ within the inner algebra, that symmetry
breaks down when we adduce the original generators $L\!o$ or $L\!oo$\,: indeed,
$Li$ is a target for both $L\!o$ and $L\!oo$, but its counterpart $Lii$ is 
a target for neither.\footnote{at least, under {\it strict alien derivation}\,: this doesn't
stand in contradiction to the fact that under {\it lateral} continuation
(upper or lower) of $L\!o$ or $L\!oo$ along the real axis, singularities
$\pm 4\,i Lii$ can be ``seen" over the point $\zeta_3$. See  \S9.8.3 below.}
Altogether, we get the six resurgence algebras depicted below,
with the inner algebra as their common core\,:
\[\begin{array}{llllllll}
\mi{inner\;algebra}&&&&&
\\[1.ex]
\;\;\{\mi{Li},\mi{Lii}\}
&\subset&\{\mi{Li},\mi{Lii},\mi{Lu}\}
&\subset&\{\mi{Li},\mi{Lii},\mi{Lu},\mi{Lo}\}&&
\\[1.ex]
\;\;\{\mi{Li},\mi{Lii}\}
&\subset&\{\mi{Li},\mi{Lii},\mi{Luu}\}
&\subset&\{\mi{Li},\mi{Lii},\mi{Luu},\mi{Loo},\}&&
\\[1.ex]
\;\;\{\mi{Li},\mi{Lii}\}
&\subset&\{\mi{Li},\mi{Lii},\mi{Le}\}
&&&&
\end{array}\]
Next, we list the points $\zeta_i$ where the singularities occur in the $zeta$-plane,
and their real logarithmic counterparts $\nu_i$ in the $\nu$-plane.
\begin{eqnarray*}
\nu_{0}&:=& -\infty
\\[0.ex]
\nu_{1}&:=& \int_0^{1/6} f(x)dx=
-\frac{Li_2(e^{2\pi i/6})-Li_2(e^{-2\pi i/6})}{2\pi i}=-0.3230659470\dots
\\[1.ex]
\nu_{2}&:=& 0
\\[0.ex]
\nu_{3}&:=& \int_0^{5/6} f(x)dx=
+\frac{Li_2(e^{2\pi i/6})-Li_2(e^{-2\pi i/6})}{2\pi i}=+0.3230659470\dots
\\[1.5 ex]
\zeta_{0}&:=& 0 
\\[0.3ex]
\zeta_{1}&:=& \exp({\nu_{1}})=0.723926112\dots =1/\zeta_{3} 
\\[0.5 ex]
\zeta_{2}&:=& 1 
\\
\zeta_{3}&:=& \exp({\nu_{3}})=1.381356444\dots
\end{eqnarray*}

The assignment of {\it generators} to {\it singular points} goes like this\,:
\footnote{in the $\zeta$-plane, for definiteness.} 
\[\begin{array}{llllllllllll}
\imi{L\!o}&\mi{and}& \imi{L\!oo}&& && \mi{at}& \zeta_0 &\quad&&&
\\[1.5 ex] 
\imi{Li}&& &&&& \mi{at}& \zeta_1 &\quad ;\quad&
\underline{\imi{Li}}& \mi{at}& \underline{\zeta_1}
\\[1.5 ex] 
\imi{Lu}&\mi{and}& \imi{Luu}&\mi{and}&\imi{Le} && \mi{at}& \zeta_2 &\quad ;\quad&
\underline{\imi{Luu}}& \mi{at}& \underline{\zeta_2}
\\[1.5 ex] 
\imi{Lii}&& &&&& \mi{at}& \zeta_3 &\quad ;\quad&
\underline{\imi{Lii}}& \mi{at}& \underline{\zeta_3}
\end{array}\]
with 
$$
\underline{\zeta_1}:=-\zeta_1\;\; ,\;\; 
\underline{\zeta_2}:=-\zeta_2\;\; ,\; \;
\underline{\zeta_3}:=-\zeta_3\;\;\; \; 
$$
and
$$ 
\underline{\imi{Li}}(\zeta):=\imi{Li}(-\zeta)\;\; ,\; \;
\underline{\imi{Lii}}(\zeta):=\imi{Lii}(-\zeta)\;\; ,\; \;
\underline{\imi{Luu}}(\zeta):=\imi{Luu}(-\zeta)\;\;\; \; 
$$
The correspondence between singularities in the $\zeta$- and $\nu$-planes is as follows\,:
\[\begin{array}{llllllll}
\mi{minors}&&\mi{minors}
&\quad&
\mi{majors}&&\mi{majors}
\\[0.5 ex]
\zeta\;\mi{plane}&&\nu\;\mi{plane}
&\quad&
\zeta\;\mi{plane}&&\nu\;\mi{plane}
\\[1.5 ex]
\mr{\imi{Li}}(\zeta)&\!=\!&\mr{\imi{li}}\big(\log(1+{\zeta}/{\zeta_{i}})\big)
&\quad&
\mr{\ima{Li}}(\zeta)&\!=\!&\mr{\ima{li}}\big(-\log(1-{\zeta}/{\zeta_{i}})\big)
\\[1.5 ex]
\mr{\imi{Lii}}(\zeta)&\!=\!&\mr{\imi{lii}}\big(\log(1+{\zeta}/{\zeta_{ii}})\big)
&\quad&
\mr{\ima{Lii}}(\zeta)&\!=\!&\mr{\ima{lii}}\big(-\log(1-{\zeta}/{\zeta_{ii}})\big)
\\[1.5 ex]
\mr{\imi{Lu}}(\zeta)&\!=\!&\mr{\imi{lu}}\big(\log(1+{\zeta})\big)
&\quad&
\mr{\ima{Lu}}(\zeta)&\!=\!&\mr{\ima{lu}}\big(-\log(1-{\zeta})\big)
\\[1.5 ex]
\mr{\imi{Luu}}(\zeta)&\!=\!&\mr{\imi{luu}}\big(\log(1+{\zeta})\big)
&\quad&
\mr{\ima{Luu}}(\zeta)&\!=\!&\mr{\ima{luu}}\big(-\log(1-{\zeta})\big)
\\[1.5 ex]
\mr{\imi{Le}}(\zeta)&\!=\!&\mr{\imi{le}}\big(\log(1+{\zeta})\big)
&\quad&
\end{array}\]
With all these notations and definitions out of the way, we are 
now in a position 
to write down the resurgence equations
connecting the various generators\,:
\\
\\ \noindent
{\bf Resurgence algebra generated by $Lo $}.
\[\begin{array}{llllllllllll}
\Delta_{\zeta_1}{\imia{Lo}}&=& 2\,{\imia{Li}}
&\quad\;&\Delta_{\zeta_3-\zeta_2}{\imia{Lu}}&=&\frac{2}{2\pi}\,{\imia{Lii}}
&\quad\;&\Delta_{\zeta_3-\zeta_1}{\imia{Li}}&=&\frac{3}{2\pi}\,{\imia{Lii}}
\\[1.5 ex]
\Delta_{\zeta_2}{\imia{Lo}}&=&1\,{\imia{Lu}}
&\quad\;&\Delta_{\zeta_1-\zeta_2}{\imia{Lu}}&=&\frac{2}{2\pi}\,{\imia{Li}}
&\quad\;&\Delta_{\zeta_1-\zeta_3}{\imia{Lii}}&=&\frac{3}{2\pi}\,{\imia{Li}}
\\[1.5 ex]
\Delta_{\zeta_3}{\imia{Lo}}&=&0\,{\imia{Lii}}
&\quad\;&
\end{array}\]
{\bf Resurgence algebra generated by $Loo $}.
\[\begin{array}{llllllllllll}
\Delta_{\zeta_1}{\imia{Loo}}&=& 2\,{\imia{Li}}
&\quad\;&\Delta_{\zeta_3-\zeta_2}{\imia{Luu}}&=&\frac{2}{2\pi}\,{\imia{Lii}}
&\quad\;&\Delta_{\zeta_3-\zeta_1}{\imia{Li}}&=&\frac{3}{2\pi}\,{\imia{Lii}}
\\[1.5 ex]
\Delta_{\zeta_2}{\imia{Loo}}&=&1\,{\imia{Luu}}
&\quad\;&\Delta_{\zeta_1-\zeta_2}{\imia{Luu}}&=&\!\!-\frac{2}{2\pi}\,{\imia{Li}}
&\quad\;&\Delta_{\zeta_1-\zeta_3}{\imia{Lii}}&=&\frac{3}{2\pi}\,{\imia{Li}}
\\[1.5 ex]
\Delta_{\zeta_3}{\imia{Loo}}&=&0\,{\imia{Lii}}
&\quad\;&
\\[1.5 ex]
\Delta_{\underline{\zeta_1}}{\imia{Loo}}&=& 0\,\underline{{\imia{Li}}}
&\quad\;&\Delta_{\underline{\zeta_3}-\underline{\zeta_2}}\underline{{\imia{Luu}}}&=&\frac{2}{2\pi}\,
\underline{{\imia{Lii}}}
&\quad\;&\Delta_{\underline{\zeta_3}-\underline{\zeta_1}}\underline{{\imia{Li}}}&=&\frac{3}{2\pi}\,
\underline{{\imia{Lii}}}
\\[1.5 ex]
\Delta_{\underline{\zeta_2}}{\imia{Loo}}&=&\!\!\!-2\,\underline{{\imia{Luu}}}
&\quad\;&\Delta_{\underline{\zeta_1}-\underline{\zeta_2}}\underline{{\imia{Luu}}}
&=&\!\!-\frac{2}{2\pi}\,
\underline{{\imia{Li}}}
&\quad\;&\Delta_{\underline{\zeta_1}-\underline{\zeta_3}}\underline{{\imia{Lii}}}&=&\frac{3}{2\pi}\,
\underline{{\imia{Li}}}
\\[1.5 ex]
\Delta_{\underline{\zeta_3}}{\imia{Loo}}&=&0\,\underline{{\imia{Lii}}}
&\quad\;&
\end{array}\]
{\bf Resurgence algebra generated by $Le $}.
\[\begin{array}{llllllllllll}
&& 
&\quad\;&\Delta_{\zeta_3-\zeta_2}{\imia{Le}}&=&\frac{2}{2\pi}\,{\imia{Lii}}
&\quad\;&\Delta_{\zeta_3-\zeta_1}{\imia{Li}}&=&\frac{3}{2\pi}\,{\imia{Lii}}
\\[1.5 ex]
&&
&\quad\;&\Delta_{\zeta_1-\zeta_2}{\imia{Le}}&=&\!\!-\frac{2}{2\pi}\,{\imia{Li}}
&\quad\;&\Delta_{\zeta_1-\zeta_3}{\imia{Lii}}&=&\frac{3}{2\pi}\,{\imia{Li}}
\end{array}\]
%%%%%%%%%%%%%%%%%%%%%%%%%%%%%%%%%%%%%%%%%%%%%%%%%%%%%%%%%%%%%%%%%%%%%%%%%%%%%%%%%%%%%%%%%%%%
%%%%%%%%%%%%%%%%%%%%%%%%%%%%%%%%%%%%%%%%%%%%%%%%%%%%%%%%%%%%%%%%%%%%%%%%%%%%%%%%%%%%%%%%%%%%
\subsection{Computational verifications.}
In order to check numerically our dozen or so resurgence equations,
we shall make systematic use of the method of \S2.3 which describes {\it
singularities} in terms of {\it Taylor coefficient asymptotics}. Three situations,
however, may present themselves\,:
\\

\noindent
(i) The singularity under investigation is closest to zero. This is the most
favourable situation, as it makes for a straightforward application of \S2.3.
\\

\noindent
(ii) The singularity under investigation is not closest to zero, but becomes so
after an {\it origin-preserving} conformal transform, after which we can
once again resort to \S2.3. This is no serious complication, 
because such conformal transforms
don't diminish the accuracy with which Taylor coefficients of a given rank are computed.
\\

\noindent
(iii) The singularity under investigation is not closest to zero, nor can it
be made so under a reasonably simple, origin-preserving conformal transform.
We must then take recourse to {\it origin-changing} conformal transforms,
the simplest instances of which are {\it shifts}. This is the least favourable
case, because origin-changing conformal transforms \--- and be they simple
shifts \--- entail a steep loss of numerical accuracy and demand
great attention to the propriety of the truncations being performed.\footnote{
indeed, inept truncations can all too easily lead to meaningless results.}
\\

Fortunately, this third, least favourable situation shall occur but once
(in \S9.8.3 , when investigating the arrow $L\!o \rightarrow Lii $) and even there
we will manage the confirm the theoretical prediction with reasonable accuracy 
(up to 7 places).
In all other instances, we shall achieve  truly remarkable numerical accuracy,
often with up to 50 or 60 exact digits.

%%%%%%%%%%%%%%%%%%%%%%%%%%%%%%%%%%%%%%%%%%%%%%%%%%%%%%%%%%%%%%%%%%%%%%%%%%%%%
%%%%%%%%%%%%%%%%%%%%%%%%%%%%%%%%%%%%%%%%%%%%%%%%%%%%%%%%%%%%%%%%%%%%%%%%%%%%%
\subsubsection{\bf From $\mi{Li}$ to $\mi{Lii}$ and back ({\it inner} to {\it inner}).}
Since the theory predicts that $Li$ and $Lii$ generate each other
under alien differentiation, but that neither of them generates $L\!o$ nor $Lu$,
we may directly solve the system $ \doS_{i\,,\,ii}^{\bf{n},\bf{m}} $\,:
\begin{equation*}
\tilde{\ell i}_{-\frac{1}{2}+n}
= 3\,\nu_{1,3}^{\frac{1}{2}-n}\!\!\sum_{0\leq m < \bf{m}} 
\tilde{riis}_{\frac{1}{2}+m}\,(-\frac{1}{2}+n)^{-\frac{1}{2}-m}
\quad ; \quad \forall n\in ]\bf{n}-\bf{m},\bf{n}]
\end{equation*}
with $\bf{n}$ equations and $\bf{m}$ unknowns $\tilde{riis}_{\frac{1}{2}+m} $.
Then we may form\,:
\begin{eqnarray*}
\imi{riis}(\rho)&:=&\sum_{0\leq m <\bf{m} } \;\;\tilde{riis}_{\frac{1}{2}+m}\;\;
\frac{\rho^{-\frac{1}{2}+m}}{(-\frac{1}{2}+m)!}
\\[1.5 ex]
\smi{liis}(\nu)&:=&\imi{riis}\Big(\log(1+\frac{\nu}{\nu_{1,3}})\Big)
\,=\,
\sum_{0\leq m <\bf{m} } \;\;\smi{liis}_{-\frac{1}{2}+m}\;\;
\nu^{-\frac{1}{2}+m}
\end{eqnarray*}
and check that the ratios 
$\mi{rat}_{-\frac{1}{2}+m}:=\frac{\smi{liis}_{-\frac{1}{2}+m}}{\smi{lii}_{-\frac{1}{2}+m}}$
are indeed $\sim 1$.
For instance, with the coefficients $ \smi{liis}_{-\frac{1}{2}+m} $ computed from 
$\doS_{i\,,\,ii}^{150,45}$,
 we already get a high degree of accuracy\,:
$$ |1-\mi{rat}_{-\frac{1}{2}}|< 10^{-58},\dots,
 |1-\mi{rat}_{\frac{15}{2}}|< 10^{-40},\dots,
 |1-\mi{rat}_{\frac{31}{2}}|< 10^{-24},\dots
$$
This confirms the (equivalent) pairs of resurgence equations
\begin{eqnarray*}
\Delta_{\nu_3-\nu_1}{\smia{\ell i}}=\frac{3}{2\pi}\,{\smia{\ell ii}}
&\; ;\; & 
\Delta_{\nu_1-\nu_3}{\smia{\ell ii}}=\frac{3}{2\pi}\,{\smia{\ell i}}
\\
\Delta_{\nu_3-\nu_1}{\imia{\ell i}}=\frac{3}{2\pi}\,{\imia{\ell ii}}
&\; ;\; & 
\Delta_{\nu_1-\nu_3}{\imia{\ell ii}}=\frac{3}{2\pi}\,{\imia{\ell i}}
\end{eqnarray*}
in the $\nu$-plane, which in turn imply
$$
\Delta_{\zeta_3-\zeta_1}{\imia{Li}}=\frac{3}{2\pi}\,{\imia{Lii}}
\quad ;\quad 
\Delta_{\zeta_1-\zeta_3}{\imia{Lii}}=\frac{3}{2\pi}\,{\imia{Li}}
$$
in the $\zeta$-plane.
%%%%%%%%%%%%%%%%%%%%%%%%%%%%%%%%%%%%%%%%%%%%%%%%%%%%%%%%%%%%%%%%%%%%%%%%%%%%%
%%%%%%%%%%%%%%%%%%%%%%%%%%%%%%%%%%%%%%%%%%%%%%%%%%%%%%%%%%%%%%%%%%%%%%%%%%%%%
\subsubsection{\bf From $\mi{Lo}$ to $\mi{Li}$ ({\it original} to {\it close-inner}).}
Since $\zeta_1$ is closest to 0,
we solve the system $ \doS_{o\,,\,i}^{\bf{n},\bf{m}} $\,:
\begin{equation*}
\imi{Lo}_{n}
= 2\,\zeta_{1}^{-n}\!\!\sum_{0\leq m < \bf{m}} 
\tilde{lis}_{-\frac{1}{2}+m}\;\,n^{\frac{1}{2}-m}
\quad ; \quad \forall n\in ]\bf{n}-\bf{m},\bf{n}]
\end{equation*}
with $\bf{n}$ equations and $\bf{m}$ unknowns $\tilde{lis}_{-\frac{1}{2}+m} $.
Then we  check that the ratios 
$\mi{rat}_{-\frac{1}{2}+m}:=\frac{\smi{lis}_{-\frac{1}{2}+m}}{\smi{li}_{-\frac{1}{2}+m}}$
are indeed $\sim 1$.
For instance, with the coefficients $ \smi{lis}_{-\frac{1}{2}+m} $ computed from 
$\doS_{o\,,\,i}^{700,50}$,
 we get this sort of accuracy\,:
$$ |1-\mi{rat}_{-\frac{1}{2}}|< 10^{-54},\dots,
 |1-\mi{rat}_{\frac{15}{2}}|< 10^{-29},\dots,
 |1-\mi{rat}_{\frac{31}{2}}|< 10^{-6},\dots
$$
This confirms the resurgence equations
$
\Delta_{\zeta_1}{\imia{Lo}}=2\;\,{\imia{Li}}
$
in the $\zeta$-plane.
%%%%%%%%%%%%%%%%%%%%%%%%%%%%%%%%%%%%%%%%%%%%%%%%%%%%%%%%%%%%%%%%%%%%%%%%%%%%%
%%%%%%%%%%%%%%%%%%%%%%%%%%%%%%%%%%%%%%%%%%%%%%%%%%%%%%%%%%%%%%%%%%%%%%%%%%%%%
\subsubsection{\bf From $\mi{Lo}$ to $\mi{Lii}$ ({\it original} to {\it distant-inner}).}
The singular point $\zeta_3$ being farthest from 0, we first resort to an origin-preserving
conformal transform $\zeta\rightarrow \xi$\,:
\begin{eqnarray*}
 h_{\zeta,\xi}  &:& \xi \mapsto \zeta:=\zeta_1-\big(\zeta_1^{1/4}-\xi \big)^4
\hspace{19. ex} \forall \xi
\\
 h_{\xi,\zeta}  &:& \zeta \mapsto \xi:=\zeta_1^{1/4}-\big(\zeta_1-\zeta \big)^{1/4}
\hspace{16. ex} \forall \zeta \in [0,\zeta_1]
\\
 h^{+}_{\xi,\zeta}  &:& \zeta \mapsto \xi:=
\zeta_1^{1/4}-\big(\zeta-\zeta_1 \big)^{1/4}\,e^{-i\pi/4}
\hspace{14. ex} \forall \zeta \in [\zeta_1,\infty]
\\
 h^{-}_{\xi,\zeta}  &:& \zeta \mapsto \xi:=
\zeta_1^{1/4}-\big(\zeta-\zeta_1 \big)^{1/4}\,e^{+i\pi/4}
\hspace{14. ex} \forall \zeta \in [\zeta_1,\infty]
\end{eqnarray*}
\[\begin{array}{lllllllll}
 h_{\xi,\zeta}  &\!:\!& \zeta_1 \mapsto \xi_1 &\!\!=\!\!& 0.9224\dots
&;& |\xi_1|\,=\,0.9224\dots & \mi{(farthest)}
\\
 h^{\pm}_{\xi,\zeta}  &\!:\!& \zeta_2 \mapsto \xi_2^{\pm}&\!\!=\!\!&
0.4098\pm0.5126 \, i\dots 
&;& |\xi_2^{\pm}|=0.6563\dots & \mi{(closest)}
\\
 h^{\pm}_{\xi,\zeta}  &\!:\!& \zeta_3 \mapsto \xi_3^{\pm}&\!\!=\!\!&
0.2857\pm0.6367\, i\dots
&;& |\xi_3^{\pm}|=0.6979 \dots & \mi{(middling)}
\end{array}\]
Since the images $\xi_3^{\pm}$ are closer, but not closest, to 0, we must
perform an additional shift $\xi \rightarrow \tau$\,:
\[\begin{array}{lllllllll}
 h_{\tau,\xi}  &\!\!:\!\!& \xi \mapsto \tau &\!\!:=\!\!& \xi-\frac{i}{2}
&  h_{\xi,\tau} \; :\; \tau \mapsto \xi:=\tau+\frac{i}{2}
\\[1.5 ex]
 h_{\tau,\xi}  &\!\!:\!\!& \xi_1 \mapsto \tau_1&\!\!=\!\!& 0.9224-0.5000\,i\dots
&|\tau_1|=1.0492\dots \mi{(farthest)}
\\
 h_{\tau,\xi}  &\!\!:\!\!& \xi_2^{+} \mapsto \tau_2^{+}&\!\!=\!\!&0.4098+0.0125\,i\dots
&|\tau^{+}_2|=0.4100\dots \mi{(middling)}
\\
 h_{\tau,\xi}  &\!\!:\!\!& \xi_3^{+} \mapsto \tau_3^{+}&\!\!=\!\!&0.2857+0.1367\,i\dots
&|\tau^{+}_3|=0.3167\dots \mi{(closest)}
\end{array}\]
The image $\tau^+_3$ at last is closest, and we can now go through the usual motions. We
form successively\,:
\begin{eqnarray*}
\imi{Lo}_{\#}\!(\zeta) &\!\!:=\!\!& \sum_{0<n<\bf{n}} \imi{Lo}_n\;\zeta^n
\hspace{32.ex} \mi{(truncation)}
\\[1.5 ex]
\imi{Lo}_{\#\#}\!(\xi) &\!\!:=\!\!&\imi{Lo}_{\#}\!(h_{\zeta,\xi}(\xi))
\hspace{31.ex} \mi{(conf.\, transf.)}
\\[1.5 ex]
\imi{Lo}_{\#\#\#}\!(\xi) &\!\!:=\!\!&\imi{Lo}_{\#\#}\!(h_{\xi,\tau}(\tau))
= \sum_{0<n<\bf{n}} {L}_n\;\tau^n + (\dots)
\hspace{8.ex} \mi{(simple\, shift.)}
\end{eqnarray*}
We then solve the system $\doS^{\bf{n},\bf{m}}_{o\,,\,ii} $\,:
$$
L_n=4\,i\,(\tau_3^+)^{-n}\, \sum_{0\leq m < \bf{m}}P_{-\frac{1}{2}+m}\;\; n^{\frac{1}{2}-m}
\hspace{10.ex} \big(n\,\in\, ]\bf{n}-\bf{m},\bf{n}]\big)
$$
with $\bf{m}$ equations and  $\bf{m}$ unknowns $P_{-\frac{1}{2}+m}$.
\begin{eqnarray*}
\imi{P}(\nu)&:=& \sum_{0\leq m < \bf{m}} P_{-\frac{1}{2}+m}\;\;
\frac{\nu^{-\frac{3}{2}+m}}{(-\frac{3}{2}+m)!}+(\dots)
\\
\imi{R}(\tau)&:=&\imi{P}(\log(1+\frac{\tau}{\tau^+_3}))\;=
\;\sum_{0\leq m < \bf{m}} R_{-\frac{3}{2}+m}\;\;\tau^{-\frac{3}{2}+m} +(\dots)
\end{eqnarray*}
 
Next, for comparison, we form series that carry the expected singularity
$Lii$ successively in the $\nu$, $\zeta$ and $\tau$-planes\,:
\begin{eqnarray*}
\imi{\ell\,ii}(\nu)&:=& \sum_{0\leq m < \bf{m}} \tilde{\ell\,ii}_{-\frac{1}{2}+m}\;
\frac{\nu^{-\frac{3}{2}+m}}{(-\frac{3}{2}+m)!}+(\dots)
\\
\imi{Lii}(\zeta)&:=& \imi{\ell\,ii}\big(\log(1+\frac{\zeta}{\zeta_3})\big)
\\
\imi{Q}(\tau)&:=&\imi{Lii}(dh_{\zeta,\tau}(\tau))\;=
\;\sum_{0\leq m < \bf{m}} Q_{-\frac{3}{2}+m}\;\;\tau^{-\frac{3}{2}+m} +(\dots)
\end{eqnarray*} 
Lastly, we form the ratios
$\mi{rat}_{-\frac{3}{2}+m}:=\frac{R_{-\frac{3}{2}+m}}{Q_{-\frac{3}{2}+m}} $
of homologous coefficients $P,Q$ and check that these ratios are $\sim 1$.
With the data derived from the linear system
$\doS_{o\,,\,ii}^{\,800,4}$
and with truncation at order $\bf{n^\ast}=20 $ in the computation of 
$\imi{Lo}_{\#\#\#} $, 
we get the following, admittedly poor degree\footnote{this is because
of the recourse to the {\it shift}  $\tau:=\xi+\frac{i}{2} $ whereas
in all the other computations we handled less disruptive 
{\it origin-preserving} conformal
transforms $\zeta\rightarrow \xi $.
} of accuracy\,:
$$
|1-\mi{rat}_{-3/2}|<10^{-7},
|1-\mi{rat}_{-1/2}|<10^{-3},
|1-\mi{rat}_{+1/2}|<10^{-2},\dots
$$
To compound the poor numerical accuracy, the theoretical interpretation
is also rather roundabout in this case. By itself, the above results only
show that\,:
\begin{equation}  \label{a199}
\Delta^{\pm}_{\zeta_3} \imia{Lo}\; =\pm 4\,i \imia{Lii}
\end{equation} 
with the one-path lateral operators $\Delta_\omega^\pm$ of \S2.3 which, unlike
the multi-path averages $\Delta_\omega$, are {\it not} alien derivations. To 
infer from (\ref{a199}) the expected resurgence equation\,:
\begin{equation}  \label{a200}
\Delta_{\zeta_3} \imia{Lo}\; =0 \imia{Lii}
\end{equation} 
we must apply the basic identity (5) of \S2.3 to $\imia{L\!o}$\,:
\begin{equation}  \label{a201}
\Big(1+\sum_{0<\omega}\Delta^{+}_{\omega} \Big)\imia{L\!o}\; =\;
 \Big(\exp\big(2\pi i\sum_{0<\omega}\Delta_{\omega} \big )\Big)\imia{L\!o}
\end{equation}
and then equate the sole term coming from the left-hand side,
namely $\Delta^{\pm}_{\zeta_3} \imia{Lo} $, with the 4 possible terms
coming from the right-hand side, namely\,:
\begin{eqnarray}  \label{a202}
2\pi i\;\Delta_{\zeta_3} \imia{Lo}  &=& \mi{unknown}
\\  \label{a203}
\frac{(2\pi i)^2}{2} \Delta_{\zeta_3-\zeta_1}\Delta_{\zeta_1} \imia{Lo}
&=& 1\;\imia{Lii}
\\  \label{a204}
\frac{(2\pi i)^2}{2} \Delta_{\zeta_3-\zeta_2}\Delta_{\zeta_2} \imia{Lo}
&=& 3\;\imia{Lii}
\\  \label{a205}
\frac{(2\pi i)^3}{6} \Delta_{\zeta_3-\zeta_2}\Delta_{\zeta_2-\zeta_1}\Delta_{\zeta_1}
 \imia{Lo}
&=& 0\;\imia{Lii}
\end{eqnarray}
Equating the terms in the left and right clusters, we find that the
sole unknown term (\ref{a202}) does indeed vanish, as required by the theory.  
%%%%%%%%%%%%%%%%%%%%%%%%%%%%%%%%%%%%%%%%%%%%%%%%%%%%%%%%%%%%%%%%%%%%%%%%%%%%%
%%%%%%%%%%%%%%%%%%%%%%%%%%%%%%%%%%%%%%%%%%%%%%%%%%%%%%%%%%%%%%%%%%%%%%%%%%%%%
\subsubsection{\bf From $\mi{Lo}$ to $\mi{Lu}$ ({\it original} to {\it outer}).}
A single, origin-preserving conformal transform $\zeta\rightarrow \xi$
takes the singular point $\zeta_2$  to middling position $\xi_2^\pm$\,:
\begin{eqnarray*}
 h_{\zeta,\xi}  &:& \xi \mapsto \zeta:=\zeta_1-\big(\zeta_1^{1/2}-\xi \big)^2
\hspace{19. ex} \forall \xi
\\
 h_{\xi,\zeta}  &:& \zeta \mapsto \xi:=\zeta_1^{1/2}-\big(\zeta_1-\zeta \big)^{1/2}
\hspace{16. ex} \forall \zeta \in [0,\zeta_1]
\\
 h^{+}_{\xi,\zeta}  &:& \zeta \mapsto \xi:=
\zeta_1^{1/2}+i\,\big(\zeta-\zeta_1 \big)^{1/2}
\hspace{14. ex} \forall \zeta \in [\zeta_1,\infty]
\\
 h^{-}_{\xi,\zeta}  &:& \zeta \mapsto \xi:=
\zeta_1^{1/2}-i\,\big(\zeta-\zeta_1 \big)^{1/2}
\hspace{14. ex} \forall \zeta \in [\zeta_1,\infty]
\end{eqnarray*}
\[\begin{array}{lllllllll}
 h_{\xi,\zeta}  &\!:\!& \zeta_1 \mapsto \xi_1 &\!\!=\!\!& 0.8508\dots
&;& |\xi_1|\,=\,0.8508\dots & \mi{(closest)}
\\
 h^{\pm}_{\xi,\zeta}  &\!:\!& \zeta_2 \mapsto \xi_2^{\pm}&\!\!=\!\!&
0.8508\pm 0.5254\, i\dots 
&;& |\xi_2^{\pm}|= 1.0000\dots & \mi{(middling)}
\\
 h^{\pm}_{\xi,\zeta}  &\!:\!& \zeta_3 \mapsto \xi_3^{\pm}&\!\!=\!\!&
0.8508\pm 0.8108\, i\dots
&;& |\xi_3^{\pm}|= 1.7573\dots & \mi{(farthest)}
\end{array}\]
Then we form\,:
\begin{eqnarray*}
\imi{Lo}_{\#}\!(\zeta) &\!\!:=\!\!& \sum_{0<n<\bf{n}} \imi{Lo}_n\;\zeta^n
\hspace{32.ex} \mi{(truncation)}
\\[0.7 ex]
\imi{Lo}_{\#\#}\!(\xi) &\!\!:=\!\!&\imi{Lo}_{\#}\!(h_{\zeta,\xi}(\xi))
\hspace{31.ex} \mi{(conf.\, transf.)}
\\[1.5 ex]
\imi{Lo}_{\#\#\#}\!(\xi) &\!\!:=\!\!&\imi{Lo}_{\#\#}\!(\xi)\; (\xi_1-\xi)^3
= \sum_{0<n<\bf{n}} {L}_n\;\xi^n + (\dots)
\hspace{4.ex} \mi{(sing.\, remov.)}
\end{eqnarray*}
Since $\zeta_1 -\zeta=(\xi_1-\xi)^2$, all the semi-integral powers
$ (\zeta_1-\zeta)^{n/2}$ present  in $\imi{L\!o}_{\#}(\zeta)$ at
$\zeta\sim \zeta_2$ vanish from $\imi{L\!o}_{\#\#}(\xi)$, except
for the first two terms\,:
$$\mi{C_{-3}}\;(\xi_1 -\xi)^{-3} +\mi{C_{-1}}\;(\xi_1 -\xi)^{-1}  $$
but even these two vanish from $\imi{L\!o}_{\#\#\#}(\xi)$ due to multiplication
by $(\xi_1 -\xi)^{3} $. So the points $\xi_2^{\pm}$ now carry the closest singularities
of $\imi{L\!o}_{\#\#\#}(\xi)$, and we can apply the usual Taylor coefficient
asymptotics.

 For comparison with the expected singularity $\imi{Lu}$, we construct
 a new triplet $\{\imi{Ro}_{\#},\imi{Ro}_{\#\#},\imi{Ro}_{\#\#\#}\}$,
but with a more severe truncation ($\bf{n^\ast}\prec \bf{n} $) and with coefficients
$\imi{Lo}_n$ replaced by the $\imi{Ro}_n$ defined as follows\,: 
$$
\imi{Ro}_n:=\frac{1}{n}\mi{SP}^F(\frac{1}{n})\quad\quad\quad \mi{with}\quad \quad\quad
\mi{SP}^{F}(x):=
\sum_{1 \leq m \leq \bf{n^\ast}} 
\prod_{1\leq k \leq m}F(k\,x) 
$$
\begin{eqnarray*}
\imi{Ro}_{\#}\!(\zeta) &\!\!:=\!\!& \sum_{0<n<\bf{n^\ast}} \imi{Ro}_n\;\zeta^n
\hspace{32.ex} \mi{(truncation)}
\\[0. ex]
\imi{Ro}_{\#\#}\!(\xi) &\!\!:=\!\!&\imi{Ro}_{\#}\!(h_{\zeta,\xi}(\xi))
\hspace{31.ex} \mi{(conf.\, transf.)}
\\[0.5 ex]
\imi{Ro}_{\#\#\#}\!(\xi) &\!\!:=\!\!&\imi{Ro}_{\#\#}\!(\xi)\; (\xi_1-\xi)^3
= \sum_{0<n<\bf{n}} {R}_n\;\xi^n + (\dots)
\hspace{4.ex} \mi{(sing.\, remov.)}
\end{eqnarray*}
Then, we solve the two parallel systems $\overline{\doS}_{o\,,\,u}^{\bf{n},\bf{m}}$ and 
$\underline{\doS}_{\,o\,,\,u}^{\bf{n},\bf{m}}$\,:
\begin{eqnarray*}
L_n=\sum_{\epsilon=\pm}\;\;
({\xi_2^{\epsilon}})^{-n}\!\!\sum_{1\leq m\leq  \bf{m}}L^{\epsilon}_m\;n^{-k}
\hspace{15.ex}\big(n\;\in \;]\bf{n}-2\,\bf{m},\bf{n}]\;\big)
\\
R_n=\sum_{\epsilon=\pm}\;\;
({\xi_2^{\epsilon}})^{-n}\!\!\sum_{1\leq m\leq  \bf{m}}R^{\epsilon}_m\;n^{-k}
\hspace{15.ex}\big(n\;\in \;]\bf{n}-2\,\bf{m},\bf{n}]\;\big)
\end{eqnarray*}
 each with $2\,\bf{m}$ equations and $2\,\bf{m}$ unknowns, $L^\epsilon_m$
or $R^\epsilon_m$ respectively. We then check that the ratios 
$ \mi{rat}^{\epsilon}_n:=\frac{L^{\epsilon}_{m}}{R^{\epsilon}_{m}} $ are $\sim 1$.
With the data obtained from the systems
$\overline{\doS}_{o\,,\,u}^{\,495,7}$ and $\underline{\doS}_{\,o\,,\,u}^{\,495,7}$ 
and with truncation at order $\bf{n^\ast}=30$ in the $\imi{Ro}$ triplet, we get
the following degree of accuracy\,:
$$
|1-\mi{rat}^{\pm}_1|<10^{-17},\dots,
|1-\mi{rat}^{\pm}_3|<10^{-13},\dots,
|1-\mi{rat}^{\pm}_6|<10^{-10},\dots,
$$
The immediate implication is $\Delta_{\zeta_2}^{+}\imia{L\!o}\;= 2\pi i \,\imia{Lu} $.
To translate this into a statement about 
$\Delta_{\zeta_2}\imia{L\!o}$, the argument is the same as in \S9.8.3, only
much simpler. Indeed, the only term coming from the left-hand side of (\ref{a201})
is now $ \Delta_{\zeta_2}^{+}\imia{L\!o}$ and the only two possible terms 
coming from the right-hand side are\,:
\begin{equation}  \label{a206}
2\pi i\;\Delta_{\zeta_2}\imia{L\!o}\;=\mr{unknown}\quad \mi{and}
\quad
\frac{(2\pi i)^2}{2}\Delta_{\zeta_2-\zeta_1}\Delta_1 \imia{L\!o}= 0
\end{equation}
Equating both sides, we find 
$ \Delta_{\zeta_2}\imia{L\!o}\;=\imia{Lu} $, as required by the theory.
%%%%%%%%%%%%%%%%%%%%%%%%%%%%%%%%%%%%%%%%%%%%%%%%%%%%%%%%%%%%%%%%%%%%%%%%%%%%%
%%%%%%%%%%%%%%%%%%%%%%%%%%%%%%%%%%%%%%%%%%%%%%%%%%%%%%%%%%%%%%%%%%%%%%%%%%%%%
\subsubsection{\bf From $\mi{Lu}$ to $\mi{Li}$ and $\mi{Lii}$ 
({\it outer} to {\it inner}).}
The singular points under investigation being closest, the investigation is straightforward.
We form the linear system $\doS_{u,i/ii}^{\bf{n},\bf{m}}$\,:
$$
\imi{\ell u}_n=2\;(\nu_3^{-n}+(-\nu_3)^{-n}\big)\sum_{0\leq m < \bf{2\,m}}
\tilde{riis}_{-\frac{1}{2}+m}\;n^{\frac{1}{2}-m}
\quad\quad \Big(n\;\in\;]\bf{n}\!-\!2\bf{m},\bf{n}] \Big)
$$
with $\bf{m}$ effective equations (for even values of $n$) and $\bf{m}$
unknowns $\tilde{riis}_{-\frac{1}{2}+m}$.
We then form\,:
\begin{eqnarray*}
\imi{riis}(\xi)&:=&\sum_{0 \leq m< \bf{m}} \tilde{riis}_{-\frac{1}{2}+m}\;\;
\frac{\xi^{-\frac{3}{2}+m}}{(-\frac{3}{2}+m)!}
\\[1.0 ex]
\imi{liis}(\xi)&:=&\imi{riis}\big(\log(1+\frac{\nu}{\nu_3})\big)
\;=\;
\sum_{0 \leq m< \bf{m}} \imi{liis}_{-\frac{3}{2}+m}\;\;
\nu^{-\frac{3}{2}+m}
\end{eqnarray*}
and check that the ratios 
$\mi{rat}_{-\frac{3}{2}+m}:=\frac{\imi{liis}_{-\frac{3}{2}+m}}{\imi{lii}_{-\frac{3}{2}+m}} $
are indeed $\sim 1$. 
With the data obtained from the system 
${\doS}_{u\,,\,i/ii}^{\,300,40}$ , we get
this high degree of accuracy\,:
\begin{eqnarray*}
|1-\mi{rat}^{\pm}_{-3/2}|<10^{-81},&\dots&,
|1-\mi{rat}^{\pm}_{17/2}|<10^{-39},\;\dots\;,
\\[1.5 ex]
|1-\mi{rat}^{\pm}_{37/2}|<10^{-21},&\dots&,
|1-\mi{rat}^{\pm}_{57/2}|<10^{-10},\;\dots
\end{eqnarray*}
This confirms, via the $\nu$-plane, the expected resurgence equations in the $\zeta$-plane,
namely\,:
$$
\Delta_{\zeta_3-\zeta_2}{\imia{Lu}}=\frac{2}{2\pi}\;{\imia{Lii}}
\quad ;\quad 
\Delta_{\zeta_1-\zeta_2}{\imia{Lu}}=\frac{2}{2\pi}\;{\imia{Li}}
$$
%%%%%%%%%%%%%%%%%%%%%%%%%%%%%%%%%%%%%%%%%%%%%%%%%%%%%%%%%%%%%%%%%%%%%%%%%%%%%
%%%%%%%%%%%%%%%%%%%%%%%%%%%%%%%%%%%%%%%%%%%%%%%%%%%%%%%%%%%%%%%%%%%%%%%%%%%%%
\subsubsection{\bf From $\mi{Loo}$ to $\mi{Li}$
 ({\it original} to {\it close-inner}).}
We proceed exactly as in \S9.8.2. 
We form the linear system $\doS_{oo\,,\,i}^{\bf{n},\bf{m}}$\,:
$$
\imi{Loo}_n=\zeta_1^{-n}\sum_{0\leq m < \bf{m}}
\tilde{lis}_{-\frac{1}{2}+m}\;n^{\frac{1}{2}-m}
\quad\quad \Big(n\;\in\;]\bf{n}\!-\!2\bf{m},\bf{n}] \Big)
$$
with $\bf{m}$  equations  and $\bf{m}$
unknowns $\tilde{lis}_{-\frac{1}{2}+m}$.
We then  check that the ratios 
$\mi{rat}_{-\frac{3}{2}+m}:=\frac{\imi{lis}_{-\frac{3}{2}+m}}{\imi{li}_{-\frac{3}{2}+m}} $
are indeed $\sim 1$. 
With the data obtained from the system 
${\doS}_{oo\,,\,i}^{\,800,30}$ , we get
this degree of accuracy\,:
$$
|1-\mi{rat}^{\pm}_{-1/2}|<10^{-51},\dots,
|1-\mi{rat}^{\pm}_{19/2}|<10^{-22},\dots,
|1-\mi{rat}^{\pm}_{39/2}|<10^{-7},\dots
$$
This confirms the expected resurgence equations in the $\zeta$-plane\,:
$$
\Delta_{\zeta_1}{\imia{Lo}}=2\;\,{\imia{Li}}
\quad ;\quad 
\Delta_{\underline{\zeta_1}}{\imia{Lo}}=0\;\,\underline{{\imia{Li}}}
$$
An alternative method  would to check that $\imi{Lo}\!(\zeta) -\!\!\imi{Loo}\!(\zeta)$
has radius of convergence 1, which means that $\imi{Lo}$
and $\imi{Loo}$ have the same singularity at $\zeta_1$, namely
$\imia{Li} $\,: see \S 9.8.2. With that method, too, the numerical confirmation is excellent.
%%%%%%%%%%%%%%%%%%%%%%%%%%%%%%%%%%%%%%%%%%%%%%%%%%%%%%%%%%%%%%%%%%%%%%%%%%%%%
%%%%%%%%%%%%%%%%%%%%%%%%%%%%%%%%%%%%%%%%%%%%%%%%%%%%%%%%%%%%%%%%%%%%%%%%%%%%%
\subsubsection{\bf From $\mi{Loo}$ to $\mi{Lii}$
 ({\it original} to {\it distant-inner}).}
The verication hasn't been done yet. The theory, however, predicts a vanishing alien
derivative 
$\Delta_{\zeta_3}(\imia{Loo})=0 $ just as with $\imia{Lo} $. Therefore,
the upper/lower lateral singularity seen at $\zeta_3$ when continuing $\imia{Loo}(\zeta) $
should be $\pm 4\,i\!\imia{Lii}$, just as was the case with the lateral continuations
of $\imia{Lo}(\zeta) $.
%%%%%%%%%%%%%%%%%%%%%%%%%%%%%%%%%%%%%%%%%%%%%%%%%%%%%%%%%%%%%%%%%%%%%%%%%%%%%
%%%%%%%%%%%%%%%%%%%%%%%%%%%%%%%%%%%%%%%%%%%%%%%%%%%%%%%%%%%%%%%%%%%%%%%%%%%%%
\subsubsection{\bf From $\mi{Loo}$ to $\mi{Luu}$
 ({\it original} to {\it outer}).}
We form the linear system $\doS_{oo\,,\,uu}^{\bf{n},\bf{m}}$\,:
\begin{eqnarray*}
\imi{Loo}_n-\imi{Lo}_n
&=&-\zeta_2^{-n}\sum_{1\leq m \leq \bf{m}}
\tilde{lu}_{1+2m}\;n^{-1-2\,m}
\\[1.5 ex]
&&+\zeta_2^{-n}\sum_{1\leq m \leq \bf{m}}
\tilde{luus}_{\,2m}\;\;n^{-2\,m}
\\[1.5 ex]
&-2&\;({-\zeta_2})^{-n}\sum_{1\leq m \leq \bf{m}}
\tilde{luus}_{\,2m}\;\;n^{-2\,m}
\end{eqnarray*}
with $\bf{m}$  equations\footnote{  with $ n$ ranging
through the interval $]\bf{n}\!-\!2\bf{m},\bf{n}] $. } and
$\bf{m}$ unknowns 
$\tilde{luus}_{2\,m}$.\;\footnote{the coefficients 
$\tilde{lu}_{1+2m} $ are already known, from \S9.8.4. }
We then  check that the ratios 
$\mi{rat}_{2m}:=\frac{\imi{luus}_{2m}}{\imi{luu}_{2m}} $
are indeed $\sim 1$. 
With the data obtained from the system 
${\doS}_{oo\,,\,uu}^{\,600,30}$ , we get
this level of accuracy\,:
$$
|1-\mi{rat}_{2}|<10^{-48},\dots,
|1-\mi{rat}_{12}|<10^{-28},\dots,
|1-\mi{rat}_{24}|<10^{-15},\dots
$$
This directly confirms the expected resurgence equations\,:
$$ \Delta_{\zeta_2}\imia{L\!oo}\;=\imia{Luu}
\quad \quad ;\quad\quad 
\Delta_{\underline{\zeta_2}}\imia{L\!oo}\;=-2\,\underline{\imia{Luu}} 
$$
%%%%%%%%%%%%%%%%%%%%%%%%%%%%%%%%%%%%%%%%%%%%%%%%%%%%%%%%%%%%%%%%%%%%%%%%%%%%%
%%%%%%%%%%%%%%%%%%%%%%%%%%%%%%%%%%%%%%%%%%%%%%%%%%%%%%%%%%%%%%%%%%%%%%%%%%%%%
\subsubsection{\bf From $\mi{Luu}$ to $\mi{Li}$ and $\mi{Lii}$
({\it outer} to {\it inner}).}
We proceed exactly as in \S9.8.8.
We solve the linear system $\doS^{\bf{n},\bf{m}}_{uu\,,\,i/ii}$ \,:
$$
\imi{luu}_n\,=\,\frac{2}{2\pi}\,(\nu_3^{-n}-(-\nu_3)^{-n})\sum_{0\leq m<\bf{m}}
\; \tilde{riis}_{-\frac{1}{2}+m}\;n^{\frac{1}{2}-m}
\quad \quad\quad\big( n\,\in\, ]\bf{n}\!-\!2\bf{m}\,,\,\bf{n}]\big)
$$
with $\bf{m}$ effective equations (for $n$ odd) and $\bf{m}$ unknowns
$\tilde{riis}_{-\frac{1}{2}+m}$. Then we form\,:
\begin{eqnarray*}
\imi{riis}(\xi)&:=& \sum_{0\leq m < \bf{m}}
 \tilde{riis}_{-\frac{1}{2}+m}\;\frac{\xi^{-\frac{3}{2}+m}}{(-\frac{3}{2}+m)!}+(\dots)
\\[1. ex]
\imi{liis}(\xi)&:=& \imi{riis}\big(\log(1+\frac{\nu}{\nu_3})\big)
\;=:\;\sum_{0\leq m < \bf{m}}
 \imi{liis}_{-\frac{3}{2}+m}\;{\nu^{-\frac{3}{2}+m}}+(\dots)
\end{eqnarray*}
and check that the ratios 
$\mi{rat}_{-\frac{3}{2}+m}:=\frac{\imi{liis}_{-\frac{3}{2}+m}}{\imi{lii}_{-\frac{3}{2}+m}}$
of homologous coefficients are indeed $\sim 1$. For the data corresponding to 
$\doS^{300,40}_{uu\,,\,i/ii}$, we find this excellent level of accuracy\,:
$$
|1-\mi{rat}_{-\frac{3}{2}}|<10^{-74},\dots,
|1-\mi{rat}_{\frac{21}{2}}|<10^{-37},\dots,
|1-\mi{rat}_{\frac{45}{2}}|<10^{-18},\dots
$$
This confirms, via the $\nu$-plane, the expected resurgence equations in the $\zeta$-plane,
namely\,:
$$
\Delta_{\zeta_3-\zeta_2}{\imia{Luu}}=\frac{2}{2\pi}\;{\imia{Lii}}
\quad ;\quad 
\Delta_{\zeta_1-\zeta_2}{\imia{Luu}}=-\frac{2}{2\pi}\;{\imia{Li}}
$$
and also, by mirror symmetry\,:
$$
\Delta_{\underline{\zeta_3}-\underline{\zeta_2}}{\underline{\imia{Luu}}}=
\frac{2}{2\pi}\;{\underline{\imia{Lii}}}
\quad ;\quad 
\Delta_{\underline{\zeta_1}-\underline{\zeta_2}}{\underline{\imia{Luu}}}=
-\frac{2}{2\pi}\;{\underline{\imia{Li}}}
$$
%%%%%%%%%%%%%%%%%%%%%%%%%%%%%%%%%%%%%%%%%%%%%%%%%%%%%%%%%%%%%%%%%%%%%%%%%%%%%
%%%%%%%%%%%%%%%%%%%%%%%%%%%%%%%%%%%%%%%%%%%%%%%%%%%%%%%%%%%%%%%%%%%%%%%%%%%%%
\subsubsection{\bf From $\mi{Le}$ to $\mi{Li}$ and $\mi{Lii}$
({\it exceptional} to {\it inner}).}
As in the preceding subsection,
we solve the linear system $\doS^{\bf{n},\bf{m}}_{e\,,\,i/ii}$ \,:
$$
\imi{le}_n\,=\,\frac{2}{2\pi}\,(\nu_3^{-n}-(-\nu_3)^{-n})\sum_{0\leq m<\bf{m}}
\; \tilde{riis}_{-\frac{1}{2}+m}\;n^{\frac{1}{2}-m}
\quad \quad\quad\big( n\,\in\, ]\bf{n}\!-\!2\bf{m}\,,\,\bf{n}]\big)
$$
with $\bf{m}$ effective equations (for $n$ odd) and $\bf{m}$ unknowns
$\tilde{riis}_{-\frac{1}{2}+m}$. Then we form\,:
\begin{eqnarray*}
\imi{riis}(\xi)&:=& \sum_{0\leq m < \bf{m}}
 \tilde{riis}_{-\frac{1}{2}+m}\;\frac{\xi^{-\frac{3}{2}+m}}{(-\frac{3}{2}+m)!}+(\dots)
\\[1.0 ex]
\imi{liis}(\xi)&:=& \imi{riis}\big(\log(1+\frac{\nu}{\nu_3})\big)
\;=\;\sum_{0\leq m < \bf{m}}
 \imi{liis}_{-\frac{3}{2}+m}\;{\nu^{-\frac{3}{2}+m}}+(\dots)
\end{eqnarray*}
and check that the ratios 
$\mi{rat}_{-\frac{3}{2}+m}:=\frac{\imi{liis}_{-\frac{3}{2}+m}}{\imi{lii}_{-\frac{3}{2}+m}}$
of homologous coefficients are indeed $\sim 1$. For the data corresponding to 
$\doS^{300,40}_{e\,,\,i/ii}$, we find this excellent level of accuracy\,:
$$
|1-\mi{rat}_{-\frac{3}{2}}|<10^{-73},\dots,
|1-\mi{rat}_{\frac{21}{2}}|<10^{-38},\dots,
|1-\mi{rat}_{\frac{45}{2}}|<10^{-17},\dots
$$
This confirms, via the $\nu$-plane, the expected resurgence equations in the $\zeta$-plane,
to wit\,:
$$
\Delta_{\zeta_3-\zeta_2}{\imia{Le}}=\frac{2}{2\pi}\;{\imia{Lii}}
\quad ;\quad 
\Delta_{\zeta_1-\zeta_2}{\imia{Le}}=-\frac{2}{2\pi}\;{\imia{Li}}
$$
An alternative method  is to check that $\imi{\ell e}\!\!(\nu)-\imi{\ell uu}\!\!(\nu)$
has a radius of convergence larger than $|\nu_1|=|\nu_3|$, which implies that $\imi{\ell e}$
and $\imi{\ell uu}$ have the same singularity at $\nu_1$ and $\nu_3$, namely
$\frac{2}{2\pi}\imia{Li} $ and $\frac{2}{2\pi}\imia{Lii} $: see \S 9.8.9. Here too,
the numerical accuracy is excellent.
%%%%%%%%%%%%%%%%%%%%%%%%%%%%%%%%%%%%%%%%%%%%%%%%%%%%%%%%%%%%%%%%%%%%%%%%%%%%%%%%%%%%%%%%%%%%
%%%%%%%%%%%%%%%%%%%%%%%%%%%%%%%%%%%%%%%%%%%%%%%%%%%%%%%%%%%%%%%%%%%%%%%%%%%%%%%%%%%%%%%%%%%%

%\end{document}
%%%%%%%%%%%%%%%%%%%%%%%%%%%%%%%%%%%%%%%%%%%%%%%%%%%%%%%%%%%%%%%%%%%%%%%%%%%%%%%%%%%%%%%%%%%%
%%%%%%%%%%%%%%%%%%%%%%%%%%%%%%%%%%%%%%%%%%%%%%%%%%%%%%%%%%%%%%%%%%%%%%%%%%%%%%%%%%%%%%%%%%%%
%%%%%%%%%%%%%%%%%%%%%%%%%%%%%%%%%%%%%%%%%%%%%%%%%%%%%%%%%%%%%%%%%%%%%%%%%%%%%%%%%%%%%%%%%%%%
%%%%%%%%%%%%%%%%%%%%%%%%%%%%%%%%%%%%%%%%%%%%%%%%%%%%%%%%%%%%%%%%%%%%%%%%%%%%%%%%%%%%%%%%%%%%

%

%%%%%%%%%%%%%%%%%%%%%%%%%%%%%%%%%%%%%%%%%%%%%%%%%%%%%%%%%%
%% Some resurgence properties of knot-related functions.
%%%%%%%%%%%%%%%%%%%%%%%%%%%%%%%%%%%%%%%%%%%%%%%%%%%%%%%%%%

%\documentclass[12pt,a4paper]{article}\input{SP_commands}\begin{document}
 
%%%%%%%%%%%%%%%%%%%%%%%%%%%%%%%%%%%%%%%%%%%%%%%%%%%%%%%%%%%%%%%%%%%%%%%%%%%%%%%%%%%%%%%%%%%%
%%%%%%%%%%%%%%%%%%%%%%%%%%%%%%%%%%%%%%%%%%%%%%%%%%%%%%%%%%%%%%%%%%%%%%%%%%%%%%%%%%%%%%%%%%%%
%%%%%%%%%%%%%%%%%%%%%%%%%%%%%%%%%%%%%%%%%%%%%%%%%%%%%%%%%%%%%%%%%%%%%%%%%%%%%%%%%%%%%%%%%%%%

\section{General tables.}
%%%%%%%%%%%%%%%%%%%%%%%%%%%%%%%%%%%%%%%%%%%%%%%%%%%%%%%%%%%%%%%%%%%%%%%%%%%%%%%%%%%%%%%%%%%%
%%%%%%%%%%%%%%%%%%%%%%%%%%%%%%%%%%%%%%%%%%%%%%%%%%%%%%%%%%%%%%%%%%%%%%%%%%%%%%%%%%%%%%%%%%%%
%%%%%%%%%%%%%%%%%%%%%%%%%%%%%%%%%%%%%%%%%%%%%%%%%%%%%%%%%%%%%%%%%%%%%%%%%%%%%%%%%%%%%%%%%%%%
%%%%%%%%%%%%%%%%%%%%%%%%%%%%%%%%%%%%%%%%%%%%%%%%%%%%%%%%%%%%%%%%%%%%%%%%%%%%%%%%%%%%%%%%%%%%
\subsection{Main formulas.} 
%%%%%%%%%%%%%%%%%%%%%%%%%%%%%%%%%%%%%%%%%%%%%%%%%%%%%%%%%%%%%%%%%%%%%%%%%%%%%%%%%%%%%%%%%%%%
%%%%%%%%%%%%%%%%%%%%%%%%%%%%%%%%%%%%%%%%%%%%%%%%%%%%%%%%%%%%%%%%%%%%%%%%%%%%%%%%%%%%%%%%%%%%
\subsubsection{Functional transforms.}
\[\begin{array}{llllll}
\mi{standard \;case}&& 
\beta(\tau):=\frac{1}{e^{\tau/2)}-e^{-\tau/2}}
&&
\beta^\dagger(\tau):=\frac{1}{e^{\tau/2)}-e^{-\tau/2}}-\frac{1}{\tau}
\\[1.ex]
\mi{free\;\beta\; case}&&
\beta(\tau):={\tau}^{-1}+\sum_{1\leq k}\beta_k\,\tau^k
&&
\beta^\dagger(\tau):=\sum_{1\leq k}\beta_k\,\tau^k
\end{array}\]
\dotfill
\\ 
\noindent
{\it mir}-transform\,:\;\; $\gbar:=1/g\mapsto \hbar:=1/h $\; with
\begin{equation}  \label{a207}
\frac{1}{\hbar(\nu)}=
\Big[\frac{1}{
\gbar(\nu)}\exp\Big(-\beta^\dagger\!\big(I\gbar(\nu)\partial_\nu\big)\gbar(\nu) \Big) 
 \Big]_{I=\partial_\nu^{-1}}
\end{equation}
\dotfill
\\
\noindent
{\it nir}-transform\,:\;\; $f \mapsto h $\; with
\begin{equation}  \label{a208}
h(\nu)=\frac{1}{2\pi i}\int_{c-i\infty}^{c+i\infty}e^{n\nu}\frac{dn}{n}
\int_0^{+\infty}\exp^{\#}
\!\Big(-\beta(\partial_\tau)f(\frac{\tau}{n}) \Big)\,d\tau
\end{equation} 
\dotfill
\\
\noindent  
{\it nir}-translocation\,:\;\; 
$f \mapsto \nabla h :=(
\mi{nir}-e^{-\eta\partial_\nu}\mi{nir}\;e^{\epsilon\partial_x})(h)$\;
 with
\begin{equation}  \label{a209}
\nabla h(\epsilon,\nu)=\frac{1}{2\pi i}\int_{c-i\infty}^{c+i\infty}e^{n\nu}\frac{dn}{n}
\int_0^{\epsilon\,
n}\exp_{\#}\!\Big(-\beta(\partial_\tau)f(\frac{\tau}{n})
\Big)\,d\tau
\end{equation} 
\dotfill
\\
\noindent
 {\it nur}-transform\,:\;\; $f \mapsto h $\; with
\begin{equation}  \label{a210}
h(\nu)=\frac{1}{2\pi i}\int_{c-i\infty}^{c+i\infty}e^{n\nu}\frac{dn}{n}
\sum_{\tau\in\,\frac{1}{2}+\doN}\exp^{\#}
\!\Big(-\beta(\partial_\tau)f(\frac{\tau}{n})
\Big)\,d\tau
\end{equation} 
\dotfill
\\
\noindent
{\it nur} in terms of {\it nir}\;:
\begin{equation}  \label{a211}
\mr{nur}(f)=\sum_{k\in\doZ}(-1)^k\,\mr{nir}(k\,2\pi i+f) 
\end{equation} 
\dotfill
\\
For the interpretation of $exp^{\#}, exp_{\#}$ see \S4.3.
%%%%%%%%%%%%%%%%%%%%%%%%%%%%%%%%%%%%%%%%%%%%%%%%%%%%%%%%%%%%%%%%%%%%%%%%%%%%%%%%%%%%%%%%%%%%
%%%%%%%%%%%%%%%%%%%%%%%%%%%%%%%%%%%%%%%%%%%%%%%%%%%%%%%%%%%%%%%%%%%%%%%%%%%%%%%%%%%%%%%%%%%%
\subsubsection{SP coefficients and SP series.}
Basic data\,: 
$\;\;F=\exp(-f)\;\;\;\;,\;\;\;\; 
\eta_F:=\int_0^1f(x)dx\;\;\;\;,\;\;\;\;
\omega_F=e^{-\eta_F}$
\[\begin{array}{llccllll}
\mi{asymptotic\; series}     &&\mi{funct.\; germs}
\\[1.ex]
\tilde{I\!g}_{_{F}}(n)=
\exp\Big(\!-\frac{1}{2}f(0)+\sum_{1\leq s \mi{odd}} \frac{\frak{b}_s}{n^s}f^{(s)}(0)\Big)
&&{I\!g}_{_{F}}(n)&\mi{ingress\;factor}
\\[1.ex]
\tilde{E\!g}_{_{F}}(n)=
\exp\Big(\!-\frac{1}{2}f(1)-\sum_{1\leq s \mi{odd}} \frac{\frak{b}_s}{n^s}f^{(s)}(1)\Big)
&&{E\!g}_{_{F}}(n)&\mi{egress\;factor}
\end{array}\]
\dotfill
\begin{eqnarray*}
\mi{``raw"}\hspace{24. ex} & \mi{``cleansed"}\hspace{14. ex}
\\[1.5 ex]
P_F(n):=\prod_{0\leq k\leq n}F(\frac{k}{n})
\hspace{8. ex} &\hspace{7.ex}
P^\#_F(n):= (\omega_F)^n=\frac{P_F(n)}{{I\!g}_{_{F}}(n){E\!g}_{_{F}}(n) }
\\[1. ex]
J_F(n):=\sum_{0\leq m<n}\prod_{0\leq k\leq m}F(\frac{k}{n})
\quad &
J^\#_F(n):=J_F(n)/{I\!g}_{_{F}}(n)
\\[1.ex]
j_F(\zeta):=\sum_{0\leq n}J_F(n)\,\zeta^n
\hspace{7.5 ex}  &
j^\#_F(\zeta):=\sum_{0\leq n}J^\#_F(n)\,\zeta^n
\end{eqnarray*}
%%%%%%%%%%%%%%%%%%%%%%%%%%%%%%%%%%%%%%%%%%%%%%%%%%%%%%%%%%%%%%%%%%%%%%%%%%%%%%%%%%%%%%%%%%%%
%%%%%%%%%%%%%%%%%%%%%%%%%%%%%%%%%%%%%%%%%%%%%%%%%%%%%%%%%%%%%%%%%%%%%%%%%%%%%%%%%%%%%%%%%%%%
\subsubsection{Parity relations.}
\begin{eqnarray*}
F^{\,\vdash}(x):=1/F(1\!-\!x) \hspace{7.ex}&\Longrightarrow& 
\\[1.5 ex]
1
=\tilde{I\!g}_{_{F}}(n)\,\tilde{E\!g}_{_{F^{\vdash}}}(n)
=\tilde{I\!g}_{_{F^{\vdash}}}(n)\,\tilde{E\!g}_{_{F}}(n)
&&
\\
J_{F^{\vdash}}(n)=J_F(n)/P_F(n) 
\hspace{4 ex}&\mi{and}&\hspace{3 ex}
J^\#_{F^{\vdash}}(n)=J^\#_F(n)/P^\#_F(n)
\\
j_{F^{\vdash}}(\zeta)\not=j_F(\zeta/\omega_F) 
\hspace{8 ex}&\mi{but}&\hspace{3.2 ex}
j^\#_{F^{\vdash}}(\zeta)=j^\#_F(\zeta/\omega_F) 
\end{eqnarray*}
\dotfill
\begin{eqnarray*}
F^{\,\bot}(x):=1/F(-x)\;\;\;,\;\;\;f^{\,\bot}(x):=-f(-x)  \hspace{7.ex}
&&\Longrightarrow 
\\[1.5 ex]
\mr{nur}(f^\bot)(\nu)=-\mr{nur}(f)(\nu)
\hspace{23.ex}
&& (\mi{tangency}\;\, \kappa=0)
\\
\mr{nir}(f^\bot)(\nu)=-\mr{nir}(f)(\nu)
\hspace{24.ex} 
&&(\mi{tangency}\;\, \kappa=0)
\\
\mr{nir}(f^\bot)\;\; \mi{and}\;\;\mr{nir}(f)\;\;\mi{unrelated}
\hspace{18.ex}
&& (\mi{tangency}\;\, \kappa\;\mi{even}\geq 2)
\\
\mr{nir}(f^\bot)(\nu)=-\mr{nur}(f)(\epsilon_\kappa\nu)\;\;\; \mi{with}\;\;\;
\epsilon_\kappa^{\frac{1}{\kappa+1}}=-1
\hspace{3.3 ex}
&& (\mi{tangency}\;\, \kappa\;\mi{odd}\geq 1)
\\
\Rightarrow h^{\bot}_{\frac{k}{\kappa+1}}=(-1)^{k-1}h_{\frac{k}{\kappa+1}}
\;\;\mi{with}\;\; \;:(f,f^\bot)\stackrel{\mr{nir}}{\mapsto} (h,h^\bot)
\hspace{2.ex}
&& (\mi{tangency}\;\, \kappa\;\mi{odd}\geq 1)
\end{eqnarray*}
\dotfill
%%%%%%%%%%%%%%%%%%%%%%%%%%%%%%%%%%%%%%%%%%%%%%%%%%%%%%%%%%%%%%%%%%%%%%%%%%%%%%%%%%%%%%%%%%%%
%%%%%%%%%%%%%%%%%%%%%%%%%%%%%%%%%%%%%%%%%%%%%%%%%%%%%%%%%%%%%%%%%%%%%%%%%%%%%%%%%%%%%%%%%%%%
%%%%%%%%%%%%%%%%%%%%%%%%%%%%%%%%%%%%%%%%%%%%%%%%%%%%%%%%%%%%%%%%%%%%%%%%%%%%%%%%%%%%%%%%%%%%
%%%%%%%%%%%%%%%%%%%%%%%%%%%%%%%%%%%%%%%%%%%%%%%%%%%%%%%%%%%%%%%%%%%%%%%%%%%%%%%%%%%%%%%%%%%%
%%%%%%%%%%%%%%%%%%%%%%%%%%%%%%%%%%%%%%%%%%%%%%%%%%%%%%%%%%%%%%%%%%%%%%%%%%%%%%%%%%%%%%%%%%%%
%%%%%%%%%%%%%%%%%%%%%%%%%%%%%%%%%%%%%%%%%%%%%%%%%%%%%%%%%%%%%%%%%%%%%%%%%%%%%%%%%%%%%%%%%%%%
\subsection{The {\it Mir} mould.}
\subsubsection{Layered form.}
$\mathbf{length:}\;\; r=1\;,\;\mathbf{order:}\;\; d=1 \;,\;\mathbf{factor:}\;\;c_{1,1}=1/24 $
\[\begin{array}{lrrrrl}
\Mir[0,1] &=& 1\;c_{3,2} &=&  1/24
 &\hspace{46.ex}
\end{array}\]
$\mathbf{length:}\;\; r=3\;,\;\mathbf{order:}\;\; d=2 \;,\;\mathbf{factor:}\;\;c_{3,2}= 1/1152 $
\[\begin{array}{lrrrrl}
\Mir[1,2,0,0] &=& 1\;c_{3,2}&=& 1/1152
 &\hspace{36.ex}
\end{array}\]
$\mathbf{length:}\;\; r=3\;,\;\mathbf{order:}\;\; d=3 \;,\;\mathbf{factor:}\;\;c_{3,3}=7/5760 $
\[\begin{array}{lrrrrl}
\Mir[0,3,0,0] &=& -1\;c_{3,3}&=& -7/5760\\ 
\Mir[1,1,1,0] &=& -4\;c_{3,3}&=&  -7/1440\\ 
\Mir[2,0,0,1] &=& -1\;c_{3,3}&=&  -7/5760

 &\hspace{36.ex}
\end{array}\]
$\mathbf{length:}\;\; r=5\;,\;\mathbf{order:}\;\; d=3 \;,\;\mathbf{factor:}\;\;c_{5,3}= 1/82944 $
\[\begin{array}{lrrrrl}
\Mir[2, 3, 0, 0, 0, 0] &=& 1\;c_{5,3}&=&  1/82944
 &\hspace{36.ex}
\end{array}\]
$\mathbf{length:}\;\; r=5\;,\;\mathbf{order:}\;\; d=4 \;,\;\mathbf{factor:}\;\;c_{5,4}=7/138240 $
\[\begin{array}{lrrrrl}
\Mir[1,4,0,0,0,0] &=&  -1\;c_{5,4}&=& -7/138240\\ 
\Mir[2,2,1,0,0,0] &=&  -4\;c_{5,4}&=& -7/34560\\ 
\Mir[3,1,0,1,0,0] &=&  -1\;c_{5,4}&=& -7/138240
 &\hspace{36.ex}
\end{array}\]
$\mathbf{length:}\;\; r=5\;,\;\mathbf{order:}\;\; d=5 \;,\;\mathbf{factor:}\;\;c_{5,5}=31/967680
$
\[\begin{array}{lrrrrl}
\Mir[0,5,0,0,0,0] &=& 1\;c_{5,5} &=& 31/967680 \\
\Mir[1,3,1,0,0,0] &=& 26\;c_{5,5} &=& 403/483840 \\
\Mir[2,1,2,0,0,0] &=& 34\;c_{5,5} &=& 527/483840\\ 
\Mir[2,2,0,1,0,0] &=& 32\;c_{5,5} &=& 31/30240\\
\Mir[3,0,1,1,0,0] &=& 15\;c_{5,5} &=& 31/64512\\ 
\Mir[3,1,0,0,1,0] &=& 11\;c_{5,5} &=& 341/967680\\ 
\Mir[4,0,0,0,0,1] &=& 1\;c_{5,5} &=& 31/967680
 &\hspace{36.ex}
\end{array}\]

%%%%%%%%%%%%%%%%%%%%%%%%%%%%%%%%%%%%%%%%%%%%%%%%%%%%%%%%%%%%%%%%%%%%%%%%%%%%%%%%%%%%%%%%%%%%
%%%%%%%%%%%%%%%%%%%%%%%%%%%%%%%%%%%%%%%%%%%%%%%%%%%%%%%%%%%%%%%%%%%%%%%%%%%%%%%%%%%%%%%%%%%%
\subsubsection{Compact form.}
$\mathbf{length:}\;\; r=1\;\;, \;\;\mathbf{gcd:}\;\;d_3=24 $
\[\begin{array}{lrrrrl}
\mrMir[1] &=& 1/d_1&=&  1/24\hspace{56.ex}
\end{array}\]
$\mathbf{length:}\;\; r=3\;\;, \;\;\mathbf{gcd:}\;\;d_3=5760 $
\[\begin{array}{lrrrrl}
\mrMir[1, 2, 0] &=& -2/d_3 &=& -1/2880&\\
\mrMir[2, 0, 1] &=& -7/d_3&=&  -7/5760&\hspace{36.ex}
\end{array}\]
$\mathbf{length:}\;\; r=5\;\;, \;\;\mathbf{gcd:}\;\;d_5= 2903040 $
\[\begin{array}{lrrrrl}
\mrMir[1, 4, 0, 0, 0] &=& 16/d_5 &=& 1/181440\\ 
\mrMir[2, 2, 1, 0, 0] &=& 540/d_5 &=& 1/ 5376\\
\mrMir[3, 0, 2, 0, 0] &=& 372/d_5 &=& 31/ 241920\\ 
\mrMir[3, 1, 0, 1, 0] &=& 504/d_5 &=& 1/ 5760\\
\mrMir[4, 0, 0, 0, 1] &=& 93/d_5 &=& 31/ 967680 &\hspace{36.ex}
\end{array}\]
$\mathbf{length:}\;\; r=7\;\;, \;\;\mathbf{gcd:}\;\;d_7= 1393459200 $
\[\begin{array}{lrrrrl}
\mrMir[1, 6, 0, 0, 0, 0, 0] &=& -144/d_7 &=& -1/ 9676800 &\\
\mrMir[2, 4, 1, 0, 0, 0, 0] &=&-28824/d_7 &=&  -1201/ 58060800 &\\
\mrMir[3, 2, 2, 0, 0, 0, 0] &=&-141576/d_7 &=&  -5899/ 58060800 &\\
\mrMir[4, 0, 3, 0, 0, 0, 0] &=&-38862/d_7 &=&  -2159/ 77414400 &\\
\mrMir[3, 3, 0, 1, 0, 0, 0] &=&-88928/d_7 &=&  -397/ 6220800 &\\
\mrMir[4, 1, 1, 1, 0, 0, 0] &=&-186264/d_7 &=&  -2587/ 19353600 &\\
\mrMir[5, 0, 0, 2, 0, 0, 0] &=&-16116/d_7 &=&  -1343/ 116121600 &\\
\mrMir[4, 2, 0, 0, 1, 0, 0] &=&-67878/d_7 &=&  -419/ 8601600 &\\
\mrMir[5, 0, 1, 0, 1, 0, 0] &=&-29718/d_7 &=&  -1651/ 77414400 &\\
\mrMir[5, 1, 0, 0, 0, 1, 0] &=&-16428/d_7 &=&  -1369/ 116121600 &\\
\mrMir[6, 0, 0, 0, 0, 0, 1] &=&-1143/d_7 &=&  -127/ 154828800 & 
\hspace{36.ex}
\end{array}\]
\nopagebreak
$\mathbf{length:}\;\; r=9\;\;, \;\;\mathbf{gcd:}\;\;d_9= 367873228800 $
\[\begin{array}{lrrrrl}
\mrMir[1, 8, 0, 0, 0, 0, 0, 0, 0]  &=&  768/ d_9  &=&   1/ 479001600  &\\
\mrMir[2, 6, 1, 0, 0, 0, 0, 0, 0]  &=&  789504/ d_9  &=&   257/ 119750400  &\\
\mrMir[3, 4, 2, 0, 0, 0, 0, 0, 0]  &=&  13702656/ d_9  &=&   811/ 21772800  &\\
\mrMir[4, 2, 3, 0, 0, 0, 0, 0, 0]  &=&  26034672/ d_9  &=&   542389/ 7664025600  &\\
\mrMir[5, 0, 4, 0, 0, 0, 0, 0, 0]  &=&  3801840/ d_9  &=&   2263/ 218972160  &\\
\mrMir[3, 5, 0, 1, 0, 0, 0, 0, 0]  &=&  6324224/ d_9  &=&   193/ 11226600  &\\
\mrMir[4, 3, 1, 1, 0, 0, 0, 0, 0]  &=&  52597760/ d_9  &=&   10273/ 71850240  &\\
\mrMir[5, 1, 2, 1, 0, 0, 0, 0, 0]  &=&  40989024/ d_9  &=&   47441/ 425779200  &\\
\mrMir[5, 2, 0, 2, 0, 0, 0, 0, 0]  &=&  18164736/ d_9  &=&   73/ 1478400  &\\
\mrMir[6, 0, 1, 2, 0, 0, 0, 0, 0]  &=&  6350064/ d_9  &=&   18899/ 1094860800  &\\
\mrMir[4, 4, 0, 0, 1, 0, 0, 0, 0]  &=&  11628928/ d_9  &=&   90851/ 2874009600  &\\
\mrMir[5, 2, 1, 0, 1, 0, 0, 0, 0]  &=&  33372912/ d_9  &=&   695269/ 7664025600  &\\
\mrMir[6, 0, 2, 0, 1, 0, 0, 0, 0]  &=&  5886720/ d_9  &=&   73/ 4561920  &\\
\mrMir[6, 1, 0, 1, 1, 0, 0, 0, 0]  &=&  9462768/ d_9  &=&   28163/ 1094860800  &\\
\mrMir[7, 0, 0, 0, 2, 0, 0, 0, 0]  &=&  429240/ d_9  &=&   511/ 437944320  &\\
\mrMir[5, 3, 0, 0, 0, 1, 0, 0, 0]  &=&  7436800/ d_9  &=&   83/ 4105728  &\\
\mrMir[6, 1, 1, 0, 0, 1, 0, 0, 0]  &=&  7391376/ d_9  &=&   51329/ 2554675200  &\\
\mrMir[7, 0, 0, 1, 0, 1, 0, 0, 0]  &=&  736848/ d_9  &=&   731/ 364953600  &\\
\mrMir[6, 2, 0, 0, 0, 0, 1, 0, 0]  &=&  1941144/ d_9  &=&   80881/ 15328051200  &\\
\mrMir[7, 0, 1, 0, 0, 0, 1, 0, 0]  &=&  490560/ d_9  &=&   73/ 54743040  &\\
\mrMir[7, 1, 0, 0, 0, 0, 0, 1, 0]  &=&  209712/ d_9  &=&   4369/ 7664025600  &\\
\mrMir[8, 0, 0, 0, 0, 0, 0, 0, 1]  &=&  7665/ d_9  &=&   73/ 3503554560 &
\hspace{36.ex}
\end{array}\]

%%%%%%%%%%%%%%%%%%%%%%%%%%%%%%%%%%%%%%%%%%%%%%%%%%%%%%%%%%%%%%%%%%%%%%%%%%%%%%%%%%%%%%%%%%%%
%%%%%%%%%%%%%%%%%%%%%%%%%%%%%%%%%%%%%%%%%%%%%%%%%%%%%%%%%%%%%%%%%%%%%%%%%%%%%%%%%%%%%%%%%%%%
%%%%%%%%%%%%%%%%%%%%%%%%%%%%%%%%%%%%%%%%%%%%%%%%%%%%%%%%%%%%%%%%%%%%%%%%%%%%%%%%%%%%%%%%%%%%
%%%%%%%%%%%%%%%%%%%%%%%%%%%%%%%%%%%%%%%%%%%%%%%%%%%%%%%%%%%%%%%%%%%%%%%%%%%%%%%%%%%%%%%%%%%%
\subsection{The {\it mir} transform\,: from $\gbar$  to  $\hbar$.}
%%%%%%%%%%%%%%%%%%%%%%%%%%%%%%%%%%%%%%%%%%%%%%%%%%%%%%%%%%%%%%%%%%%%%%%%%%%%%%%%%%%%%%%%%%%%

\subsubsection{Tangency 0, ramification 1.} 
Recall that $\gbar=1/g $ and  $ \hbar=1/h$ .
\[\begin{array}{lll}
h_{1}-g_{1}&\!=\!&
\frac{1}{24}{\gbar}_{1}
\\[1.5 ex]
h_{2}-g_{2}&\!=\!&
\frac{1}{24}{\gbar}_{2}+ \frac{1}{2304}{\gbar}_{0}{\gbar}_{1}^2
\\[1.5 ex]
h_{{3}}-g_{{3}}&\!=\!&
\frac{1}{24}\,{\gbar}_{{3}}-{\frac {1}{17280}}\,{{\gbar}_{{1}}}^{3}-{\frac {1}{
960}}\,{\gbar}_{{0}}{\gbar}_{{1}}{\gbar}_{{2}}-{\frac {7}{5760}}\,{{
\gbar}_{{0}}}^{2}{\gbar}_{{3}}+{\frac {1}{497664}}\,{{\gbar}_{{0}}}^{2
}{{\gbar}_{{1}}}^{3}
\\[1.5 ex]
h_{{4}}-g_{{4}}&\!=\!&
{\frac{1}{24}}\,{\gbar}_{{4}}-{\frac {1}{1920}}\,{\gbar}_{{0}}{{\gbar}_{{2}}}^{2
}-{\frac {1}{2880}}\,{{\gbar}_{{1}}}^{2}{\gbar}_{{2}}
-{\frac {1}{720}}\,{
\gbar}_{{0}}{\gbar}_{{1}}{\gbar}_{{3}}-{\frac {7}{5760}}\,{{\gbar}_
{{0}}}^{2}{\gbar}_{{4}}
\\[1. ex]&&\!\!\!
-{\frac {23}{1658880}}\,{{\gbar}_{{0}}}^{2}
{{\gbar}_{{1}}}^{2}{\gbar}_{{2}}
-{\frac {11}{9953280}}\,{\gbar}_{{0}}{{\gbar}_{{1}}}^{4}
-{\frac {7}{552960}}\,{{\gbar}_{{0}}}^{3}{\gbar}_{{1}}{\gbar}_{{3}}
+{\frac {1}{191102976}}\,{{\gbar}_{{0}}}^{3}{{\gbar}_{{1}}}^{4}
\\[1.5 ex]
h_{{5}}-g_{{5}}&\!=\!&
{\frac{1}{24}}\,{\gbar}_{{5}}
-{\frac {11}{28800}}\,{\gbar}_{{1}}{{\gbar}_{{2}}}^{2}
-{\frac {23}{14400}}\,{\gbar}_{{0}}{\gbar}_{{1}}{\gbar}_{{4}}
-{\frac {7}{14400}}\,{{\gbar}_{{1}}}^{2}{\gbar}_{{3}}
-{\frac {17}{14400}}\,
{\gbar}_{{0}}{\gbar}_{{2}}{\gbar}_{{3}}
\\[1.ex] &&\!\!\!
-{\frac {7}{5760}}\,{{\gbar}
_{{0}}}^{2}{\gbar}_{{5}}
+{\frac {143}{21772800}}\,{\gbar}_{{0}}{{\gbar}_{{1}}}^{3}{\gbar}_{{2}}
+{\frac {61}{2419200}}\,{{\gbar}_{{0}}}^{2}{\gbar}_{{1}}{{\gbar}_{{2}}}^{2}
+{\frac {73}{1209600}}\,{{\gbar}_{{0}}}^
{3}{\gbar}_{{1}}{\gbar}_{{4}}
\\[1.ex] &&\!\!\!
+{\frac {19}{537600}}\,{{\gbar}_{{0}}}^{2}{{\gbar}_{{1}}}^{2}{\gbar}_{{3}}
+{\frac {13}{302400}}\,{{\gbar}_{{0}}}^{3}{\gbar}_{{2}}{\gbar}_{{3}}
+{\frac {31}{967680}}\,{{\gbar}_{{0}}}^
{4}{\gbar}_{{5}}+{\frac {1}{21772800}}\,{{\gbar}_{{1}}}^{5}
\\[1.ex] &&\!\!\!
-{\frac {37}{597196800}}\,{{\gbar}_{{0}}}^{3}{{\gbar}_{{1}}}^{3}{\gbar}_{{2}}
-{\frac {1}{176947200}}\,{{\gbar}_{{0}}}^{2}{{\gbar}_{{1}}}^{5}
-{\frac {7}{132710400}}\,{{\gbar}_{{0}}}^{4}{{\gbar}_{{1}}}^{2}{\gbar}_{{3}}
\\[1.ex] &&\!\!\!
+{\frac {1}{114661785600}}\,{{\gbar}_{{0}}}^{4}{{\gbar}_{{1}}}^{5}
\\[1.5 ex]
h_{{6}}-g_{{6}}&\!=\!&

+{\frac {1}{24}}\,{\gbar}_{{6}}

-{\frac {13}{14400}}\,{\gbar}_{{1}}{\gbar}_{{2}}{\gbar}_{{3}}
-{\frac {17}{28800}}\,{{\gbar}_{{1}}}^{2}{\gbar}_{{4}}
-{\frac {19}{14400}}\,{\gbar}_{{0}}{\gbar}_{{2}}{\gbar}_{{4}}
-{\frac {1}{576}}\,{\gbar}_{{0}}{\gbar}_{{1}}{\gbar}_{{5}}
\\[1.ex] &&\!\!\!
-{\frac {7}{5760}}\,{{\gbar}_{{0}}}^{2}{\gbar}_{{6}}
-{\frac {11}{86400}}\,{{\gbar}_{{2}}}^{3}
-{\frac {17}{28800}}\,{\gbar}_{{0}}{{\gbar}_{{3}}}^{2}
 
+{\frac {17}{241920}}\,{{\gbar}_{{0}}}^{2}{\gbar}_{{1}}{\gbar}_{{2}}{\gbar}_{{3}}
\\[1.ex] &&\!\!\!
+{\frac {163}{10886400}}\,{\gbar}_{{0}}{{\gbar}_{{1}}}^{2}{{\gbar}_{{2}}}^{2}
+{\frac {41}{2721600}}\,{\gbar}_{{0}}{{\gbar}_{{1}}}^{3}{\gbar}_{{3}}
+{\frac {13}{241920}}\,{{\gbar}_{{0}}}^{2}{{\gbar}_{{1}}}^{2}{\gbar}_{{4}}
\\[1.ex] &&\!\!\!
+{\frac {59}{1209600}}\,{{\gbar}_{{0}}}^{3}{\gbar}_{{2}}{\gbar}_{{4}}
+{\frac {17}{14515200}}\,{\gbar}_{{2}}{{\gbar}_{{1}}}^{4}
+{\frac {13}{181440}}\,{{\gbar}_{{0}}}^{3}{\gbar}_{{1}}{\gbar}_{{5}}
+{\frac {31}{967680}}\,{{\gbar}_{{0}}}^{4}{\gbar}_{{6}}
\\[1.ex] &&\!\!\!
+{\frac {61}{7257600}}\,{{\gbar}_{{0}}}^{2}{{\gbar}_{{2}}}^{3}
+{\frac {13}{604800}}\,{{\gbar}_{{0}}}^{3}{{\gbar}_{{3}}}^{2}

+{\frac {461}{1161216000}}\,{{\gbar}_{{0}}}^{4}{\gbar}_{{1}}{\gbar}_{{2}}{\gbar}_{{3}}
\\[1.ex] &&\!\!\!
+{\frac {5347}{20901888000}}\,{{\gbar}_{{0}}}^{3}{{\gbar}_{{1}}}^{2}{{\gbar}_{{2}}}^{2}
+{\frac {25259}{83607552000}}\,{{\gbar}_{{0}}}^{3}{{\gbar}_{{1}}}^{3}{\gbar}_{{3}}
+{\frac {211}{464486400}}\,{{\gbar}_{{0}}}^{4}{{\gbar}_{{1}}}^{2}{\gbar}_{{4}}
\\[1.ex] &&\!\!\!
+{\frac {12109}{167215104000}}\,{{\gbar}_{{0}}}^{2}{{\gbar}_{{1}}}^{4}{\gbar}_{{2}}
+{\frac {31}{139345920}}\,{{\gbar}_{{0}}}^{5}{\gbar}_{{1}}{\gbar}_{{5}}
+{\frac {37}{37158912000}}\,{\gbar}_{{0}}{{\gbar}_{{1}}}^{6}
\\[1.ex] &&\!\!\!
+{\frac {49}{1327104000}}\,{{\gbar}_{{0}}}^{5}{{\gbar}_{{3}}}^{2}

-{\frac {17}{114661785600}}\,{{\gbar}_{{0}}}^{4}{{\gbar}_{{1}}}^{4}{\gbar}_{{2}}
-{\frac {1}{68797071360}}\,{{\gbar}_{{0}}}^{3}{{\gbar}_{{1}}}^{6}
\\[1.ex] &&\!\!\!
-{\frac {7}{57330892800}}\,{{\gbar}_{{0}}}^{5}{{\gbar}_{{1}}}^{3}{\gbar}_{{3}}

+{\frac {1}{99067782758400}}\,{{\gbar}_{{0}}}^{5}{{\gbar}_{{1}}}^{6}
\end{array}\]

%%%%%%%%%%%%%%%%%%%%%%%%%%%%%%%%%%%%%%%%%%%%%%%%%%%%%%%%%%%%%%%%%%%%%%%%%%%%%%%%%%%%%%%%%%%%
%%%%%%%%%%%%%%%%%%%%%%%%%%%%%%%%%%%%%%%%%%%%%%%%%%%%%%%%%%%%%%%%%%%%%%%%%%%%%%%%%%%%%%%%%%%%
\subsubsection{Tangency 1, ramification 2.}
\[\begin{array}{cll}
h_{1/2}-g_{1/2}&\!=\!&
{\frac{1}{24}}\,{\gbar}_{{1/2}}

\\[1.5 ex]
h_{1}-g_{1}&\!=\!&
{\frac{1}{24}}\,{\gbar}_{{1}}

\\[1.5 ex]
h_{3/2}-g_{3/2}&\!=\!&
{\frac{1}{24}}\,{\gbar}_{{3/2}}
+{\frac {1}{3456}}\,{{\gbar}_{{1/2}}}^{3}

\\[1.5 ex]
\ell_{2}-g_{2}&\!=\!&
{\frac{1}{24}}\,{\gbar}_{{2}}
+{\frac {5}{9216}}\,{{\gbar}_{{1/2}}}^{2}{\gbar}_{{1}}

\\[1.5 ex]
h_{5/2}-g_{5/2}&\!=\!&
{\frac{1}{24}}\,{\gbar}_{{5/2}}
-{\frac {1}{9600}}\,{\gbar}_{{1/2}}{{\gbar}_{{1}}}^{2}
-{\frac {7}{28800}}\,{{\gbar}_{{1/2}}}^{2}{\gbar}_{{3/2}}
+{\frac {1}{1244160}}\,{{\gbar}_{{1/2}}}^{5}

\\[1.5 ex]
h_{3}-g_{3}&\!=\!&
{\frac{1}{24}}\,{\gbar}_{{3}}

-{\frac {1}{17280}}\,{{\gbar}_{{1}}}^{3}
-{\frac {1}{1920}}\,{\gbar}_{{1}}{\gbar}_{{3/2}}{\gbar}_{{1/2}}
-{\frac {1}{2304}}\,{{\gbar}_{{1/2}}}^{2}{\gbar}_{{2}}

+{\frac {1}{497664}}\,{{\gbar}_{{1/2}}}^{4}{\gbar}_{{1}}

\\[1.5 ex]
h_{7/2}-g_{7/2}&\!=\!&
{\frac{1}{24}}\,{\gbar}_{{7/2}}

-{\frac {53}{201600}}\,{{\gbar}_{{1}}}^{2}{\gbar}_{{3/2}}
-{\frac {13}{40320}}\,{\gbar}_{{1/2}}{{\gbar}_{{3/2}}}^{2}
-{\frac {11}{14400}}\,{\gbar}_{{1/2}}{\gbar}_{{1}}{\gbar}_{{2}}
\\[1. ex]&&\!\!\!
-{\frac {113}{201600}}\,{{\gbar}_{{1/2}}}^{2}{\gbar}_{{5/2}}
-{\frac {11}{21772800}}\,{{\gbar}_{{1/2}}}^{3}{{\gbar}_{{1}}}^{2}
-{\frac {113}{43545600}}\,{{\gbar}_{{1/2}}}^{4}{\gbar}_{{3/2}}
\\[1. ex] &&\!\!\!
+{\frac {1}{836075520}}\,{{\gbar}_{{1/2}}}^{7}

\\[1.5 ex]
h_{4}-g_{4}&\!=\!&
{\frac{1}{24}}\,{\gbar}_{{4}}
-{\frac {1}{1280}}\,{\gbar}_{{1/2}}{\gbar}_{{3/2}}{\gbar}_{{2}}
-{\frac {7}{23040}}\,{\gbar}_{{1}}{{\gbar}_{{3/2}}}^{2}
-{\frac {11}{11520}}\,{\gbar}_{{1/2}}{\gbar}_{{1}}{\gbar}_{{5/2}}
\\[1. ex] &&\!\!\!
-{\frac {1}{2880}}\,{{\gbar}_{{1}}}^{2}{\gbar}_{{2}}
-{\frac {1}{1536}}\,{{\gbar}_{{1/2}}}^{2}{\gbar}_{{3}}
-{\frac {229}{159252480}}\,{{\gbar}_{{1/2}}}^{2}{{\gbar}_{{1}}}^{3}
\\[1. ex] &&\!\!\!
-{\frac {229}{39813120}}\,{{\gbar}_{{1/2}}}^{3}{\gbar}_{{1}}{\gbar}_{{3/2}}
-{\frac {49}{15925248}}\,{{\gbar}_{{1/2}}}^{4}{\gbar}_{{2}}
+{\frac {11}{3057647616}}\,{{\gbar}_{{1/2}}}^{6}{\gbar}_{{1}}

\\[1.5 ex]
h_{9/2}-g_{9/2}&\!=\!&
{\frac{1}{24}}\,{\gbar}_{{9/2}}
-{\frac {43}{60480}}\,{\gbar}_{{1}}{\gbar}_{{3/2}}{\gbar}_{{2}}
-{\frac {67}{60480}}\,{\gbar}_{{1/2}}{\gbar}_{{1}}{\gbar}_{{3}}
-{\frac {17}{40320}}\,{\gbar}_{{1/2}}{{\gbar}_{{2}}}^{2}
\\[1. ex] &&\!\!\!
-{\frac {13}{120960}}\,{{\gbar}_{{3/2}}}^{3}
-{\frac {29}{40320}}\,{{\gbar}_{{1/2}}}^{2}{\gbar}_{{7/2}}
-{\frac {17}{40320}}\,{{\gbar}_{{1}}}^{2}{\gbar}_{{5/2}}
-{\frac {11}{12096}}\,{\gbar}_{{1/2}}{\gbar}_{{3/2}}{\gbar}_{{5/2}}
\\[1. ex] &&\!\!\!
+{\frac {907}{457228800}}\,{{\gbar}_{{1/2}}}^{2}{{\gbar}_{{1}}}^{2}{\gbar}_{{3/2}}
+{\frac {757}{228614400}}\,{{\gbar}_{{1/2}}}^{3}{\gbar}_{{1}}{\gbar}_{{2}}
+{\frac {79}{914457600}}\,{\gbar}_{{1/2}}{{\gbar}_{{1}}}^{4}
\\[1. ex] &&\!\!\!
+{\frac {659}{457228800}}\,{{\gbar}_{{1/2}}}^{3}{{\gbar}_{{3/2}}}^{2}
+{\frac {1531}{914457600}}\,{{\gbar}_{{1/2}}}^{4}{\gbar}_{{5/2}}
-{\frac {13}{12541132800}}\,{{\gbar}_{{1/2}}}^{5}{{\gbar}_{{1}}}^{2}
\\[1. ex] &&\!\!\!
-{\frac {29}{4180377600}}\,{{\gbar}_{{1/2}}}^{6}{\gbar}_{{3/2}}
+{\frac {1}{902961561600}}\,{{\gbar}_{{1/2}}}^{9}

\\[1.5 ex]
h_{5}-g_{5}&\!=\!&
{\frac{1}{24}}\,{\gbar}_{{5}}
-{\frac {11}{28800}}\,{\gbar}_{{1}}{{\gbar}_{{2}}}^{2}
-{\frac {7}{14400}}\,{{\gbar}_{{1}}}^{2}{\gbar}_{{3}}
-{\frac {47}{57600}}\,{\gbar}_{{1}}{\gbar}_{{3/2}}{\gbar}_{{5/2}}
\\[1. ex] &&\!\!\!
-{\frac {71}{57600}}\,{\gbar}_{{1/2}}{\gbar}_{{1}}{\gbar}_{{7/2}}
-{\frac {53}{57600}}\,{\gbar}_{{2}}{\gbar}_{{1/2}}{\gbar}_{{5/2}}
-{\frac {59}{57600}}\,{\gbar}_{{1/2}}{\gbar}_{{3/2}}{\gbar}_{{3}}
\\[1. ex] &&\!\!\!
-{\frac {41}{115200}}\,{{\gbar}_{{3/2}}}^{2}{\gbar}_{{2}}
-{\frac {89}{115200}}\,{{\gbar}_{{1/2}}}^{2}{\gbar}_{{4}}
+{\frac {1}{21772800}}\,{{\gbar}_{{1}}}^{5}
\\[1. ex] &&\!\!\!
+{\frac {383}{43545600}}\,{{\gbar}_{{1/2}}}^{3}{\gbar}_{{1}}{\gbar}_{{5/2}}
+{\frac {797}{116121600}}\,{{\gbar}_{{1/2}}}^{2}{{\gbar}_{{1}}}^{2}{\gbar}_{{2}}
\\[1. ex] &&\!\!\!
+{\frac {361}{174182400}}\,{\gbar}_{{1/2}}{{\gbar}_{{1}}}^{3}{\gbar}_{{3/2}}
+{\frac {41}{6220800}}\,{{\gbar}_{{1/2}}}^{3}{\gbar}_{{3/2}}{\gbar}_{{2}}
\\[1. ex] &&\!\!\!
+{\frac {17}{3225600}}\,{{\gbar}_{{1/2}}}^{2}{\gbar}_{{1}}{{\gbar}_{{3/2}}}^{2}
+{\frac {19}{4976640}}\,{{\gbar}_{{1/2}}}^{4}{\gbar}_{{3}}
-{\frac {1}{188743680}}\,{{\gbar}_{{1/2}}}^{4}{{\gbar}_{{1}}}^{3}
\\[1. ex] &&\!\!\!
-{\frac {47}{2548039680}}\,{{\gbar}_{{1/2}}}^{5}{\gbar}_{{1}}{\gbar}_{{3/2}}
-{\frac {107}{15288238080}}\,{{\gbar}_{{1/2}}}^{6}{\gbar}_{{2}}
+{\frac {7}{1834588569600}}\,{{\gbar}_{{1/2}}}^{8}{\gbar}_{{1}}
\end{array}\]
%%%%%%%%%%%%%%%%%%%%%%%%%%%%%%%%%%%%%%%%%%%%%%%%%%%%%%%%%%%%%%%%%%%%%%%%%%%%%%%%%%%%%%%%%%%%
%%%%%%%%%%%%%%%%%%%%%%%%%%%%%%%%%%%%%%%%%%%%%%%%%%%%%%%%%%%%%%%%%%%%%%%%%%%%%%%%%%%%%%%%%%%%
\subsubsection{Tangency 2, ramification 3.}
\[\begin{array}{cll}
h_{{2/3}}-g_{{2/3}}\!&\!\!=\!\!&\!
{\frac {1}{24}}\,{\gbar}_{{2/3}}
\\[1.5 ex]
h_{{1}}-g_{{1}}\!&\!\!=\!\!&\!
{\frac {1}{24}}\,{\gbar}_{{1}}
\\[1.5 ex]
h_{{4/3}}-g_{{4/3}}\!&\!\!=\!\!&\!
{\frac {1}{24}}\,{\gbar}_{{4/3}}
\\[1.5 ex]
\ell_{{5/3}}-g_{{5/3}}\!&\!\!=\!\!&\!
{\frac {1}{24}}\,{\gbar}_{{5/3}}
\\[1.5 ex]
h_{{2}}-g_{{2}}\!&\!\!=\!\!&\!
{\frac {1}{24}}\,{\gbar}_{{2}} 
+{\frac {1}{5184}}\,{{\gbar}_{{2/3}}}^{3}

\\[1.5 ex]
h_{{7/3}}-g_{{7/3}}\!&\!\!=\!\!&\!
{\frac {1}{24}}\,{\gbar}_{{7/3}}
-{\frac {1}{40320}}\,{{\gbar}_{{2/3}}}^{2}{\gbar}_{{1}}

\\[1.5 ex]
h_{{8/3}}-g_{{8/3}}\!&\!\!=\!\!&\!
{\frac {1}{24}}\,{\gbar}_{{8/3}}
-{\frac {7}{57600}}\,{\gbar}_{{2/3}}{{\gbar}_{{1}}}^{2}
-{\frac {1}{5760}}\,{{\gbar}_{{2/3}}}^{2}{\gbar}_{{4/3}}

\\[1.5 ex]
h_{{3}}-g_{{3}}\!&\!\!=\!\!&\!
{\frac {1}{24}}\,{\gbar}_{{3}}
-{\frac {1}{17280}}\,{{\gbar}_{{1}}}^{3}
-{\frac {11}{25920}}\,{\gbar}_{{2/3}}{\gbar}_{{1}}{\gbar}_{{4/3}}
-{\frac {1}{3456}}\,{{\gbar}_{{2/3}}}^{2}{\gbar}_{{5/3}}

\\[1.5 ex]
h_{{10/3}}-g_{{10/3}}\!&\!\!=\!\!&\!
{\frac {1}{24}}\,{\gbar}_{{10/3}} 

-{\frac {53}{201600}}\,{\gbar}_{{2/3}}{{\gbar}_{{4/3}}}^{2}
-{\frac {47}{201600}}\,{{\gbar}_{{1}}}^{2}{\gbar}_{{4/3}}
\\[1. ex] &&\!\!\!
-{\frac {59}{100800}}\,{\gbar}_{{2/3}}{\gbar}_{{1}}{\gbar}_{{5/3}}
-{\frac {11}{28800}}\,{{\gbar}_{{2/3}}}^{2}{\gbar}_{{2}}
+{\frac {1}{2903040}}\,{{\gbar}_{{2/3}}}^{5}

\\[1.5 ex]
h_{{11/3}}-g_{{11/3}}\!&\!\!=\!\!&\!
{\frac {1}{24}}\,{\gbar}_{{11/3}}  

-{\frac {23}{31680}}\,{\gbar}_{{2/3}}{\gbar}_{{1}}{\gbar}_{{2}}
-{\frac {29}{63360}}\,{{\gbar}_{{2/3}}}^{2}{\gbar}_{{7/3}}
-{\frac {17}{63360}}\,{\gbar}_{{1}}{{\gbar}_{{4/3}}}^{2}
\\[1. ex] &&\!\!\!
-{\frac {37}{126720}}\,{{\gbar}_{{1}}}^{2}{\gbar}_{{5/3}}
-{\frac {1}{1584}}\,{\gbar}_{{2/3}}{\gbar}_{{4/3}}{\gbar}_{{5/3}}
+{\frac {1}{22809600}}\,{{\gbar}_{{2/3}}}^{4}{\gbar}_{{1}}

\\[1.5 ex]
h_{{4}}-g_{{4}}\!&\!\!=\!\!&\!
{\frac {1}{24}}\,{\gbar}_{{4}} 
-{\frac {1}{2880}}\,{{\gbar}_{{1}}}^{2}{\gbar}_{{2}}
-{\frac {1}{1620}}\,{\gbar}_{{1}}{\gbar}_{{4/3}}{\gbar}_{{5/3}}
-{\frac {11}{12960}}\,{\gbar}_{{2/3}}{\gbar}_{{1}}{\gbar}_{{7/3}}
\\[1. ex] &&\!\!\!
-{\frac {1}{2880}}\,{\gbar}_{{2/3}}{{\gbar}_{{5/3}}}^{2}
-{\frac {1}{10368}}\,{{\gbar}_{{4/3}}}^{3}
-{\frac {19}{25920}}\,{\gbar}_{{2/3}}{\gbar}_{{4/3}}{\gbar}_{{2}}
-{\frac {1}{1920}}\,{{\gbar}_{{2/3}}}^{2}{\gbar}_{{8/3}}
\\[1. ex] &&\!\!\!
-{\frac {47}{44789760}}\,{{\gbar}_{{2/3}}}^{3}{{\gbar}_{{1}}}^{2}
-{\frac {1}{1119744}}\,{{\gbar}_{{2/3}}}^{4}{\gbar}_{{4/3}}

\\[1.5 ex]
h_{{13/3}}-g_{{13/3}}\!&\!\!=\!\!&\!
{\frac {1}{24}}\,{\gbar}_{{13/3}} 
-{\frac {131}{187200}}\,{\gbar}_{{1}}{\gbar}_{{4/3}}{\gbar}_{{2}}
-{\frac {179}{187200}}\,{\gbar}_{{2/3}}{\gbar}_{{1}}{\gbar}_{{8/3}}
-{\frac {43}{74880}}\,{{\gbar}_{{2/3}}}^{2}{\gbar}_{{3}}
\\[1. ex] &&\!\!\!
-{\frac {119}{374400}}\,{{\gbar}_{{4/3}}}^{2}{\gbar}_{{5/3}}
-{\frac {5}{14976}}\,{\gbar}_{{1}}{{\gbar}_{{5/3}}}^{2}
-{\frac {149}{374400}}\,{{\gbar}_{{1}}}^{2}{\gbar}_{{7/3}}
\\[1. ex] &&\!\!\!
-{\frac {31}{37440}}\,{\gbar}_{{2/3}}{\gbar}_{{4/3}}{\gbar}_{{7/3}}
-{\frac {11}{14400}}\,{\gbar}_{{2/3}}{\gbar}_{{5/3}}{\gbar}_{{2}}
-{\frac {1}{22014720}}\,{{\gbar}_{{2/3}}}^{2}{{\gbar}_{{1}}}^{3}
\\[1. ex] &&\!\!\!
+{\frac {37}{330220800}}\,{{\gbar}_{{2/3}}}^{3}{\gbar}_{{1}}{\gbar}_{{4/3}}
+{\frac {23}{264176640}}\,{{\gbar}_{{2/3}}}^{4}{\gbar}_{{5/3}}

\\[1.5 ex]
h_{{14/3}}-g_{{14/3}}\!&\!\!=\!\!&\!
{\frac {1}{24}}\,{\gbar}_{{14/3}} 
-{\frac {173}{221760}}\,{\gbar}_{{1}}{\gbar}_{{4/3}}{\gbar}_{{7/3}}
-{\frac {23}{31680}}\,{\gbar}_{{1}}{\gbar}_{{5/3}}{\gbar}_{{2}}
-{\frac {233}{221760}}\,{\gbar}_{{2/3}}{\gbar}_{{1}}{\gbar}_{{3}}
\\[1. ex] &&\!\!\!
-{\frac {179}{443520}}\,{\gbar}_{{2/3}}{{\gbar}_{{2}}}^{2}
-{\frac {197}{443520}}\,{{\gbar}_{{1}}}^{2}{\gbar}_{{8/3}}
-{\frac {29}{31680}}\,{\gbar}_{{2/3}}{\gbar}_{{4/3}}{\gbar}_{{8/3}}
\\[1. ex] &&\!\!\!
-{\frac {31}{88704}}\,{{\gbar}_{{4/3}}}^{2}{\gbar}_{{2}}
-{\frac {5}{8064}}\,{{\gbar}_{{2/3}}}^{2}{\gbar}_{{10/3}}
-{\frac {149}{443520}}\,{\gbar}_{{4/3}}{{\gbar}_{{5/3}}}^{2}
\\[1. ex] &&\!\!\!
-{\frac {37}{44352}}\,{\gbar}_{{2/3}}{\gbar}_{{5/3}}{\gbar}_{{7/3}}
 
+{\frac {1}{9313920}}\,{\gbar}_{{2/3}}{{\gbar}_{{1}}}^{4}
+{\frac {1}{776160}}\,{{\gbar}_{{2/3}}}^{2}{{\gbar}_{{1}}}^{2}{\gbar}_{{4/3}}
\\[1. ex] &&\!\!\!
+{\frac {1}{1451520}}\,{{\gbar}_{{2/3}}}^{3}{{\gbar}_{{4/3}}}^{2}
+{\frac {17}{10160640}}\,{{\gbar}_{{2/3}}}^{3}{\gbar}_{{1}}{\gbar}_{{5/3}}
+{\frac {23}{31933440}}\,{{\gbar}_{{2/3}}}^{4}{\gbar}_{{2}}
\\[1. ex] &&\!\!\!
+{\frac {1}{3065610240}}\,{{\gbar}_{{2/3}}}^{7}
\end{array}\]

%%%%%%%%%%%%%%%%%%%%%%%%%%%%%%%%%%%%%%%%%%%%%%%%%%%%%%%%%%%%%%%%%%%%%%%%%%%%%%%%%%%%%%%%%%%%
%%%%%%%%%%%%%%%%%%%%%%%%%%%%%%%%%%%%%%%%%%%%%%%%%%%%%%%%%%%%%%%%%%%%%%%%%%%%%%%%%%%%%%%%%%%%
\subsubsection{Tangency 3, ramification 4.}
\[\begin{array}{cll}
           h_{{3/4}}-g_{{3/4}}\!&\!\!=\!\!&\!  
{\frac{1}{24}}\,{\gbar}_{{3/4}}

\\[1.5 ex] h_{{1}}-g_{{1}}\!&\!\!=\!\!&\!  
{\frac{1}{24}}\,{\gbar}_{{1}}

\\[1.5 ex] h_{{5/4}}-g_{{5/4}}\!&\!\!=\!\!&\!  
{\frac{1}{24}}\,{\gbar}_{{5/4}}

\\[1.5 ex] h_{{3/2}}-g_{{3/2}}\!&\!\!=\!\!&\! 
{\frac{1}{24}}\,{\gbar}_{{3/2}}

\\[1.5 ex] h_{{7/4}}-g_{{7/4}}\!&\!\!=\!\!&\! 
{\frac{1}{24}}\,{\gbar}_{{7/4}}

\\[1.5 ex] h_{{2}}-g_{{2}}\!&\!\!=\!\!&\! 
{\frac{1}{24}}\,{\gbar}_{{2}}

\\[1.5 ex] h_{{9/4}}-g_{{9/4}}\!&\!\!=\!\!&\! 
{\frac{1}{24}}\,{\gbar}_{{9/4}}
+{\frac {1}{86400}}\,{{\gbar}_{{3/4}}}^{3}

\\[1.5 ex] h_{{5/2}}-g_{{5/2}}\!&\!\!=\!\!&\!  
{\frac{1}{24}}\,{\gbar}_{{5/2}}
-{\frac {1}{14400}}\,{{\gbar}_{{3/4}}}^{2}{\gbar}_{{1}}

\\[1.5 ex] h_{{11/4}}-g_{{11/4}}\!&\!\!=\!\!&\!
{\frac{1}{24}}\,{\gbar}_{{11/4}}
-{\frac {59}{443520}}\,{\gbar}_{{3/4}}{{\gbar}_{{1}}}^{2}
-{\frac {71}{443520}}\,{{\gbar}_{{3/4}}}^{2}{\gbar}_{{5/4}}

\\[1.5 ex] h_{{3}}-g_{{3}}\!&\!\!=\!\!&\!  
{\frac{1}{24}}\,{\gbar}_{{3}}
-{\frac {1}{17280}}\,{{\gbar}_{{1}}}^{3}
-{\frac {1}{2560}}\,{\gbar}_{{3/4}}{\gbar}_{{1}}{\gbar}_{{5/4}}
-{\frac {11}{46080}}\,{{\gbar}_{{3/4}}}^{2}{\gbar}_{{3/2}}

\\[1.5 ex] h_{{13/4}}-g_{{13/4}}\!&\!\!=\!\!&\!  
{\frac{1}{24}}\,{\gbar}_{{13/4}}
-{\frac {49}{224640}}\,{{\gbar}_{{1}}}^{2}{\gbar}_{{5/4}}
-{\frac {19}{37440}}\,{\gbar}_{{3/4}}{\gbar}_{{1}}{\gbar}_{{3/2}}
\\[1. ex] &&\!\!\!
-{\frac {53}{224640}}\,{\gbar}_{{3/4}}{{\gbar}_{{5/4}}}^{2}
-{\frac {23}{74880}}\,{{\gbar}_{{3/4}}}^{2}{\gbar}_{{7/4}}

\\[1.5 ex] h_{{7/2}}-g_{{7/2}}\!&\!\!=\!\!&\! 
{\frac{1}{24}}\,{\gbar}_{{7/2}} 
-{\frac {53}{201600}}\,{{\gbar}_{{1}}}^{2}{\gbar}_{{3/2}}
-{\frac {31}{50400}}\,{\gbar}_{{3/4}}{\gbar}_{{1}}{\gbar}_{{7/4}}
 \\[1. ex] &&\!\!\!
-{\frac {1}{1800}}\,{\gbar}_{{3/4}}{\gbar}_{{5/4}}{\gbar}_{{3/2}}
-{\frac {37}{100800}}\,{{\gbar}_{{3/4}}}^{2}{\gbar}_{{2}}
-{\frac {1}{4032}}\,{\gbar}_{{1}}{{\gbar}_{{5/4}}}^{2}

\\[1.5 ex] h_{{15/4}}-g_{{15/4}}\!&\!\!=\!\!&\! 
{\frac{1}{24}}\,{\gbar}_{{15/4}}
-{\frac {113}{158400}}\,{\gbar}_{{3/4}}{\gbar}_{{1}}{\gbar}_{{2}}
-{\frac {89}{158400}}\,{\gbar}_{{3/2}}{\gbar}_{{1}}{\gbar}_{{5/4}}
-{\frac {17}{190080}}\,{{\gbar}_{{5/4}}}^{3}
 \\[1. ex] &&\!\!\!
-{\frac {97}{316800}}\,{\gbar}_{{3/4}}{{\gbar}_{{3/2}}}^{2}
-{\frac {97}{316800}}\,{{\gbar}_{{1}}}^{2}{\gbar}_{{7/4}}
-{\frac {101}{158400}}\,{\gbar}_{{3/4}}{\gbar}_{{5/4}}{\gbar}_{{7/4}}
 \\[1. ex] &&\!\!\!
-{\frac {133}{316800}}\,{{\gbar}_{{3/4}}}^{2}{\gbar}_{{9/4}}
+{\frac {1}{53222400}}\,{{\gbar}_{{3/4}}}^{5}

\\[1.5 ex] h_{{4}}-g_{{4}}\!&\!\!=\!\!&\!  
{\frac{1}{24}}\,{\gbar}_{{4}}
-{\frac {1}{2880}}\,{{\gbar}_{{1}}}^{2}{\gbar}_{{2}}
-{\frac {7}{23040}}\,{\gbar}_{{1}}{{\gbar}_{{3/2}}}^{2}
-{\frac {37}{46080}}\,{\gbar}_{{3/4}}{\gbar}_{{1}}{\gbar}_{{9/4}}
\\[1. ex] &&\!\!\!
-{\frac {29}{46080}}\,{\gbar}_{{1}}{\gbar}_{{5/4}}{\gbar}_{{7/4}}
-{\frac {43}{92160}}\,{{\gbar}_{{3/4}}}^{2}{\gbar}_{{5/2}}
-{\frac {11}{15360}}\,{\gbar}_{{3/4}}{\gbar}_{{5/4}}{\gbar}_{{2}}
\\[1. ex] &&\!\!\!
-{\frac {3}{10240}}\,{{\gbar}_{{5/4}}}^{2}{\gbar}_{{3/2}}
-{\frac {31}{46080}}\,{\gbar}_{{3/4}}{\gbar}_{{3/2}}{\gbar}_{{7/4}}
-{\frac {29}{70778880}}\,{{\gbar}_{{3/4}}}^{4}{\gbar}_{{1}}

\\[1.5 ex] h_{{17/4}}-g_{{17/4}}\!&\!\!=\!\!&\! 
{\frac{1}{24}}\,{\gbar}_{{17/4}}
-{\frac {563}{636480}}\,{\gbar}_{{3/4}}{\gbar}_{{1}}{\gbar}_{{5/2}}
-{\frac {419}{636480}}\,{\gbar}_{{3/2}}{\gbar}_{{1}}{\gbar}_{{7/4}}
 \\[1. ex] &&\!\!\!
-{\frac {443}{636480}}\,{\gbar}_{{1}}{\gbar}_{{5/4}}{\gbar}_{{2}}
-{\frac {7}{19584}}\,{\gbar}_{{3/4}}{{\gbar}_{{7/4}}}^{2}
-{\frac {407}{1272960}}\,{{\gbar}_{{5/4}}}^{2}{\gbar}_{{7/4}}
 \\[1. ex] &&\!\!\!
-{\frac {467}{636480}}\,{\gbar}_{{3/4}}{\gbar}_{{3/2}}{\gbar}_{{2}}
-{\frac {79}{254592}}\,{\gbar}_{{5/4}}{{\gbar}_{{3/2}}}^{2}
-{\frac {647}{1272960}}\,{{\gbar}_{{3/4}}}^{2}{\gbar}_{{11/4}}
 \\[1. ex] &&\!\!\!
-{\frac {491}{1272960}}\,{{\gbar}_{{1}}}^{2}{\gbar}_{{9/4}}
-{\frac {503}{636480}}\,{\gbar}_{{3/4}}{\gbar}_{{5/4}}{\gbar}_{{9/4}}
-{\frac {1}{7128576}}\,{{\gbar}_{{3/4}}}^{3}{{\gbar}_{{1}}}^{2}
 \\[1. ex] &&\!\!\!
-{\frac {43}{641571840}}\,{{\gbar}_{{3/4}}}^{4}{\gbar}_{{5/4}}
\end{array}\]
%%%%%%%%%%%%%%%%%%%%%%%%%%%%%%%%%%%%%%%%%%%%%%%%%%%%%%%%%%%%%%%%%%%%%%%%%%%%%%%%%%%%%%%%%%%%
%%%%%%%%%%%%%%%%%%%%%%%%%%%%%%%%%%%%%%%%%%%%%%%%%%%%%%%%%%%%%%%%%%%%%%%%%%%%%%%%%%%%%%%%%%%%
%%%%%%%%%%%%%%%%%%%%%%%%%%%%%%%%%%%%%%%%%%%%%%%%%%%%%%%%%%%%%%%%%%%%%%%%%%%%%%%%%%%%%%%%%%%%
%%%%%%%%%%%%%%%%%%%%%%%%%%%%%%%%%%%%%%%%%%%%%%%%%%%%%%%%%%%%%%%%%%%%%%%%%%%%%%%%%%%%%%%%%%%%
\subsection{The {\it nir} transform\,: from $f$  to  $h$.}
%%%%%%%%%%%%%%%%%%%%%%%%%%%%%%%%%%%%%%%%%%%%%%%%%%%%%%%%%%%%%%%%%%%%%%%%%%%%%%%%%%%%%%%%%%%%
%%%%%%%%%%%%%%%%%%%%%%%%%%%%%%%%%%%%%%%%%%%%%%%%%%%%%%%%%%%%%%%%%%%%%%%%%%%%%%%%%%%%%%%%%%%%
\subsubsection{Tangency 0, ramification 1.}
\[\begin{array}{lllll}
{h}_{{0}}\!\!&=& \!\!
{{\it f}_{{0}}}^{-1}
&\hspace{30.ex}\\[1.5 ex]
{h}_{{1}}\!\!&=&\!\!
-{{\it f}_{{0}}}^{-3}f_{{1}}
\;\;\hesi
+{\frac {1}{24}}\,{\it f}_{{0}}^{-1}f_{{1}}

&\hspace{30.ex}\\[1.5 ex]
{h}_{{2}}\!\!&=&\!\!
+{\frac {3}{2}}\,{{\it f}_{{0}}}^{-5}{f_{{1}}}^{2}
-                {{\it f}_{{0}}}^{-4}f_{{2}}
\;\;\hesi
-{\frac {1}{48}}\,{{\it f}_{{0}}}^{-3}{f_{{1}}}^{2}
+{\frac {1}{24}}\,{{\it f}_{{0}}}^{-2}f_{{2}}
\;\;\hesi
+{\frac {1}{2304}}\,{{\it f}_{{0}}}^{-1}{f_{{1}}}^{2}

&\hspace{30.ex}\\[1.5 ex]
{h}_{{3}}\!\!&=&\!\!
-{\frac {5}{2}}\,    {{\it f}_{{0}}}^{-7}{f_{{1}}}^{3}
+{\frac {10}{3}}\,   {{\it f}_{{0}}}^{-6}f_{{1}}f_{{2}}
-                    {{\it f}_{{0}}}^{-5}f_{{3}}
\;\;\hesi
+{\frac {1}{48}}\,   {{\it f}_{{0}}}^{-5}{f_{{1}}}^{3}
-{\frac {1}{18}}\,   {{\it f}_{{0}}}^{-4}f_{{1}}f_{{2}}
+{\frac {1}{24}}\,   {{\it f}_{{0}}}^{-3}f_{{3}}
\\[1.ex]&&
\;\;\hesi
-{\frac {1}{6912}}\, {{\it f}_{{0}}}^{-3}{f_{{1}}}^{3}
+{\frac {1}{1728}}\, {{\it f}_{{0}}}^{-2}f_{{1}}f_{{2}}
-{\frac {7}{5760}}\, {{\it f}_{{0}}}^{-1}f_{{3}}
\;\;\hesi
+{\frac {1}{497664}}\,{{\it f}_{{0}}}^{-1}{f_{{1}}}^{3}

&\hspace{30.ex}\\[1.5 ex]
{h}_{{4}}\!\!&=&\!\!
 {\frac {35}{8}}\,    {{\it f}_{{0}}}^{-9}{f_{{1}}}^{4}
-{\frac {35}{4}}\,    {{\it f}_{{0}}}^{-8}{f_{{1}}}^{2}f_{{2}}
+{\frac {5}{3}}\,     {{\it f}_{{0}}}^{-7}{f_{{2}}}^{2}
+{\frac {15}{4}}\,    {{\it f}_{{0}}}^{-7}f_{{1}}f_{{3}}
-                     {{\it f}_{{0}}}^{-6}f_{{4}}
\;\;\hesi
\\[1.ex]&&
-{\frac {5}{192}}\,   {{\it f}_{{0}}}^{-7}{f_{{1}}}^{4}
+{\frac {25}{288}}\,  {{\it f}_{{0}}}^{-6}{f_{{1}}}^{2}f_{{2}}
-{\frac {7}{96}}\,    {{\it f}_{{0}}}^{-5}f_{{1}}f_{{3}}
-{\frac {1}{36}}\,    {{\it f}_{{0}}}^{-5}{f_{{2}}}^{2}
+{\frac {1}{24}}\,    {{\it f}_{{0}}}^{-4}f_{{4}}
\;\;\hesi
\\[1.ex]&&
+{\frac {1}{9216}}\,  {{\it f}_{{0}}}^{-5}{f_{{1}}}^{4}
-{\frac {7}{13824}}\, {{\it f}_{{0}}}^{-4}{f_{{1}}}^{2}f_{{2}}
+{\frac {17}{23040}}\,{{\it f}_{{0}}}^{-3}f_{{1}}f_{{3}}
+{\frac {1}{3456}}\,  {{\it f}_{{0}}}^{-3}{f_{{2}}}^{2}
-{\frac {7}{5760}}\,  {{\it f}_{{0}}}^{-2}f_{{4}}
\\[1.ex]&&
\hesi\!
-\!{\frac {1}{1990656}}\,{{\it f}_{{0}}}^{-3}{f_{{1}}}^{4}
+{\frac {1}{331776}}\, {{\it f}_{{0}}}^{-2}{f_{{1}}}^{2}f_{{2}}
-{\frac {7}{552960}}\, {{\it f}_{{0}}}^{-1}f_{{1}}f_{{3}}
\hesi\!\!
+{\frac {1}{191102976}}\,{{\it f}_{{0}}}^{-1}{f_{{1}}}^{4}

&\hspace{30.ex}\\[1.5 ex]
{h}_{{5}}\!\!&=&\!\!
-{\frac {63}{8}}\,       {{\it f}_{{0}}}^{-11}{f_{{1}}}^{5}
+21\,                    {{\it f}_{{0}}}^{-10}{f_{{1}}}^{3}f_{{2}}
-{\frac {28}{3}}\,       {{\it f}_{{0}}}^{-9}f_{{1}}{f_{{2}}}^{2}
-{\frac {21}{2}}\,       {{\it f}_{{0}}}^{-9}{f_{{1}}}^{2}f_{{3}}
+{\frac {7}{2}}\,        {{\it f}_{{0}}}^{-8}f_{{2}}f_{{3}}
\\[1.ex]&&
+{\frac {21}{5}}\,       {{\it f}_{{0}}}^{-8}f_{{1}}f_{{4}}
-                        {{\it f}_{{0}}}^{-7}f_{{5}}
\;\;\hesi
+{\frac {7}{192}}\,      {{\it f}_{{0}}}^{-9}{f_{{1}}}^{5}
-{\frac {7}{48}}\,       {{\it f}_{{0}}}^{-8}{f_{{1}}}^{3}f_{{2}}
+{\frac {1}{8}}\,        {{\it f}_{{0}}}^{-7}{f_{{1}}}^{2}f_{{3}}
\\[1.ex]&&
+{\frac {7}{72}}\,       {{\it f}_{{0}}}^{-7}f_{{1}}{f_{{2}}}^{2}
-{\frac {1}{16}}\,       {{\it f}_{{0}}}^{-6}f_{{2}}f_{{3}}
-{\frac {11}{120}}\,     {{\it f}_{{0}}}^{-6}f_{{1}}f_{{4}}
+{\frac {1}{24}}\,       {{\it f}_{{0}}}^{-5}f_{{5}}
\;\;\hesi
-{\frac {1}{9216}}\,     {{\it f}_{{0}}}^{-7}{f_{{1}}}^{5}
\\[1.ex]&&
+{\frac {1}{1728}}\,     {{\it f}_{{0}}}^{-6}{f_{{1}}}^{3}f_{{2}}
-{\frac {43}{57600}}\,   {{\it f}_{{0}}}^{-5}{f_{{1}}}^{2}f_{{3}}
-{\frac {1}{1728}}\,     {{\it f}_{{0}}}^{-5}f_{{1}}{f_{{2}}}^{2}
+{\frac {37}{57600}}\,   {{\it f}_{{0}}}^{-4}f_{{2}}f_{{3}}
\\[1.ex]&&
+{\frac {31}{28800}}\,   {{\it f}_{{0}}}^{-4}f_{{1}}f_{{4}}
-{\frac {7}{5760}}\,     {{\it f}_{{0}}}^{-3}f_{{5}}
\;\;\hesi
+{\frac {1}{3317760}}\,  {{\it f}_{{0}}}^{-5}{f_{{1}}}^{5}
-{\frac {1}{497664}}\,   {{\it f}_{{0}}}^{-4}{f_{{1}}}^{3}f_{{2}}
\\[1.ex]&&
+{\frac {1}{230400}}\,   {{\it f}_{{0}}}^{-3}{f_{{1}}}^{2}f_{{3}}
+{\frac {1}{414720}}\,   {{\it f}_{{0}}}^{-3}f_{{1}}{f_{{2}}}^{2}
-{\frac {7}{1382400}}\,  {{\it f}_{{0}}}^{-2}f_{{2}}f_{{3}}
-{\frac {7}{691200}}\,   {{\it f}_{{0}}}^{-2}f_{{1}}f_{{4}}
\\[1.ex]&&
+{\frac {31}{967680}}\,  {{\it f}_{{0}}}^{-1}f_{{5}}
\;\;\hesi
-{\frac {1}{955514880}}\,{{\it f}_{{0}}}^{-3}{f_{{1}}}^{5}
+{\frac {1}{119439360}}\,{{\it f}_{{0}}}^{-2}{f_{{1}}}^{3}f_{{2}}
\\[1.ex]&&
-{\frac {7}{132710400}}\,{{\it f}_{{0}}}^{-1}{f_{{1}}}^{2}f_{{3}}
\;\;\hesi
+{\frac {1}{114661785600}}\,{{\it f}_{{0}}}^{-1}{f_{{1}}}^{5}
\end{array}\]
%%%%%%%%%%%%%%%%%%%%%%%%%%%%%%%%%%%%%%%%%%%%%%%%%%%%%%%%%%%%%%%%%%%%%%%%%%%%%%%%%%%%%%%%%%%%
%%%%%%%%%%%%%%%%%%%%%%%%%%%%%%%%%%%%%%%%%%%%%%%%%%%%%%%%%%%%%%%%%%%%%%%%%%%%%%%%%%%%%%%%%%%%
\subsubsection{Tangency 1, ramification 2.}
\[\begin{array}{lllll}

{h}_{{-1/2}}\!\!&=&\!\!
{2}^{-1/2}
\{
{{\it f}_{{1}}}^{-1/2}
\}
&\hspace{30.ex}\\[1.5 ex]

{h}_{{0}}\!\!&=& \!\!
\;\;\;\;\{-{\frac {2}{3}}{{\it f}_{{1}}}^{-2}f_{{2}} \}
&\hspace{30.ex}\\[1.5 ex]

{h}_{{1/2}}\!\!&=&\!\!
{2}^{1/2}
\{
 {\frac {5}{6}}\,{{\it f}_{{1}}}^{-7/2}{f_{{2}}}^{2}
-{\frac {3}{4}}\,{{\it f}_{{1}}}^{-5/2}f_{{3}}
\;\;\hesi
+{\frac {1}{24}}\, {{\it f}_{{1}}}^{1/2}
\}
&\hspace{30.ex}\\[1.5 ex]

{h}_{{1}}\!\!&=&\!\!
2\;\;\;
\{
-{\frac {4}{5}}\,  {{\it f}_{{1}}}^{-3}f_{{4}}
-{\frac {32}{27}}\,{{\it f}_{{1}}}^{-5}{f_{{2}}}^{3}
+2\,               {{\it f}_{{1}}}^{-4}f_{{2}}f_{{3}}
\;\;\hesi
+{\frac {1}{36}}\, {{\it f}_{{1}}}^{-1}f_{{2}}
\}
&\hspace{30.ex}\\[1.5 ex]

{h}_{{3/2}}\!\!&=&\!\!
{2}^{3/2}
\{
-{\frac {5}{6}}\,    {{\it f}_{{1}}}^{-7/2}f_{{5}}
+{\frac {7}{3}}\,    {{\it f}_{{1}}}^{-9/2}f_{{2}}f_{{4}}
-{\frac {35}{8}}\,   {{\it f}_{{1}}}^{-11/2}{f_{{2}}}^{2}f_{{3}}
+{\frac {385}{216}}\,{{\it f}_{{1}}}^{-13/2}{f_{{2}}}^{4}
\\[1. ex] &&\;\;\;\;\;\;
+{\frac {35}{32}}\,  {{\it f}_{{1}}}^{-9/2}{f_{{3}}}^{2}
\;\;\hesi
-{\frac {7}{432}}\,  {{\it f}_{{1}}}^{-5/2}{f_{{2}}}^{2}
+{\frac {1}{32}}\,   {{\it f}_{{1}}}^{-3/2}f_{{3}}
\;\;\hesi
+{\frac {1}{3456}}\, {{\it f}_{{1}}}^{3/2}
\}
&\hspace{30.ex}\\[1.5 ex]

{h}_{{2}}\!\!&=&\!\!
{2}^{2}\;\;
\{
-{\frac {16}{3}}\,  {{\it f}_{{1}}}^{-6}{f_{{2}}}^{2}f_{{4}}
-5\,                {{\it f}_{{1}}}^{-6}f_{{2}}{f_{{3}}}^{2}
+{\frac {80}{9}}\,  {{\it f}_{{1}}}^{-7}{f_{{2}}}^{3}f_{{3}}
+{\frac {8}{3}}\,   {{\it f}_{{1}}}^{-5}f_{{2}}f_{{5}}
\\[1. ex] &&\;\;\;\;\;\;
+{\frac {12}{5}}\,  {{\it f}_{{1}}}^{-5}f_{{3}}f_{{4}}
-{\frac {6}{7}}\,   {{\it f}_{{1}}}^{-4}f_{{6}}
-{\frac {224}{81}}\,{{\it f}_{{1}}}^{-8}{f_{{2}}}^{5}
\;\;\hesi
+{\frac {1}{30}}\,  {{\it f}_{{1}}}^{-2}f_{{4}}
\\[1. ex] &&\;\;\;\;\;\;
+{\frac {5}{324}}\, {{\it f}_{{1}}}^{-4}{f_{{2}}}^{3}
-{\frac {1}{24}}\,  {{\it f}_{{1}}}^{-3}f_{{2}}f_{{3}}
\;\;\hesi
+{\frac {5}{13824}}\,f_{{2}}
\}
&\hspace{30.ex}\\[1.5 ex]

{h}_{{5/2}}\!\!&=&\!\!
{2}^{5/2}
\{
-{\frac {231}{20}}\,     {{\it f}_{{1}}}^{-13/2}f_{{2}}f_{{3}}f_{{4}}
-{\frac {7}{8}}\,        {{\it f}_{{1}}}^{-9/2}f_{{7}}
-{\frac {5005}{288}}\,   {{\it f}_{{1}}}^{-17/2}{f_{{2}}}^{4}f_{{3}}
\\[1. ex] &&\;\;\;\;\;\;
+{\frac {1001}{90}}\,    {{\it f}_{{1}}}^{-15/2}{f_{{2}}}^{3}f_{{4}}
-{\frac {231}{128}}\,    {{\it f}_{{1}}}^{-13/2}{f_{{3}}}^{3}
+{\frac {17017}{3888}}\, {{\it f}_{{1}}}^{-19/2}{f_{{2}}}^{6}
\\[1. ex] &&\;\;\;\;\;\;
-{\frac {77}{12}}\,      {{\it f}_{{1}}}^{-13/2}f_{{5}}{f_{{2}}}^{2}
+{\frac {21}{8}}\,       {{\it f}_{{1}}}^{-11/2}f_{{3}}f_{{5}}
+{\frac {63}{50}}\,      {{\it f}_{{1}}}^{-11/2}{f_{{4}}}^{2}
\\[1. ex] &&\;\;\;\;\;\;
+{\frac {1001}{64}}\,    {{\it f}_{{1}}}^{-15/2}{f_{{2}}}^{2}{f_{{3}}}^{2}
+3\,                     {{\it f}_{{1}}}^{-11/2}f_{{2}}f_{{6}}
\;\;\hesi
-{\frac {91}{5184}}\,    {{\it f}_{{1}}}^{-11/2}{f_{{2}}}^{4}
\\[1. ex] &&\;\;\;\;\;\;
+{\frac {5}{144}}\,      {{\it f}_{{1}}}^{-5/2}f_{{5}}
-{\frac {17}{768}}\,     {{\it f}_{{1}}}^{-7/2}{f_{{3}}}^{2}
+{\frac {35}{576}}\,     {{\it f}_{{1}}}^{-9/2}{f_{{2}}}^{2}f_{{3}}
-{\frac {19}{360}}\,     {{\it f}_{{1}}}^{-7/2}f_{{2}}f_{{4}}
\\[1. ex] &&\;\;\;\;\;\;
\;\;\hesi
-{\frac {7}{38400}}\,    {{\it f}_{{1}}}^{-1/2}f_{{3}}
+{\frac {1}{20736}}\,    {{\it f}_{{1}}}^{-3/2}{f_{{2}}}^{2}
\;\;\hesi
{\frac {1}{1244160}}\,   {{\it f}_{{1}}}^{5/2}
\}
&\hspace{30.ex}\\[1.5 ex]
{h}_{{3}}\!\!&=&\!\!
{2}^{3}\;\;
\{
{\frac {112}{3}}\,   {{\it f}_{{1}}}^{-8}{f_{{2}}}^{2}f_{{3}}f_{{4}}
-{\frac {40}{3}}\,    {{\it f}_{{1}}}^{-7}f_{{2}}f_{{3}}f_{{5}}
+{\frac {10}{3}}\,    {{\it f}_{{1}}}^{-6}f_{{2}}f_{{7}}
+{\frac {20}{7}}\,    {{\it f}_{{1}}}^{-6}f_{{3}}f_{{6}}
\\[1. ex] &&\;\;\;\;\;\;
-{\frac {160}{21}}\,  {{\it f}_{{1}}}^{-7}{f_{{2}}}^{2}f_{{6}}
+{\frac {8}{3}}\,     {{\it f}_{{1}}}^{-6}f_{{4}}f_{{5}}
+{\frac {1120}{81}}\, {{\it f}_{{1}}}^{-8}{f_{{2}}}^{3}f_{{5}}
-{\frac {32}{5}}\,    {{\it f}_{{1}}}^{-7}f_{{2}}{f_{{4}}}^{2}
\\[1. ex] &&\;\;\;\;\;\;
-{\frac {8}{9}}\,     {{\it f}_{{1}}}^{-5}f_{{8}}
-6\,                  {{\it f}_{{1}}}^{-7}{f_{{3}}}^{2}f_{{4}}
-{\frac {1792}{81}}\, {{\it f}_{{1}}}^{-9}{f_{{2}}}^{4}f_{{4}}
-{\frac {1120}{27}}\, {{\it f}_{{1}}}^{-9}{f_{{2}}}^{3}{f_{{3}}}^{2}
\\[1. ex] &&\;\;\;\;\;\;
+{\frac {35}{3}}\,    {{\it f}_{{1}}}^{-8}f_{{2}}{f_{{3}}}^{3}
-{\frac {5120}{729}}\,{{\it f}_{{1}}}^{-11}{f_{{2}}}^{7}
+{\frac {896}{27}}\,  {{\it f}_{{1}}}^{-10}{f_{{2}}}^{5}f_{{3}}
\;\;\hesi
+{\frac {1}{28}}\,    {{\it f}_{{1}}}^{-3}f_{{6}}
\\[1. ex] &&\;\;\;\;\;\;
-{\frac {7}{108}}\,   {{\it f}_{{1}}}^{-4}f_{{2}}f_{{5}}
-{\frac {5}{54}}\,    {{\it f}_{{1}}}^{-6}{f_{{2}}}^{3}f_{{3}}
+{\frac {16}{729}}\,  {{\it f}_{{1}}}^{-7}{f_{{2}}}^{5}
+{\frac {11}{135}}\,  {{\it f}_{{1}}}^{-5}{f_{{2}}}^{2}f_{{4}}
\\[1. ex] &&\;\;\;\;\;\;
+{\frac {5}{72}}\,    {{\it f}_{{1}}}^{-5}f_{{2}}{f_{{3}}}^{2}
-{\frac {1}{20}}\,    {{\it f}_{{1}}}^{-4}f_{{3}}f_{{4}}
\;\;\hesi
-{\frac {1}{23328}}\, {{\it f}_{{1}}}^{-3}{f_{{2}}}^{3}
-{\frac {1}{2880}}\,  {{\it f}_{{1}}}^{-1}f_{{4}}
\\[1. ex] &&\;\;\;\;\;\;
+{\frac {1}{5760}}\,  {{\it f}_{{1}}}^{-2}f_{{2}}f_{{3}}
\;\;\hesi
+{\frac {1}{746496}}\,{\it f}_{{1}}f_{{2}}
\}
&\hspace{30.ex} \end{array}\] 
%%%%%%%%%%%%%%%%%%%%%%%%%%%%%%%%%%%%%%%%%%%%%%%%%%%%%%%%%%%%%%%%%%%%%%%%%%%%%%%%%%%%%%%%%%%%
%%%%%%%%%%%%%%%%%%%%%%%%%%%%%%%%%%%%%%%%%%%%%%%%%%%%%%%%%%%%%%%%%%%%%%%%%%%%%%%%%%%%%%%%%%%%
\subsubsection{Tangency 2, ramification 3.}
\[\begin{array}{lllll}
{h}_{{-2/3}} \!\!&=&\!\!  
{3}^{-2/3}
\;\{
{{\it f}_{{2}}}^{-1/3}
\}
\\ [1.5 ex] 
{h}_{{-1/3}} \!\!&=&\!\! 
{3}^{-1/3}\{
-{\frac {1}{2}}\,{{\it f}_{{2}}}^{-5/3}f_{{3}}\}
\\ [1.5 ex] 
{h}_{{0}} \!\!&=&\!\! 
\;\;\;\;\;\{
{\frac {9}{16}}{{\it f}_{{2}}}^{-3}f_{{3}}^{{2}}
-{\frac {3}{5}}{{\it f}_{{2}}}^{-2}f_{{4}}
\}
\\ [1.5 ex] 
{h}_{{1/3}} \!\!&=&\!\! 
{3}^{1/3}
\,\{
-{\frac {2}{3}}\,{{\it f}_{{2}}}^{-7/3}f_{{5}}
-{\frac {35}{48}}\,{{\it f}_{{2}}}^{-13/3}{f_{{3}}}^{3}
+{\frac {7}{5}}\,{{\it f}_{{2}}}^{-10/3}f_{{3}}f_{{4}}
\} 
&\hspace{30.ex} 
\\ [1.5 ex] 
{h}_{{2/3}} \!\!&=&\!\! 
{3}^{2/3}
\,\{
{\frac {385}{384}}\,{{\it f}_{{2}}}^{{- {17}/{3}}}{f_{{3}}}^{4}
+{\frac {5}{3}}\,{{\it f}_{{2}}}^{-11/3}f_{{3}}f_{{5}}
+{\frac {4}{5}}\,{{\it f}_{{2}}}^{-11/3}{f_{{4}}}^{2}
-{\frac {5}{7}}\,{{\it f}_{{2}}}^{-8/3}f_{{6}}
\\[1.ex]&&\;\;\;\;
-{\frac {11}{4}}\,{{\it f}_{{2}}}^{-14/3}{f_{{3}}}^{2}f_{{4}}
\;\;\hesi
+{\frac {1}{24}}\,{{\it f}_{{2}}}^{1/3}
\}
\\ [1.5 ex]
{h}_{{1}} \!\!&=&\!\! 
3
\,\{
{\frac {9}{5}}\,{{\it f}_{{2}}}^{-4}f_{{4}}f_{{5}}
+{\frac {81}{16}}\,{{\it f}_{{2}}}^{-6}f_{{4}}{f_{{3}}}^{3}
-{\frac {3}{4}}\,{{\it f}_{{2}}}^{-3}f_{{7}}
+{\frac {27}{14}}\,{{\it f}_{{2}}}^{-4}f_{{3}}f_{{6}}
-{\frac {27}{8}}\,{{\it f}_{{2}}}^{-5}{f_{{3}}}^{2}f_{{5}}
\\[1.ex] &&\;\;\;\;
-{\frac {81}{25}}\,{{\it f}_{{2}}}^{-5}f_{{3}}{f_{{4}}}^{2}
-{\frac {729}{512}}\,{{\it f}_{{2}}}^{-7}{f_{{3}}}^{5}
\hesi
+{\frac {1}{48}}\,f_{{3}}{{\it f}_{{2}}}^{-1}
\}
&\hspace{30.ex}\\ [1.5 ex]
{h}_{{4/3}} \!\!&=&\!\! 
{3}^{4/3}
\,\{
-{\frac {91}{75}}\,{{\it f}_{{2}}}^{-16/3}{f_{{4}}}^{3}
-{\frac {91}{12}}\,{{\it f}_{{2}}}^{-16/3}f_{{3}}f_{{4}}f_{{5}}
-{\frac {7}{9}}\,{{\it f}_{{2}}}^{-10/3}f_{{8}}
-{\frac {65}{16}}\,{{\it f}_{{2}}}^{-16/3}{f_{{3}}}^{2}f_{{6}}
\\[1.ex] &&\;\;\;
+{\frac {35}{16}}\,{{\it f}_{{2}}}^{-13/3}f_{{3}}f_{{7}}
+               2\,{{\it f}_{{2}}}^{-13/3}f_{{4}}f_{{6}}
-{\frac {1729}{192}}\,{{\it f}_{{2}}}^{{- {22}/{3}}}{f_{{3}}}^{4}f_{{4}}
+{\frac {19019}{9216}}\,{{\it f}_{{2}}}^{{-{25}/{3}}}{f_{{3}}}^{6}
\\[1.ex]&&\;\;\;
+{\frac {455}{72}}\,{{\it f}_{{2}}}^{{- {19}/{3}}}{f_{{3}}}^{3}f_{{5}}
+{\frac {35}{36}}\,{{\it f}_{{2}}}^{-13/3}{f_{{5}}}^{2}
+{\frac {91}{10}}\,{{\it f}_{{2}}}^{{- {19}/{3}}}{f_{{3}}}^{2}{f_{{4}}}^{2}
\\[1.ex]&&\;\;\;
\hesi
+{\frac {1}{40}}\,{{\it f}_{{2}}}^{-4/3}f_{{4}}
-{\frac {5}{384}}\,{{\it f}_{{2}}}^{-7/3}{f_{{3}}}^{2}
\}
&\hspace{30.ex}\\ [1.5 ex]
{h}_{{5/3}} \!\!&=&\!\! 
{3}^{5/3}
\,\{
-{\frac {4}{5}}\,{{\it f}_{{2}}}^{-11/3}f_{{9}}
-{\frac {44}{5}}\,{{\it f}_{{2}}}^{{- {17}/{3}}}f_{{4}}f_{{6}}f_{{3}}
-{\frac {77}{16}}\,{{\it f}_{{2}}}^{{- {17}/{3}}}{f_{{3}}}^{2}f_{{7}}
-{\frac {1309}{60}}\,{{\it f}_{{2}}}^{{- {23}/{3}}}{f_{{3}}}^{3}{f_{{4}}}^{2}
\\[1.ex]&&
-{\frac {6545}{576}}\,{{\it f}_{{2}}}^{{- {23}/{3}}}{f_{{3}}}^{4}f_{{5}}
+{\frac {187}{24}}\,{{\it f}_{{2}}}^{{- {20}/{3}}}{f_{{3}}}^{3}f_{{6}}
+{\frac {30107}{1920}}\,{{\it f}_{{2}}}^{{- {26}/{3}}}{f_{{3}}}^{5}f_{{4}}
-{\frac {55913}{18432}}\,{{\it f}_{{2}}}^{{- {29}/{3}}}{f_{{3}}}^{7}
\\[1.ex]&&
+{\frac {1309}{60}}\,{{\it f}_{{2}}}^{{- {20}/{3}}}{f_{{3}}}^{2}f_{{4}}f_{{5}}
+{\frac {44}{21}}\,{{\it f}_{{2}}}^{-14/3}f_{{6}}f_{{5}}
+{\frac {2618}{375}}\,{{\it f}_{{2}}}^{{- {20}/{3}}}f_{{3}}{f_{{4}}}^{3}
+{\frac {22}{9}}\,{{\it f}_{{2}}}^{-14/3}f_{{3}}f_{{8}}
\\[1.ex]&&
+{\frac {11}{5}}\,{{\it f}_{{2}}}^{-14/3}f_{{7}}f_{{4}}
-{\frac {308}{75}}\,{{\it f}_{{2}}}^{{- {17}/{3}}}{f_{{4}}}^{2}f_{{5}}
-{\frac {77}{18}}\,{{\it f}_{{2}}}^{{- {17}/{3}}}{f_{{5}}}^{2}f_{{3}}
\hesi
-{\frac {1}{30}}\,{{\it f}_{{2}}}^{-8/3}f_{{3}}f_{{4}}
\\[1.ex]&&
+{\frac {7}{576}}\,{{\it f}_{{2}}}^{-11/3}{f_{{3}}}^{3}
+{\frac {}{}}1/36\,{{\it f}_{{2}}}^{-5/3}f_{{5}}
\}
 &\hspace{30.ex}\\ [1.5 ex]
{h}_{{2}} \!\!&=&\!\! 
{3}^{2}
\,\{
{\frac {729}{28}}\,{{\it f}_{{2}}}^{-7}{f_{{3}}}^{2}f_{{4}}f_{{6}}
-{\frac {135}{14}}\,{{\it f}_{{2}}}^{-6}f_{{3}}f_{{5}}f_{{6}}
-{\frac {1701}{32}}\,{{\it f}_{{2}}}^{-8}{f_{{3}}}^{3}f_{{4}}f_{{5}}
+{\frac {243}{10}}\,{{\it f}_{{2}}}^{-7}f_{{3}}{f_{{4}}}^{2}f_{{5}}
\\[1.ex]&&
-{\frac {9}{2}}\,{{\it f}_{{2}}}^{-6}{f_{{5}}}^{2}f_{{4}}
-{\frac {81}{8}}\,{{\it f}_{{2}}}^{-6}f_{{3}}f_{{4}}f_{{7}}
+{\frac {295245}{65536}}\,{{\it f}_{{2}}}^{-11}{f_{{3}}}^{8}
+{\frac {27}{10}}\,{{\it f}_{{2}}}^{-5}f_{{3}}f_{{9}}
\\[1.ex]&&
-{\frac {3645}{256}}\,{{\it f}_{{2}}}^{-8}{f_{{3}}}^{4}f_{{6}}
+{\frac {243}{125}}\,{{\it f}_{{2}}}^{-7}{f_{{4}}}^{4}
-{\frac {9}{11}}\,{{\it f}_{{2}}}^{-4}f_{{10}}
-{\frac {137781}{5120}}\,{{\it f}_{{2}}}^{-10}{f_{{3}}}^{6}f_{{4}}
\\[1.ex]&&
+{\frac {54}{49}}\,{{\it f}_{{2}}}^{-5}{f_{{6}}}^{2}
-{\frac {45}{8}}\,{{\it f}_{{2}}}^{-6}{f_{{3}}}^{2}f_{{8}}
-{\frac {5103}{200}}\,{{\it f}_{{2}}}^{-8}{f_{{3}}}^{2}{f_{{4}}}^{3}
+{\frac {5103}{256}}\,{{\it f}_{{2}}}^{-9}{f_{{3}}}^{5}f_{{5}}
\\[1.ex]&&
+{\frac {15309}{320}}\,{{\it f}_{{2}}}^{-9}{f_{{3}}}^{4}{f_{{4}}}^{2}
+{\frac {12}{5}}\,{{\it f}_{{2}}}^{-5}f_{{8}}f_{{4}}
+{\frac {405}{32}}\,{{\it f}_{{2}}}^{-7}{f_{{3}}}^{2}{f_{{5}}}^{2}
-{\frac {162}{35}}\,{{\it f}_{{2}}}^{-6}{f_{{4}}}^{2}f_{{6}}
\\[1.ex]&&
+{\frac {1215}{128}}\,{{\it f}_{{2}}}^{-7}{f_{{3}}}^{3}f_{{7}}
+{\frac {9}{4}}\,{{\it f}_{{2}}}^{-5}f_{{5}}f_{{7}}
\;\;\hesi
-{\frac {27}{2048}}\,{{\it f}_{{2}}}^{-5}{f_{{3}}}^{4}
-{\frac {1}{24}}\,{{\it f}_{{2}}}^{-3}f_{{3}}f_{{5}}
\\[1.ex]&&
+{\frac {5}{168}}\,{{\it f}_{{2}}}^{-2}f_{{6}}
-{\frac {11}{600}}\,{{\it f}_{{2}}}^{-3}{f_{{4}}}^{2}
+{\frac {3}{64}}\,{{\it f}_{{2}}}^{-4}{f_{{3}}}^{2}f_{{4}}
\;\;\hesi
+{\frac {1}{5184}}\,{{\it f}_{{2}}}^{-1}
\}
&\hspace{30.ex}
\end{array}\]
%%%%%%%%%%%%%%%%%%%%%%%%%%%%%%%%%%%%%%%%%%%%%%%%%%%%%%%%%%%%%%%%%%%%%%%%%%%%%%%%%%%%%%%%%%%%
%%%%%%%%%%%%%%%%%%%%%%%%%%%%%%%%%%%%%%%%%%%%%%%%%%%%%%%%%%%%%%%%%%%%%%%%%%%%%%%%%%%%%%%%%%%%
\subsubsection{Tangency 3, ramification 4.}
\[\begin{array}{lllll}

{h}_{{-3/4}}\!\!&=&\!\!
{4}^{-3/4}
\{
{{\it f}_{{3}}}^{-3/4}
\}
&\hspace{30.ex}\\[1.5 ex]

{h}_{{-1/2}}\!\!&=&\!\!
{4}^{-1/2}
\{
-{\frac {2}{5}}\,{{\it f}_{{3}}}^{-3/2}f_{{4}}
\}
&\hspace{30.ex}\\[1.5 ex]

{h}_{{-1/4}}\!\!&=&\!\!
{4}^{-1/4}
\{
{\frac {21}{50}}\,{{\it f}_{{3}}}^{-11/4}{f_{{4}}}^{2}
-{\frac {1}{2}}\, {{\it f}_{{3}}}^{-7/4}f_{{5}}
\}

&\hspace{30.ex}\\[1.5 ex]
{h}_{{0}}\!\!&=& \!\!
\;\;\;\;\;\{
-{\frac {64}{125}}{{\it f}_{{3}}}^{-4}f_{{4}}^{{3}}
+{\frac {16}{15}}{{\it f}_{{3}}}^{-3}f_{{4}}f_{{5}}
-{\frac {4}{7}}{{\it f}_{{3}}}^{-2}f_{{6}}
\}

&\hspace{30.ex}\\[1.5 ex]

{h}_{{1/4}}\!\!&=&\!\!
{4}^{1/4}
\{
-{\frac {39}{20}}\,   {{\it f}_{{3}}}^{{- {17}/{4}}}{f_{{4}}}^{2}f_{{5}}
+{\frac {9}{7}}\,     {{\it f}_{{3}}}^{{- {13}/{4}}}f_{{4}}f_{{6}}
+{\frac {5}{8}}\,     {{\it f}_{{3}}}^{{- {13}/{4}}}{f_{{5}}}^{2}
\\[1.ex]&& \;\;\;\;\;\;
+{\frac {663}{1000}}\,{{\it f}_{{3}}}^{{- {21}/{4}}}{f_{{4}}}^{4}
-{\frac {5}{8}}\,     {{\it f}_{{3}}}^{-9/4}f_{{7}}
\}

&\hspace{30.ex}\\[1.5 ex]
{h}_{{1/2}}\!\!&=&\!\!
{4}^{1/2}
\{
{\frac {84}{25}}\,     {{\it f}_{{3}}}^{-11/2}{f_{{4}}}^{3}f_{{5}}
-{\frac {2772}{3125}}\,{{\it f}_{{3}}}^{-13/2}{f_{{4}}}^{5}
-{\frac {12}{5}}\,     {{\it f}_{{3}}}^{-9/2}{f_{{4}}}^{2}f_{{6}}
-{\frac {7}{3}}\,      {{\it f}_{{3}}}^{-9/2}f_{{4}}{f_{{5}}}^{2}
\\[1.ex]&& \;\;\;\;\;\;
+{\frac {3}{2}}\,      {{\it f}_{{3}}}^{-7/2}f_{{4}}f_{{7}}
+{\frac {10}{7}}\,     {{\it f}_{{3}}}^{-7/2}f_{{5}}f_{{6}}
-{\frac {2}{3}}\,      {{\it f}_{{3}}}^{-5/2}f_{{8}}
\}

&\hspace{30.ex}\\[1.5 ex]

{h}_{{3/4}}\!\!&=&\!\!
{4}^{3/4}
\{
-{\frac {231}{80}}\,  {{\it f}_{{3}}}^{{- {19}/{4}}}{f_{{4}}}^{2}f_{{7}}
-{\frac {11}{2}}\,    {{\it f}_{{3}}}^{{- {19}/{4}}}f_{{4}}f_{{5}}f_{{6}}
-{\frac {385}{432}}\, {{\it f}_{{3}}}^{{- {19}/{4}}}{f_{{5}}}^{3}
\\[1.ex]&& \;\;\;
+{\frac {209}{50}}\,  {{\it f}_{{3}}}^{{- {23}/{4}}}{f_{{4}}}^{3}f_{{6}}
+{\frac {1463}{240}}\,{{\it f}_{{3}}}^{{- {23}/{4}}}{f_{{4}}}^{2}{f_{{5}}}^{2}
-{\frac {7}{10}}\,    {{\it f}_{{3}}}^{-11/4}f_{{9}}
+{\frac {11}{14}}\,   {{\it f}_{{3}}}^{{- {15}/{4}}}{f_{{6}}}^{2}
\\[1.ex]&& \;\;\;
+{\frac {302841}{250000}}\,{{\it f}_{{3}}}^{{- {31}/{4}}}{f_{{4}}}^{6}
+{\frac {77}{45}}\,   {{\it f}_{{3}}}^{{- {15}/{4}}}f_{{4}}f_{{8}}
+{\frac {77}{48}}\,   {{\it f}_{{3}}}^{{- {15}/{4}}}f_{{7}}f_{{5}}
\\[1.ex]&& \;\;\;\;\;\;
-{\frac {33649}{6000}}\,{{\it f}_{{3}}}^{{- {27}/{4}}}f_{{5}}{f_{{4}}}^{4}
\;\;\hesi
+{\frac {1}{24}}\,{{\it f}_{{3}}}^{1/4}
\}

&\hspace{30.ex}\\[1.5 ex]
{h}_{{1}}\!\!&=&\!\!
4\;\;
\{
+{\frac {28672}{3125}}\,  {{\it f}_{{3}}}^{-8}{f_{{4}}}^{5}f_{{5}}
-{\frac {131072}{78125}}\,{{\it f}_{{3}}}^{-9}{f_{{4}}}^{7}
-{\frac {8}{11}}\,        {{\it f}_{{3}}}^{-3}f_{{10}}
-{\frac {32}{5}}\,        {{\it f}_{{3}}}^{-5}f_{{4}}f_{{5}}f_{{7}}
\\[1.ex]&& \;\;\;
-{\frac {256}{75}}\,      {{\it f}_{{3}}}^{-5}{f_{{4}}}^{2}f_{{8}}
-{\frac {768}{245}}\,     {{\it f}_{{3}}}^{-5}f_{{4}}{f_{{6}}}^{2}
-{\frac {64}{21}}\,       {{\it f}_{{3}}}^{-5}{f_{{5}}}^{2}f_{{6}}
+{\frac {48}{25}}\,       {{\it f}_{{3}}}^{-4}f_{{4}}f_{{9}}
\\[1.ex]&& \;\;\;
+{\frac {16}{9}}\,        {{\it f}_{{3}}}^{-4}f_{{5}}f_{{8}}
+{\frac {12}{7}}\,        {{\it f}_{{3}}}^{-4}f_{{6}}f_{{7}}
-{\frac {6144}{875}}\,    {{\it f}_{{3}}}^{-7}{f_{{4}}}^{4}f_{{6}}
-{\frac {1024}{75}}\,     {{\it f}_{{3}}}^{-7}{f_{{4}}}^{3}{f_{{5}}}^{2}
\\[1.ex]&& \;\;\;
+{\frac {512}{35}}\,      {{\it f}_{{3}}}^{-6}{f_{{4}}}^{2}f_{{5}}f_{{6}}
+{\frac {128}{25}}\,      {{\it f}_{{3}}}^{-6}{f_{{4}}}^{3}f_{{7}}
+{\frac {128}{27}}\,      {{\it f}_{{3}}}^{-6}f_{{4}}{f_{{5}}}^{3}
\;\;\hesi
+{\frac {1}{60}}\,        {\it f}_{{3}}^{-1}f_{{4}}
\}
&\hspace{30.ex} 

&\hspace{30.ex}\\[1.5 ex]
{h}_{{5/4}}\!\!&=&\!\!
{4}^{5/4}
\{
 {\frac {23476167}{10000000}}\,{{\it f}_{{3}}}^{{- {41}/{4}}}{f_{{4}}}^{8}
+{\frac {13}{7}}\,             {{\it f}_{{3}}}^{{- {17}/{4}}}f_{{6}}f_{{8}}
+{\frac {57681}{5000}}\,       {{\it f}_{{3}}}^{{- {33}/{4}}}{f_{{4}}}^{5}f_{{6}}
\\[1.ex]&& \;\;\;
+{\frac {1547}{1152}}\,        {{\it f}_{{3}}}^{{- {25}/{4}}}{f_{{5}}}^{4}
-{\frac {1989}{500}}\,         {{\it f}_{{3}}}^{{- {21}/{4}}}{f_{{4}}}^{2}f_{{9}}
-{\frac {221}{64}}\,           {{\it f}_{{3}}}^{{- {21}/{4}}}{f_{{5}}}^{2}f_{{7}}
\\[1.ex]&& \;\;\;
-{\frac {663}{196}}\,          {{\it f}_{{3}}}^{{- {21}/{4}}}f_{{5}}{f_{{6}}}^{2}
-{\frac {1480479}{100000}}\,   {{\it f}_{{3}}}^{{- {37}/{4}}}{f_{{4}}}^{6}f_{{5}}
+{\frac {5967}{700}}\,         {{\it f}_{{3}}}^{{- {25}/{4}}}{f_{{4}}}^{2}{f_{{6}}}^{2}
\\[1.ex]&& \;\;\;
-{\frac {13923}{1600}}\,       {{\it f}_{{3}}}^{{- {29}/{4}}}{f_{{4}}}^{4}f_{{7}}
-{\frac {1547}{96}}\,          {{\it f}_{{3}}}^{{- {29}/{4}}}{f_{{4}}}^{2}{f_{{5}}}^{3}
+{\frac {117}{55}}\,           {{\it f}_{{3}}}^{{- {17}/{4}}}f_{{4}}f_{{10}}
\\[1.ex]&& \;\;\;
+{\frac {39}{20}}\,            {{\it f}_{{3}}}^{{- {17}/{4}}}f_{{9}}f_{{5}}
+{\frac {1547}{250}}\,         {{\it f}_{{3}}}^{{- {25}/{4}}}{f_{{4}}}^{3}f_{{8}}
+{\frac {117}{128}}\,          {{\it f}_{{3}}}^{{- {17}/{4}}}{f_{{7}}}^{2}
-{\frac {3}{4}}\,              {{\it f}_{{3}}}^{{- {13}/{4}}}f_{{11}}
\\[1.ex]&& \;\;\;
-{\frac {663}{20}}\,           {{\it f}_{{3}}}^{{- {29}/{4}}}{f_{{4}}}^{3}f_{{5}}f_{{6}}
+{\frac {663}{40}}\,           {{\it f}_{{3}}}^{{- {25}/{4}}}{f_{{5}}}^{2}f_{{6}}f_{{4}}
+{\frac {44863}{1600}}\,       {{\it f}_{{3}}}^{{- {33}/{4}}}{f_{{4}}}^{4}{f_{{5}}}^{2}
\\[1.ex]&& \;\;\;
-{\frac {1989}{280}}\,         {{\it f}_{{3}}}^{{- {21}/{4}}}f_{{4}}f_{{6}}f_{{7}}
+{\frac {13923}{800}}\,        {{\it f}_{{3}}}^{{- {25}/{4}}}{f_{{4}}}^{2}f_{{5}}f_{{7}}
-{\frac {221}{30}}\,           {{\it f}_{{3}}}^{{- {21}/{4}}}f_{{4}}f_{{5}}f_{{8}}
\\[1.ex]&& \;\;\;
\;\;\hesi
+{\frac {1}{48}}\,             {{\it f}_{{3}}}^{-5/4}f_{{5}}
-{\frac {13}{1200}}\,          {{\it f}_{{3}}}^{-9/4}{f_{{4}}}^{2}
\}
\end{array}\] 

%%%%%%%%%%%%%%%%%%%%%%%%%%%%%%%%%%%%%%%%%%%%%%%%%%%%%%%%%%%%%%%%%%%%%%%%%%%%%%%%%%%%%%%%%%%%
%%%%%%%%%%%%%%%%%%%%%%%%%%%%%%%%%%%%%%%%%%%%%%%%%%%%%%%%%%%%%%%%%%%%%%%%%%%%%%%%%%%%%%%%%%%%
%%%%%%%%%%%%%%%%%%%%%%%%%%%%%%%%%%%%%%%%%%%%%%%%%%%%%%%%%%%%%%%%%%%%%%%%%%%%%%%%%%%%%%%%%%%%
%%%%%%%%%%%%%%%%%%%%%%%%%%%%%%%%%%%%%%%%%%%%%%%%%%%%%%%%%%%%%%%%%%%%%%%%%%%%%%%%%%%%%%%%%%%%
\subsection{The {\it nur} transform\,: from $f$  to  $h$.}
%%%%%%%%%%%%%%%%%%%%%%%%%%%%%%%%%%%%%%%%%%%%%%%%%%%%%%%%%%%%%%%%%%%%%%%%%%%%%%%%%%%%%%%%%%%%
%%%%%%%%%%%%%%%%%%%%%%%%%%%%%%%%%%%%%%%%%%%%%%%%%%%%%%%%%%%%%%%%%%%%%%%%%%%%%%%%%%%%%%%%%%%%
\subsubsection{Tangency 0.}
\[\begin{array}{lllll}
{h}_{{0}}\!\!&=& \!\!
{{\it f}_{{0}}}^{<-1>}
&\hspace{30.ex}\\[1.5 ex]
{h}_{{1}}\!\!&=&\!\!
-{{\it f}_{{0}}}^{<-3>}f_{{1}}
\;\;\hesi
+{\frac {1}{24}}\,{\it f}_{{0}}^{<-1>}f_{{1}}

&\hspace{30.ex}\\[1.5 ex]
{h}_{{2}}\!\!&=&\!\!
+{\frac {3}{2}}\,{{\it f}_{{0}}}^{<-5>}{f_{{1}}}^{2}
-                {{\it f}_{{0}}}^{<-4>}f_{{2}}
\;\;\hesi
-{\frac {1}{48}}\,{{\it f}_{{0}}}^{<-3>}{f_{{1}}}^{2}
+{\frac {1}{24}}\,{{\it f}_{{0}}}^{<-2>}f_{{2}}
\;\;\hesi
+{\frac {1}{2304}}\,{{\it f}_{{0}}}^{<-1>}{f_{{1}}}^{2}
\end{array}\]
\dotfill
$$
\alpha^{<-k>}:=
\sum_{n\,\in\,\doZ}\frac{(-1)^n}{(\alpha+2\pi\,i\,n)^k}
=\frac{1}{(k\!-\!1)!}\,\Big[\partial_\sigma^{k-1}\frac{2}{\sinh(\frac{\alpha-\sigma}{2})}
\Big]_{\sigma=0}
\equiv (-1)^k\,(-\alpha)^{<-k>}
$$
\[\begin{array}{llllll}
\alpha^{<-1>}\!&\!=\!&\!
{\frac {\sqrt {a}}{a-1}}\quad\quad\quad \mi{with}\quad a:=e^\alpha
\\[1.5 ex] 
\alpha^{<-2>}\!&\!=\!&\!
{\frac {\sqrt {a} \left( a+1 \right) }{ 2\;\left( a-1 \right) ^{2}}}
\\[1.5 ex] 
\alpha^{<-3>}\!&\!=\!&\!
{\frac {\sqrt {a} \left( {a}^{2}+6\,a+1 \right) }
{ 8\;\left( a-1 \right) ^{3}}}
\\[1.5 ex] 
\alpha^{<-4>}\!&\!=\!&\!
{\frac {\sqrt {a} \left( {a}^{3}+23\,{a}^{2}+23\,a+1 \right) }
{ 48\;\left( a-1 \right) ^{4}}}
\\[1.5 ex] 
\alpha^{<-5>}\!&\!=\!&\!
{\frac {\sqrt {a} \left( {a}^{4}+76\,{a}^{3}+230\,{a}^{2}+76\,a+1
 \right) }{ 384\;\left( a-1 \right) ^{5}}}
\\[1.5 ex] 
\alpha^{<-6>}\!&\!=\!&\!
{\frac {\sqrt {a} \left( {a}^{5}+237\,{a}^{4}+1682\,{a}^{3}+1682\,{a}^{2
}+237\,a+1 \right) }{ 3840\;\left( a-1 \right) ^{6}}}
\\[1.5 ex] 
\alpha^{<-7>}\!&\!=\!&\!
{\frac {\sqrt {a} \left( {a}^{6}+722\,{a}^{5}+10543\,{a}^{4}+23548\,{a}
^{3}+10543\,{a}^{2}+722\,a+1 \right) }{46080\; \left( a-1 \right) ^{7}}}
\\[1.5 ex] 
\alpha^{<-8>}\!&\!=\!&\!
{\frac {\sqrt {a} \left( {a}^{7}+2179\,{a}^{6}+60657\,{a}^{5}+259723\,
{a}^{4}+259723\,{a}^{3}+60657\,{a}^{2}+2179\,a+1 \right) }{645120\; \left( a-1 \right) ^{8}}}
\\[1.5 ex] 
\alpha^{<-9>}\!&\!=\!&\!
{\frac {\sqrt {a} \left( {a}^{8}+6552\,{a}^{7}+331612\,{a}^{6}+
2485288\,{a}^{5}+4675014\,{a}^{4}+2485288\,{a}^{3}+331612\,{a}^{2}+6552\,a+1 \right) }{
 10321920\;\left( a-1 \right) ^{9}}}
\end{array}\]
\dotfill
%%%%%%%%%%%%%%%%%%%%%%%%%%%%%%%%%%%%%%%%%%%%%%%%%%%%%%%%%%%%%%%%%%%%%%%%%%%%%%%%%%%%%%%%%%%%
%%%%%%%%%%%%%%%%%%%%%%%%%%%%%%%%%%%%%%%%%%%%%%%%%%%%%%%%%%%%%%%%%%%%%%%%%%%%%%%%%%%%%%%%%%%%
%\end{document}

\subsubsection{Tangency $\mg{>}$ 0.}
We set $ Dh:=\mi{nur}(f)-\mi{nir}(f)$
and to calculate the general $Dh_n$, we take $h_n$ in the {\it nir}-table
\S10.4 and perform the substitution\,:
$$
f_0^{-k}\longrightarrow \sum_{n\,\in\doZ^\ast}\frac{1}{(2\pi\,i)^k}
=-\beta_{k-1}=\frac{1}{(k-1)!}
\Big[\partial_\sigma^{k-1}\frac{2}{\sinh(-\frac{\sigma}{2})}
\Big]_{\sigma=0}
$$ 
\[\begin{array}{lllll}
{Dh}_{{0}}\!\!&=& {Dh}_{{1}} =\, 0
&\hspace{30.ex}

\\[1.5 ex] 
Dh_2  \!&\!=\!&\!
+\frac{17}{5760}\,f_2
\quad\quad\quad \equiv
\big({{\it \beta}}_{3}-{\frac {1}{24}}\,{{\it \beta}}_{1}\big)\,f_{{2}}
\\[1.5 ex] 
Dh_3  \!&\!=\!&\!
+\frac{1}{5040}\,f_2 f_1
\quad\quad\quad \equiv
\big(
-{\frac {10}{3}}\,   {{\it \beta}}_{5}
+{\frac {1}{18}}\,   {{\it \beta}}_{3}
-{\frac {1}{1728}}\, {{\it \beta}}_{1}\big)f_{{1}}f_{{2}}
\\[1.5 ex] 
Dh_4  \!&\!=\!&\!
-\frac{43}{322560}\,f_4
+\frac{11927}{1114767360}\,f_2 f_1^2

\\[1.5 ex] 
Dh_5  \!&\!=\!&\!
-\frac{18839}{2322432000}\,f_4 f_1
-\frac{27241}{4644864000}\,f_3 f_2

+\frac{28283}{200658124800}\,f_2 f_1^3

\\[1.5 ex] 
Dh_6  \!&\!=\!&\!
+\frac{769}{154828800}\,f_6

-\frac{143}{746496000}\,f_4 f_1^2
-\frac{28709}{111476736000}\,f_3 f_2 f_1
-\frac{677}{20065812480}\,f_2^3
\\[1.0 ex]&&
+\frac{1657}{2407897497600}\,f_2 f_1^4

\\[1.5 ex] 
Dh_7  \!&\!=\!&\!
+\frac{319}{1560674304}\,f_6 f_1
+\frac{13277}{117050572800}\,f_5 f_2
+\frac{373}{4335206400}\,f_4 f_3
-\frac{41539}{28092137472000}\,f_4 f_1^3
\\[1.0 ex]&&
-\frac{156581}{56184274944000}\,f_3 f_2 f_1^2
-\frac{667}{936404582400}\,f_2^3 f_1

+\frac{6133}{3467372396544000}\,f_2 f_1^5

\\[1.5 ex] 
Dh_8  \!&\!=\!&\!
-\frac{163}{1114767360}\,f_8

+\frac{32843}{10487731322880}\,f_6 f_1^2
+\frac{27151}{8739776102400}\,f_5 f_2 f_1
\\[1.0 ex]&&
+\frac{207481}{93640458240000}\,f_4 f_3 f_1
+\frac{2693}{3121348608000}\,f_4 f_2^2
+\frac{269489}{374561832960000}\,f_3^2 f_2
\\[1.0 ex]&&
-\frac{32447}{6472428473548800}\,f_4 f_1^4
-\frac{6953}{577895399424000}\,f_3 f_2 f_1^3
\\[1.0 ex]&&
-\frac{7157}{1618107118387200}\,f_2^3 f_1^2
+\frac{633697}{223687128045846528000}\,f_2 f_1^6

\end{array}\]
\dotfill

%%%%%%%%%%%%%%%%%%%%%%%%%%%%%%%%%%%%%%%%%%%%%%%%%%%%%%%%%%%%%%%%%%%%%%%%%%%%%%%%%%%%%%%%%%%%
%%%%%%%%%%%%%%%%%%%%%%%%%%%%%%%%%%%%%%%%%%%%%%%%%%%%%%%%%%%%%%%%%%%%%%%%%%%%%%%%%%%%%%%%%%%%
%%%%%%%%%%%%%%%%%%%%%%%%%%%%%%%%%%%%%%%%%%%%%%%%%%%%%%%%%%%%%%%%%%%%%%%%%%%%%%%%%%%%%%%%%%%%
%%%%%%%%%%%%%%%%%%%%%%%%%%%%%%%%%%%%%%%%%%%%%%%%%%%%%%%%%%%%%%%%%%%%%%%%%%%%%%%%%%%%%%%%%%%%
\subsection{Translocation of {\it nir}.}
%%%%%%%%%%%%%%%%%%%%%%%%%%%%%%%%%%%%%%%%%%%%%%%%%%%%%%%%%%%%%%%%%%%%%%%%%%%%%%%%%%%%%%%%%%%%
%%%%%%%%%%%%%%%%%%%%%%%%%%%%%%%%%%%%%%%%%%%%%%%%%%%%%%%%%%%%%%%%%%%%%%%%%%%%%%%%%%%%%%%%%%%%
\subsubsection{Standard case.}
\[\begin{array}{llll}
%%%%%%%%%%%%%%%%% lprint(ssarsuumo(1,0));
(\delta_1\,h)_0 \!&\!=\!&\!
+\frac{1}{24} f_1
%%%%%%%%%%%%%%%%% lprint(ssarsuumo(1,1));
\\[1.ex]
(\delta_1\,h)_1 \!&\!=\!&\!
+\frac{1}{1152} f_1^2

%%%%%%%%%%%%%%%%% lprint(ssarsuumo(1,2));
\\[1.ex]
(\delta_1\,h)_2 \!&\!=\!&\!
-\frac{7}{1920} f_3

+\frac{1}{165888} f_1^3

%%%%%%%%%%%%%%%%% lprint(ssarsuumo(1,3));
\\[1.ex]
(\delta_1\,h)_3 \!&\!=\!&\!
-\frac{7}{138240} f_1 f_3

 +\frac{1}{47775744} f_1^4

%%%%%%%%%%%%%%%%% lprint(ssarsuumo(1,4));
\\[1.ex]
(\delta_1\,h)_4 \!&\!=\!&\!
+\frac{31}{193536} f_5

- \frac{7}{26542080} f_3 f_1^2

+\frac{1}{22932357120} f_1^5

%%%%%%%%%%%%%%%%% lprint(ssarsuumo(1,5));
\\[1.ex]
(\delta_1\,h)_5 \!&\!=\!&\!
 +\frac{49}{221184000} f_3^2
+\frac{31}{23224320} f_1 f_5

 -\frac{7}{9555148800} f_1^3 f_3

 +\frac{1}{16511297126400} f_1^6

%%%%%%%%%%%%%%%%% lprint(ssarsuumo(1,6));
\\[1.ex]
(\delta_1\,h)_6 \!&\!=\!&\!
-\frac{127}{22118400} f_7

+\frac{31}{6688604160} f_1^2 f_5
+\frac{49}{31850496000} f_1 f_3^2

- \frac{7}{5503765708800} f_1^4 f_3
\\[1.ex] &&\!
 +\frac{1}{16643387503411200} f_1^7

%%%%%%%%%%%%%%%%%%%%%%%%%%%%%%%%%%%%%%%%%%%%%%%%%%%%%%%%%%%%%%%%%%%%%5
%%%%%%%%%%%%%%%%%  lprint(ssarsuumo(2,0));
\\[1.ex]
(\delta_2\,h)_0 \!&\!=\!&\!
+\frac{1}{24} f_2

-\frac{1}{2304} f_0 f_1^2

%%%%%%%%%%%%%%%%%  lprint(ssarsuumo(2,1));
\\[1.ex]
(\delta_2\,h)_1 \!&\!=\!&\!
+\frac{1}{576} f_1 f_2
+\frac{7}{1920} f_0 f_3

- \frac{1}{165888} f_0 f_1^3

%%%%%%%%%%%%%%%%%  lprint(ssarsuumo(2,2));
\\[1.ex]
(\delta_2\,h)_2 \!&\!=\!&\!
-\frac{7}{960} f_4

+\frac{7}{92160} f_0 f_1 f_3
+\frac{1}{55296} f_1^2 f_2

-\frac{1}{31850496} f_0 f_1^4

%%%%%%%%%%%%%%%%%  lprint(ssarsuumo(2,3));
\\[1.ex]
(\delta_2\,h)_3 \!&\!=\!&\!
- \frac{7}{138240} f_3 f_2
-\frac{7}{69120} f_4 f_1
-\frac{31}{96768} f_0 f_5

+\frac{7}{13271040} f_0 f_3 f_1^2
\\[1.ex] &&\!
+\frac{1}{11943936} f_1^3 f_2

-\frac{1}{11466178560} f_0 f_1^5

%%%%%%%%%%%%%%%%%  lprint(ssarsuumo(2,4));
\\[1.ex]
(\delta_2\,h)_4 \!&\!=\!&\!
+\frac{31}{64512} f_6

- \frac{7}{13271040} f_3 f_1 f_2
-\frac{31}{9289728} f_0 f_1 f_5
-\frac{49}{88473600} f_0 f_3^2
\\[1.ex] &&\!
-\frac{7}{13271040} f_1^2 f_4

+\frac{7}{3822059520} f_0 f_1^3 f_3
+\frac{1}{4586471424} f_1^4 f_2

- \frac{1}{6604518850560} f_0 f_1^6

%%%%%%%%%%%%%%%%%  lprint(ssarsuumo(2,5));
\\[1.ex]
(\delta_2\,h)_5 \!&\!=\!&\!
+\frac{31}{7741440} f_6 f_1
+\frac{31}{23224320} f_2 f_5
+\frac{127}{7372800} f_0 f_7
+\frac{49}{55296000} f_3 f_4 
\\[1.ex] &&\!
-\frac{7}{3185049600} f_1^2 f_2 f_3
-\frac{49}{10616832000} f_0 f_3^2 f_1
-\frac{31}{2229534720} f_0 f_1^2 f_5
\\[1.ex] &&\!
-\frac{7}{4777574400} f_1^3 f_4

+\frac{7}{1834588569600} f_0 f_1^4 f_3
+\frac{1}{2751882854400} f_1^5 f_2
\\[1.ex] &&\!

-\frac{1}{5547795834470400} f_0 f_1^7

%%%%%%%%%%%%%%%%%%%%%%%%%%%%%%%%%%%%%%%%%%%%%%%%%%%%%%%%%%%%%%%%%%%%%%%%%%%%%%%%
%%%%%%%%%%%%%%%%%  lprint(ssarsuumo(3,0));
\\[1.ex]
(\delta_3\,h)_0 \!&\!=\!&\!
+\frac{1}{24} f_3

- \frac{1}{864} f_0 f_1 f_2
-\frac{7}{5760} f_3 f_0^2
-\frac{1}{6912} f_1^3

+\frac{1}{497664} f_0^2 f_1^3

%%%%%%%%%%%%%%%%%  lprint(ssarsuumo(3,1));
\\[1.ex]
(\delta_3\,h)_1 \!&\!=\!&\!
+\frac{7}{720} f_0 f_4
+\frac{17}{5760} f_3 f_1
+\frac{1}{864} f_2^2

-\frac{7}{138240} f_3 f_0^2 f_1
-\frac{1}{41472} f_0 f_1^2 f_2
\\[1.ex] &&\!
-\frac{1}{497664} f_1^4 

+\frac{1}{47775744} f_0^2 f_1^4

%%%%%%%%%%%%%%%%%  lprint(ssarsuumo(3,2));
\\[1.ex]
(\delta_3\,h)_2 \!&\!=\!&\!
-\frac{7}{576} f_5

+\frac{7}{34560} f_0 f_4 f_1
+\frac{7}{69120} f_0 f_3 f_2
+\frac{1}{41472} f_1 f_2^2
+\frac{1}{23040} f_3 f_1^2
\\[1.ex] &&\!
+\frac{31}{96768} f_0^2 f_5

-\frac{7}{13271040} f_3 f_0^2 f_1^2
-\frac{1}{5971968} f_0 f_1^3 f_2
-\frac{1}{95551488} f_1^5 
\\[1.ex] &&\!
+\frac{1}{11466178560} f_0^2 f_1^5

%%%%%%%%%%%%%%%%% lprint(ssarsuumo(3,3));
\\[1.ex]
(\delta_3\,h)_3 \!&\!=\!&\!
-\frac{31}{24192} f_0 f_6
-\frac{7}{51840} f_4 f_2
-\frac{5}{18144} f_1 f_5
-\frac{7}{138240} f_3^2

+\frac{31}{6967296} f_1 f_0^2 f_5
\\[1.ex] &&\!
+\frac{7}{4976640} f_0 f_1^2 f_4
+\frac{1}{5971968} f_1^2 f_2^2
+\frac{49}{66355200} f_0^2 f_3^2
+\frac{31}{119439360} f_1^3 f_3
\\[1.ex] &&\!
+\frac{7}{4976640} f_0 f_3 f_1 f_2

-\frac{7}{2866544640} f_0^2 f_1^3 f_3
-\frac{1}{1719926784} f_0 f_1^4 f_2
\\[1.ex] &&\!
-\frac{1}{34398535680} f_1^6

+\frac{1}{4953389137920} f_0^2 f_1^6

\end{array}\]

%%%%%%%%%%%%%%%%%%%%%%%%%%%%%%%%%%%%%%%%%%%%%%%%%%%%%%%%%%%%%%%%%%%%%%%%%%%%%%%%%%%%%%%%%%%%
%%%%%%%%%%%%%%%%%%%%%%%%%%%%%%%%%%%%%%%%%%%%%%%%%%%%%%%%%%%%%%%%%%%%%%%%%%%%%%%%%%%%%%%%%%%%
\subsubsection{Free-$\beta$ case.}
%%%%%%%%%%%%%%%%%%%%%%%%%%%%%%%%%%%%%%%%%%%%%%%%%%%%%%%%%%%%%%%%%%%%%%%%%%%%%%%%%%%%%%%%%%%%
%%%%%%%%%%%%%%%%%%%%%%%%%%%%%%%%%%%%%%%%%%%%%%%%%%%%%%%%%%%%%%%%%%%%%%%%%%%%%%%%%%%%%%%%%%%%
\[\begin{array}{llll}
%% lprint(ssarsuumo(1,0));
(\delta_1\,h)_0 \!&\!=\!&\!
-f_1 \beta_1

%%lprint(ssarsuumo(1,1));
\\[1.ex]
(\delta_1\,h)_1 \!&\!=\!&\!
-2 f_2 \beta_2

+\frac{1}{2} f_1^2 \beta_1^2

%% lprint(ssarsuumo(1,2));
\\[1.ex]
(\delta_1\,h)_2 \!&\!=\!&\!
-3 f_3 \beta_3

+ f_1 f_2 \beta_1 \beta_2 

-\frac{1}{12} f_1^3 \beta_1^3

%% lprint(ssarsuumo(1,3));
\\[1.ex]
(\delta_1\,h)_3 \!&\!=\!&\!
-4 f_4 \beta_4

-4 f_4 \beta_4 

+f_1 f_3 \beta_1 \beta_3  
+\frac{1}{3} f_2^2 \beta_2^2

-\frac{1}{6} f_1^2 f_2 \beta_1^2 \beta_2 

+\frac{1}{144} f_1^4 \beta_1^4

%% lprint(ssarsuumo(1,4));
\\[1.ex]
(\delta_1\,h)_4 \!&\!=\!&\!
-5 f_5 \beta_5 

+  f_1 f_4 \beta_1 \beta_4 
+\frac{1}{2} f_2 f_3 \beta_2  \beta_3

-\frac{1}{8}  f_1^2 f_3 \beta_1^2\beta_3 
-\frac{1}{12} f_1 f_2^2 \beta_1\beta_2^2
\\[1.ex] &&\!
+\frac{1}{72} f_1^3 f_2 \beta_1^3 \beta_2

-\frac{1}{2880} f_1^5 \beta_1^5

%% lprint(ssarsuumo(1,5));
\\[1.ex]
(\delta_1\,h)_5 \!&\!=\!&\!
- 6 f_6 \beta_6 

+ f_1  f_5 \beta_1\beta_5 
+ \frac{2}{5} f_2 f_4  \beta_2 \beta_4 
+ \frac{3}{20} f_3^2 \beta_3^2 

-\frac{1}{10 } f_1^2  f_4 \beta_1^2\beta_4
\\[1.ex] &&\!
-\frac{1}{90}f_2^3 \beta_2^3

-\frac{1}{10} f_1  f_2 f_3 \beta_1\beta_2  \beta_3

+\frac{1}{120} f_1^3 f_3 \beta_1^3 \beta_3 
+ \frac{1}{120} f_1^2  f_2^2 \beta_1^2 \beta_2^2 
\\[1.ex] &&\!
-\frac{1}{1440} f_1^4  f_2 \beta_1^4 \beta_2 

+\frac{1}{86400} f_1^6 \beta_1^6

%%%%%%%%%%%%%%%%%%%%%%%%%%%%%%%%%%%%%%%%%%%%%%
%% lprint(ssarsuumo(2,0));
\\[1.ex]
(\delta_2\,h)_0 \!&\!=\!&\!
- f_2 \beta_1 

+  f_0 f_2\beta_2

-\frac{1}{4} f_0 f_1^2 \beta_1^2

%% lprint(ssarsuumo(2,1));
\\[1.ex]
(\delta_2\,h)_1 \!&\!=\!&\!
-3 f_3 \beta_2 

+3 f_0 f_3 \beta_3
+ f_1  f_2 \beta_1^2

-  f_0 f_1  f_2 \beta_1 \beta_2

+\frac{1}{12} f_0 f_1^3 \beta_1^3

%% lprint(ssarsuumo(2,2));
\\[1.ex]
(\delta_2\,h)_2 \!&\!=\!&\!
-6 f_4 \beta_3 

+\frac{3}{2} f_1 f_3 \beta_1 \beta_2
+6 f_0 f_4 \beta_4
+ f_2^2 \beta_1 \beta_2 

-\frac{1}{2} f_0 f_2^2 \beta_2^2
-\frac{1}{4} f_1^2  f_2 \beta_1^3
\\[1.ex] &&\!
-\frac{3}{2} f_0 f_1 f_3  \beta_1 \beta_3

+\frac{1}{4} f_0 f_1^2  f_2 \beta_1^2\beta_2

-\frac{1}{96} f_0 f_1^4 \beta_1^4

%% lprint(ssarsuumo(2,3));
\\[1.ex]
(\delta_2\,h)_3 \!&\!=\!&\!
-10 f_5 \beta_4 

+2  f_1 f_4 \beta_1 \beta_3
+  f_2 f_3 \beta_1 \beta_3
+10 f_0 f_5 \beta_5
+  f_2 f_3 \beta_2^2
\\[1.ex] &&\!
-2 f_0  f_1 f_4 \beta_1 \beta_4

-\frac{1}{4} f_1^2  f_3 \beta_1^2\beta_2
-\frac{1}{3} f_1  f_2^2 \beta_1^2\beta_2
- f_0  f_2 f_3 \beta_2  \beta_3

+\frac{1}{36} f_1^3  f_2\beta_1^4
\\[1.ex] &&\!
+\frac{1}{4} f_0 f_1^2  f_3 \beta_1^2 \beta_3

+\frac{1}{6} f_0  f_1 f_2^2 \beta_1 \beta_2^2

-\frac{1}{36} f_0 f_1^3  f_2 \beta_1^3 \beta_2

+\frac{1}{1440} f_0 f_1^5 \beta_1^5

%%%%%%%%%%%%%%%%%%%%%%%%%%%%%%%%%%%%%%%%%%%%%%%%%%%%%%%%%%%%%%%%%%%%%%%
%% lprint(ssarsuumo(3,0));
\\[1.ex]
(\delta_3\,h)_0 \!&\!=\!&\!
 -f_3 \beta_1 

+\frac{1}{3} f_1 f_2 \beta_2
+2 f_0 f_3 \beta_2
 
-\frac{1}{12} f_1^3 \beta_1^2
-  f_0^2 f_3\beta_3
-\frac{2}{3} f_0 f_1  f_2 \beta_1^2
\\[1.ex] &&\!
+\frac{1}{3} f_0^2 f_1  f_2  \beta_1 \beta_2

-\frac{1}{36} f_0^2 f_1^3 \beta_1^3

%% lprint(ssarsuumo(3,1));
\\[1.ex]
(\delta_3\,h)_1 \!&\!=\!&\!
- 4 f_4 \beta_2 

+  f_1 f_3 \beta_1^2
+  f_1 f_3 \beta_3
+\frac{2}{3} f_2^2 \beta_1^2
+8 f_0 f_4 \beta_3 

-4 f_0^2 f_4 \beta_4
\\[1.ex] &&\!
-\frac{1}{3} f_1^2  f_2 \beta_1\beta_2
-\frac{4}{3} f_0 f_2^2 \beta_1 \beta_2
-2 f_0 f_1  f_3 \beta_1\beta_2

+\frac{1}{36 } f_1^4 \beta_1^3
+\frac{1}{3} f_0^2  f_2^2\beta_2^2
\\[1.ex] &&\!
+\frac{1}{3} f_0 f_1^2 f_2 \beta_1^3
+f_0^2 f_1 f_3 \beta_1  \beta_3

-\frac{1}{6} f_0^2 f_1^2 f_2 \beta_1^2 \beta_2 

+\frac{1}{144} f_0^2 f_1^4 \beta_1^4

%% lprint(ssarsuumo(3,2));
\\[1.ex]
(\delta_3\,h)_2 \!&\!=\!&\!
-10 f_5 \beta_3 

+2  f_1 f_4 \beta_1 \beta_2
+20 f_0 f_5 \beta_4
+2  f_1 f_4\beta_4
+3  f_2 f_3 \beta_1 \beta_2 

-10 f_0^2 f_5 \beta_5
\\[1.ex] &&\!
-\frac{1}{4 } f_1^2 f_3 \beta_1^3
-\frac{1}{3} f_1  f_2^2 \beta_1^3
-\frac{1}{6} f_1 f_2^2 \beta_2^2
-\frac{1}{2} f_1^2 f_3\beta_1 \beta_3 
-2 f_0  f_2 f_3 \beta_2^2
\\[1.ex] &&\!
-2 f_0  f_2 f_3 \beta_1 \beta_3

-4 f_0  f_1 f_4\beta_1 \beta_3

+\frac{1}{12} f_1^3  f_2 \beta_1^2 \beta_2
+\frac{2}{3} f_0  f_1 f_2^2 \beta_1^2 \beta_2
\\[1.ex] &&\!
+\frac{1}{2}  f_0  f_1^2 f_3\beta_1^2 \beta_2
+f_0^2 f_2 f_3 \beta_2  \beta_3

+2 f_0^2  f_1 f_4 \beta_1  \beta_4

-\frac{1}{288} f_1^5 \beta_1^4
-\frac{1}{18} f_0 f_1^3  f_2 \beta_1^4
\\[1.ex] &&\!
-\frac{1}{6} f_0^2  f_1 f_2^2\beta_1 \beta_2^2
-\frac{1}{4} f_0^2  f_1^2  f_3\beta_1^2  \beta_3

+\frac{1}{36} f_0^2 f_1^3  f_2 \beta_1^3\beta_2

-\frac{1}{1440} f_0^2 f_1^5 \beta_1^5
\end{array}\]

% \end{document}
%%%%%%%%%%%%%%%%%%%%%%%%%%%%%%%%%%%%%%%%%%%%%%%%%%%%%%%%%%%%%%%%%%%%%%%%%%%%%%%%%%%%%%%%%%%%
%%%%%%%%%%%%%%%%%%%%%%%%%%%%%%%%%%%%%%%%%%%%%%%%%%%%%%%%%%%%%%%%%%%%%%%%%%%%%%%%%%%%%%%%%%%%
%%%%%%%%%%%%%%%%%%%%%%%%%%%%%%%%%%%%%%%%%%%%%%%%%%%%%%%%%%%%%%%%%%%%%%%%%%%%%%%%%%%%%%%%%%%%
%%%%%%%%%%%%%%%%%%%%%%%%%%%%%%%%%%%%%%%%%%%%%%%%%%%%%%%%%%%%%%%%%%%%%%%%%%%%%%%%%%%%%%%%%%%%
%%%%%%%%%%%%%%%%%%%%%%%%%%%%%%%%%%%%%%%%%%%%%%%%%%%%%%%%%%%%%%%%%%%%%%%%%%%%%%%%%%%%%%%%%%%%
%%%%%%%%%%%%%%%%%%%%%%%%%%%%%%%%%%%%%%%%%%%%%%%%%%%%%%%%%%%%%%%%%%%%%%%%%%%%%%%%%%%%%%%%%%%%
%%%%%%%%%%%%%%%%%%%%%%%%%%%%%%%%%%%%%%%%%%%%%%%%%%%%%%%%%%%%%%%%%%%%%%%%%%%%%%%%%%%%%%%%%%%%
%%%%%%%%%%%%%%%%%%%%%%%%%%%%%%%%%%%%%%%%%%%%%%%%%%%%%%%%%%%%%%%%%%%%%%%%%%%%%%%%%%%%%%%%%%%%
%%%%%%%%%%%%%%%%%%%%%%%%%%%%%%%%%%%%%%%%%%%%%%%%%%%%%%%%%%%%%%%%%%%%%%%%%%%%%%%%%%%%%%%%%%%%

%

%%%%%%%%%%%%%%%%%%%%%%%%%%%%%%%%%%%%%%%%%%%%%%%%%%%%%%%%%%
%% Some resurgence properties of knot-related functions.
%%%%%%%%%%%%%%%%%%%%%%%%%%%%%%%%%%%%%%%%%%%%%%%%%%%%%%%%%%

%%%%%%% \newcommand{\kos}[2]{[\frac{#1}{#2}] }
%\documentclass[12pt,a4paper]{article}\input{SP_commands}\begin{document}

%%%%%%%%%%%%%%%%%%%%%%%%%%%%%%%%%%%%%%%%%%%%%%%%%%%%%%%%%%%%%%%%%%%%%%%%%%%%%%%%%%%%%%%%%%%%
%%%%%%%%%%%%%%%%%%%%%%%%%%%%%%%%%%%%%%%%%%%%%%%%%%%%%%%%%%%%%%%%%%%%%%%%%%%%%%%%%%%%%%%%%%%%
%%%%%%%%%%%%%%%%%%%%%%%%%%%%%%%%%%%%%%%%%%%%%%%%%%%%%%%%%%%%%%%%%%%%%%%%%%%%%%%%%%%%%%%%%%%%

\section{Tables relative to the $4_1$ knot.}
%%%%%%%%%%%%%%%%%%%%%%%%%%%%%%%%%%%%%%%%%%%%%%%%%%%%%%%%%%%%%%%%%%%%%%%%%%%%%%%%%%%%%%%%%%%%
%%%%%%%%%%%%%%%%%%%%%%%%%%%%%%%%%%%%%%%%%%%%%%%%%%%%%%%%%%%%%%%%%%%%%%%%%%%%%%%%%%%%%%%%%%%%

\subsection{The original generators $\mi{Lo}$ and $\mi{Loo}$ .}
\[
\begin{array}{llll}
\kos{p}{q}&\!\!:=\!\!&\cos(\pi \frac{p}{q}) 
\\[1.7 ex]
\imi{\mi{Jo}}_2 &\!\!=\!\!& 4\,\kos{0}{2} 
\;\;\;;\;\;\;
\imi{\mi{Jo}}_4 \,=\, 26\,\kos{0}{4}
\;\;\;;\;\;\;
\imi{\mi{Jo}}_6 \,=\, 60\,\kos{0}{6}+56\,\kos{2}{6}
\\ [0.2 ex]
\imi{\mi{Jo}}_8 &\!\!=\!\!& 186\,\kos{0}{8}+90\,\kos{2}{8}
\;\;\;;\;\;\;
\imi{\mi{Jo}}_{10} \,=\, 348\,\kos{0}{10}+366\,\kos{2}{10}+22\,\kos{4}{10}
\\ [0.2 ex]
\imi{\mi{Jo}}_{12} &\!\!=\!\!& 650\,\kos{0}{12}+748\,\kos{2}{12}+624\,\kos{4}{12}
\\[0.2 ex]
\imi{\mi{Jo}}_{14} &\!\!=\!\!& 1396\,\kos{0}{14}+1854\,\kos{2}{14}+1030\,\kos{4}{14}
        +568\,\kos{6}{14}
\\ [0.2 ex]
\imi{\mi{Jo}}_{16} &\!\!=\!\!& 2776\,\kos{0}{16}+3804\,\kos{2}{16}+2816\,\kos{4}{16}
+1570\,\kos{6}{16}
\\ [0.2 ex]
\imi{\mi{Jo}}_{18} &\!\!=\!\!&
4862\,\kos{0}{18}+8078\,\kos{2}{18}+6550\,\kos{4}{18}+4802\,\kos{6}{18}+1472\,\kos{8}{18}
\\[0.2 ex]
\imi{\mi{Jo}}_{20} &\!\!=\!\!& 9864\,\kos{0}{20}
+16588\,\kos{2}{20}+14484\,\kos{4}{20}+10242\,\kos{6}{20}+4646\,\kos{8}{20}
\\
\end{array}\]

\[\begin{array}{llll}
\imi{\mi{Joo}}_{2} &\!\!=\!\!& 4\,\kos{0}{2}
 \\
\imi{\mi{Joo}}_{4} &\!\!=\!\!& 12\,\kos{0}{4}
 \\
\imi{\mi{Joo}}_{6} &\!\!=\!\!& 20\,\kos{0}{6}+16\,\kos{2}{6}
 \\
\imi{\mi{Joo}}_{8} &\!\!=\!\!& 44\,\kos{0}{8}+24\,\kos{2}{8}
 \\
\imi{\mi{Joo}}_{10} &\!\!=\!\!& 68\,\kos{0}{10}+72\,\kos{2}{10}+8\,\kos{4}{10}
 \\
\imi{\mi{Joo}}_{12} &\!\!=\!\!& 
108\,\kos{0}{12}+128\,\kos{2}{12}+96\,\kos{4}{12}
 \\
\imi{\mi{Joo}}_{14} &\!\!=\!\!& 
196\,\kos{0}{14}+264\,\kos{2}{14}+152\,\kos{4}{14}+80\,\kos{6}{14}
 \\
\imi{\mi{Joo}}_{16}  &\!\!=\!\!&  
340\,\kos{0}{16}+480\,\kos{2}{16}+352\,\kos{4}{16}+200\,\kos{6}{16} 
 \\
\end{array}\]

\[\begin{array}{llll}
 \imi{\mi{D\!Jo}}_{n}&\!\!\!:=\!\!& \imi{\mi{Jo}}_{n} -\frac{1}{2}\,\imi{\mi{Joo}}_{n} \;=\;
n\,(\imi{\mi{Lo}}_n - \imi{\mi{Loo}}_n)
\\[-0.4 ex] 
\dots &&\dots\dots\dots\dots\dots\dots\dots\dots\dots
\\[-0.4 ex]
\imi{\mi{D\!Jo}}_{2}  &\!\!=\!\!&   0 
 \\
\imi{\mi{D\!Jo}}_{4}  &\!\!=\!\!&  2\,\kos{0}{4}
 \\
\imi{\mi{D\!Jo}}_{6}  &\!\!=\!\!&  8\,\kos{2}{6}
 \\
\imi{\mi{D\!Jo}}_{8}  &\!\!=\!\!&  10\,\kos{0}{8}-6\,\kos{2}{8}
 \\
\imi{\mi{D\!Jo}}_{10}  &\!\!=\!\!&  8\,\kos{0}{10}+6\,\kos{2}{10}-18\,\kos{4}{10}
 \\
\imi{\mi{D\!Jo}}_{12}  &\!\!=\!\!&  2\,\kos{0}{12}-20\,\kos{2}{12}+48\,\kos{4}{12} 
 \\
\end{array}\]

%%%%%%%%%%%%%%%%%%%%%%%%%%%%%%%%%%%%%%%%%%%%%%%%%%%%%%%%%%%%%%%%%%%%%%%%%%%%%%%%%%%%%%%%%%%%
%%%%%%%%%%%%%%%%%%%%%%%%%%%%%%%%%%%%%%%%%%%%%%%%%%%%%%%%%%%%%%%%%%%%%%%%%%%%%%%%%%%%%%%%%%%%
 
\subsection{The outer generators $\mi{Lu}$ and $\mi{Luu}$.} 

\[\begin{array}{llll}

\smi{\mi{\ell u}}_{3}&\!\!=\!\!& \frac{
2}{3}\;\pi^{2}
\\[1.5 ex]
\smi{\mi{\ell u}}_{5}&\!\!=\!\!& \frac{ 
47}{90}\;\pi^{4}
\\[1.5 ex] 
\smi{\mi{\ell u}}_{7}&\!\!=\!\!& \frac{ 
12361}{28350}\;\pi^{6}
\\[1.5 ex] 
\smi{\mi{\ell u}}_{9}&\!\!=\!\!& \frac{
10771487}{28576800}\;\pi^{8}
\\[1.5 ex] 
\smi{\mi{\ell u}}_{11}&\!\!=\!\!& \frac{ 
23554574521}{70727580000}\;\pi^{10}
\\[1.5 ex] 
\smi{\mi{\ell u}}_{13}&\!\!=\!\!& \frac{
108745677770927}{364105581840000}\;\pi^{12}
\\[1.5 ex] 
\smi{\mi{\ell u}}_{15}&\!\!=\!\!& \frac{
941991271620333481}{
3479028834481200000}\;\pi^{14}
\\[1.5 ex] 
\smi{\mi{\ell u}}_{17}&\!\!=\!\!& \frac{
14052251175701893474367}{56777750578733184000000}\;\pi^{16}
\\
\end{array}\]

\[\begin{array}{llll}
 
\smi{\mi{\ell uu}}_{2}&\!\!=\!\!& \frac{ 
20}{81}\;\pi^{2}
\\[1.5 ex] 
\smi{\mi{\ell uu}}_{4}&\!\!=\!\!& \frac{  
1219}{6561}\;\pi^{4}
\\[1.5 ex] 
\smi{\mi{\ell uu}}_{6}&\!\!=\!\!& \frac{ 
401353}{2657205}\;\pi^{6}
\\[1.5 ex]  
\smi{\mi{\ell uu}}_{8}&\!\!=\!\!& \frac{  
36170973257}{281238577200}\;\pi^{8}
\\[1.5 ex]  
\smi{\mi{\ell uu}}_{10}&\!\!=\!\!& \frac{  
7690394022421}{68340974259600}\;\pi^{10}
\\[1.5 ex]  
\smi{\mi{\ell uu}}_{12}&\!\!=\!\!& \frac{ 
6713640059454013219}{66980988871833960000}\;\pi^{12}
\\[1.5 ex]  
\smi{\mi{\ell uu}}_{14}&\!\!=\!\!& \frac{  
270805989843610341382811}{2995215671777347922904000*Pi^14}
\\[1.5 ex]  
\smi{\mi{\ell uu}}_{16}&\!\!=\!\!& \frac{  
159786141600838021397411486857}{1940899755311721454041792000000}\;\pi^{16}
\\
\end{array}\]
%%%%%%%%%%%%%%%%%%%%%%%%%%%%%%%%%%%%%%%%%%%%%%%%%%%%%%%%%%%%%%%%%%%%%%%%%%%%%%%%%%%%%%%%%%%%
%%%%%%%%%%%%%%%%%%%%%%%%%%%%%%%%%%%%%%%%%%%%%%%%%%%%%%%%%%%%%%%%%%%%%%%%%%%%%%%%%%%%%%%%%%%%

\subsection{The inner generators $\mi{Li}$ and $\mi{Lii}$.}
\[\begin{array}{llll}
\\
\smi{\ell i}_{n}\!&\!\!:=\!&\! (\sqrt{3}\,\pi)^n\; \smi{\elliast}_{n} \quad\quad
{with} \quad\quad \smi{\elliast}_{n}\in\doQ^{+} \quad\quad \forall n \in -\frac{1}{2}+\doN
\\
\dots\dots&&\dots\dots\dots\dots\dots\dots\dots\dots\dots\dots\dots\dots
\\ 
\smi{{\elliast}}_{-1/2}\!&\!=\!&\!{\frac {1}{2}}
\quad\quad\quad
\smi{{\elliast}}_{{1/2}}={\frac {11}{108}}
\quad\quad\quad
\smi{{\elliast}}_{{3/2}}={\frac {697}{34992}}
\quad\quad\quad
\smi{{\elliast}}_{{5/2}}={\frac {724351}{141717600}}
\\[1.ex]
\smi{{\elliast}}_{{7/2}}\!&\!=\!&\!{\frac {278392949}{214277011200}}
\\[1.ex]
\smi{{\elliast}}_{{9/2}}\!&\!=\!&\!{\frac {244284791741}{728970392102400}}
\\[1.ex]
\smi{{\elliast}}_{{11/2}}\!&\!=\!&\!{\frac {1140363907117019}{12990252387264768000}}
\\[1.ex]
\smi{{\elliast}}_{{13/2}}\!&\!=\!&\!{\frac {212114205337147471}{9119157175859867136000}}
\\
\end{array}\]

%%%%%%%%%%%%%%%%%%%%%%%%%%%%%%%%%%%%%%%%%%%%%%%%%%%%%%%%%%%%%%%%%%%%%%%%%%%%%%%%%%%%%%%%%%%%
%%%%%%%%%%%%%%%%%%%%%%%%%%%%%%%%%%%%%%%%%%%%%%%%%%%%%%%%%%%%%%%%%%%%%%%%%%%%%%%%%%%%%%%%%%%%

\subsection{The exceptional generator $\mi{Le}$.}

\[\begin{array}{lllll}
\smi{\mi{\ell\!e}}_{2n+1}\!\!&\!\!\equiv\!\!&\! 0 
\\[1.ex]

\smi{\mi{\ell\!e}}_{2n}\!\!&\!\!\equiv\!\!&\!\big(\frac{\pi}{4}\big)^{2n}\,\lambda^{-3n-2}
\;\smi{\mi{\ell\!e}^\ast}_{2n} 
\quad\mi{with}\quad \quad\lambda:=\log2 \quad\quad \mi{and}
\\ [1.5 ex]

\smi{\mi{\ell\!e}^\ast}_{2n}
\!\!&\!\!\equiv\!\!&
\sum_{k=1}^{k=3n}\smi{\mi{\ell\!e}^\ast}_{2n,k}\;\lambda^k
\;\;,\;\;\smi{\mi{\ell\!e}^\ast}_{2n,k}\in\doQ\quad
\\\dotfill&\dots&\dotfill
\\[1. ex]

\smi{\mi{\ell\!e}^\ast}_0  &\!\!=\!\!& 
\frac{1}{2}\lambda 

\\[1.9 ex]
\smi{\mi{\ell\!e}^\ast}_2  &\!\!=\!\!& 
+\lambda
-\frac{1}{6}\lambda^3
  
\\[1.9 ex]
\smi{\mi{\ell\!e}^\ast}_4  &\!\!=\!\!& 
+\frac{10}{3}\lambda
+\frac{2}{3}\lambda^2
-\frac{2}{9}\lambda^3 
-\frac{1}{9}\lambda^4
+\frac{1}{108}\lambda^5
+\frac{7}{540}\lambda^6

\\[1.9 ex]
\smi{\mi{\ell\!e}^\ast}_6  &\!\!=\!\!& 

+\frac{112}{9}\lambda
+\frac{224}{45}\lambda^2
+\frac{26}{135}\lambda^3
-\frac{52}{135}\lambda^4
-\frac{43}{405}\lambda^5
+\frac{34}{2025}\lambda^6
+\frac{331}{24300}\lambda^7
\\
\end{array}\]
%%%%%%%%%%%%%%%%%%%%%%%%%%%%%%%%%%%%%%%%%%%%%%%%%%%%%%%%%%%%%%%%%%%%%%%%%%%%%%%%%%%%%%%%%%%%
%%%%%%%%%%%%%%%%%%%%%%%%%%%%%%%%%%%%%%%%%%%%%%%%%%%%%%%%%%%%%%%%%%%%%%%%%%%%%%%%%%%%%%%%%%%%

% \end{document}
%%%%%%%%%%%%%%%%%%%%%%%%%%%%%%%%%%%%%%%%%%%%%%%%%%%%%%%%%%%%%%%%%%%%%%%%%%%%%%%%%%%%%%%%%%%%
%%%%%%%%%%%%%%%%%%%%%%%%%%%%%%%%%%%%%%%%%%%%%%%%%%%%%%%%%%%%%%%%%%%%%%%%%%%%%%%%%%%%%%%%%%%%
%%%%%%%%%%%%%%%%%%%%%%%%%%%%%%%%%%%%%%%%%%%%%%%%%%%%%%%%%%%%%%%%%%%%%%%%%%%%%%%%%%%%%%%%%%%%
%%%%%%%%%%%%%%%%%%%%%%%%%%%%%%%%%%%%%%%%%%%%%%%%%%%%%%%%%%%%%%%%%%%%%%%%%%%%%%%%%%%%%%%%%%%%
%%%%%%%%%%%%%%%%%%%%%%%%%%%%%%%%%%%%%%%%%%%%%%%%%%%%%%%%%%%%%%%%%%%%%%%%%%%%%%%%%%%%%%%%%%%%
%%%%%%%%%%%%%%%%%%%%%%%%%%%%%%%%%%%%%%%%%%%%%%%%%%%%%%%%%%%%%%%%%%%%%%%%%%%%%%%%%%%%%%%%%%%%
%%%%%%%%%%%%%%%%%%%%%%%%%%%%%%%%%%%%%%%%%%%%%%%%%%%%%%%%%%%%%%%%%%%%%%%%%%%%%%%%%%%%%%%%%%%%
%%%%%%%%%%%%%%%%%%%%%%%%%%%%%%%%%%%%%%%%%%%%%%%%%%%%%%%%%%%%%%%%%%%%%%%%%%%%%%%%%%%%%%%%%%%%
%%%%%%%%%%%%%%%%%%%%%%%%%%%%%%%%%%%%%%%%%%%%%%%%%%%%%%%%%%%%%%%%%%%%%%%%%%%%%%%%%%%%%%%%%%%%

%

%%%%%%%%%%%%%%%%%%%%%%%%%%%%%%%%%%%%%%%%%%%%%%%%%%%%%%%%%%
%% Some resurgence properties of knot-related functions.
%%%%%%%%%%%%%%%%%%%%%%%%%%%%%%%%%%%%%%%%%%%%%%%%%%%%%%%%%%

% \documentclass[12pt,a4paper]{article}\input{SP_commands}\begin{document}

%%%%%%%%%%%%%%%%%%%%%%%%%%%%%%%%%%%%%%%%%%%%%%%%%%%%%%%%%%%%%%%%%%%%%%%%%%%%%%%%%%%%%%%%%%%%
%%%%%%%%%%%%%%%%%%%%%%%%%%%%%%%%%%%%%%%%%%%%%%%%%%%%%%%%%%%%%%%%%%%%%%%%%%%%%%%%%%%%%%%%%%%%
%%%%%%%%%%%%%%%%%%%%%%%%%%%%%%%%%%%%%%%%%%%%%%%%%%%%%%%%%%%%%%%%%%%%%%%%%%%%%%%%%%%%%%%%%%%%
%\setcounter{section}{11}
\section{Acknowledgments and references.}
%%%%%%%%%%%%%%%%%%%%%%%%%%%%%%%%%%%%%%%%%%%%%%%%%%%%%%%%%%%%%%%%%%%%%%%%%%%%%%%%%%%%%%%%%%%%
%%%%%%%%%%%%%%%%%%%%%%%%%%%%%%%%%%%%%%%%%%%%%%%%%%%%%%%%%%%%%%%%%%%%%%%%%%%%%%%%%%%%%%%%%%%%

\subsection{Acknowledgments.}
Since our interest in knot-connected power series (i.e. the series
$G^{NP}_\mathcal{K}$ and $G^{P}_\mathcal{K}$ associated with a knot $\mathcal{K}$\,: cf \S9.1)
 and the closely related notion of SP-series (a natural and conceptually more appealing generalisation, in terms of which we chose to reframe the problem)
was first awakened after the 2006 visit to Orsay of Stavros Garoufalidis and his pioneering
joint work with Ovidiu Costin and since, despite tackling
the problem from very different angles, we have been keeping in touch
for about one year,  comparing methods and results, we feel we owe it to the reader to outline the
main differences between our two approaches \--- to justify, as it were,
their parallel existence.
 
The very first step is the same in both cases\,: we all rely on
a quite natural method \footnote{see [O.C.] and \S2.3 of the present paper.}
for deducing the shape of a function's closest singularity, or singularities,
from the exact asymptotics of its Taylor coefficients at 0. 
\footnote{this convergence is hardly surprising\,: the functions on hand (knot-related or SP)
tend to verify no useable equations, whether differential or functional,
that might give us a handle on their analytic properties, and so the Taylor
coefficients are all we have to go by. Two of us (O.C. in [C2] and
 J.E. in a 1993 letter to prof. G.K. Immink) hit independently on the 
same method \--- which must also have occurred, time and again, more or less explicitely, 
to many an analyst grappling with singularities.}
 
But then comes the question of handling the other singularities \--- those farther afield
\--- and this is where our approaches start diverging. In [C.G1.]-[C.G4],
the idea  is to re-write the functions under investigation in the form of multiple
integrals amenable to the Riemann-Hilbert theory and then use the well-oiled machinery
that goes with that theory. In this approach, the global picture (exact location of the
singularities on the various Riemann leaves,
rough nature of these singularities etc) emerges first, and the exact description
of each singularity, while also achievable at the cost of some extra work, comes
second. 

Our own approach reverses this sequence\,: the local aspect takes precedence,
and we then piece the global picture together from the local data. To that end, we distinguish
three types of ``resurgence generators" (i.e. local singularities that generate the resurgence
algebra under alien derivation\,) : the actually occuring {\it inner} and {\it outer}
generators\footnote{while there are only two
{\it outer} generators (which may coalesce into one), there can be any number
of {\it inner} generators.}, and the  auxiliary {\it exceptional} or {\it
movable} generators. The basic object here is the {\it inner resurgence algebra}, spanned by the
{\it inner} generators, which recur indefinitely under alien derivation. The {\it outer}
generators, on the other hand, produce only {\it inner} ones under alien derivation.\footnote{ 
which is only natural, since the {\it outer} generators can be interpreted as infinite sums
of (self-reproducing) {\it inner} generators.} We give exact descriptions of
both the {\it inner} and {\it outer} generators by means of special
integro-differential functionals of infinite order\,: {\it nir}, {\it mir} and {\it nur}, {\it mur}.

So much for the local aspect. To arrive at the global picture,
we  resort to an auxiliary construct, the so-called {\it exceptional} or {\it movable} generators,
which are very useful on account of three features\,:
\\
(i) they depend on a arbitrary {\it base point}, which can be taken 
as close as we wish to any particular singularity we want to zoom in onto 
\footnote{thus bringing it within the purview of the method of Taylor coefficient
asymptotics (see above).}
\\
(ii) their own set of singularities include all the {\it inner}
generators of the SP function 
\\
(iii) they may also possess {\it parasitical} singularities\footnote{
this is  the case {\it iff} the driving function $F$ has at least
one zero or one pole.} (i.e. singularities
other than the above), but these 
always lie farther away from the
{\it base point} than the closest {\it inner} generators.
\\ Thus, by moving the base point around, we can reduce the global
investigation to a local, or should we say, semi-local one, and derive the full picture,
beginning with the crucial {\it inner algebra}.

A further difference between our approaches is this:
while O. Costin and S. Garoufalidis are more directly concerned with the knot-related
series $G^{NP}_\mathcal{K}$ and $G^{P}_\mathcal{K}$ and the so-called
{\it volume conjecture} which looms ominously over
the whole field, the framework we have chosen for our
investigation is that of SP-series, i.e. general Taylor series with  coefficients that are syntactically of
sum-product type.
But this latter difference might well be less than appears, since each of the two
methods would seem, in principle, to be capable of extension in both directions.

%%%%%%%%%%%%%%%%%%%%%%%%%%%%%%%%%%%%%%%%%%%%%%%%%%%%%%%%%%%%%%%%%%%%%%%%%%%%%%%%%%%%%%%%%%%%
%%%%%%%%%%%%%%%%%%%%%%%%%%%%%%%%%%%%%%%%%%%%%%%%%%%%%%%%%%%%%%%%%%%%%%%%%%%%%%%%%%%%%%%%%%%%
%%%%%%%%%%%%%%%%%%%%%%%%%%%%%%%%%%%%%%%%%%%%%%%%%%%%%%%%%%%%%%%%%%%%%%%%%%%%%%%%%%%%%%%%%%%%
%%%%%%%%%%%%%%%%%%%%%%%%%%%%%%%%%%%%%%%%%%%%%%%%%%%%%%%%%%%%%%%%%%%%%%%%%%%%%%%%%%%%%%%%%%%%
\subsection{References.}
{\bf [C1]} O. Costin, {\it Asymptotics and Borel summability,}  2008, CRC Press.
\\
{\bf [C2]} O. Costin, {\it Global reconstruction of analytic functions from local expansions,}
Preprint 2007, Ohio State University.
\\
{\bf [CG1]} O. Costin and S. Garoufalidis,
{\it Resurgence of the Kontsevich-Zagier power series, }
preprint 2006 math. GT/0609619.
\\
{\bf [CG2]} O. Costin and S. Garoufalidis,
{\it Resurgence of the fractional polylogarithms, }
preprint 2007, math. CA/0701743.
\\
{\bf [CG3]} O. Costin and S. Garoufalidis,
{\it Resurgence of the Euler-MacLaurin summation formula, }
preprint 2007, math. CA/0703641.
\\
{\bf [CG4]} O. Costin and S. Garoufalidis,(forthcoming paper on general knot-related series)
\\
{\bf [E1]} J. Ecalle, {\it Les fonctions resurgentes: Vol.1}, Publ. Math. Orsay, 1981\\
{\bf [E2]} J. Ecalle, {\it Les fonctions resurgentes: Vol.2}, Publ. Math. Orsay, 1981\\
{\bf [E3]} J. Ecalle, {\it Les fonctions resurgentes: Vol.3}, Publ. Math. Orsay, 1985\\
{\bf [E4]} J. Ecalle., {\it Six lectures on transseries etc...,}
in: Bifurcations and Periodic Orbits etc, D. Schlomiuk ed., 1993, Kluwer, p 75-184.
\\
{\bf [E5]} J. Ecalle, {\it Recent advances in the analysis of divergence and singularities},
in: Normal Forms etc, Y.Ilyashenko and C.Rousseau eds., 2004, Kluwer,p 87-186.
 \\
{\bf [ES]} J. Ecalle and Sh. Sharma, {\it On a remarkable resurgence algebra linked to sum-product Taylor series.},
(forthcoming).
\\
{\bf [LR]} J.-L. Loday and M.O. Ronco, {\it Hopf algebra of a permutation and of planar binary trees.}
J. Alg. Combinatorics, {\bf 15}, 2002, pp 253-270.
\\
{\bf [HNT]} F. Hivert, J.-C. Novelli, J.-Y. Thibon, {\it The algebra of binary search trees}.
 Theoret. Computer. Sc.. {\bf 339}, 2005, pp 129-165.
 \\
{\bf [SS1]} Sh. Sharma, {\it Some examples of SP-series\,: a theoritical and numerical investigation},
(forthcoming).
\\
{\bf [SS2]} Sh. Sharma, {\it SP-series and the interplay between their resurgence and differential properties},
(forthcoming).
%%%%%%%%%%%%%%%%%%%%%%%%%%%%%%%%%%%%%%%%%%%%%%%%%%%%%%%%%%%%%%%%%%%%%%%%%%%%%%%%%%%%%%%%%%%%
%%%%%%%%%%%the leaves of B%%%%%%%%%%%%%%%%%%%%%%

%%%%%%%%%%%%%%%%%%%%%%%%%%%%%%%%%%%%%%%%%%%%%%%%%%%%%%%%%%%%%%%%%%%%%%%%%%%%%%%%%%%%%%%%%%%%
%\end{document}
%%%%%%%%%%%%%%%%%%%%%%%%%%%%%%%%%%%%%%%%%%%%%%%%%%%%%%%%%%%%%%%%%%%%%%%%%%%%%%%%%%%%%%%%%%%%
%%%%%%%%%%%%%%%%%%%%%%%%%%%%%%%%%%%%%%%%%%%%%%%%%%%%%%%%%%%%%%%%%%%%%%%%%%%%%%%%%%%%%%%%%%%%
%%%%%%%%%%%%%%%%%%%%%%%%%%%%%%%%%%%%%%%%%%%%%%%%%%%%%%%%%%%%%%%%%%%%%%%%%%%%%%%%%%%%%%%%%%%%
%%%%%%%%%%%%%%%%%%%%%%%%%%%%%%%%%%%%%%%%%%%%%%%%%%%%%%%%%%%%%%%%%%%%%%%%%%%%%%%%%%%%%%%%%%%%
%%%%%%%%%%%%%%%%%%%%%%%%%%%%%%%%%%%%%%%%%%%%%%%%%%%%%%%%%%%%%%%%%%%%%%%%%%%%%%%%%%%%%%%%%%%%
%%%%%%%%%%%%%%%%%%%%%%%%%%%%%%%%%%%%%%%%%%%%%%%%%%%%%%%%%%%%%%%%%%%%%%%%%%%%%%%%%%%%%%%%%%%%
%%%%%%%%%%%%%%%%%%%%%%%%%%%%%%%%%%%%%%%%%%%%%%%%%%%%%%%%%%%%%%%%%%%%%%%%%%%%%%%%%%%%%%%%%%%%
%%%%%%%%%%%%%%%%%%%%%%%%%%%%%%%%%%%%%%%%%%%%%%%%%%%%%%%%%%%%%%%%%%%%%%%%%%%%%%%%%%%%%%%%%%%%
%%%%%%%%%%%%%%%%%%%%%%%%%%%%%%%%%%%%%%%%%%%%%%%%%%%%%%%%%%%%%%%%%%%%%%%%%%%%%%%%%%%%%%%%%%%%

%
\end{document}